\documentclass{article}
\usepackage{amsmath,amssymb}

\numberwithin{equation}{section}

\newtheorem{Thm}{Theorem}[section]
\newtheorem{Defi}[Thm]{Definition}
\newtheorem{Cor}[Thm]{Corollary}
\newtheorem{Lemma}[Thm]{Lemma}
\newtheorem{Prop}[Thm]{Proposition}
\newtheorem{Rem}[Thm]{Remark}

\newenvironment{thm}[0]{\par\vspace{2 mm}\Thm\sl}%
{\rm\par\vspace{2 mm}}
\newenvironment{defi}[0]{\par\vspace{2 mm}\Defi\sl}%
{\par\vspace{2 mm}}
\newenvironment{cor}[0]{\par\vspace{2 mm}\Cor\sl}%
{\rm\par\vspace{2 mm}}
\newenvironment{lemma}[0]{\par\vspace{2 mm}%
\Lemma\sl}{\rm\par\vspace{2 mm}}
\newenvironment{prop}[0]{\par\vspace{2 mm}\Prop\sl}%
{\rm\par\vspace{2 mm}}
\newenvironment{rem}[0]{\par\vspace{2 mm}\Rem\rm}{\par\vspace{2 mm}}
{\rm\par\vspace{2 mm}}
{\par\vspace{2 mm}}
\def\medno{\medbreak\noindent}

\def\qed{~\hfill$\square$\medbreak\noindent\medbreak}
\def\proof{\par\noindent{\bf Proof:}{\ }{\ }}

\def\naam#1{\label{#1}\footnote{{\bf LAB:} #1}}
\def\refer#1{\ref{#1}\footnote{{\bf REF:} #1}}
\def\bib#1{{\rm [#1]}}

\def\text#1{\;\;\;\;{\rm \hbox{#1}}\;\;\;\;}
\def\qquad{\quad\quad}

\def\minspace{\vspace{-2mm}}

\def\itema{\vspace{-1mm}\item[{\rm (a)}]}
\def\itemb{\minspace\item[{\rm (b)}]}
\def\itemc{\minspace\item[{\rm (c)}]}

\def\msy#1{{\mathbb #1}}
\def\C{{\msy C}}
\def\N{{\msy N}}
\def\Z{{\msy Z}}
\def\R{{\msy R}}
\def\D{{\msy D}}

\def\ga{\alpha}
\def\gb{\beta}
\def\gd{\delta}

\def\gf{\varphi}
\def\gg{\gamma}

\def\gl{\lambda}

\def\gs{\sigma}

\def\gD{\Delta}

\def\gS{\Sigma}

\def\fa{{\mathfrak a}}
\def\fb{{\mathfrak b}}
\def\fg{{\mathfrak g}}
\def\fh{{\mathfrak h}}

\def\fk{{\mathfrak k}}

\def\fm{{\mathfrak m}}
\def\fn{{\mathfrak n}}

\def\fp{{\mathfrak p}}
\def\fq{{\mathfrak q}}

\def\implies{\Rightarrow}
\def\to{\rightarrow}
\def\Re{{\rm Re}\,}

\def\inp#1#2{\langle#1\,,\,#2\rangle}
\def\hinp#1#2{\langle#1\,|\,#2\rangle}

\def\Ad{{\rm Ad}}
\def\End{{\rm End}}
\def\Hom{{\rm Hom}}

\def\ad{{\rm ad}}
\def\after{\,{\scriptstyle\circ}\,}

\def\pr{{\rm pr}}

\def\tensor{\otimes}
\def\iM{{\rm M}}

\def\ik{{\rm k}}
\def\ih{{\rm h}}
\def\iq{{\rm q}}

\def\iC{{\scriptscriptstyle \C}}
\def\biC{\C}

\def\cA{{\cal A}}

\def\cC{{\cal C}}
\def\cD{{\cal D}}
\def\cE{{\cal E}}
\def\cF{{\cal F}}
\def\cH{{\cal H}}

\def\cK{{\cal K}}
\def\cL{{\cal L}}
\def\cM{{\cal M}}

\def\cO{{\cal O}}
\def\cP{{\cal P}}
\def\cR{{\cal R}}

\def\cW{{\cal W}}

\def\coleq{\mathchar"303A=}%defines := with : as a relation
\mathcode`:="603A%defines : as a punctuation instead of a relation
\def\col{\,:\,}%redefines col as a colon, in light of above
\def\cO{{\cal O}}
\def\sphXvp{\Ci(\Xvp \col \tau)}
\def\Ci{C^\infty}

\def\sphXp{\Ci(\spXp \col \tau)}
\def\sphXpI{\Ci(\spXp \col \tau \col I)}
\def\CepXpI{\Ci(\spXp \col \tau \col I)}
\def\oCtau{{}^\circ\cC(\tau)}

\def\oC{{}^\circ {\cal C}}

\def\glob{{\rm glob}}
\def\FQYOmegad{\cE_{Q,Y}(\spXp\col \tau\col \Omega\col \gdmap)}
\def\FQYOmegadP{{\FQYOmegad}_\glob}
\def\FQOmegad{\cE_Q(\spXp\col\tau\col\Omega\col \gdmap)}
\def\FQOmegadP{{\FQOmegad}_\glob}
\def\cEQY{\cE_{Q,Y}}
\def\cEQ{\cE_Q}

\def\hyp{{\rm hyp}}

\def\cEzerohyp{\cE_0^{\rm hyp}}
\def\cEstarzero{\cEzerohyp}

\def\cEhypgL{\cEzerohyp(\spXp\col \tau\col \gL)}

\def\specfamgL{\cEhypgL}

\def\ep{{\rm ep}}
\def\Cep{C^\ep}
\def\ExppolXptau{\Cep(\spXp\col \tau)}

\def\Exppol{\Cep}
\def\CepQY{C^\ep_{Q,Y}}
\def\anfamQY{\CepQY(\spXp\col \tau \col \Omega)}
\def\cFep{\cF^\ep}
\def\hglob{{\rm hglob}}
\def\cEhypQYgd{\cE^\hyp_{Q,Y}(\spXp\col \tau\col \gd)}
\def\cEhypQgd{\cE^\hyp_{Q}(\spXp\col \tau\col \gd)}
\def\cEhypQgdhglob{\cE^\hyp_{Q}(\spXp\col \tau\col \gd)_{\rm hglob}}
\def\cEhypQYgdhglob{\cE^\hyp_{Q,Y}(\spXp\col \tau\col \gd)_{\rm hglob}}

\def\Cephyp{C^{\ep,{\rm hyp}}_0(\spXp\col \tau)}
\def\CephypQY{C^{\ep, \hyp}_{Q,Y}(\spXp\col \tau)}

\def\Vtau{V_\tau}
\def\VtauuvH{V_\tau^{\KM \cap uvHv^{-1}u^{-1}}}
\def\tauQ{\tau_Q}
\def\tauP{\tau_P}
\def\tauM{\tau_{\iM}}

\def\DGH{{\msy D}(\spGH)}
\def\DXv{{\msy D}(\spXv)}
\def\Igdgl{I_{\gd, \gl}}

\def\spGH{\spX}
\def\spX{{\rm X}}
\def\spXgawp{\spX_{\ga,w,+}}
\def\Xvp{\spX_{v,+}}
\def\XQv{\spX_{Q,v}}
\def\XQvp{\spX_{Q,v,+}}
\def\spXgaw{\spX_{\ga,w}}
\def\spXgawp{\spX_{\ga,w,+}}
\def\spXgaw{\spX_{\ga,w}}
\def\spXQv{\spX_{Q,v}}
\def\spXQvp{\spX_{Q,v,+}}

\def\spXp{\spX_+}

\def\MoneQ{M_{1Q}}

\def\starP{{}^*P}

\def\vH{vHv^{-1}}

\def\laur{{\rm laur}}

\def\DP{{\Delta(P)}}

\def\DrQ{\Delta_r(Q)}
\def\DsubQP{\Delta_Q(P)}

\def\fadc{\fa_\iC^*}
\def\faq{\fa_\iq}
\def\faqc{{\fa_{\iq\iC}}}
\def\faqd{\fa_{\iq}^*}
\def\faqdc{\fa_{\iq\iC}^*}
\def\faQq{\fa_{Q\iq}}
\def\faQqd{\fa_{Q\iq}^*}
\def\faQqdc{\fa_{Q\iq\iC}^*}
\def\faQqp{\fa_{Q\iq}^+}
\def\fagaq{{\fa_{\ga\iq}}}

\def\fbk{\fb_\ik}

\def\fbdc{\fb_\iC^*}

\def\fbkdc{\fb^*_{\ik\iC}}

\def\staQq{{}^*\fa_{Q\iq}}
\def\staQqp{{}^*\fa_{Q\iq}^+}
\def\staQqdc{{}^*\fa_{Q\iq\iC}^*}
\def\stAqp{{}^*\! A_{Q\iq}^+}

\def\stbQ{{}^*\fb_Q}
\def\stbQdc{{}^*\fb_{Q\iC}^*}
\def\starP{{}^*P}

\def\Aqp{A_\iq^+}
\def\Aq{A_\iq}
\def\AQqp{A_{Q\iq}^+}
\def\AQq{A_{Q\iq}}
\def\AQqp{A_{Q\iq}^+}

\def\minparabs{\cP_\gs^{\rm min}}
\def\allparabs{\cP_\gs}

\def\NKQaq{N_{K_Q}(\faq)}
\def\NKaq{N_K(\faq)}

\def\QW{{}^Q\cW}

\def\nC{C^\circ}
\def\nE{{E^\circ}}
\def\nEp{{E_+}}

\def\precD{\preceq_\Delta}
\def\simD{\sim_\Delta}
\def\setmid{\mid}

\def\cWQv{\cW_{Q,v}}
\def\WKH{W_{K \cap H}}

\def\gL{\Lambda}
\def\uz{{\underline z}}
\def\dotvar{\,\cdot\,}

\def\tsj#1{\c{#1}}

\def\bsl{\backslash}

\def\displaycup{\bigcup}

\def\Exp{{\rm Exp}}
\def\ExpL{{\rm Exp}_{\rm L}}

\def\pr{{\rm pr}}
\def\embed{\hookrightarrow}

\def\sing{{\rm sing}}
\def\reg{{\rm reg}}

\def\supp{{\rm supp}\,}

\def\dK{{\widehat K}}
\def\types{\vartheta}

\def\spXv{\spX_v}
\def\spXzerov{\spX_{0, v}}

\def\KM{K_\iM}
\def\asmid{\,|\,}
\def\DX{\D(\spX)}
\def\faQ{\fa_Q}
\def\SrQ{\Sigma_r(Q)}
\def\WKvH{W_{K\cap \vH}}

\def\Hyp{\cH}
\def\Mer{\cM}

\def\Eps{E_{+,s}}
\def\Lau{\cL}
\def\Laustar{\cL_*}

\def\DrQ{\Delta_r(Q)}

\def\dpr{d'}
\def\Hyppr{\Hyp'}
\def\stbQdc{{}^*\fb_{Q\iC}^*}
\def\staQqd{{}^*\fa_{Q\iq}^*}
\def\faQqdczero{\fa_{Q\iq\iC}^{*\circ}}

\def\rhotauu{\rho_{\tau,u}}
\def\embeds{\hookrightarrow}
\def\HM{H_\iM}
\def\HMone{H_{\iM_1}}

\def\Xov{\spX_{0,v}}

\def\DMoneH{\D(M_1/\HMone)}

\def\ev{{\rm ev}}

\def\CiXptau{\Ci(\spXp\col \tau)}

\def\ofkQ{\fk{\scriptstyle(Q)}}
\def\barfn{\bar \fn}
\def\Cartan{\theta}
\def\fmoneQ{\fm_{1Q}}
\def\HoneQ{H_{1Q}}
\def\barfnQ{\barfn_Q}

\def\MoneQpr{M_{1Q}'}
\def\detQ{\gd_Q}
\def\ringQ{\cR_Q}

\def\congruent{\equiv}
\def\radQalg{\ringQ \otimes \End(\Vtau) \otimes U(\fmoneQ)}
\def\MoneQppr{{M_{1Q,+}'}}
\def\MoneQp{M_{1Q,+}}
\def\radQalgplus{\cR_Q^+\otimes\End(\Vtau) \otimes U(\fmoneQ)}
\def\ringQplus{\cR_Q^+}

\def\restQ{{\rm R}_Q}

\def\wH{wHw^{-1}}
\def\TdownPw{T^\downarrow_{P,w}}
\def\TdownPv{T^\downarrow_{P,v}}
\def\TdownPminw{T^\downarrow_{\Pmin,w}}
\def\TdownQw{T^\downarrow_{Q,w}}
\def\VtauKwH{V_\tau^{\KM\cap w H w^{-1}}}
\def\VtauKvH{V_\tau^{\KM\cap v H v^{-1}}}
\def\VtauMKwH{V_\tau^{\KM\cap w H w^{-1}}}
\def\KMwH{\KM \cap w H w^{-1}}

\def\TdownPcW{T^\downarrow_{P,\cW}}

\def\Ap{A^+}
\def\spXzerow{\spX_{0,w}}
\def\uq{{\underline q}}
\def\stAQq{{}^*\!A_{Q\iq}}
\def\stAQqp{{}^*\!A_{Q\iq}^+}
\def\stfaQq{{}^*\fa_{Q\iq}}
\def\diffop{\partial}
\def\stdiffop{\partial}
\def\MQ{M_Q}
\def\fv{\mathfrak v}

\def\bs{\backslash}
\def\spXPvp{\spX_{P,v,+}}
\def\DrP{\Delta_r(P)}
\def\APqp{A_{P\iq}^+}
\def\stAPqp{{}^*\!A_{P\iq}^+}
\def\DPmin{\gD(\Pmin)}
\def\NKPaq{N_{K_P}(\faq)}
\def\spXPv{\spX_{P,v}}
\def\Hyp{\cH}
\def\Mer{\cM}
\def\simQP{\sim_{Q|P}}
\def\simPQ{\sim_{P|Q}}
\def\faPq{\fa_{P\iq}}

\def\faPqdc{\fa_{P\iq\iC}^*}
\def\WQP{W/\!\sim_{Q|P}}
\def\WPQ{W/\!\sim_{P|Q}}

\def\faQqdzero{\fa_{Q\iq\iC}^{*\circ}}

\def\bfn{\bar\fn}
\def\faPq{\fa_{P\iq}}
\def\DrP{\Delta_r(P)}
\def\DQP{\Delta_Q(P)}

\def\gdmap{\gd}

\def\Ep{E_+}

\def\Dr{\Delta_r}

\def\ppfaqdc{\bpr\bpr\faqdc}

\def\VtauMKH{V_\tau^{M \cap K \cap H}}

\def\PiSRaq{\Pi_{\Sigma, \R}}

\def\Ept{E_{+,t}}

\def\prQv{\pr_{Q,v}}

\def\Vc{V_\biC}
\def\natX{\N^X}

\def\SprojVX{S_\leftarrow(V,X)}
\def\ad{{\rm ad}\,}

\def\Sproj{S_\leftarrow}
\def\Laur{{\rm Laur}\,}
\def\subXr{X_r^0}
\def\XL{X_L}

\def\subXbL{X^0(L)}
\def\subXL{X_L^0}
\def\aff{\cA}
\def\VLperpc{V^\perp_{L\biC}}
\def\VLperp{V^\perp_L}

\def\HypbLX{\Hyp(L,X)}
\def\natHypbLX{\N^{\Hyp(L,X)}}

\def\faQqdc{\fa_{Q\iq\iC}^*}

\def\expset{\Xi}
\def\fagaqdc{{\fa_{\ga \iq \iC}^*}}

\def\subX{X^0}

\def\discdual{\widehat\spX_d}
\def\dK{\widehat K}

\def\stagaq{{}^*\fagaq}

\def\stagaq{{}^*\fagaq}

\def\cF{{\cal F}}
\def\fmQgs{\fm_{Q\gs}}
\def\MQgs{M_{Q\gs}}
\def\MPgs{M_{P\gs}}
\def\spXoneQv{\spX_{1Q,v}}
\def\spXoneQvp{\spX_{1Q,v,+}}
\def\cDQgs{\cD_{Q\gs}}
\def\cDoneQ{\cD_{1Q}}
\def\cDoneQgs{{}^\circ\cD_{1Q}}
\def\spXQ{\spX_Q}

\def\bpr{{}^{\scriptscriptstyle\backprime}}

\def\MQgsp{M_{Q\gs,+}}
\def\spXQep{\spX_{Q,e,+}}

\def\uplus{u_+}
\def\parone{\cP_\gs^1}

\def\ZKaq{Z_K(\faq)}

\def\faqp{\fa_\iq^+}
\def\simeqarrow{\;\;{\buildrel \simeq \over \longrightarrow} \;\;}

\def\MQgsb#1{\MQgs[#1]}
\def\MoneQb#1{\MoneQ[#1]}
\def\fhoneQ{\fh_{1Q}}
\def\XQvb#1{\XQv[#1]}
\def\XQvpb#1{\XQvp[#1]}
\def\XQepb#1{\spX_{Q,e,+}[#1]}
\def\br#1{[#1]}
\def\lbr#1{_{[#1]}}
\def\MQgsb#1{M_{Q\gs}[#1]}
\def\cDoneQplus{\cD_{1Q}^+}
\def\spXu{\spX_{u}}
\def\spXup{\spX_{u,+}}
\def\spXoneQv{\spX_{1Q,v}}
\def\allparabsv{\cP_{\gs^v}}
\def\AQ{A_Q}

\def\tauuQ{\tau_{uQu^{-1}}}
\def\spher{\varsigma}
\def\tauQga{\tau_{\stPga}}

\def\Vtypes{{V}_\types}

\def\DQmaps{{\rm D}_Q}

\def\hide#1{{}}
\def\refappa{A}

\def\Pmin{P_1}%%
\def\starPmin{{}^*P_1}%%
\def\rmD{{\rm D}}

\def\oCQv{\oC(Q,v)}

\def\stPga{P_\ga}
\def\stPo{P_0}
\def\spXP{\spX_P}
\def\spXPp{\spX_{P,+}}
\def\Ps{{P,s}}
\def\starstPo{{}^*P_0}
\def\stagaqdc{{}^*\fa_{\ga \iq\iC}^*}
\def\cpta{\cK'}
\def\cptb{\cK}
\def\WQQ{W/\sim_{Q|Q}}

\def\dega{{\rm deg}_{\rm a}}
\def\rega{{\rm reg}_{\rm a}}

\def\iQv{{\rm i}_{Q,v}}
\def\rmi{{\rm i}}

\def\naam{\label}

\def\refer{\ref}
\def\bib#1{\cite{#1}}

\begin{document}
%%%%%%%%%%%%%%%%%%%%%%%%%%%%%%%%%%%%%%%%%%%%%%%%
\title{Analytic families of eigenfunctions\\
on a reductive symmetric space}
\author{E.P. van den Ban and H. Schlichtkrull}
\date{}%24/11-2000, revised 7/2-2001
\maketitle
\begin{abstract}
The asymptotic behavior of holomorphic families of generalized
eigenfunctions on a reductive symmetric space is studied.
The family parameter is a complex character on the split component of a 
parabolic subgroup. The main result asserts 
that the family vanishes if a particular asymptotic coefficient 
does. This allows an induction of relations between families
that will be applied in forthcoming 
work on the Plancherel and the Paley-Wiener theorem.
\end{abstract}
%%%%%%%%%%%%%%%%%%%%%%%%%%%%%%%%%%%%%%%%%%%%%%%%
\tableofcontents
\section*{Introduction}
In harmonic analysis on a reductive symmetric space 
$\spX$ an important role is played
by families of generalized eigenfunctions for the algebra 
$\DX$ of invariant differential 
operators. Such families arise for instance as matrix coefficients of
 representations
that come in series, such as the (generalized) principal series. 
In particular, relations between such families are of great interest. 
We recall that a real reductive group $G,$ equipped with the left times
right multiplication action, is a reductive symmetric space. 
In the case of the group, examples of the mentioned relations are
functional equations for Eisenstein integrals, 
see \bib{HCEis} and \bib{HCIII}, or Arthur-Campoli relations
for Eisenstein integrals, see \bib{Arthur}, \bib{Campoli}.
In this paper we develop a general tool to 
establish relations of this kind. 
We show that they can be derived from similar relations satisfied 
by the family of functions obtained by taking 
one particular coefficient in a
certain asymptotic expansion. Since the functions in the  
family so obtained are eigenfunctions on symmetric spaces of lower split rank, 
this yields a powerful inductive method; we call it 
{\em induction of relations.} In the case of the group, a closely related
lifting theorem by Casselman was used by Arthur in the proof of
the Paley-Wiener theorem, see \bib{Arthur}, Thm.\ II.4.1. However,
no proof seems yet to have appeared of Casselman's theorem.

The tools developed in this paper are used in \bib{BSfi},
and they will also be applied 
in the forthcoming papers \bib{BSpl} and \bib{BSpw}.
For example, it is the induction   
of relations that allows us to establish symmetry properties
of certain integral kernels appearing in a Fourier inversion formula in 
\bib{BSfi}. Also in \bib{BSfi}, the induction of relations is used
to define generalized Eisenstein integrals corresponding to
non-minimal principal series. In \bib{BSpl}, the results of this paper will
be applied to identify these `formal' Eisenstein integrals 
with those defined in Delorme \bib{DelEis}. This is a key 
step towards the Plancherel decomposition.
The results will also be applied to 
establish functional equations for the Eisenstein integrals. 
Applied in this manner our technique serves as a replacement
for the use of the Maass-Selberg relations as in Harish-Chandra \bib{HCIII}
and \bib{DelEis}. On the other hand, in \bib{BSpw}
we apply our tool to show that Arthur-Campoli relations
satisfied by normalized Eisenstein integrals of spaces of lower split
rank induce similar relations for normalized 
Eisenstein integrals of $\spX.$ This result is then used to prove
a Paley-Wiener theorem for $\spX$ that generalizes Arthur's theorem
for the group. In particular, the missing proof of Casselman's
theorem will then be circumvented by means of a technique of 
the present paper. 

It should be mentioned that in the case of 
the group, induction of Arthur-Campoli relations for unnormalized 
Eisenstein integrals is easily derived from their 
integral representations (see \bib{Arthur}, p.\ 77, proof of Lemma 2.3). 
For normalized Eisenstein integrals, which are
not representable by integrals, the result seems to be much 
deeper, also in the group case. 

One of the interesting 
features of the theory is that it also deals with families of functions
that are not necessarily globally defined on the space 
$\spX$ but on a suitable open dense subset.

Asymptotic behavior of eigenfunctions on a symmetric space has been 
studied at many other places in the literature. The following
papers hold results that are related to some of the ideas
of the present paper
\bib{HCsf}, \bib{Gang}, \bib{TroVar}, \bib{HCI}, \bib{HCIII},   \bib{Helg}, 
\bib{KKMOOT}, \bib{OS}, \bib{CM},  
 \bib{Wal}, \bib{Arthur}, \bib{Osh}, \bib{BS}, \bib{Carmona}.

The core results of this paper were found and announced 
in the fall of 1995, when both authors were guests at the 
Mittag-Leffler Institute. In the same period Delorme announced 
his proof of the Plancherel 
theorem, which has now appeared in \bib{Delpl}.

We shall now explain the contents of this paper in more detail.
The space $\spX$ is of the form $G/H,$ with $G$ a real reductive
Lie group of Harish-Chandra's class, and $H$ an open subgroup 
of the set of fixed points for an involution $\gs$ of~$G.$ 

The group $G$ has a $\gs$-stable maximal compact subgroup $K,$ 
let $\Cartan$ be the associated Cartan 
involution of $G$. 
Let $P_0 = M_0 A_0 N_0$ be a fixed minimal 
$\gs \after \Cartan$-invariant parabolic subgroup of $G,$ 
with the indicated Langlands decomposition.
The Lie algebra
$\fa_0$ of $A_0$ is invariant 
under the infinitesimal involution 
$\gs;$ we denote the associated $-1$ eigenspace
in $\fa_0$ by $\faq.$ Its dimension is called 
the split rank of $\spX.$ 
Let $\Aq$ be 
the  vectorial subgroup of $G$ with Lie algebra $\faq$ 
and let $\Aq^\reg$ be the set of regular points relative to the adjoint
action of $\Aq$ in $\fg$. Then $\spXp\coleq K \Aq^\reg H$ is a 
$K$-invariant dense open subset of $\spX$. Let
$\Aqp$ 
be the open chamber in $\Aq$ determined by $P_0.$ 
Then
$\spXp$ is a finite union of disjoint sets of the form
$K\Aqp vH,$ with $v$ in the normalizer of $\faq$ in $K.$ 
In this introduction we assume, for simplicity of exposition,
that $\spXp = K\Aqp H.$ This assumption is actually fulfilled in the 
case that $\spX$ is a group.

Let $(\tau, \Vtau)$ be a finite
dimensional continuous representation of $K.$ 
Then by $\Ci(\spXp \col \tau)$ 
we denote the space of smooth functions 
$f: \spXp \to \Vtau$ that are $\tau$-spherical,
i.e., $f(kx) = \tau(k)f(x),$ for all $x \in \spXp$ and $k \in K.$

Let $\allparabs$ denote the (finite) set of $\gs\after \Cartan$-invariant
parabolic subgroups of $G$ containing $\Aq.$ 
Let $Q = M_Q A_Q N_Q$ be an element of $\allparabs.$ 
Then $\gs$ restricts to an involution
of $\faQ,$ the Lie algebra of $A_Q;$ we denote its $-1$ eigenspace 
by $\faQq.$ 
In the first part of the paper we study 
a family $f$  of the following type (cf.\ Definition
\refer{d: anfamQY newer}).
The family is a smooth map of the form
$$ 
f: \Omega \times \spXp \to \Vtau,
$$
with $\Omega$ an open subset of $\faQqdc,$ 
the complexified linear dual of $\faQq.$ It is assumed 
that  $f$ is holomorphic in its first variable.
Moreover, for every $\gl \in \Omega$ the function 
$f_\gl:=  f(\gl, \dotvar)$ belongs 
to $\Ci(\spXp \col \tau).$ It is furthermore assumed that the functions
$f_\gl$ allow suitable exponential polynomial
expansions along $\Aqp.$ More precisely, we assume,
for $m \in M_0$ and $a \in \Aqp,$ that 
\begin{equation}
\naam{e: prototype expansion}
f_\gl (ma) = \sum_{s \in W/W_Q} 
a^{s\gl - \rho_{P_0}} \sum_{\xi \in - sW_QY + \N\Sigma(P_0)} a^{-\xi} q_{s,\xi}
(\gl, \log a, m).
\end{equation} 
Here $W$ is the Weyl group of $\Sigma = \Sigma(\fg, \faq)$ and 
$W_Q$ is the centralizer of $\faQq$ in $W.$ 
Moreover, $\Sigma(P_0)$ denotes the collection of 
roots from $\Sigma$  occurring 
in $N_0$ and $Y$ is a finite subset of $\staQqdc,$ the annihilator 
of $\faQq$ in $\faqdc.$  
Finally, the
$q_{s, \xi}$ are smooth functions, holomorphic in the 
first and polynomial in the second variable. Thus, we impose a 
limitation on the set of exponents and assume that the 
coefficients depend holomorphically on the parameter~$\gl.$  
The type of convergence that we impose on the expansion
(\refer{e: prototype expansion}) is described in general terms in the 
preliminary Section \refer{s: exp pol series}.

We show that the functions
$f_\gl$ actually allow exponential polynomial expansions 
similar to (\refer{e: prototype expansion}) along 
any (possibly non-minimal) $P \in \allparabs.$ 
These expansions are investigated in detail in Sections  \refer{s: asymp walls} 
and \refer{s: analytic families}. Their coefficients are families of 
$\tau|_{M_P\cap K}$-spherical functions on $\spXPp,$ the analogue of 
$\spXp$ for 
the lower split rank symmetric space $\spXP:= M_P/M_P \cap H.$ 

The operators from $\DX$ do also allow expansions along every
 $P\in \allparabs.$ 
In Section \refer{s: diff op along walls}
this is shown by investigating a radial decomposition that
reflects the decomposition
$G = K M_P A_{Pq} H.$ 
It is of importance that the coefficients in these expansions
are globally defined smooth functions on $M_P,$ 
see Prop.\ \refer{p: radial deco with MQ}
and Cor.\ \refer{c: cor on ringQ}.
{}From the expansions we derive that the algebra $\DX$ acts on the space 
of families of the above type, see Prop.\ \refer{p: D on families}.

In Section \refer{s: asymptotic globality} we 
introduce the notion of asymptotic 
$s$-globality of a family  along $P.$ 
Losely speaking, it means that the coefficients $q_{s, \xi}(\gl, \log a , \dotvar)$ of the expansion along $P$ 
extend smoothly from $\spXPp$ to the full space $\spXP,$ for every 
$\xi \in (sW_Q Y - \N\Sigma(P))|_{\faPq}.$ This notion is proved to be stable under the action of $\DX.$ 

In Section \refer{s: vanishing thm} we impose three other conditions on the family.
The first is that each 
member satisfies a system of differential equations of the form
$$ 
Df_\gl = 0\qquad (D \in \Igdgl).
$$ 
Here $\Igdgl$ is a certain cofinite ideal in the algebra $\DX$  
depending polynomially on $\gl\in \faQqdc$ in a suitable way. 
Accordingly, $\gl$ is called
the spectral parameter of the family.
The second condition imposed is a suitable condition of asymptotic globality along certain
parabolic subgroups $P$ with $\dim(\faq/\faPq) =1.$
Thirdly, it is required that the domain $\Omega$ for the
parameter $\gl$ is unbounded in certain directions
(see Defn.~\refer{d: Q-distinguished}).

The first main result of the paper is then the following vanishing theorem, see 
Theorem \refer{t: vanishing theorem new}.
\medbreak\noindent
{\bf The vanishing theorem.\ }\sl
Let $f$ be a family as above, and assume that the coefficient of $\gl - \rho_Q$ 
in the expansion along $Q$ vanishes for $\gl$ in a non-empty open subset of $\Omega.$ 
Then the family $f$ is identically zero.\rm
\medbreak
In the proof the globality assumption is needed to link suitably many asymptotic coefficients
together; the vanishing of one of them then inductively causes the vanishing of others.
In the induction step a key role is played by the observation that a symmetric space 
cannot have a continuum of discrete series.
 
The importance of the vanishing theorem is that it applies to 
many families that naturally 
arise in representation theory. 
In the present paper we show that this is so for Eisenstein integrals
associated with the minimal principal series for $\spX;$ in \bib{BSpl} 
we will show that Eisenstein integrals obtained by parabolic induction from discrete
series form a family of the above type. The idea is that the latter Eisenstein
integrals can be obtained from those associated with the minimal principal
series by the application of residual operators with respect
to the spectral parameter. Such residual operators occur
in our papers \bib{BSres} and \bib{BSfi}. 

A suitable class of operators
containing the residual operators is formed by the Laurent operators.
In the second half of the  paper we study the application of 
them 
to suitable families
of eigenfunctions, with respect to the spectral parameter.  The Laurent operators
are best described by means of Laurent functionals, see Sections \refer{s: Laurent functionals}
and \refer{s: Laurent operators}. 

In Section
\refer{s: spec fam new} we introduce a special type of families 
$g$ of eigenfunctions.
It is of the above type, with $\Omega$ dense in $\faPqdc,$ 
$P$ a minimal parabolic subgroup in $\allparabs,$ and satisfies some additional
requirements, see Definition \refer{d: defi cEhyp Q Y gd}. One of these is that
the family and its asymptotic expansions should depend meromorphically on the spectral 
parameter $\gl \in \faPqdc$ with singularities along translated root hyperplanes.
This allows the application of Laurent functionals with respect to the spectral parameter.
More precisely, let $Q \in \allparabs$ contain $P,$ 
and let $\Lau$ be a Laurent functional on $\staQqdc.$ 
{}From the family $g$ 
a new family $f = \cL_* g,$ with a spectral 
parameter from $\faQqdc,$ 
is obtained by the application of
$\cL$  to the $\staQqdc$-component
of the spectral parameter.
In Theorem \refer{t: source of functions by Lau new} it is shown that 
the resulting family $\cL_* g$ satisfies the requirements
of the vanishing theorem, provided the special family $g$ satisfies certain holomorphic asymptotic
globality conditions.

In Section \refer{s: partial Eisenstein integrals} we introduce partial Eisenstein 
integrals associated with a minimal parabolic subgroup $P$ from $\allparabs.$ 
The partial Eisenstein integrals 
are spherical generalized eigenfunctions on $\spXp$ obtained from the normalized
Eisenstein integral $\nE(P\col \gl),$ ($\gl \in \faqdc$ generic), by splitting it
according to its exponential polynomial expansion along $P.$ 
More precisely, the exponents
of $\nE(P\col \gl)$ are contained in $W\gl - \rho_P -\N\Sigma(P);$ the partial Eisenstein
integrals $\Eps(P\col \gl),$ for $s \in W,$ are the smooth spherical functions on $\spXp$ 
determined by the requirements that
$$
\nE(P\col \gl) = \sum_{s \in W} \Eps(P\col \gl)
$$ 
and the set of exponents of $\Eps(P \col \gl)$ along $P$ should be contained in 
$s\gl - \rho_P -\N\Sigma(P).$ 
It is then  shown that the partial Eisenstein integrals yield examples 
of the special families mentioned above. Moreover, if $Q\in \allparabs,$ $Q\supset P,$ 
let $W^Q$ be the collection of minimal length 
(with respect to  $\Sigma(P)$) coset representatives for $W/W_Q$ in $W.$ 
Then it is shown that for each $t\in W_Q$ the family 
\begin{equation}
\naam{e: g as sum of partial Eis}
f_t = \sum_{s \in W^Q} E_{+,st}(P\col \dotvar)
\end{equation}
satisfies the additional holomorphic asympotic globality property 
guaranteeing that $\cL_* f_t$ satisfies the hypothesis of 
the vanishing theorem, for $\Lau$ a Laurent 
functional on $\staQqdc.$ 

In Section \refer{s: asymptotics of partial Eisenstein integrals} 
the asymptotic behavior of $\cL_* f_t$ is investigated, and the 
coefficient of $a^{\gl - \rho_Q}$ in the expansion along $Q$ is expressed
in terms of partial Eisenstein integrals of $\spXQ.$ 

The above preparations pave the way for the induction of relations
in Section \refer{s: induction of relations}. The idea is as follows.
Let $f_t$  be the family defined by 
(\refer{e: g as sum of partial Eis}), and
let a Laurent functional $\cL_t$ on $\staQqdc$ be given for each $t\in W_Q$. 
Then by the vanishing theorem 
a relation of the form $\sum_t \cL_t f_t=0$ is valid if a similar 
relation is valid for the $(\gl - \rho_Q)$-coefficients along $Q;$ 
this in turn may be expressed as a similar relation between partial Eisenstein 
integrals for the lower split rank space $\spXQ.$ 
In this setting, taking the $(\gl -\rho_Q)$-coefficient along $Q$ 
essentially inverts the procedure 
of parabolic induction from $Q$ to $G.$ This motivaties our choice
of terminology.  The precise result is formulated in Theorem
\refer{t: induction of relations, new}. An equivalent result, closer
to the formulation of Casselman's theorem in \bib{Arthur} is stated
at the end of the section.

\section{Exponential polynomial series}
\naam{s: exp pol series}
Let $A$ be a vectorial group and $\fa$ 
its Lie algebra. The exponential map 
$\exp: \fa \to A$ is a diffeomorphism;
we denote its inverse by $\log.$ 
If $\xi$ belongs to $\fadc,$ the complexified linear dual of $\fa,$ 
then we define the function
$e^\xi: a \mapsto a^\xi$ on $A$ by $a^\xi = e^{\xi(\log a)}.$
Let $P(\fa)$ denote the algebra of  
polynomial functions $\fa \to \C.$ If $d \in \N,$ let 
$P_d(\fa)$ denote the (finite dimensional) subspace
of polynomials of degree at most $d.$ Let $\Delta$ be
a set of linearly independent vectors in $\fa$ 
(we do not require
this set to span $\fa$).
We put 
$$
\fa^+ = \fa^+(\Delta) := \{ X \in \fa \mid \ga(X) > 0, \;\;\;\forall\;\ga \in \Delta\},
$$
and $A^+ = A^+(\Delta) = \exp (\fa^+).$ 
We define $\N \gD$ to be the $\N$-span of $\gD;$ if $\gD = \emptyset$ then 
$\N\gD = \{0\}.$ 
Moreover, if $X$ is a subset of $\fadc,$ we denote by $X- \N \gD$ the 
vectorial sum of $X$ and $\N \gD.$ 

Let $V$ be a complete locally convex space; here and in the following we
will always assume such a space to be Hausdorff. 
If $\xi \in \fadc,$ then by a $V$-valued $\xi$-exponential polynomial function on $A$
we mean a function $A \to V$ of the form $a \mapsto a^\xi q(\log a),$ with 
$q \in P(\fa)\otimes V.$ 

\begin{defi}
\naam{d: exp pol series}
By a $\Delta$-exponential polynomial series on $A$ 
with coefficients in $V$ we mean a  formal series $F$ 
of exponential polynomial functions of the form
\begin{equation}
\naam{e: intro exp pol series}
\sum_{\xi \in \fadc} a^\xi \,q_\xi (\log a),
\end{equation}
with $\xi \mapsto q_\xi$ a map $\fadc \to P(\fa) \otimes V,$  such that
\begin{enumerate}
\itema
there exists a finite subset $X \subset \fadc$ such that $q_\xi = 0$ for $\xi \notin X - \N \gD;$ 
\itemb
there exists a constant $d \in \N$ such that $q_\xi \in P_d(\fa) \otimes V$ 
for all $\xi\in \fadc.$ 
\end{enumerate}
The smallest $d \in \N$ with property (b) will be called the polynomial degree of the series;
this number is denoted by $\deg(F).$  

The collection of all  $\gD$-exponential polynomial series 
with coefficients in $V$ is denoted by 
$\cFep(A,V) = \cFep_\gD(A,V).$ 
\end{defi}

If $F \in \cFep(A,V)$ is an expansion of the form (\refer{e: intro exp pol series}) 
then, for every $\xi \in \fadc,$ 
we write $q_\xi(F)$ for $q_\xi.$ Moreover, we write $q(F)$ for the map 
$\xi \mapsto q_\xi(F)$ 
from $\fadc$ to $P_d(\fa)\otimes V.$ Then 
$F \mapsto q(F)$ 
defines a bijection from $\cFep(A,V)$ onto a linear subspace of 
$(P_d(\fa)\otimes V)^{\fadc},$ the space of maps $\faqdc \to P_d(\fa) \otimes V.$ 
Via this bijection we equip $\cFep(A,V)$ with the structure of a linear space.

If $F \in \cFep(A,V),$ then
$$
\Exp(F): = \{ \xi \in \fadc \mid q_\xi(F) \not= 0 \}
$$ 
is called the set 
of exponents of  $F.$ If $F_1, F_2 \in \cFep(A,V),$ 
we call $F_1$ a subseries of $F_2$ if 
$q_\xi(F_2)  = q_\xi(F_1)$ 
for all $\xi \in \Exp(F_1).$ 

The series (\refer{e: intro exp pol series}) 
is said to converge absolutely 
in a fixed point $a_0 \in A$ 
if  the series 
$$
\sum_{\xi \in \Exp(F)} a_0^{\xi} q_\xi(\log a_0)
$$
with coefficients in $V$ 
converges absolutely. It is said to converge absolutely on a subset $\Omega \subset A$ 
if it converges absolutely in every point $a_0 \in \Omega.$ In this case 
pointwise summation of the series defines a function $\Omega \to V.$ 
 
We will also need a more special type of convergence for the series
(\refer{e: intro exp pol series}).
 
\begin{defi}
\naam{d: neat convergence}
The series (\refer{e: intro exp pol series}) 
is said to converge neatly at a fixed point $a_0 \in A$ 
if for every continuous seminorm $s$ on $P_d(\fa) \tensor V,$ where $d = \deg(F),$ 
the series 
$$
\sum_{\xi \in \Exp(F)} s(q_\xi) a_0^{\Re \xi}
$$
converges. 

The series (\refer{e: intro exp pol series}) is said to converge neatly on a
subset $\Omega$ of $A$ if it converges neatly at every
point of $\Omega.$
\end{defi}

\begin{rem}
If the series (\refer{e: intro exp pol series}) converges neatly at a point 
$a_0 \in A,$ then so does every subseries. Moreover, neat convergence at $a_0$ 
implies absolute convergence in $a_0.$ However, we should warn the reader
that neat convergence at $a_0$ cannot be seen from the 
series with coefficients in $V$ arising from (\refer{e: intro exp pol series}) by evaluation
of its terms at $a = a_0,$ since this type of convergence involves
the global behavior of the polynomials $q_\xi.$ In particular, it is 
possible that the series (\refer{e: intro exp pol series}) 
does not converge neatly at $a_0,$ whereas 
its evaluation in $a_0$ is identically zero.

The motivation for the definition of neat convergence is provided 
later by Lemmas \refer{l: neat conv of exp pol series} and
\refer{l: formal application of Ufa}, which
express that neat convergence of the series (\refer{e: intro exp pol series})
on an open subset $\Omega \subset A$ guarantees that 
(a) the function $f: \Omega \to V$ 
defined by   (\refer{e: intro exp pol series}) 
is real analytic on $\Omega;$ (b) its derivatives
are given by series obtained by termwise differentiation 
from (\refer{e: intro exp pol series}).
\end{rem}

By a $\Delta$-power series on $A,$ with coefficients 
in  $V,$ 
we mean a $\gD$-exponential polynomial series $F$ with $\deg F = 0$ 
and $\Exp(F) \subset - \N \gD,$ i.e.,
\begin{equation}
\naam{e: Delta power series}
F = \sum_{\xi \in -\N\Delta} a^\xi c_\xi,
\end{equation}
with $c_\xi \in V,$ for $\xi \in -\N\Delta.$ 
Note that for a $\Delta$-power series the notion of neat convergence 
at a point $a_0 \in A$ 
coincides with the notion of absolute convergence in the point $a_0.$

The terminology `power series' is motivated
by the following consideration. If $\mu \in \N\gD,$ we put $\mu = \sum_{\ga \in \gD} \mu_\ga \ga,$ 
with $\mu_\ga \in \N.$ 
For $z \in \C^{\Delta},$ we write 
$$
z^\mu = \prod_{\ga \in \Delta} z_\ga^{\mu_\ga}.
$$
Finally, to the series (\refer{e: Delta power series}) we associate the power series
\begin{equation}
\naam{e: power series in z}
\sum_{\mu \in \N \Delta} z^\mu c_{-\mu}
\end{equation} 
with coefficients in $V.$ 

Let $\uz: A \to \C^\Delta$ be the map defined by
$\uz(a)_\ga = a^{-\ga}.$ Then it is obvious that the series (\refer{e: Delta power series})
converges  with sum $S$ for $a = a_0$ if and only if the power series 
(\refer{e: power series in z}) converges with sum $S$ for $z = \uz(a_0).$ 
If $r \in \,]\,0,\infty\,[^\Delta$ we write $D(0,r)$ for the polydisc in $\C^\Delta$ 
consisting of the points $z$ with $|z_\ga| < r_\ga$ for all $\ga \in \Delta.$
Note that the preimage of this set in $A$ under the map $\uz$ is given
by 
$$
A^+(\Delta, r) := \{ a \in A \mid a^{-\ga} < r_\ga, \;\; \forall \; \ga \in \Delta\}.
$$
If $R >0,$ we also agree to write $A^+(\gD, R)$ for $A^+(\gD,r)$ with $r$ defined
by $r_\ga = R$ for all $\ga \in \gD.$ 
Finally, if $a_0 \in A,$ we write $A^+(\gD,a_0) := A^+(\gD, \uz(a_0)).$ Thus,
\begin{equation}
\naam{e: Ap Delta a zero}
A^+(\Delta, a_0) := \{ a \in A \mid a^{\ga} > a_0^\ga, \;\; \forall \; \ga \in \Delta\}
= A^+ a_0.
\end{equation}
We now note that if (\refer{e: Delta power series}) converges absolutely for $a=a_0,$ 
then the power series (\refer{e: power series in z}) 
converges absolutely for $z = z(a_0),$ hence uniformly absolutely
on the closure of the polydisc $D(0,\uz(a_0)).$ It follows that the series 
(\refer{e: Delta power series}) then converges uniformly absolutely on 
the closure of $A^+(\gD, a_0).$ 

Let $a_0 \in A.$ By $\cO(A^+(\Delta, a_0), V)$ we denote the space 
of functions $f: A^+(\Delta, a_0) \to V$ that are given by an absolutely converging
series of the form (\refer{e: Delta power series}).
For such a function the associated power series (\refer{e: power series in z})
converges absolutely on the polydisc $D(0,r),$ with $r = \uz(a_0);$
let $\tilde f: D(0,r) \to V$ be the holomorphic function defined by it. 
Then obviously  
$$
f(a) = \tilde f(\uz(a)),\qquad (a \in A^+(\gD, a_0)).
$$ 
We see that the $\gD$-power series representing $f \in \cO(A^+(\gD,a_0)$ is
unique. Moreover, let $\cO(D(0,r), V)$ denote the space of holomorphic 
functions $D(0,r) \to V,$ then it follows that the map 
$$
f \mapsto \tilde f,\quad \cO(A^+(\gD, r),V) \to \cO(D(0,r),V) 
$$ 
is a linear isomorphism.

In particular it follows that every $f\in \cO(A^+(\gD,r),V)$ 
is real analytic  on $A^+(\gD, r).$ Moreover, its $\gD$-power series
converges uniformly absolutely on every set of the form
$A^+(\gD, \rho),$ where $\rho \in \,]\,0, \infty\,[^\gD,$ $\rho_\ga < r_\ga$ 
for all $\ga \in \gD.$

If $\fv$ is a real linear space, then by $S(\fv)$ we denote the symmetric 
algebra of its complexification $\fv_\iC.$ 
Via the right regular 
action we  identify $S(\fa)$ with the algebra of invariant differential 
operators on $A.$ If $f\in \cO(A^+(\gD,r),V)$ and $u \in S(\fa),$ 
then  $uf$ belongs to $\cO(A^+(\gD, r), V)$ 
again; its series may be obtained from the series of $f$ by termwise application of 
$u.$ 

We now return to the more general exponential polynomial series 
(\refer{e: intro exp pol series})
with 
coefficients in $V.$ 
Let $d \geq \deg(F).$ 
Fix a basis $\gL$ of $\fa.$ For $m \in \N\gL$ we write
$m = \sum_{\gl \in \gL} m_\gl \gl$ and $|m| = \sum_{\gl} m_\gl.$ 
For such $m$ we define the polynomial function
$X \mapsto X^m$ on $\fa$ by 
$$
X^m = \prod_{\gl \in \gL} \gl(X)^{m_\gl}.
$$
These polynomial functions with $|m|\leq d$ constitute a basis 
for $P_d(\fa).$ Accordingly, we may write:
\begin{equation}
\naam{e: q and coefficients}
q_\xi(X) = \sum_{|m| \leq d} X^m c_{\xi,m},
\end{equation}
with $c_{\xi,m} \in V.$ 

\begin{lemma}
\naam{l: neat convergence each m}
The series (\refer{e: intro exp pol series})
converges neatly 
on a set $\Omega\subset A$ if and only if for  every $m \in \N\Lambda$
with $|m| \leq d$ the series 
$$
\sum_{\xi \in \Exp(F)} a^\xi c_{\xi,m}
$$
with coefficients in $V$
converges absolutely 
for all $a \in \Omega.$
\end{lemma}
\proof
This is a straightforward consequence of the definition of neat 
convergence and the finite dimensionality of the space $P_d(\fa).$ 
\qed
  
We define a partial ordering $\preceq_\Delta$ on $\fadc$ 
by
\begin{equation}
\naam{e: partial ordering preceq gD}
\xi_1 \precD \xi_2 \iff \xi_2 - \xi_1 \in \N \Delta.
\end{equation}
Moreover, we define the relation of $\Delta$-integral equivalence on $\fadc$ by
$$
\xi_1 \simD \xi_2 \iff \xi_2 - \xi_1 \in \Z\Delta.
$$
Let $F \in \cFep(A,V)$ be as in (\refer{e: intro exp pol series}) and 
have polynomial degree at most $d.$ 
In view of condition (a) of Definition \refer{d: exp pol series},
the restriction of the relation $\simD$ to the set $\Exp(F)$ induces a
finite partition of it.  Every class $\omega$ in this
partition has a least $\precD$-upper bound $s(\omega)$ in $\fadc.$  
Let $S = S_F$ be the set of these upper bounds. 
For every $s \in S$ and every $m \in \N\gL$ with $|m| \leq d$ 
we define the $\gD$-power series
\begin{equation}
\naam{e: series for f s m} 
f_{s,m}(a) = \sum_{\mu \in  \N \Delta}  a^{- \mu} c_{s -\mu ,m},
\end{equation}
with coefficients determined by (\refer{e: q and coefficients}).

\begin{lemma}
\naam{l: neat conv of exp pol series}
Let the series (\refer{e: intro exp pol series}) be 
neatly convergent at the point $a_0 \in A.$ 
Then the series (\refer{e: intro exp pol series}) and,
for every
 $s \in S= S_F$ and $m \in \N\gL$ with $|m|\leq d,$ the series
(\refer{e: series for f s m}) is  neatly convergent
on the closure of the set $A^+(\gD, a_0).$
The functions $f_{s,m},$ defined by (\refer{e: series for f s m}),
belong to $\cO(A^+(\gD, a_0), V).$ Moreover, let $f: A^+(\gD, a_0) \to V$ be the function 
defined by the summation of (\refer{e: intro exp pol series}). Then 
\begin{equation}
\naam{e: f as sum fsm}
f(a) = \sum_{ s \in S \atop |m| \leq d } 
a^s (\log a)^m f_{s,m}(a),\qquad 
(a \in A^+(\Delta, a_0)).
\end{equation}
In particular, the function $f: A^+(\Delta,a_0) \to V$ is real analytic.
\end{lemma}
\proof
{}From the neat convergence of  (\refer{e: intro exp pol series})  at $a_0$ 
it follows by Lemma \refer{l: neat convergence each m} that for every 
$s$ and $m$ the series $\sum_{\mu \in \N\Delta} a^{s -\mu}c_{s - \mu, m}$  
converges absolutely for $a = a_0.$ This implies that 
the $\Delta$-power series (\refer{e: series for f s m}) converges 
absolutely for $a = a_0.$ Hence it converges (uniformly) absolutely 
on the closure of $A^+(\gD, a_0);$ in particular it converges neatly on that set.
It follows from this that $f_{s,m} \in \cO(A^+(\gD, a_0),V),$ for $s \in S$ and $m \in \N \gL$ with
$|m| \leq d.$ Moreover, 
\begin{equation}
\naam{e: series for a s log a m times f s m}
a^s (\log a)^m \,f_{s,m} = \sum_{\xi \in s - \N\gD} a^\xi (\log a)^m c_{\xi, m}
\end{equation} 
where the $\gD$-exponential polynomial series on the right-hand side converges
neatly on the closure of $A^+(\gD, a_0).$ The 
series (\refer{e: series for a s log a m times f s m}), for $s \in S$ and 
$m\in \N \gL$ with $|m|\leq d$ add up to the series  
(\refer{e: intro exp pol series}), which is therefore neatly convergent as well.
Moreover, (\refer{e: f as sum fsm}) follows. This in turn implies the real 
analyticity of the function $f.$ 
\qed
\begin{rem}
\naam{r: 1.5 bis}
Let $\fa_\gD: = \cap_{\ga \in \gD} \ker \ga$ and $A_\gD: = \exp (\fa_\gD).$ 
Then the functions $f_{s,m},$ defined by 
(\refer{e: series for f s m}) satisfy $f_{s,m}(a a_\gD) = f_{s,m}(a)$ 
for all $a \in A,\, a_\gD \in A_\gD.$ In particular, the function $f$ of 
(\refer{e: f as sum fsm}) generates
a finite dimensional $A_\gD$-module with respect to the right regular action.
Thus, if $\gD = \emptyset,$ then $f$ is an exponential polynomial function.
\end{rem}

\begin{lemma}
{\rm (Uniqueness of asymptotics)\ }
\naam{l: uniqueness of asymp}
Let  $a_0 \in A,$ and assume that the $\Delta$-exponential polynomial 
series  
(\refer{e: intro exp pol series})
converges neatly on $\Ap(\gD,a_0).$  
If the sum of the series is zero  for all $a \in \Ap(\gD,a_0),$ 
then $q_\xi = 0$ for all $\xi \in \fadc.$ 
\end{lemma}

\proof
Let $f: \Ap(\gD,a_0) \to V$ be defined by summation 
of the series (\refer{e: intro exp pol series}). 
Then it follows from Lemma \refer{l: neat conv of exp pol series}
that the series (\refer{e: intro exp pol series}) is 
an asymptotic expansion for $f$ in the sense
of \bib{BS}, Sect.\ 3.
Hence, if $f=0,$ then by uniqueness of asymptotics, see \bib{HCsf}, p.\ 305, Cor.\ and
\bib{BS}, Prop.\ 3.1, it follows that the series vanishes identically.
\qed

\begin{defi}
\naam{d: defi Cep A a zero}
Let $a_0 \in A.$ By $\Cep(\Ap(\gD,a_0), V)$ 
we denote the space of  functions 
$f: A^+(\Delta,a_0) \to V$ that are given by the summation of a 
(necessarily unique)
neatly converging
$\Delta$-exponential 
polynomial series  of the form (\refer{e: intro exp pol series}).

If $f \in \Cep(\Ap(\gD,a_0), V),$ then by  $\ep(f)$  we denote 
the unique series from $\cFep(A,V)$ whose summation gives $f.$ 
Moreover, the asymptotic degree of $f$  
is defined to be the number
$$
\dega(f): = \deg (\ep(f)).
$$
\end{defi}

Note that the map 
$$
\ep: \Cep(A^+(\gD, a_0),V) \to \cFep(A,V),
$$
defined above, is a linear embedding. 

Let $f \in \Cep(A^+(\gD, a_0),V).$ 
We briefly write $\Exp(f)$ for the set 
$\Exp(\ep(f));$ its elements are called the exponents 
of $f.$ 
We put $q_\xi(f,\dotvar):= q_\xi(\ep(f),\dotvar),$ for $\xi \in \fadc.$ 
Then $\xi \in \Exp(f) \iff q_\xi(f)\neq 0.$ 
 
The $\precD$-maximal elements in $\Exp(f)$ are called
the ($\Delta$-)leading exponents of $f$ (or of the expansion). The set
of these is denoted by $\ExpL(f).$ 

By the formal application of $S(\fa)$ to $\cF^\ep(A, V)$ we shall mean the linear
map 
$$
S(\fa) \otimes \cF^\ep(A, V) \to \cF^\ep(A, V)
$$ 
induced by termwise differentiation
(recall that $S(\fa)$ acts on $C^\infty(A)$ via the right regular action). The image of
an element $u\otimes F$ under this map will be denoted by $uF.$ 

\begin{lemma}
\naam{l: formal application of Ufa}
Let $a_0 \in A$ and let $f\in C^\ep(A^+(\gD, a_0),V).$ If $u \in S(\fa)$ then the function
$uf: a \mapsto R_uf(a)$ belongs to $C^\ep(A^+(\gD, a_0),V).$ Moreover, 
$$
\ep (uf) = u \, \ep(f).
$$ 
\end{lemma}
\proof
We may assume that $u \in \fa.$ 
Express $f$ as in (\refer{e: f as sum fsm}). For  
each $s, m$ the function
$u f_{s,m}$ belongs to $\cO(A^+(\gD, a_0), V);$ its expansion is obtained 
from $\ep(f_{s,m})$ by termwise application of $u,$ hence by the formal 
application of $u.$ 
\qed

We shall also need a second type of formal application.
Suppose that complete locally convex spaces $U$ and $W$ are given, 
and a continous bilinear map $U \times V \to W,$ denoted by $(u,v) \mapsto uv.$ 
By the formal application of $\cF^\ep(A, U)$ 
to $\cF^\ep(A, V)$ we mean the linear map  
$$ 
\cF^\ep(A, U) \otimes \cF^\ep(A, V) \to \cF^\ep(A, W),
$$ 
given by 
\begin{equation}
\naam{e: defi second formal appl} 
\sum_{\xi\in \fadc}  a^\xi p_\xi(\log a)  \otimes 
\sum_{\eta\in \faqd} a^\eta q_\eta(\log a)
\mapsto 
\sum_{\nu \in \fadc} 
a^\nu \sum_{\xi + \eta = \nu} p_\xi(\log a) q_\eta(\log a).
\end{equation}
This map is indeed well defined. To see this, let $F$ denote the first
series and $G$ the second. Then for every $\nu \in \faqdc,$ the collection
$S_\nu$ of $(\xi, \eta) \in \Exp(F)\times \Exp(G)$ with $\xi + \eta = \nu$ is finite.
Hence the $W$-valued polynomial function
$$
r_\nu: X\mapsto  \sum_{(\xi, \eta)\in S_\nu} p_\xi(X)q_\eta(X)
$$ 
has degree at most $\deg(F) + \deg(G).$ 
Moreover, let $X_1, X_2 \subset \fadc$ be finite subsets such that 
$\Exp(F) \subset X_1 - \N\gD$ 
and $\Exp(G) \subset X_2 - \N\gD$ and put $X = X_1 + X_2.$ Then for 
$\nu \in \fadc \setminus [X - \N\gD]$ the collection $S_\nu$ is empty, 
hence $r_\nu = 0.$ 
Therefore, the formal series on the right-hand side
of (\refer{e: defi second formal appl}) satisfies 
the conditions of Definition \refer{d: exp pol series}.

The image of an element $F \otimes G$ under the  map 
(\refer{e: defi second formal appl}) 
is denoted
by $FG.$ 
Again we have a lemma relating the formal
application with neat convergence.

\begin{lemma}
\naam{l: formal application of hom valued series}
Let $U \times V \to W, (u,v) \mapsto uv$ be a continuous bilinear map 
of complete locally convex spaces.
Let $a_0 \in A$ and let $f\in C^\ep(A(\gD, a_0),U)$ and
$g \in C^\ep(A(\gD, a_0), V).$ 
Then the function $fg: a \mapsto f(a)g(a)$ belongs
to $C^\ep(A(\gD, a_0),W).$ 
Moreover, its $\gD$-exponential polynomial expansion is given by 
$$ 
\ep(fg) = \ep(f)\, \ep(g).
$$ 
\end{lemma}

\proof
This follows by a straightforward application of 
Lemma \refer{l: neat conv of exp pol series}.
\qed

\section{Basic notation, spherical functions}
In this section we study spherical functions  that are defined 
on a certain open dense 
subset $\spXp$ of the symmetric space $\spX,$ 
and are (radially) given by 
exponential polynomial series. This class of functions 
will play an important role in the paper. Later
we will see that $\DGH$-finite spherical funtions
belong to this class.

Throughout this paper, we assume that $\spX$ is a reductive symmetric space of 
Harish-Chandra's class, i.e., $\spX= G/H$ with $G$ a real reductive group of Harish-Chandra's
class and $H$ an open subgroup of $G^\gs,$ 
the group of fixed points for an involution $\gs$ of $G.$  
There exists a Cartan involution $\Cartan$ of $G,$ 
commuting with $\gs.$ The associated fixed point group $K$ is a $\gs$-stable 
maximal compact subgroup. 

We adopt the usual convention to denote Lie groups by Roman capitals
and their Lie algebras  by the corresponding Gothic lower cases.
The infinitesimal involutions $\Cartan$ and  $\gs$  of $\fg$ commute;  
let 
\begin{equation}
\naam{e: Cartan decos}
\fg = \fk \oplus \fp = \fh \oplus \fq 
\end{equation} 
be the associated decompositions into $+1$ and $-1$ eigenspaces
for $\Cartan$ and $\gs$, respectively.  
We equip $\fg$ with a positive definite inner product $\inp{\dotvar}{\dotvar}$ 
that is invariant under the compact group of automorphisms generated by 
$\Ad(K),$ $e^{i \ad(\fp)},$ $\Cartan$ and $\gs.$ Then the decompositions 
(\refer{e: Cartan decos}) are orthogonal.

Let $\faq$ be a maximal abelian subspace of $\fp \cap \fq.$
We equip $\faq$ with the restricted inner product $\inp{\dotvar}{\dotvar}$ 
and its dual $\faqd$ with the dual inner product. The latter is 
extended to a complex bilinear form, also denoted  $\inp{\dotvar}{\dotvar},$ 
on the  complexified dual $\faqdc.$ 

The exponential 
map is a diffeomorphism from $\faq$ onto a vectorial subgroup $\Aq$ of $G.$ 
We recall that $G = K\Aq H.$ 
Let $\Sigma$ be the restricted root system of $\faq$ in $\fg;$ 
we recall that the associated Weyl group $W$ is naturally isomorphic
to $\NKaq/ \ZKaq,$ the normalizer modulo the centralizer of $\faq$ in $K.$ 
Let $\faq^\reg$ denote the associated set of regular elements in $\faq,$ 
i.e., the complement of the union of the root hyperplanes $\ker \ga,$ 
as $\ga \in \Sigma.$ We put $\Aq^\reg\coleq\exp(\faq^\reg)$ and define
the dense subset $\spXp$ of $\spX$ by
$$\spXp=K \Aq^\reg H.$$

If $Q$ is a parabolic subgroup of $G,$ we denote its Langlands decomposition
by $Q = M_Q A_Q N_Q.$ 
By a $\gs$-parabolic subgroup of $G$ we mean a parabolic subgroup
that is invariant under the composition $\gs \after \Cartan.$ 
It follows from \bib{Bprser1}, Lemmas 2.5 and 2.6, that 
the collection $\allparabs$ of $\gs$-parabolic 
subgroups of $G$ containing $\Aq$ is finite. 

If $Q$ is a $\gs$-parabolic subgroup then 
the Lie algebra $\faQ$ of its split component 
is $\gs$-stable, hence decomposes as
$\faQ = \fa_{Q\ih} \oplus \faQq,$ the vector sum of the associated 
$+1$ and $-1$ eigenspaces
of $\gs|_{\faQ},$ respectively. 
We write 
$\AQq := \exp \faQq$ and $\MQgs:= M_Q (A_Q\cap H);$ 
the decomposition 
$$
Q = \MQgs \AQq N_Q 
$$ 
is called the $\gs$-Langlands decomposition of $Q.$ 
If  $Q \in \allparabs,$ then $\MoneQ = Q \cap \Cartan(Q)$ contains $\Aq.$ Hence 
$\faQq$ is contained in $\fp \cap \fq$ and centralizes $\faq;$ it follows
that $\faQq \subset \faq.$ By $\Sigma(Q)$ we denote the set of roots of $\Sigma$ occurring 
in $\fn_Q.$ Obviously, 
$$ 
\fn_Q = \oplus_{\ga \in \Sigma(Q)} \;\fg_\ga.
$$

Let $\minparabs$ denote the collection of elements of $\allparabs$ that
are minimal with respect to inclusion. An element $P \in \allparabs$ 
belongs to $\minparabs$ if and only if $\faPq = \faq,$ see
\bib{Bprser1}, Cor.\ 2.7. 
This implies that the 
associated groups $M_P$ and $A_P$ are independent of $P \in \minparabs.$ 
We denote them by $M$ and $A,$ respectively.
 From the maximality of
$\faq$ in $\fp\cap \fq$ it follows that $\fm \cap \fp \subset \fh.$ 
Thus, if $\KM:= K \cap M$ and $\HM:= H \cap M,$ 
then the inclusion map $\KM \to M$ 
induces a diffeomorphism
\begin{equation}
\naam{e: isomorphism for M mod HM}
\KM/\KM\cap H\simeqarrow M/\HM. 
\end{equation} 
In particular, the symmetric space $M/\HM$ is compact.

According to \bib{Bprser1}, Lemma 2.8, the map $P \mapsto \Sigma(P)$ 
induces a bijective map from $\minparabs$ onto the collection of positive
systems for $\Sigma.$ If  $\Phi$ is a positive system for $\Sigma,$ 
then the associated element $P \in \minparabs$ is given by the following 
characterization of its Lie algebra:
${\rm Lie}\,(P) = \fm + \fa + \sum_{\ga \in \Phi}\fg_\ga.$ From this we see that 
$\NKaq$ acts on $\minparabs$ by conjugation; moreover, the action
commutes with the map $P \mapsto \Sigma(P).$ Accordingly, the action
factors to a free transitive action of $W$ on $\minparabs,$ 
see also \bib{Bprser1}, Lemma 2.8.

If $P \in \minparabs,$ 
then the collection of simple roots for the 
positive system $\Sigma(P)$ is denoted by $\DP;$ 
the associated positive chamber
in $\faq$ is denoted by $\faqp(P)$  and 
the corresponding chamber in $\Aq$  by $\Aqp(P).$ 
Thus, we see that  $\Aq^\reg$ 
is the disjoint union of the chambers
$\Aqp(P),$ as $P \in \minparabs.$ 

More generally, if $Q \in \allparabs,$ we write
\begin{equation}
\naam{e: defi faQqp}
\faQqp: = \{X \in \faQq \mid \ga( X) > 0 \text{for} \ga \in \Sigma(Q) \}.
\end{equation} 
It follows from \bib{Bprser1}, Lemmas 2.5 and 2.6, that $\faQqp \not= \emptyset.$ 
Moreover, if $X \in \faQqp,$ then the parabolic subgroup $Q$ is determined
by the following characterization of its Lie algebra
\begin{equation}
\naam{e: char Lie Q}
{\rm Lie}\, (Q) = \fm \oplus \fa \oplus \bigoplus_{\ga \in \Sigma\atop \ga(X) \geq 0} \fg_\ga.
\end{equation}
Conversely, if $X$ is any element of $\faq,$ then (\refer{e: char Lie Q}) 
defines the Lie
algebra of a group $Q$ from $\allparabs;$ moreover, $X \in \faQqp.$ 
{}From this we readily see that conjugation
induces an action of $\NKaq$ on $\allparabs,$ which factors to an action of 
$W.$ 

By a straightforward calculation involving root spaces, 
it follows that the multiplication
map $K \times \Aq^\reg \to \spX$ induces a diffeomorphism 
$$ 
K \times_{\NKaq \cap H} \Aq^\reg \simeqarrow \spXp. 
$$
In particular, it follows that $\spXp$ is an open dense subset of $\spX$.
Let $\WKH$ denote the canonical image of $\NKaq \cap H$ in $W$  
and let $\cW$ be a complete set of representatives for $W/\WKH$ in $\NKaq.$ 
If $P \in \minparabs,$ then it follows that 
\begin{equation}
\naam{e: space X plus as union}
\spXp = \cup_{w \in \cW} \;\;K\Aqp(P)wH \quad \text{(disjoint union).}
\end{equation}
Moreover, for each $w \in \cW$ the multiplication map $(k,a) \mapsto kawH$ 
induces a diffeomorphism
\begin{equation}
\naam{e: diffeo onto KAqpPwH}
K \times_{\KM \cap wHw^{-1}} \Aqp(P)
\simeqarrow K\Aqp(P)wH.
\end{equation}
Here we have written $\KM = K \cap M;$ in (\refer{e: diffeo onto KAqpPwH})
the set on the right 
 is an open 
subset of $\spX.$

Let $(\tau, \Vtau)$ be a smooth representation
of  $K$ in a complete locally convex space.
For later applications it will be crucial that we allow
$\tau$ to be infinite dimensional (see
the proof of Theorem \refer{t: behavior along the walls for families new}).

By $\sphXp$ we denote the space of smooth functions 
$f: \spXp \to \Vtau$ that are $\tau$-spherical, i.e.,
\begin{equation}
\naam{e: spherical transformation rule}
f(kx) = \tau(k)f(x),
\end{equation}
for $x \in \spXp, \, k\in K.$ 
The space $\Ci(\spX\col \tau)$ of smooth $\tau$-spherical 
 functions
on $\spX$ will be identified with the subspace 
of functions in $\sphXp$ that extend smoothly to all of $\spX.$

In the following we assume that $P \in \minparabs$
is  fixed. 
If $w \in \NKaq,$ 
then by $C^\infty_{P,w}(\spXp\col \tau)$ or
$C^\infty_w(\spXp\col \tau)$ 
we denote the space of functions $f \in \Ci(\spXp\col \tau)$ 
with support contained in $K\Aqp(P)wH.$ 
{}From (\refer{e: space X plus as union}) we see that
$$
\Ci(\spXp\col \tau) = \oplus_{w\in \cW} \;\;C^\infty_w(\spXp\col \tau).
$$
Let  $w \in \NKaq$ be fixed for the moment.
For $f \in \Ci(\spXp \col \tau)$ we define the
function $\TdownPw f \in \Ci(\Aqp(P), \Vtau^{\KMwH})$ by 
$$
\TdownPw f (a) = f(awH).
$$
Since (\refer{e: diffeo onto KAqpPwH}) is a diffeomorphism, 
the restriction of $\TdownPw$ to 
$C^\infty_w(\spXp\col \tau)$ 
is an isomorphism of complete locally convex spaces 
onto the space $\Ci(\Aqp(P), \Vtau^{\KMwH}).$ 
Taking the direct sum of the maps
$\TdownPw,$ as $w \in \cW,$ 
we therefore obtain an isomorphism of complete locally convex spaces
\begin{equation}
\naam{e: the iso T down P cW}
\TdownPcW: 
\;\;\Ci(\spXp \col \tau)
\simeqarrow
\oplus_{w \in \cW}\;\; \Ci(\Aqp(P), \Vtau^{\KMwH}).
\end{equation}

\begin{defi}
\naam{d: Cep spXp tau}
We denote by $\ExppolXptau$ the space of functions
$f \in \CiXptau$ such that for every $w \in \cW$ the
function $\TdownPw f$ belongs to $\Cep(\Aqp(P),\VtauKwH),$ 
where the latter space is defined as in Definition 
\refer{d: defi Cep A a zero},
with $\fa,$ $a_0$ and $\gD$ replaced
by $\faq,$ $e$ and $\gD(P),$ respectively. 

If $f \in \ExppolXptau,$ we define its asymptotic degree to be the number
$$ 
\dega (f): = \max_{w \in \cW}\; \deg (\TdownPw f).
$$
\end{defi}
\noindent
It follows from the above definition that 
restriction of $\TdownPcW$ induces a linear isomorphism
\begin{equation}
\naam{e: isomorphism of exppol}
\ExppolXptau \simeq 
\oplus_{w \in \cW}\;\; \Cep(\Aqp(P), \Vtau^{\KMwH}).
\end{equation}
Using conjugations by elements of $\NKaq$ 
it is readily seen 
that the space $\ExppolXptau$
and the map $\dega: \ExppolXptau\to \N$  
are independent of the
particular choices of $P$ and $\cW.$ In particular, if $P \in \minparabs$ 
and $w \in \NKaq,$ then $\TdownPw f \in \Cep(\Aqp(P), \VtauKwH)$ 
and
$\deg (\TdownPw f) \leq \dega (f).$ 
We put 
$$
\Exp(P,w\asmid f) := \Exp(\TdownPw f),\quad \text{and}
\quad\ExpL(P,w\asmid f) := \ExpL(\TdownPw f).
$$ 
Moreover, for all $\xi \in \faqdc$ we define 
$\uq_\xi(P,w \asmid f) = q_\xi(\TdownPw f).$ 
Then, for every $a \in \Aqp(P),$ 
\begin{equation}
\naam{e: exp pol expression f a w}
f(aw) = \sum_{\xi \in \Exp(P,w\asmid f)} a^\xi\,\uq_\xi(P,w\asmid f, \log a),
\end{equation}
where the $\gD(P)$-exponential polynomial series on the right-hand side
neatly converges on $\Aqp(P).$ 

For $w \in \NKaq,$ we will use the notation
\begin{equation}
\naam{e: defi spXzerow}
\spXzerow := M /M \cap \wH;
\end{equation}
moreover,
we 
put $\tauM:= \tau_{\KM}$ and 
write  
$
C^\infty(\spXzerow\col \tauM)
$ 
for the space of 
$\tauM$-spherical 
$C^\infty$ functions from
$\spXzerow$ to $V_\tau,$ i.e., the space of functions $f \in \Ci(\spXzerow, \Vtau)$ satisfying the 
rule (\refer{e: spherical transformation rule})  for $k \in \KM$ and $x \in \spXzerow.$ 
 From (\refer{e: isomorphism for M mod HM}) with $wHw^{-1}$ in place of $H$ we see that 
the inclusion $\KM \to M$ induces a diffemorphism from $\KM/\KM \cap \wH$ onto $\spXzerow.$ 
Hence evaluation at the point $e(M\cap \wH)$ induces a linear isomorphism
from $C^\infty(\spXzerow \col \tauM)$ onto $\VtauMKwH.$ Thus, if 
$f \in \ExppolXptau,$ 
then for every $\xi \in \faqc$ there exists a unique 
$C^\infty(\spXzerow \col\tauM)$-valued 
polynomial function $q_\xi(P,w\asmid f)$ on $\faq$ such that 
$$
q_\xi(P,w \asmid f, X, e)= \uq_\xi(P,w\asmid f)(X)
\qquad (X \in \faq).
$$ 
Using sphericality of the function $f$ 
we obtain from (\refer{e: exp pol expression f a w})
that 
\begin{equation}
\naam{e: expansion f on PqPw with m}
f(maw) = \sum_{\xi \in \Exp(P,w\asmid f)} a^\xi q_\xi(P,w\asmid f , \log a, m),
\end{equation}
for $m \in M,\; a \in \Aqp(P).$ The series on the right-hand side
is a $\gD(P)$-exponential polynomial series in the variable
$a,$ with coefficients
in $\Ci(\spXzerow \col \tauM),$ relative to the variable $m.$ As 
such it converges neatly on $\Aqp(P).$ 

We shall now discuss a lemma whose main purpose is to enable us to 
reduce on the set of exponents in certain proofs, in order to simplify 
the exposition.

\begin{lemma}
\naam{l: splitting lemma}
Let $P \in \minparabs$ and let $\cW \subset \NKaq$ be a 
complete set of representatives
of $W/\WKH.$ Assume that
$f \in \ExppolXptau.$ 

There exists a finite set $S \subset \faqdc$ 
of mutually $\gD(P)$-integrally inequivalent elements such that 
$\Exp(P,v\asmid f) \subset S - \N\gD(P)$ 
for every $v \in \cW.$ 

If $S$ is a set as above, then there exist unique 
functions $f_s \in \ExppolXptau,$ for $s \in S,$ 
such that 
$$
f = \sum_{s \in S} f_s,
$$ 
and such that $\Exp(P,v\asmid f_s) \subset s - \N \gD(P),$ 
for every $v \in \cW.$ 
\end{lemma}
\proof
There exists a finite set $X \subset \faqdc$ such that 
$\Exp(P,v\asmid f) \subset
X - \N \gD(P)$ for all $v \in \cW.$ Obviously there exists 
a finite set $S$ as
required, such that $X - \N \gD(P) \subset S - \N \gD(P).$ 

If $S$ is such as mentioned, then for 
$s \in S$ and $v \in \cW$ we define the function 
$f_{s,v}: \Aqp(P) \to \VtauKvH$ 
by 
$$
f_{s,v}(a) = \sum_{\nu \in \N\gD(P)} a^{s -\nu} q_{s -\nu} (P,v\mid f,\log a,e);
$$ 
here the exponential polynomial series is neatly convergent,
hence $f_{s,v} $ belongs to the space $\Cep(\Aqp(P), \VtauKvH),$ for every $v \in \cW.$ 
By the isomorphism (\refer{e: isomorphism of exppol}) there exists a unique 
function $f_s \in \Cep(\spXp\col \tau)$ 
such that $f_s(av) = f_{s,v}(a)$ for $v \in \cW,\, a \in \Aqp(P).$ 
By the hypothesis on $S$ the sets $s - \N\gD(P),$ for $s \in S$, are disjoint. 
Hence $f = \sum_{s \in S} f_s$  on $\Aqp(P)v,$ for every $v \in \cW.$ 
By (\refer{e: space X plus as union}) and
sphericality this equality holds on all of $\spXp.$ 
\qed

\section{Asymptotic behavior along walls}
\naam{s: asymp walls}
In this section we study the asymptotic behavior along  walls
of functions from $\ExppolXptau;$ here $\tau$ 
is a smooth representation
in a complete locally convex  space $\Vtau.$ 

Let $P \in \minparabs$  and let 
$Q $ be a $\gs$-parabolic subgroup  with Langlands decomposition
$Q =  M_Q A_Q N_Q,$ containing $P.$ 
In addition to the notation introduced in the beginning of the previous
section, the following notation will also  be convenient. 

We agree to write
$K_Q:= K \cap M_Q$ and $H_Q: = H \cap M_Q.$ Moreover, 
$W_Q$ denotes the centralizer of $\faQq$ in $W.$ Then $W_Q \simeq N_{K_Q}(\faq)/Z_{K_Q}(\faq).$ 
On the other hand $W_Q$ is also the subgroup of $W$ generated by the reflections in the roots
from the set 
$$
\DsubQP  := \{ \ga \in \DP \mid \ga|_{\faQq} 
 = 0 \}.
$$
We note that $\Sigma(Q) = \Sigma(P)\setminus \N\Delta_Q(P).$ Moreover, let
$\SrQ$ denote the collection of $\faQq$-weights in $\fn_Q.$ 
Then  
$$
\SrQ = \{ \ga|_{\faQq} \mid \ga \in \Sigma(Q)\}.
$$ 
Let $\DrQ$ be the collection of weights from the set $\SrQ$ that cannot be written as the 
sum of two weights from that set; then one readily verifies that 
$\DrQ$ equals the set of restrictions of elements from $\Delta(P)\setminus \Delta_Q(P)$ 
to $\faQq.$ In
particular, the elements of $\DrQ$ are linearly independent.

Given $a_0 \in \AQq$ we shall briefly write $\AQqp(a_0)$ for the set
$\AQqp(\DrQ, a_0)$ defined as 
in (\refer{e: Ap Delta a zero}) with $\faQq$ and $\DrQ$ in place of $\fa$ and $\Delta,$  
respectively. 
Similarly, if $\rho \in \,]\,0,\infty\,[^{\DrQ},$ we briefly write 
$$
\AQqp(\rho):= \AQqp(\DrQ, \rho) = 
\{a \in \AQq \mid a^{-\ga} < \rho_\ga,\;\;\forall  \ga \in \DrQ \}.
$$
If $R>0,$ we write $\AQqp(R)$ for $\AQqp(\rho),$ 
where $\rho$ is defined by $\rho_\ga = R$ 
for every $\ga \in \DrQ.$  Note that
$\AQqp(1)$ equals the positive chamber $\AQqp: = \exp(\faQqp),$ 
see (\refer{e: defi faQqp}).
 
If $v \in \NKaq,$ we define
\begin{equation}
\naam{e: defi XoneQv}
\spXoneQv := \MoneQ / \MoneQ \cap vHv^{-1}.
\end{equation}
This is a symmetric space for the involution $\gs^v$ of $\MoneQ$ defined
by $\gs^v(m) = v\gs(v^{-1} m v)v^{-1}.$ Note that this involution commutes with 
the Cartan involution $\theta|_{\MoneQ}.$  
Note also that $\faq$ is a maximal
abelian subspace of $ \Ad(v)(\fp \cap \fq) = \fp \cap \Ad(v)\fq.$ Hence it
is the analogue of $\faq$ for the triple $(\MoneQ, K_Q, \MoneQ \cap v H v^{-1}).$ 
The corresponding group $\Aq$ may naturally be identified with
a subspace of $\spXoneQv.$ 

The image of $M_Q$ in $\spXoneQv$ may be identified
with  
$$
\XQv := M_Q/M_Q \cap \vH,
$$
the symmetric space for the involution $\gs^v|_{M_Q}.$ 
It follows from the characterization of $\allparabs$ expressed by  
(\refer{e: char Lie Q}) that 
\begin{equation}
\naam{e: allparabsv}
\allparabs = \allparabsv
\end{equation}
Hence $Q$ 
is a $\gs^v$-parabolic subgroup as well. Hence $\faQ \cap \Ad(v) \fq = 
\faQ \cap \faq = \faQq,$ and we deduce that the inclusion $ \AQq \to \AQ$ 
induces a diffeomorphism $\AQq \simeq \AQ/\AQ\cap vHv^{-1}.$ From this 
we conclude that the multiplication map $\MQ \times \AQq \to \MoneQ$ 
induces the decomposition
\begin{equation}
\naam{e: deco spXoneQv}
\spXoneQv \simeq \spXQv \times \AQq.
\end{equation}
\begin{rem}
\naam{r: extreme cases subspaces}
In particular,
the 
above definitions  cover the two extreme cases  that $Q$ is minimal
and that it equals $G.$ 

In the case that $Q \in \minparabs,$ we have $Q = M A N_Q,$ and $\spXQv$ equals 
the space $\spXzerov$ defined in (\refer{e: defi spXzerow}). 
Moreover, $\spXoneQv \simeq \spXzerov \times \Aq.$ 

In the other extreme case we have $\spX_{1G,v} = G/vHv^{-1}.$ This symmetric space will
also be denoted by $\spXv.$ 
Note that right multiplication by $v$ induces an isomorphism of $\spX_v$ onto $\spX.$ 
Note also that $M_G$ equals ${}^\circ G,$ 
the intersection of $\ker \chi,$ as $\chi$ ranges over the positive characters
of $G.$ Hence $\spX_{G,v} = {}^\circ G / {}^\circ G \cap vH v^{-1}.$ 
Finally, $\spX_v \simeq \spX_{G,v} \times A_{G\iq},$ 
where $A_{G\iq}$ is the image under $\exp$ of the space 
$\fa_{G\iq},$ which in turn is the intersection
of the root hyperplanes $\ker \ga$ as $\ga \in \Sigma.$ 
\end{rem} 

Let
$\bfn_Q: = \Cartan \fn_Q$ be equipped with the restriction
of the inner product $\inp{\dotvar}{\dotvar}$ from $\fg.$
If $Q \neq G$ we define the function
$R_{Q,v}: \MoneQ \to ]0,\infty[$ by 
$$%\begin{equation}
%\naam{e: defi R Q v}
R_{Q,v}(m) = \|\Ad(m \gs^v(m)^{-1})|_{\bfn_Q}\|_{\rm op}^{1/2},
$$%\end{equation}
where $\|\dotvar\|_{\rm op}$ denotes the operator norm. We 
define $R_{G,v}$ to be the constant function $1.$

The
function $R_{Q,v}$  
is right $\MoneQ \cap v H v^{-1}$-invariant.  It  
may therefore also be viewed as a function on $\spXoneQv.$  
We shall describe the function $R_{Q,v}$ in more detail below.

The orthocomplement of $\faQq$ in $\faq$ is denoted by $\staQq.$
Note that
\begin{equation}
\naam{e: char staQq}
\staQq = \fm_Q \cap \faq;
\end{equation}
hence $\staQq$ is the analogue of
$\faq$ for the triple $(M_Q, K_Q, H_Q).$ 
We recall from the text following (\refer{e: defi XoneQv}) 
that $\faq$ is maximal abelian in $\fp \cap \Ad(v)\fq$ hence is the analogue
of $\faq$ for the triple $(G, K, vHv^{-1}).$ Accordingly, $\staQq$ is also 
the analogue of $\faq$ for the triple $(M_Q, K_Q , M_Q \cap \vH).$ 

In  view of 
(\refer{e: allparabsv}), the  group 
$\starP = P \cap M_Q$ is readily seen to be 
a minimal $\gs^v$-parabolic subgroup for $M_Q;$ 
the associated positive chamber 
in $\stAQq = \exp (\stfaQq)$ is denoted by $\stAQqp(\starP).$ 

Let $\cW_{Q,v}$ be an analogue for $\XQv$ of $\cW,$ that is, $\cW_{Q,v}$ 
is a complete set of representatives in $\NKQaq$ for the quotient
$W_Q/W_{K_Q \cap \vH}.$  Let $\XQvp$ be the analogue for $\XQv$ 
of the open dense subset
$\spXp$ of $\spX.$ According to (\refer{e: space X plus as union}) this set may
be expressed as the following disjoint union of open subsets of $\XQv$ 
\begin{equation}
\naam{e: deco XQvp}
\XQvp: = \displaycup_{u \in \cWQv} K_Q \stAQqp(\starP)\, u\, (M_Q \cap \vH) 
\qquad \text{(disjoint union).}
\end{equation}
Let
$\spXoneQvp$ be the analogue of $\spXp$ for $\spXoneQv;$ then from 
(\refer{e: deco spXoneQv}) we see that 
$
\spXoneQvp \simeq \spXQvp \times \AQq.
$
In terms of this decomposition and (\refer{e: deco XQvp})
 the function $R_{Q,v}$ may be expressed as 
follows.
\begin{lemma}
\naam{l: properties RQv}
The function $R_{Q,v}: \MoneQ \to ]0,\infty [$ is continuous, and  
right $\MoneQ \cap vHv^{-1}$- and 
left $K_Q$-invariant.
Moreover, if $Q \not=G$ and if $a \in \Aq$ and $u \in \NKQaq,$ then
\begin{equation}
\naam{e: value RQv}
R_{Q,v}(au) = \max_{\ga\in \Sigma(Q)} a^{-\ga}.
\end{equation}
Finally, $R_{Q,v} \geq 1$ on $\XQv.$ 
\end{lemma}
\proof
Since $R_{G,v}$ is the constant function $1,$ we may as well assume that $Q \not= G.$ 
Continuity of the function $R_{Q,v}$ is obvious from its definition.
The group $K_Q$ is $\gs^v$ invariant and acts unitarily on $\bfn_Q;$ 
hence the left $K_Q$-invariance is obvious from the definition.
If $a \in \Aq,$ then $a\gs^v(a)^{-1} =a^2.$ 
Hence the operator norm
of $\Ad(a\gs^v(a)^{-1})$ on $\bfn_Q$ equals the maximal value of $a^{-2\ga}$ as 
$\ga \in \Sigma(Q).$  This implies (\refer{e: value RQv}) for $u = 1.$ 

The element $u \in \NKQaq$ belongs to $M_Q,$ hence $\Ad(u)$ 
normalizes
$\fn_Q.$ Therefore,  $\Ad(u)$ leaves the collection 
$\Sigma(Q)$ of $\faq$-roots in $\fn_Q$ invariant.
Put $a' = u^{-1} a u.$ Then 
$
R_{Q,v}(a u) = R_{Q,v}(a') = \max_{\ga \in \Sigma(Q)} (a')^{-\ga}.$ 
Since $\Ad(u)$ leaves $\Sigma(Q)$ invariant, (\refer{e: value RQv}) follows.

If  $\ga \in \Sigma,$ let $h_\ga$ be the element of $\faq$ determined by
$\ga(X) = \inp{h_\ga}{X},$ for $X \in \faq.$ Then the closure of ${\staQqp(\starP)}$ is
contained in the closed convex cone 
generated by the elements $h_\gb,$ for $\gb \in \Delta_Q(P).$ 
If $\ga \in \DP\setminus \Delta_Q(P),$ then 
$\ga(h_\gb) = \inp{\ga}{\gb} \leq 0,$ 
for $\gb \in \Delta_Q(P);$ hence $\ga \leq 0$ on $\staQqp(\starP).$ 
But $\DP\setminus \Delta_Q(P)\subset \Sigma(Q),$ hence 
it follows that $R_{Q,v} \geq 1$ on $\stAQqp(\starP)u,$ 
for every $u \in \cWQv.$ 
The final assertion follows from combining this observation
with (\refer{e: deco XQvp}),  the left $K_Q$-invariance of $R_{Q,v}$ 
and density of $\spXQvp$ in $\spXQv.$ 
\qed

If $1\leq R \leq \infty$ we define
\begin{equation}
\naam{e: defi spXQvb R}
\XQvb{R}:= \{m \in \XQv \mid R_{Q,v}(m) < R\}.
\end{equation}
Note that $\XQvb{1} = \emptyset$ and
$\XQvb{\infty} = \spXQv;$ moreover, 
$R_1 < R_2 \implies \XQvb{R_1} \subset \XQvb{R_2}.$
Finally, the union of the sets $\XQvb{R}$ as $1 \leq R < \infty$ equals $\XQv.$ 

In accordance
with (\refer{e: defi spXQvb R}) we define
$\XQvpb{R}: = \spXQvp \cap \XQvb{R},$  for
$1 \leq R \leq \infty.$ 
Moreover, we also put
$$
\stAQqp(\starP)\lbr{R}: = 
\{ a \in \stAQqp(\starP)\mid a^{-\ga} < R,\;\;\forall \ga \in \Sigma(Q)\}.
$$
Note that, if $\ga \in \gS(P)\setminus\gS(Q),$ then $a^{-\ga} < 1 \leq R$ for 
all $a \in \stAQqp(\starP).$ Hence 
$$
\stAQqp(\starP)\lbr{R} = \stAQqp(\starP) \cap \Aqp(P, R).
$$ 
It follows from (\refer{e: deco XQvp}) 
and Lemma \refer{l: properties RQv} that
\begin{equation}
\naam{e: deco XQvp added R}
\quad\XQvpb{R} = \displaycup_{u \in \cWQv} K_Q \stAQqp(\starP)\lbr{R}\, u\, (M_Q \cap \vH) 
\quad \text{(disjoint union).}
\end{equation}

The function $R_{Q,v}$ plays a role in the description  
of the asymptotic behavior of a function $f \in \Cep(\spXp\col \tau)$
along `the wall' $\AQqp v.$  This behaviour is described
in terms of an expansion of $f(mav)$ in the variable $a \in \AQqp,$  
for $m \in \spXQvp.$ Thus, it is of interest to know when $mavH$ belongs
to $ \spXp,$ the domain of $f.$ 
\vspace{-4mm}
\begin{lemma}\naam{l: inclusions for asymp}
\hbox{\hspace{1mm}}\hfill\break
\vspace{-4mm}
\begin{enumerate}
\itema
If $b \in \stAQqp(\starP)$ and $a \in \AQqp(R_{Q,v}(b)^{-1})$ 
then $ba\in \Aqp(P).$ 
\itemb
Let $m \in \spXQvp.$ Then $mavH \in \spXp$ for all $a \in \AQqp(R_{Q,v}(m)^{-1}).$ 
\itemc
Let $R \geq 1.$ Then $\XQvpb{R} \AQqp(R^{-1}) v H \subset \spXp.$ 
\end{enumerate}
\end{lemma}

\proof
Let $b$ and $a$ fulfill the hypotheses of (a).
If $\ga \in \gD_Q(P),$ then $(ba)^{-\ga} = b^{-\ga} <1.$ 
On the other hand we have,
for $\ga \in \gD(P)\setminus \gD_Q(P),$ that $\ga \in \Sigma(Q),$ hence 
$(ba)^{-\ga} \leq R_{Q,v}(b) a^{-\ga} < 1,$ by Lemma \refer{l: properties RQv}.
Hence $ba \in \Aqp(P),$ and (a) is proved.

Let $m$ be as in (b), and let $a \in \AQqp(R_{Q,v}(m)^{-1}).$ 
In view of (\refer{e: deco XQvp}) 
we may write $m = kbuh$ with $k \in K_Q,$ $b \in \stAQqp(\starP),$ 
$u \in \cW_{Q,v}$ and $h \in M_Q \cap vHv^{-1}.$ Now 
$mavH = kbuhavH= kbauvH.$ Thus, it suffices to show that 
$ba \in \Aqp(P).$ This follows from (a) and the observation that 
$R_{Q,v}(b) = R_{Q,v}(m),$ by Lemma \refer{l: properties RQv}.

Finally, (c) is a straightforward consequence of (b).
\qed

If
$Q \in \allparabs$
we put $\tauQ:= \tau|_{K_Q}.$ Then, for 
$v \in \NKaq,$ 
the space $\Cep(\spXQvp \col \tauQ)$ is defined as 
above (\refer{e: isomorphism of exppol}) with 
$\spXQv$ and 
$\tauQ$ 
in place of $\spX$ and $\tau,$ respectively.
  
\begin{thm}
\naam{t: expansion along the walls}
Let $f \in \ExppolXptau.$ 
Let $Q \in \allparabs$ 
and $v \in \NKaq.$ 
\begin{enumerate}
\itema
There exist a constant $k\in\N$, a
finite set $Y \subset \faQqdc,$ and
for each $\eta \in Y - \N \DrQ$ a
$C(\XQvp, V_\tau)$-valued 
polynomial function $q_\eta = q_{\eta}(Q,v\asmid f)$ on $\faQq$ of degree at most
$k,$  such that for every  $m \in \XQvp$ 
\begin{equation}
\naam{e: expansion f along Q v}
f(mav) = \sum_{\eta \in Y - \N \DrQ} a^\eta q_\eta(\log a, m),
\qquad(a \in \AQqp(R_{Q,v}(m)^{-1})),
\end{equation}
where the $\DrQ$-exponential polynomial series with coefficients in $\Vtau$ 
converges neatly  on the indicated subset of $\AQq.$
\itemb
The set 
$\Exp(Q,v\asmid f):=\{\eta \in Y - \N \DrQ \mid q_\eta \neq 0 \}$ is uniquely determined.
Moreover, the functions $q_\eta,$ where $\eta \in Y - \N \DrQ,$ are unique and
belong to $P_d(\faQq) \otimes \Exppol(\XQvp \col \tauQ),$ where $d:=\dega(f)$.
Finally,
if $R> 1,$ then the 
the series on the right-hand side of (\refer{e: expansion f along Q v})
converges neatly on
$\AQqp(R^{-1})$ as a 
$\DrQ$-exponential polynomial series with coefficients in 
$\Ci(\XQvpb{R}\col \tauQ).$ 
\end{enumerate}
\end{thm}

\proof
We will establish existence. Uniqueness then follows from 
uniqueness of asymptotics, see Lemma \refer{l: uniqueness of asymp}.

Fix $P \in \minparabs$ with $P \subset Q.$ 
Select a complete set $\cWQv \subset \NKQaq$ of representatives
for $W_Q/W_Q\cap \WKvH.$ 

The set $\cWQv v$ maps injectively into the coset space $W/\WKH.$ 
Hence it may be extended to a complete set $\cW$ of representatives 
in $\NKaq$ for 
$W/\WKH.$ In view of Lemma \refer{l: splitting lemma} we may therefore decompose 
$f,$ if necessary, so that we arrive in the situation
that there exists a $s \in \faqdc$ such that $\Exp(P,uv\asmid f) \subset 
s - \N \gD(P),$ for all $u \in \cWQv.$ 
We put $s_Q = s|\faQq.$

Let $u \in \cWQv.$ Then the function
$f_{uv}: a \mapsto f(auv)$ has a (unique) $\Delta(P)$-exponential
polynomial expansion on $\Aqp(P)$
of the following type:
\begin{equation}
\naam{e: expansion with p u xi}
f_{uv}(a) = f(auv) = \sum_{\xi \in s - \N\gD(P)} q_{u,\xi}(\log a) a^\xi.
\end{equation}
Here $q_{u,\xi}(\dotvar) = q_\xi(P,uv\asmid f, \dotvar,e)$ 
belongs to $P_d(\faq)\otimes \VtauuvH$. 

Let $\diffop \in S(\faq).$ Then according to Lemma
\refer{l: formal application of Ufa}, the function  $\diffop f_{uv}$ is given on $\Aqp(P)$ 
by a neatly convergent $\gD(P)$-exponential 
polynomial series that is obtained from (\refer{e: expansion with p u xi})
 by term by term application
of $\diffop.$ That is,
\begin{equation}
\naam{e: series with diffop}
\diffop f_{uv} (a) =  \sum_{\xi \in s - \N\gD(P)} q_{\diffop, u,\xi}(\log a) a^\xi,
\end{equation}
where  $q_{\diffop, u, \xi}$ is the $\VtauuvH$-valued polynomial function 
on $\faq$ of degree at most $d$  given by 
$$
q_{\diffop, u, \xi}(X) =
 e^{-\xi(X)} \diffop [e^{\xi(\dotvar)} q_{u,\xi}](X)\qquad (X \in \faq).
$$
Let now $R >1$ and let $\cptb$ and $\cpta$ be compact subsets of $\stAqp(\starP)_{[R]}$ 
and $\AQqp(R^{-1}),$ respectively. Then $\cpta \cptb$ is a compact
subset of $\Aqp(P),$ by Lemma \refer{l: inclusions for asymp} (a).
Thus, if $a \in \cpta$ and $b \in \cptb,$ then 
the series in 
(\refer{e: series with diffop}) with $ba$ in place of $a$ converges 
absolutely, and may be rearranged as follows:
\begin{equation}
\naam{e: rearranged expansion for f}
\diffop f_{uv}(ab) = 
\sum_{\eta \in s_Q - \N\DrQ} 
a^\eta \sum_{\xi \in s - \N \DP\atop {\xi |{\faQq} = \eta} }
b^\xi\, q_{\diffop, u,\xi}(\log b + \log a).
\end{equation}
In view of Lemma \refer{l: neat conv of exp pol series}, the convergence is 
absolutely uniformly for $(a,b) \in \cpta \times \cptb.$ 

By a similar reasoning it  follows from the neat convergence of the series
(\refer{e: series with diffop})
that, for any continuous seminorm $\gs_0$ on $P_d(\faq) \otimes \Vtau,$ 
the series 
\begin{equation}
\naam{e: series TAG}
\sum_{\eta\in s_Q - \N \DrQ}
a^{\Re \eta}
\sum_{\xi \in s- \N \DP \atop \xi|{\faQq} = \eta}
b^{\Re \xi} \,
\gs_0(q_{\partial, u, \xi}) 
\end{equation} 
converges uniformly for $a \in \cpta$ and $b \in \cptb.$ 

Let now $\eta \in s_Q - \N \Delta_r(Q)$ and let $b \in \stAQqp(\starP)$ 
and $a \in \AQqp(R_{Q,v}(b)^{-1}).$ 
Then there exists a $R>1$ such that $b \in \stAqp(\starP)_{[R]}$ 
and $a \in \AQqp(R^{-1}).$ Hence the series (\refer{e: series TAG}) converges, and by positivity 
of all of its terms we infer that the series 
\begin{equation}
\naam{e: series with gs zero and diffop}
\sum_{\xi \in s - \N \DP\atop {\xi | \faQq = \eta} }
 b^{\Re \xi} \,\gs_0(q_{\diffop, u,\xi})
\end{equation}
converges for every continuous seminorm $\gs_0$ on
$P_{d}(\faq) \tensor \Vtau,$ for every $b \in \stAQqp(\starP).$ 

We now specialize to $\diffop = 1$ and note that $q_{1,u,\xi} = q_{u, \xi}.$ 
Let $X \in \faQq.$ 
We define the linear endomorphism $T_X$ of $P_d(\faq) \tensor \Vtau$ 
by $T_Xp(H) = p(X + H).$ This endomorphism is continuous linear by finite dimensionality.
Combining this with the convergence of 
(\refer{e: series with gs zero and diffop}) we infer, 
for every $X \in \faQq,$ that 
\begin{equation}
\naam{e: expansion q eta}
q_{Q,u,\eta}(X, b) := \sum_{\xi \in s - \N \DP\atop {\xi |\faQq = \eta} }
b^\xi \,T_X (q_{u,\xi})(\log b)
\end{equation}
is a function of $b$ defined by a neatly convergent $\Delta_Q(P)$-exponential polynomial series
on $\stAQqp(\starP).$ It is polynomial in $X$ of degree at most $d,$ 
and real analytic in $b \in \stAqp(\starP).$
Moreover, its values are in the space $\VtauuvH.$ Thus 
$q_{Q,u, \eta} \in P_d(\faQq) \otimes \Cep(\stAQqp(\starP), \VtauuvH).$ 
In view of the isomorphism (\refer{e: isomorphism of exppol}) for 
$\XQvp$, $\tauQ$, $\cWQv$ in 
place of $\spXp$, $\tau$, $\cW,$ see also the decomposition 
(\refer{e: deco XQvp}),
there exists a unique polynomial function 
$q_\eta = q_{\eta}(Q,v\asmid f)$ on $\faQq$ with values in $\Cep(\XQvp \col \tauQ)$
such that 
\begin{equation}
\naam{e: defi q eta in proof}
q_\eta(X, bu) = q_{Q, u,\eta}(X, b), 
\qquad (X \in \faQq, u \in \cWQv, b \in \stAQqp(\starP)).
\end{equation}
The degree of $q_\eta$ as a polynomial function
on $\faQq$ is at most $d.$ 
Combining this 
with (\refer{e: expansion q eta}) and (\refer{e: rearranged expansion for f}) 
and using that $R_{Q,v}(bu) = R_{Q,v}(b),$ we arrive at the 
expansion (\refer{e: expansion f along Q v}) for $m = bu$ and 
$a \in \AQqp(R_{Q,v}(m)^{-1}).$ 
Using the left $K_Q$-invariance of $R_{Q,v}$ and the sphericality 
of $f$ and the functions $m \mapsto q_\eta(\log a, m),$ we now obtain
(\refer{e: expansion f along Q v}) with absolute convergence;
the first two assertions of (b) follow as well. The assertion
of neat convergence in (a) is a consequence of the final assertion in (b),
which we will now proceed to establish.

Let  $u \in \cWQv$ and
 $R > 1$ be fixed. 
Then in view of the union (\refer{e: deco XQvp added R}) it 
suffices to prove the neat 
convergence of
the series (\refer{e: expansion f along Q v}) as a 
$\DrQ$-exponential polynomial series with coefficients
in $C^\infty(K_Q \stAQqp(\starP)\lbr{R} \,u(\MQ \cap \vH)\col \tauQ).$ The 
map $(k, a) \mapsto kau(M_Q\cap \vH)$ induces a diffeomorphism from
$K_Q/(K_Q \cap \vH) \times \stAQqp(\starP)\lbr{R}$ onto the open subset
$K_Q \stAQqp(\starP)\lbr{R}\,u(\MQ \cap \vH)$ of $\XQvp.$ 
By sphericality of the coefficients of the series 
(\refer{e: expansion f along Q v}) we see
that it suffices to prove that
$$
\sum_{\eta \in s_Q - \N \DrQ} a^\eta \gs_1(q_{Q,u,\eta})
$$ 
converges absolutely, for $a \in \AQqp(R^{-1})$ and for 
$\gs_1$ any continuous seminorm on $P_d(\faQq) \otimes \Ci(\stAQqp(\starP)\lbr{R}, \VtauuvH).$ 

Fix $X \in \faQq,$ $\stdiffop \in S(\staQq),$ 
$a \in \AQqp(R^{-1})$ and  $\cK\subset  \stAQqp(\starP)\lbr{R}$ 
a compact subset.
Then it suffices to prove that 
\begin{equation}
\naam{e: series with stdiffop to be estimated}
\sum_{\eta \in s_Q - \N \DrQ} a^\eta \sup_\cK \| \stdiffop ( q_{Q,u,\eta}(X, \dotvar))\|
\end{equation}
converges absolutely.

{}From the neat convergence of the series (\refer{e: expansion q eta}),
for $b \in \stAQq(\starP),$  
it follows 
that term by term differentiation is allowed. Since $\stdiffop \in S(\staQq),$
whereas $X \in \faQq,$  
we have
$$
b^{-\xi}\stdiffop ( b^\xi T_X (q_{u,\xi})(\log b))
=
q_{\stdiffop, u, \xi}(X + \log b).
$$ 
Hence, for every $\eta \in s_Q - \N\DrQ,$   
\begin{equation}
\naam{e: series for stdiffop q} 
\stdiffop ( q_{Q,u,\eta}(X, \dotvar))(b) = 
\sum_{\xi \in s - \N \DP\atop \xi|\faQq = \eta}
b^\xi q_{\stdiffop, u, \xi}(X + \log b).
\end{equation}

There exists a continuous seminorm $\gs_2$ on
$P_{d}(\faq) \otimes \Vtau,$ 
such that, for every $b \in \cK$ and all
$q \in P_{d}(\faq) \otimes \Vtau,$ 
$$ 
\|q(X + \log b) \| \leq \gs_2(q).
$$ 
In particular, this implies that 
\begin{equation}
\naam{e: estimate q on gs two}
\| q_{\stdiffop, u, \xi}(X + \log b)\| \leq 
 \gs_2(q_{\stdiffop, u, \xi}),
\end{equation}
for every $b \in \cK.$

Combining  
(\refer{e: series for stdiffop q})
with (\refer{e: estimate q on gs two})
we now obtain
$$
| a^\eta |  \sup_{\cK} \| \stdiffop(q_{Q,u,\eta})(X, \dotvar)\|
\leq 
\sum_{\xi \in s - \N\DP\atop \xi|\faQq = \eta}
a^{\Re \eta} b^{\Re \xi}\, \gs_2(q_{\stdiffop, u, \xi}).
$$
Thus, the absolute convergence of 
(\refer{e: series with stdiffop to be estimated}) 
follows from the uniform convergence of 
(\refer{e: series TAG}), $b \in \cK.$ 
\qed
Let $f \in \ExppolXptau$
 and let 
$Q \in \allparabs$ and $v \in \NKaq.$ Moreover, 
let the set $Y \subset \faQqdc$ and the polynomials $q_\eta = q_\eta(Q,v\asmid f),$ 
for $\eta \in Y- \N\DrQ$ 
be  as in Theorem \refer{t: expansion along the walls}. As in that
theorem, we define 
$$
\Exp(Q,v\asmid f) = \{ \eta \in Y - \N\DrQ \mid q_\eta \neq 0\}
$$
and call the elements of this set the exponents
of $f$ along $(Q,v).$ 
If $\eta \in \faQqdc$ does not belong to $\Exp(Q,v\asmid f),$ we agree to write
$q_\eta(Q,v\asmid f) = 0.$ 

Let now $P \in \minparabs$ be contained in
$Q$ and put $\starP:= P \cap M_Q.$ Then, for $u \in \NKQaq,$ we define  
$$
\Exp(Q,v\asmid f)_{P, u} = \{ \eta \in \faQqdc \mid q_{\eta} \neq 0 \text{on}
\faQq \times K_Q \stAqp(\starP) u (M_Q \cap \vH) \}.
$$
The elements of this set are called the $(Q,v)$-exponents 
of $f$ on  $\stAQqp(\starP)u.$ Let $\cWQv \subset N_{K_Q}(\faq)$ be a complete
 set of
 representatives
of $W_Q / W_Q \cap W_{K \cap vHv^{-1}}.$ Then it follows from 
(\refer{e: deco XQvp})  that 
\begin{equation}
\naam{e: Exp Q v f as union over index set cWQv}
\Exp(Q,v\asmid f) = \displaycup_{u \in \cWQv} \Exp(Q,v\asmid f)_{P, u}.
\end{equation} 
We now have the following result.
\begin{thm}
\naam{t: transitivity of asymptotics}
{\rm (Transitivity of asymptotics)\ }
Let $f \in \ExppolXptau.$ 
Let $P,Q \in \allparabs,$ assume that $P$ is minimal
and $P \subset Q$ and put $\starP = P \cap M_Q.$
Then for all $v \in \NKaq$ and $u \in \NKQaq$ we 
have:
\begin{equation}
\naam{e: equality of exponents}
\Exp(Q,v \asmid f)_{P,u} = \Exp(P, uv\asmid f )\,|_{\faQq}.
\end{equation}
Moreover, if $\eta\in \Exp(P, uv\asmid f ) |_{ \faQq},$ 
then for every $b \in \stAqp(\starP),$ $X \in \faQq,$ and $m \in M,$ 
\begin{equation}
\naam{e: q eta by transitivity}
q_{\eta}(Q,v \asmid f , X, mbu) =
\sum_{\xi \in \Exp(P,uv\asmid f ) \atop {\xi |_{\faQq} = \eta }}
b^\xi q_{\xi}(P, uv \asmid f,  X + \log b, m),
\end{equation}
where the $\DQP$-exponential polynomial series (in the variable $b$) 
on the right is 
neatly convergent on $\stAqp(\starP).$ 
Furthermore, 
the series
\begin{equation}
\naam{e: series for q eta without m}
\sum_{\xi \in \Exp(P,uv\asmid f ) \atop {\xi |_{ \faQq} = \eta }}
b^\xi q_{\xi}(P, uv \asmid f,  X + \log b)
\end{equation} 
converges neatly as a $\Delta_Q(P)$-exponential polynomial series
in the variable $b \in \stAqp(\starP)$ 
with coefficients in $\Ci(\spX_{0,uv}\col \tauM).$ 
\end{thm}

\proof Let $v \in \NKaq$ and $u \in \NKQaq$ be fixed.
Fix a set $\cW_{Q,v}$ such as in the beginning of the proof 
of Theorem \refer{t: expansion along the walls}, and such that it contains $u.$ 
Moreover, we select a set $\cW$ of representatives for $W/\WKH$ in $\NKaq$ 
containing $\cW_{Q,v}v.$ As in the proof of the mentioned
theorem we may restrict ourselves to the situation 
that $\Exp(P, u'v\asmid f) \subset s - \N \DP,$ for some $s \in \faqdc$ 
and all $u' \in \cW_{Q,v}.$ In the following we may now use the notation
and results of the proof
of Theorem \refer{t: expansion along the walls}. 
 
Let $\eta \in s_Q - \N \DrQ.$ Then from  
(\refer{e: defi q eta in  proof}) and (\refer{e: expansion q eta}) we infer that, for 
every $X \in \faQq,$ 
$$
q_\eta(Q,v\asmid f, X, bu) = \sum_{\xi \in s - \N \Delta(P)\atop \xi|\faQq  = \eta} b^\xi 
q_\xi(P, uv\asmid f, X + \log b, e), \qquad (b \in \stAQqp(\starP));
$$ 
the series on the left-hand side converges neatly as a $\DQP$-exponential polynomial 
series in the variable $b \in \stAQqp(\starP).$ 
The function $m \mapsto q_\eta(Q,v\asmid f, X, mbu)$ belongs to $\Ci(\spX_{0, uv}\col \tauM),$ 
and so does the function $m \mapsto q_\xi(P, uv\asmid f, X + \log b, m),$ for 
every $\xi \in s - \N\DP.$ Evaluation at $e$ induces 
a topological linear isomorphism
$\Ci(X_{0,w}\col \tauM) \simeq V_\tau^{M \cap w H w^{-1}},$ for every 
$w \in \NKaq,$ hence in particular
for $w = uv.$  Thus, it follows from the above that   
(\refer{e: q eta by transitivity}) holds, with the asserted convergence. In addition,
it follows that the series (\refer{e: series for q eta without m}) converges as asserted. 

In the proof of Theorem \refer{t: expansion along the walls} we saw that 
$\Exp(Q,v\asmid f) \subset s_Q - \N\DrQ.$ 
It follows from the derived expansion (\refer{e: q eta by transitivity}) that
(\refer{e: equality of exponents}) holds with the inclusion `$\subset$'
in place of the equality sign.
For the converse inclusion, let $\xi_0 \in \Exp(P, uv\asmid f )$ and
put $\eta = \xi_0 |_{ \faQq}.$ We  select $X \in \faQq$ such 
that the function $b \mapsto q_{\xi_0}(P, uv \asmid f,  X + \log b, e)$ does not
vanish identically on $\stAQq.$ The equality (\refer{e: q eta by transitivity}) 
holds for all $b \in \stAQqp(\starP)$ with a $\Delta_Q(P)$-exponential polynomial series
that converges neatly on $\stAQqp(\starP).$ Any
exponent $\xi$ of this series
coincides with $\eta = \xi_0|_{\faQq}$ on $\faQq;$ if it also coincides with
$\xi_0$ on $\stfaQq,$ then $\xi = \xi_0.$ 
Therefore, the function of $b$ defined by the series 
on the right-hand side of (\refer{e: q eta by transitivity})
is non-zero. Hence $q_\eta(Q,v\asmid f)$ does not vanish identically 
on $\faQq \times \stAQqp(\starP)u$ and we conclude that 
$\eta \in \Exp(Q,v \asmid f)_{P,u}.$ 
\qed

We proceed by discussing 
some useful transformation properties
for the coefficients in the expansion (\refer{e: expansion f along Q v}).

If
 $u \in \NKaq$ it will sometimes be convenient to write $u X := \Ad(u) X$ 
for $X \in \faq.$ Similarly, we will write $u\xi : = \xi \after \Ad(u)^{-1},$ 
for $\xi \in \faqdc.$ 

If $u,v \in \NKaq$ and $Q \in \allparabs,$ then conjugation by $u$ induces
a diffeomorphism $\gg_{u}$ 
from the space $\spX_{Q,v}$  onto   $\spX_{uQu^{-1}, uv};$
we note that $\gg_u$   maps $\spXQvp$ onto $\spX_{uQu^{-1}, uv, +}.$ 
It is easily seen that $R_{uQu^{-1},uv}(\gg_u(m)) = R_{Q,v}(m),$ 
for $m \in \spXQv.$ 

For
$\gf \in \Ci(\spX_{Q,v,+}\col \tauQ),$ we define the function
$\rhotauu\gf: \spX_{uQu^{-1}, uv, +} \to \Vtau$ by 
\begin{equation}
\naam{e: defi rho tau u}
\rhotauu\gf(x) = \tau(u) \gf(\gg_u^{-1}(x)).
\end{equation}
Then $\rhotauu$ is a topological linear isomorphism from
the space
 $\Ci(\spX_{Q,v,+}\col \tauQ)$ onto 
the space
 $\Ci(\spX_{uQu^{-1},uv,+} \col \tauuQ).$ 
Likewise, by similar definitions we obtain a topological linear isomorphism 
from
$\Ci(\spX_{1Q,v,+}\col \tauQ)$ onto
 $\Ci(\spX_{1\,uQu^{-1},uv,+} \col \tauuQ),$
also denoted by $\rhotauu.$ 
\begin{lemma}
\naam{l: transformation of coeffs}
Let $f \in \Exppol(\spXp\col\tau),$ 
let $Q\in\allparabs$ and $u,v \in \NKaq.$
Then 
$$
\Exp(uQu^{-1},uv\asmid f) = u\, \Exp(Q,v\asmid f).
$$ 
Moreover, for every $\eta \in  \Exp(Q,v\asmid f),$
$$
q_{u\eta} (uQu^{-1} ,uv\asmid f) = [\Ad(u^{-1})^* \otimes \rhotauu]\, q_\eta(Q,v\asmid f).
$$ 
\end{lemma}
\proof
Put $Q' = uQu^{-1}.$ 
Let $m \in \spX_{Q', uv,+}.$ 
Then, by Theorem \refer{t: expansion along the walls},
\begin{equation}
\naam{e: first series transf rule}
f(ma uv) = \sum_{\eta \in \Exp(Q',uv\asmid f)} 
a^\eta q_\eta(Q',uv\asmid f)(\log a, m), 
\end{equation} 
for $a \in A_{Q' \iq}^+(R_{Q',uv}(m)^{-1}),$ 
where the series on the right-hand side is neatly convergent. 
On the other hand, from $f(ma uv) = \tau(u) f(\gg_u^{-1}(m )  u^{-1}au \, v)$ 
we see, using Theorem \refer{t: expansion along the walls} again,  that
\begin{equation}
\naam{e: second series transf rule}
f(mauv) = \tau(u) \sum_{\zeta \in \Exp(Q,v\asmid f)} 
                 a^{u\zeta} q_\zeta (Q,v\asmid f)(\Ad(u)^{-1}\log a, \gg_u^{-1}(m)),
\end{equation}
for $u^{-1} a u \in \AQqp(R_{Q,v}(\gg_u^{-1}(m))^{-1}).$ 
We now note that the latter condition is equivalent to 
$$
a \in  A_{Q'\iq}^+(R_{Q,v}(\gg_{u}^{-1}(m))^{-1}) 
= A_{Q'\iq}^+(R_{Q', uv}(m)^{-1}).
$$ 
Hence the series (\refer{e: first series transf rule}) and 
(\refer{e: second series transf rule}) both converge neatly for 
$a \in A_{Q'\iq}(R_{Q', uv}(m)^{-1}).$ 
All assertions now follow by uniqueness of asymptotics.
\qed

For later purposes,
we also need another type of transformation property.
Recall from Remark \refer{r: extreme cases subspaces}
that for $u \in \NKaq$ we write $\spXu= \spX_{1G,u} = G/uHu^{-1};$
let $\spXup$ denote the analogue of $\spXp$ for this symmetric space.
We note that right multiplication by $u$ induces a diffeomorphism $ r_u$
from $\spXu$ onto $\spX,$ mapping $\spXup$ onto $\spXp.$ 
Hence pull-back by $r_u$  the topological 
linear isomorphism $R_u: = r_u^*$ from 
$\Ci(\spXp\col \tau)$ onto $\Ci(\spXup \col \tau);$ it is given by $R_u f(x) = f(xu).$ 
We note that the map $R_u$ coincides with the map $\rho_{\tau, u},$ 
introduced in the text above Lemma \refer{l: transformation of coeffs},
by sphericality of the functions involved.

The following result 
is now an immediate consequence of the definitions.
\begin{lemma}
\naam{l: q of Rv f}
Let $f \in C^\ep(\spXp\col \tau)$  and $u \in \NKaq.$ 
Then $R_u f \in C^\ep(\spXup \col \tau).$ Moreover, for each $Q \in \allparabs$ 
and every $v \in \NKaq,$ the set $\Exp(Q,vu \asmid f)$ equals $\Exp(Q, v \asmid R_uf).$ 
Finally, if $\xi \in \Exp(Q,vu \asmid f),$ then 
$$
q_\xi(Q,vu\asmid f) =  q_\xi(Q,v\asmid R_u f).
$$ 
\end{lemma}

%%%%%%%%%%%%%%%%%%%%%%%%%%%%%%%%%%%%%%%%%%%%%%%%%%%%%%%%%
\section{Behavior of differential operators along walls}
\naam{s: diff op along walls}
We assume that $Q\in \allparabs$ is fixed.
The purpose of this section is to study a $Q$-radial decomposition
of invariant differential operators on $\spX.$ This leads to 
a series expansion of such operators along $(Q,e),$
with coefficients
that turn out to be globally defined on the group $\MQgs.$  
This will be of crucial importance for the applications 
later on (see Proposition \refer{p: stability of globality new}).

The involution $\Cartan\sigma$ fixes $\faq$ pointwise, hence leaves every root space
$\fg_\ga,$ for $\ga \in \Sigma,$ invariant. We denote the 
 associated 
eigenspaces 
 of $\Cartan\gs|_{\fg_\ga}$ 
for the eigenvalues $+1$ and $-1$ by $\fg_\ga^+$ and $\fg_\ga^-,$ 
respectively. Moreover, we put $m_\ga^\pm := \dim \fg_\ga^\pm.$ 

We recall that $K_Q = K \cap M_Q$ and $H_Q = H \cap M_Q.$ 
Define
$\HoneQ := H \cap \MoneQ;$ 
then $\HoneQ = H_Q (A_Q \cap H).$ Note that $K_Q = K \cap \MoneQ.$ 
The group $\MoneQ$ admits the Cartan decomposition
$\MoneQ = K_Q \Aq H_{1Q}$  and  
normalizes the subalgebra 
$\barfn_Q.$

For $m \in \MoneQ$ we define the endomorphism $A(m) = A_Q(m)\in \End(\bfn_Q)$ 
by 
\begin{equation}
\naam{e: defi A m}
A(m):= \gs \after \Ad(m^{-1}) \after \Cartan \after\Ad(m).
\end{equation}
Moreover, we define the real analytic function $\gd = \gd_Q: \MoneQ \to \R$ 
by 
\begin{equation}
\naam{e: defi gdQ}
\gd(m) = \det(I - A(m)).
\end{equation}
Finally, we define the following subset of $\MoneQ$ 
\begin{equation}
\naam{e: defi MoneQpr} 
\MoneQpr := \MoneQ\setminus \gd^{-1}(0).
\end{equation}
 
\begin{lemma}
\naam{l: first lemma on A m}{\ }
\begin{enumerate}
\itema
Let $m \in \MoneQ,$ $k \in K_Q$ and $h \in \HoneQ.$ Then 
$A(kmh) = \Ad(h^{-1}) \after A(m) \after \Ad(h).$ 
\itemb
The endomorphism $A(m) \in \End(\bfn_Q)$ 
is diagonalizable, for every $m \in \MoneQ.$
The eigenvalues are given as follows.
Let $m = k a h,$ with $k \in K_Q,$ $a \in \Aq$ and $h \in \HoneQ.$ 
Then the eigenvalues of $A(m)$ are $\pm a^{-2\ga},$ $\ga \in \Sigma(Q),$ 
with multiplicities $m_{\ga}^\pm.$ 
\itemc
The operator norm of $A(m)$ 
is given by $\|A(m)\|_{\rm op} = R_{Q,1}(m)^2.$ 
\end{enumerate}
\end{lemma}

\proof
(a) is an immediate consequence of (\refer{e: defi A m}).
Hence, for (b) we may assume that $m = a \in \Aq.$ It is easily seen that
$A(a)|_{\fg^\pm_{-\ga}} = \pm a^{-2\ga} I$ for $\alpha\in\Sigma(Q)$. 

Finally, (c) is an immediate consequence of (b) and
 (\refer{e: value RQv})
with $v = 1.$  
\qed

\begin{cor}
If $k\in K_Q,$ $a \in \Aq,$ $h \in \HoneQ$ then  
$$%\begin{equation}
%\naam{e: formula detQ}
\gd(kah) = \prod_{\ga \in \Sigma(Q)} (1 - a^{-2\ga})^{m_\ga^+}(1 + a^{-2\ga})^{m_\ga^-}.
$$%\end{equation}
The set $\MoneQpr$ is left $K_Q$- and right $\HoneQ$-invariant, and open dense
in $\MoneQ.$ 
\end{cor}
\proof
This follows immediately from Lemma \refer{l: first lemma on A m} combined with 
(\refer{e: defi gdQ}) and (\refer{e: defi MoneQpr}).
\qed

We  define the linear subspace $\ofkQ$ of $\fk$ by 
$\ofkQ: = \fk \cap (\fn_Q + \barfn_Q).$ Then the 
map $(I +\Cartan): X \mapsto X + \Cartan X$ 
is a linear isomorphism from $\barfn_Q$ onto $\ofkQ.$ 

\begin{lemma}
\naam{l: Ad m k and h deco}
{\ }
\begin{enumerate}
\itema
If $m \in \MoneQ,$ then $\Ad(m^{-1}) \ofkQ + \fh \subset \bfn_Q + \fh.$ 
\itemb
If $m \in \MoneQpr,$ then $\Ad(m^{-1}) \ofkQ \oplus \fh = \bfn_Q + \fh.$
\end{enumerate}
\end{lemma}

\proof
(a) Since $\ofkQ\subset \bfn_Q +  \fn_Q \subset \bfn_Q + \fh,$ 
we have, for all $m \in \MoneQ,$
$$ 
\Ad(m^{-1}) \ofkQ \subset \Ad(m^{-1}) (\bfn_Q + \fn_Q) =
 \bfn_Q + \fn_Q \subset \bfn_Q + \fh.
$$ 
(b) 
The dimension of $\Ad(m^{-1}) \ofkQ$ equals that of $\ofkQ,$ 
which in turn equals that of $\bfn_Q.$ 
Hence it suffices to prove, for $m \in \MoneQpr,$ that
$\Ad(m^{-1}) \ofkQ \cap \fh = 0.$ 

Let $X \in \Ad(m^{-1})\ofkQ \cap \fh.$ Then $\Cartan \Ad(m)X = \Ad(m) X$
and $\gs X = X,$ and we see that $(I - A(m))X = 0.$ If $m \in \MoneQpr$ 
then $\det(I - \Ad(m)) = \gd(m) \not= 0 $ and it follows that $X = 0.$ 
\qed

{}From Lemma \refer{l: Ad m k and h deco}(b) we see that for $m \in \MoneQpr$ we 
may define linear maps $\Psi(m) = \Psi_Q(m) \in \Hom(\bfn_Q, \ofkQ)$ and 
$R(m) = R_Q(m) \in \Hom(\bfn_Q, \fh)$ by 
\begin{equation}
\naam{e: defi Psi and R}
X = \Ad(m^{-1}) \Psi(m) X + R(m) X.
\end{equation}
\begin{lemma}
Let $m \in \MoneQpr,$ $k \in K_Q$ and $h \in \HoneQ.$ Then 
\begin{eqnarray*}
\Psi(kmh) &= & \Ad(k)\after  \Psi(m) \after \Ad(h),\\
R(kmh) &=& \Ad(h^{-1}) \after R(m) \after \Ad(h).
\end{eqnarray*}
\end{lemma}

\proof 
This is an immediate consequence of (\refer{e: defi Psi and R}) 
combined with Lemma \refer{l: Ad m k and h deco}(b).
\qed

\begin{lemma}
\naam{l: formulas for Psi and R}
Let $m \in \MoneQpr.$ 
Then 
\begin{eqnarray*}
\Psi(m) \after (I - A(m)) &=&  (I + \Cartan) \after \Ad(m),\\
R(m) \after (I - A(m))\, &=& - (I + \gs) \after A(m).
\end{eqnarray*}
\end{lemma}
\proof
{}From (\refer{e: defi A m}) it follows that 
$$ 
I + \gs \after A(m) = \Ad(m^{-1})\after (I + \Cartan) \after \Ad(m).
$$ 
This implies in turn that 
\begin{equation}
\naam{e: expression I min A m}
I - A(m) = \Ad(m^{-1})\after (I + \Cartan) \after \Ad(m) - (I + \gs) \after A(m).
\end{equation}
Since $I + \Cartan$ and $I + \gs$ map  $\bfn_Q$ into 
$\ofkQ$ and $\fh,$ respectively,  the lemma follows from 
combining (\refer{e: expression I min A m}) with (\refer{e: defi Psi and R}).
\qed

\begin{cor}
The functions $\Psi: \MoneQpr \to \Hom(\bfn_Q, \ofkQ)$ 
and $R: \MoneQpr \to \Hom(\bfn_Q, \fh)$ are real analytic.
Moreover, the functions $\gd\,\Psi$ and $\gd \,R$ extend to 
real analytic functions on $\MoneQ.$ 
\end{cor}

\proof 
{}From (\refer{e: defi gdQ}) and  (\refer{e: defi MoneQpr}) we see 
that $I - A(m)$ is an invertible endomorphism of $\bfn_Q,$ for $m \in \MoneQpr.$ 
Since $\Ad(m)$ and $A(m)$ depend real analytically on $m \in \MoneQ,$ all statements
now follow from Lemma \refer{l: formulas for Psi and R}.
\qed

If $R >0,$ then in accordance with (\refer{e: defi spXQvb R}) we define
$$
\MoneQ\br{R}: = \{ m\in \MoneQ \mid R_{Q,1}(m) < R\}.
$$ 
Moreover, we set $\MQgs\br{R}: = \MQgs \cap \MoneQ\br{R}.$ 

\begin{lemma}{\ }
\naam{l: about MoneQpr}
\begin{enumerate}
\itema
$\MoneQb{1} \subset \MoneQpr.$
\itemb
Let $R_1,R_2 > 0.$ Then $\MQgsb{R_1}\,\AQqp(R_2) \subset \MoneQb{R_1 R_2}.$ 
\end{enumerate}
\end{lemma}

\proof
Let $m \in \MoneQ\br{1}.$ Then $\|A(m)\|_{\rm op}<1$ by Lemma 
\refer{l: first lemma on A m}(c), and hence $\delta(m)\neq0$.
This establishes (a).

Assume that  $m \in \MQgsb{R_1}$ and $a \in \AQqp(R_2).$ 
Write $m = k b h$ with $k \in K_Q,$ $b \in \stAQq$ 
and $h \in \HoneQ.$ Then $ma = k (ab) h,$ hence $R_{Q,1}(ma) = \max_{\ga \in \Sigma(Q)}
a^{-\ga} b^{-\ga} < R_2 R_{Q,1}(m)< R_1R_2.$  It follows that $ma \in \MoneQb{R_1 R_2}.$ 
\qed

\begin{prop}
\naam{p: series expansions for Psi and R}
There exist unique real analytic functions $\Psi_\mu, R_\mu: M_{Q\gs} \to \End(\bfn_Q),$ 
for $\mu \in \N \DrQ,$ such that for every $m \in M_{Q\gs}$ and every 
$a \in \AQqp(R_{Q,1}(m)^{-1}),$ 
\begin{eqnarray*}
\Psi(ma) &=& 
(1 + \Cartan) \after \sum_{\mu \in \N \DrQ} a^{-\mu} \Psi_\mu(m),          \\
R(ma) &=& (1 + \gs) \after \sum_{\mu \in \N \DrQ} a^{-\mu} R_\mu(m),
\end{eqnarray*}
with absolutely convergent series.
For every $R>1$ the above series converge neatly on $\AQqp(R^{-1})$ 
as $\DrQ$-power series with coefficients in $\Ci(\MQgs[R], \End(\bfn_Q)).$ 
\end{prop}

\proof
Let $m \in \MQgs$ and $a \in \AQqp(R_{Q,1}(m)^{-1}).$ It 
follows from Lemma \refer{l: about MoneQpr} 
that $ma \in \MoneQ[1] \subset \MoneQpr.$ 
Hence $\Psi(ma)$ and $R(ma)$ are defined.

It follows from Lemma \refer{l: first lemma on A m} that $\|A(ma)\|_{\rm op} < 1.$ 
Hence the series
$$
(I - A(ma))^{-1} = \sum_{n=0}^\infty A(ma)^n
$$ 
converges absolutely. Let $\ga \in \SrQ.$ Then $A(m)$ leaves the space
$\fg_{-\ga}$ invariant, and $A(ma)|_{\fg_{-\ga}} = a^{-2\ga} A(m)|_{\fg_{-\ga}}.$ 
Hence, in view of Lemma \refer{l: formulas for Psi and R}, 
$$ 
\Psi(ma)|_{\fg_{-\ga}} = 
(I + \Cartan) \after \Ad(m) \after \sum_{n=0}^\infty a^{-(2n +1)\ga}
A(m)^n |_{\fg_{-\ga}}
$$
and 
$$ 
R(ma)|_{\fg_{-\ga}} = - (I + \gs) \after \sum_{n=1}^\infty a^{-2n\ga}
A(m)^n |_{\fg_{-\ga}}.
$$
It is now easy to complete the proof.
\qed

We denote by $\ringQplus$ the algebra of functions on $\MoneQpr$ generated  by
the functions $\xi \after \Psi_Q,$ where $\xi \in \Hom(\barfnQ, \ofkQ)^*,$
and by the functions $\eta \after \restQ,$ where $\eta \in \Hom(\barfnQ, \fh)^*.$ 
By $\ringQ$ we denote the algebra of functions generated by $1$ and $\ringQplus.$ 
Note that $\ringQplus$ is an ideal in $\ringQ.$ 

\begin{cor}
\naam{c: cor on ringQ} 
The elements of $\ringQ$ are left $K_Q$- and right $\HoneQ$-finite functions
on $\MoneQpr.$ 

Let $\gf \in \ringQ.$ There exists a $k \in \N$ such that
$\detQ^k \gf$ extends to a real analytic function on $\MoneQ.$ 
Moreover, there exist unique real analytic functions $\gf_\xi$ 
on $\MQgs,$ for $\xi \in \N\DrQ,$ such that for every 
$m \in \MQgs$ and every $a \in \AQq(R_{Q,1}(m)^{-1}),$ 
\begin{equation}
\naam{e: series for phi}
\gf(ma) =\sum_{\xi \in \N\DrQ} a^{-\xi} \gf_\xi(m).
\end{equation}
Let $R\geq 1.$ Then the series (\refer{e: series for phi}) 
converges neatly on $\AQq(R^{-1}),$ as an exponential polynomial
series with coefficients in $\Ci(\MQgsb{R}).$ 

Finally, if $\gf \in \ringQplus,$ then (\refer{e: series for phi}) 
holds with $\gf_0 = 0.$ 
\end{cor} 

\proof
Uniqueness of the functions $\gf_\xi$ is obvious. Therefore it
suffices to prove existence and the remaining assertions.
One readily checks that it suffices to prove the  
assertions for a collection of generators of the algebra $\ringQplus.$ 
Such a collection of generators is formed by the functions of 
the form $\gf = \xi \after \Psi,$ with $\xi \in \Hom(\barfnQ, \ofkQ)^*,$
and by the functions of the form $\gf = \eta \after \restQ,$ where $\eta \in \Hom(\barfnQ, \fh)^*.$ 
For both types of generators all assertions follow immediately from 
Proposition \refer{p: series expansions for Psi and R}.
\qed

As is the previous section we  assume that
$\tau$ is a smooth representation of $K$ in a locally convex space 
$\Vtau.$ 
The space of continuous linear endomorphisms of $\Vtau$ is denoted 
by $\End(\Vtau).$ 

If  an element $u$ of the space
\begin{equation}
\naam{e: defi cDoneQ}
\cDoneQ : =  \radQalg
\end{equation}
 is of the form 
$ \gf \otimes L \otimes v,$ with  $\gf \in \ringQ,$ $L \in \End(\Vtau),$ and
$v \in U(\fmoneQ),$ then we define the differential operator
$u_*$ on $\Ci(\MoneQpr, \Vtau)$ by $u_*f = \gf L\after [R_v f];$ 
here $R$ denotes the right regular representation.
The map $u \mapsto u_*$ extends to an injective  linear map from 
$\cDoneQ$ to the
space of smooth $\End(\Vtau)$-valued differential operators of $\Ci(\MoneQpr, \Vtau).$ 
We also define the subspace 
$$
\cDoneQplus := \radQalgplus.
$$
Via the map $u \mapsto 1 \otimes I \otimes u$ we identify $U(\fmoneQ)$ with 
a subspace of $\cDoneQ.$ Then $u_* = R_u$ for $u \in U(\fmoneQ).$  

Let $\MoneQp$ be the preimage in $\MoneQ$ of the set $\spX_{1Q,1,+}$
(see below (\refer{e: deco XQvp})).
The  set 
$$
\MoneQppr:= \MoneQp \cap \MoneQpr
$$
 is an open dense subset
of $\MoneQ$ that is left $K_Q$- and right $\HoneQ$-invariant.

In view of the decomposition $\fg = \barfn_Q \oplus (\fmoneQ + \fh),$ 
there exists, 
for every $D \in U(\fg),$ an element $D_0 \in U(\fmoneQ)$ 
with $\deg(D_0) \leq \deg(D),$ such that 
\begin{equation}
\naam{e: characterization D zero} 
D - D_0 \in \barfnQ U(\fg) + U(\fg)\fh.
\end{equation}
The element $D_0$ is  uniquely determined
modulo $U(\fmoneQ)\fhoneQ.$
We recall from \bib{B91}, Sect.\ 2, see also \bib{BSft}, p.\ 548-549,
that the assignment $D \mapsto D_0$ induces an algebra homomorphism 
$\mu_Q'={}^\backprime\mu_{\bar Q}: \DX\to \D(\MoneQ/\HoneQ),$
and that the homomorphism  
$\mu_Q: \DX\to \D(\MoneQ/\HoneQ)$, defined by
$\mu_Q(D)= d_Q\after \mu_Q'(D) \after d_Q^{-1}$ with
$d_Q(m):=|\det(\Ad(m)|_{\fn_Q})|^{1/2}$ for $m\in\MoneQ$, 
only depends on $Q$ through the Levi component $\MoneQ$.

\begin{prop}
\naam{p: radial deco with MQ}
Let $D \in \DX.$ 
There exists a $\uplus \in \cDoneQplus$ of degree $\deg(\uplus) < \deg (D)$ 
such that,
for every $f \in \Ci(\spXp\col \tau),$ 
$$
Df|_\MoneQppr = [\mu_Q'(D) + u_{+*}](f|_\MoneQppr).
$$
\end{prop}

\proof 
By induction on the degree  we will first establish the following assertion
for an element $D$ of $U(\fg).$ Let $D_0 \in U(\fmoneQ)$ satisfy 
(\refer{e: characterization D zero}).
Then there exist finitely many $\gf_i \in \ringQplus,$ $u_i \in U(\fk),$ 
and $v_i \in U(\fmoneQ),$ for $1\leq i \leq n,$ such that
$\deg(u_i) + \deg(v_i) < \deg(D),$ and such that 
\begin{equation}
\naam{e: general radial expression for D}
 D -  D_0 \congruent \sum_{i=1}^n \;\gf_i(m)\, [\Ad(m)^{-1}u_i]\, v_i 
\quad \text{mod} U(\fg)\fh,
\end{equation} 
for every $m\in \MoneQpr.$ 

The assertion is trivially true for $ D$ constant. 
Thus, assume that
$D$ is not constant and that the assertion has been established for 
$D$ of strictly smaller degree. Let $ D_0 \in U(\fmoneQ)$ be as above.
Then, modulo $U(\fg)\fh,$  $D - D_0$ equals a finite sum of terms
of the form $X D_1,$ with $X \in \barfnQ$ and $D_1 \in U(\barfn_Q \oplus \fmoneQ)$ 
such that $\deg D_1 < \deg D.$ 

For $m \in \MoneQppr$ we have 
$X = \Ad(m)^{-1} \Psi(m)X + \restQ(m)X;$ hence 
$$%\begin{equation}
%\naam{e: equation for X Done}
X D_1 \congruent ( \Ad(m)^{-1} \Psi(m)X) D_1 + [\restQ(m)X, D_1] \quad \text{mod} U(\fg)\fh.
$$%\end{equation}
Now $\Ad(m)^{-1} \Psi(m)X$ is a finite sum of terms of 
the form $\gf(m) [\Ad(m)^{-1} u]$ with $u \in \ofkQ$ and  $\gf \in \ringQplus.$ 
Applying the induction hypothesis to $D_1$ we see that $[\Ad(m)^{-1}\Psi(m)X]D_1$ 
may be expressed as a sum similar to the one on the right-hand side of 
(\refer{e: general radial expression for D}). 

On the other hand, $[\restQ(m)X, D_1]$ 
is a finite sum of elements of the form $\psi(m) D_2,$ 
with $\psi \in \ringQplus$ and $D_2 \in U(\fg),$ $\deg D_2 < \deg D.$ 
Applying the induction
hypothesis to $D_2,$ we see that $[\restQ(m)X, D_1]$ may also be expressed
as a sum of the form (\refer{e: general radial expression for D}). 
This establishes the assertion involving 
(\refer{e: general radial expression for D}) of the beginning of the proof.

Let now $D \in \DX.$ By abuse of notation we use the same symbol 
$ D$ for a representative of $D$ in $U(\fg)^H,$ and
let $ D_0$ be as above. Then $\mu_Q'(D)$ equals the canonical image of
$ D_0$ in $U(\fmoneQ)^{\HoneQ}.$ Let $\gf_i, u_i, v_i$ be as above and such that
(\refer{e: general radial expression for D}) holds.
Then for every $f \in \Ci(\spXp \col \tau)$ and all $m\in \MoneQppr$ we have
\begin{eqnarray*}
Df(m) & =& \mu_Q'(D) (f|_\MoneQppr)(m) + \sum_{i=1}^n \gf_i(m) R_{\Ad(m)^{-1}u_i}R_{v_i} f(m)
\\ 
&=& \mu_Q'(D) (f|_\MoneQppr)(m) + \sum_{i=1}^n \gf_i(m) \tau(\check u_i)R_{v_i} f(m)
\end{eqnarray*}
where we have used that $R_{\Ad(m)^{-1}u_i}R_{v_i} f(m) = 
L_{\check u_i} R_{v_i} f(m) = 
\tau(u_i) R_{v_i} f(m).$ 
Thus, we obtain the desired expression with 
$ 
\uplus = \sum_{i=1}^n \gf_i \otimes \tau(u_i) \otimes v_i.
$ 
\qed
Let $U \subset M_{Q\gs}$ be an open subset.
It will be convenient to be able to refer to a `formal application' of 
elements of the space  
$\cDoneQ,$ defined in  (\refer{e: defi cDoneQ}), 
to
$\cF^\ep(\AQq, C^\infty(U, \Vtau)),$ the space of (formal) $\DrQ$-exponential polynomial series
with coefficients in  $C^\infty(U, \Vtau),$ see the definition preceding 
Lemma \refer{l: formal application of Ufa}. There is 
a natural way to define a formal application that is  
compatible with the expansions of Corollary \refer{c: cor on ringQ} and with 
the map $u \mapsto u_*,$ defined in the text following
(\refer{e: defi cDoneQ}). The motivation for the following somewhat 
tedious chain of definitions will become clear in Lemma
\refer{l: formal application of cDoneQ}.

The product decomposition $M_{1Q} \simeq M_{Q\gs} \times \AQq$ induces 
a natural isomorphism
 from $U(\fmoneQ)$ onto $U(\fmQgs) \otimes U(\faQq),$ 
by which we shall identify.
Accordingly we have  a natural isomorphism
\begin{equation}
\naam{e: tensor prod deco of cDoneQ}
\cDoneQ \simeq \cDoneQgs \otimes U(\faQq),
\end{equation}
where $\cDoneQgs:= \cR_Q \otimes \End(\Vtau) \otimes U(\fm_{Q\gs}).$
To each element $\gf \in \cR_Q$ we may associate its $\DrQ$-exponential polynomial
series of the form (\refer{e: series for phi}); this induces a linear
embedding $\cR_Q \to \cF^\ep(\AQq, C^\infty(\MQgs))$ which by identity on the 
other tensor components may be extended to a linear embedding
$$
\cDoneQgs \to \cF^\ep(\AQq, \cDQgs),
$$ 
where 
$\cDQgs := C^\infty(\MQgs) \otimes \End(\Vtau) \otimes U(\fmQgs).$ 
By identity on the second tensor component in (\refer{e: tensor prod deco of cDoneQ}) 
this embedding extends to a
linear embedding
\begin{equation}
\naam{e: the embedding ep}
\ep:\;\;\cDoneQ \to \cF^\ep(\AQq, \cDQgs) \otimes U(\faQq).
\end{equation} 
The image $\ep(u)$ of an element $u \in \cDoneQ$ under this embedding will be called
the $\DrQ$-exponential polynomial expansion of $u.$  
Via the right regular action of $U(\fmQgs)$ we may naturally 
identify $\cDQgs$  
with the space of  $C^\infty$-differential operators 
acting on  $C^\infty(\MQgs, \Vtau).$ 
Accordingly, we have a continuous bilinear pairing 
$\cD_{Q\gs} \times \Ci(U, \Vtau) \to \Ci(U, \Vtau).$ 
This induces a formal application map 
from $\cF^\ep(\AQq, \cDQgs)\otimes \cF^\ep(\AQq, C^\infty(U, \Vtau))$ to
$\cF^\ep(\AQq, C^\infty(U, \Vtau))$ in the fashion described above Lemma 
\refer{l: formal application of hom valued series}. The image of an element of the form
$u \otimes f$ under this map will be denoted by $uf.$ 
 
On the other hand,
in Lemma \refer{l: formal application of Ufa}
we described the formal 
application map
$U(\faQq) \otimes \cF^\ep(\AQq, C^\infty(U, \Vtau)) \to \cF^\ep(\AQq, C^\infty(U, \Vtau)).$ 
The image of an element of the 
form $v \otimes f$ under this map is denoted by $vf.$ 
Combination of the above formal application maps 
leads to the formal application map
$$ 
[\cF^\ep(\AQq, \cDQgs) \otimes U(\faQq)] \otimes \cF^\ep(\AQq, C^\infty(U, \Vtau)) \to 
\cF^\ep(\AQq, C^\infty(U, \Vtau)),
$$ 
given by $(u \otimes v) \otimes f \mapsto (u\otimes v)f:= u(vf),$ 
for $u \in \cF^\ep(\AQq, \cDQgs),$ $v\in U(\faQq)$ 
and $f \in \cF^\ep(\AQq, C^\infty(U, \Vtau)).$
Composing with the embedding (\refer{e: the embedding ep}) we finally obtain the linear map 
$$
\cDoneQ 
\otimes \cF^\ep(\AQq, C^\infty(U, \Vtau)) \to 
\cF^\ep(\AQq, C^\infty(U, \Vtau))
$$ 
given by $u \otimes f \mapsto uf:= \ep(u)f,$ for $u \in \cDoneQ$ 
and $f \in \cF^\ep(\AQq, C^\infty(U, \Vtau)).$ We shall call this map the 
formal application  of $\cDoneQ$ 
to $\cF^\ep(\AQq, C^\infty(U, \Vtau)).$ 

Let now $R\geq 1$ and let $U \subset \MQgsb{R}$ be an open subset. 
We use the obvious natural isomorphism to identify the space 
$C^\ep(\AQqp(R^{-1}), C^\infty(U, \Vtau))$ with 
a subspace of $C^\infty(U\AQqp(R^{-1}), \Vtau).$  If 
$u \in \cDoneQ,$ then the associated differential operator
$u_*$ induces a map from the first space into the latter.

\begin{lemma}
\naam{l: formal application of cDoneQ}
Let $u \in \cDoneQ,$ let $R\geq 1$ and let $U\subset \MQgsb{R} $ be an open
subset.
Then $u_*$ maps the space $C^\ep(\AQqp(R^{-1}), C^\infty(U, \Vtau))$ 
into itself. Moreover, if $f$ belongs to that space, 
then the $\DrQ$-exponential polynomial expansion of $u_*f$ is obtained
from the formal application of $u$ to the exponential polynomial expansion of $f.$ 
\end{lemma}

\proof
This follows from retracing the definitions of $u_*$ and of the formal application 
of $u$ given above 
and applying Corollary \refer{c: cor on ringQ} and Lemmas
 \refer{l: formal application of Ufa} and 
\refer{l: formal application of hom valued series}.
\qed

Given $v\in\NKaq$ we define
$\mu_{Q,v}: \DX\to \D(\spXoneQv)=\D(\MoneQ / \MoneQ \cap vHv^{-1})$
by 
$$\mu_{Q,v}=\Ad(v)\after\mu_{v^{-1}Qv},$$
where $\Ad(v): \D(\spX_{1v^{-1}Qv,e})\to 
\D(\spXoneQv)$ is induced by the restriction to $U(\fm_{1v^{-1}Qv})$
of $\Ad(v)$ on $U(\fg)$.
Then $\mu_{Q,v}$ depends on $Q$ only through $\MoneQ$. It is easily seen
that 
\begin{equation}
\naam{e: muQv}
\mu_{Q,v}=\mu_Q^v\after\Ad(v)
\end{equation}
where $\mu_Q^v: \DXv=\D(G/vHv^{-1})\to \D(\spXoneQv)=
\D(\MoneQ / \MoneQ \cap vHv^{-1})$ is defined similarly as $\mu_Q$, but
with $H$ replaced by $vHv^{-1}$, and where $\Ad(v): \DX\to 
\DXv$ is induced by $\Ad(v)$ on $U(\fg)$.

Let $\MQgsp=\MQgs\cap\MoneQp$ and, for $R\geq 1$, 
$\MQgsp[R]=\MQgs[R]\cap\MoneQp$.

\begin{lemma}
\naam{l: radial component applied to expansion}
Let $f\in \Cep(\spXp\col \tau)$ and let $D \in \DX$.
Then $Df \in \Cep(\spXp\col \tau).$
 
Let $Q \in \allparabs$ and let 
$\uplus \in \cR_Q^+ \otimes \End(\Vtau) \otimes U(\fmoneQ)$ 
be associated with $D$ as in Proposition \refer{p: radial deco with MQ}. 
Then the following holds.
\begin{enumerate}
\itema
The $\DrQ$-exponential expansion of $Df$ along $(Q,e)$ is 
obtained by the formal application of $\mu_Q'(D) + \uplus$ to the $\DrQ$-exponential polynomial expansion
of $f$ along $(Q,e).$ 
\itemb
Let $v \in \NKaq$, then 
$ 
\Exp(Q,v\asmid Df) \subset \Exp(Q,v\asmid f) - \N\DrQ.
$ 
\itemc
If $\xi$ is a leading exponent of $f$ along $(Q,v),$ then
\begin{equation}
\naam{e: poly of Df}
a^{\xi+\rho_Q} q_\xi(Q,v \asmid Df, \log a, m) =  [\mu_{Q,v}(D) \gf](m a),
\qquad (m \in \MQgsp,\;a\in \AQq),
\end{equation} 
where the function  $\gf: \MoneQp \to \Vtau$ is defined by 
$\gf(ma) = a^{\xi+\rho_Q} q_\xi(Q,v\asmid f, \log a, m),$ 
for $m \in \MQgsp$ and $a\in \AQq.$ 
\end{enumerate}
\end{lemma}

\proof
Let $R\geq 1$ and let $\underline f$ be the function $\AQqp(R^{-1}) \to 
\Ci(\MQgsp[R], \Vtau)$ 
defined by $\underline f (a, m) = f(ma).$ It follows from the hypothesis on $f$ 
and Theorem \refer{t: expansion along the walls} that
$\underline f(a,m)$ belongs to 
$C^\ep(\AQqp(R^{-1}), \Ci(\MQgsp[R], \Vtau)).$ 
Moreover, its $\DrQ$-exponential polynomial expansion coincides with the expansion
of $f$ along $(Q,e).$ 
Put $u = \mu_Q'(D) + \uplus.$ 
Then it follows from the previous lemma that $u_*\underline f$ belongs to 
$C^\ep(\AQqp(R^{-1}), \Ci(\MQgsp[R], \Vtau));$ its expansion is obtained from the 
formal application of $u$ to the $(Q,e)$-expansion of $f.$ It follows from 
Theorem \refer{t: expansion along the walls}
that the expansion is independent of $R$ and that its coefficients 
are functions in $C^\infty(\MQgsp, \Vtau).$ 
On the other hand, it follows from
Proposition \refer{p: radial deco with MQ} that $u_*\underline f(a, m) = Df(ma).$ This implies that
$Df$ has a $\DrQ$-exponential polynomial expansion along $(Q,e)$ with coefficients in 
$C^\infty(\MQgsp, \Vtau).$ Since $Df$ is right $H$-invariant, the coefficients
are actually functions in $C^\infty(\spXQep, \Vtau).$ Moreover, 
the expansion is independent of $R$ and converges neatly on $\AQqp(R^{-1})$ as an expansion with 
coefficients in $\Ci(\XQepb{R},\Vtau).$ In particular this holds for every 
minimal parabolic 
subgroup $Q;$ 
hence $Df \in C^\ep(\spXp\col \tau).$ 

In the above we have established assertion (a). It follows from this assertion that
(b) holds with $v =1$ for every $Q \in \allparabs.$ By Lemma \refer{l: transformation of coeffs}
it also holds for arbitrary
$Q\in \allparabs$ and $v \in \NKaq.$ 

It remains to establish (c). Assume first that $v=e$.
Fix  $\xi \in \ExpL(Q,e\asmid f).$ 
Then by (a), $a^{\xi} q_\xi(Q,e \asmid Df, \log a, m)$ is the term with 
exponent $\xi$ in the series that arises from the formal application 
of $\mu_Q'(D) + \uplus$ to the $(Q,e)$-expansion of $f.$
The exponents of the expansion $\ep(\uplus)$ of $\uplus$ all belong to 
$-[\N \DrQ]\setminus \{0\}.$ The application of 
$\uplus$ therefore 
gives rise to an expansion with exponents in $\Exp(Q,e\asmid f) - 
[\N \DrQ]\setminus \{0\}.$ 
The latter set does not contain $\xi,$  since $\xi$ is leading. 
Hence $a^{\xi} q_\xi(Q,e \asmid Df, \log a, m)$
is the term with exponent $\xi$ in the
expansion that arises from the formal application of $\mu_Q'(D)$ to the 
$(Q,e)$-expansion of $f.$ 
Now $\mu_Q'(D) \in U(\fmoneQ) \simeq U(\fmQgs) \otimes U(\faQq)$ 
and we see that the formal 
application of $\mu_Q'(D)$ to the $(Q,e)$ expansion of $f$ 
is induced by term by term differentiation
in the $\AQq$ and the $\MQgs$ variables. This implies that
$a^{\xi} q_\xi(Q,e \asmid Df, \log a, m) =  [\mu_{Q}'(D) \gf'](m a),$
where $\gf'(ma) = a^{\xi} q_\xi(Q,e\asmid f, \log a, m).$ This implies
(\refer{e: poly of Df}) for $v=e$.

Let now $v \in \NKaq$ be arbitrary, and put $f^v=R_vf$. We shall apply the
version of (\refer{e: poly of Df}) just established to the expansion
along $(Q,e)$ of the function $f^v$ on $\spXv$. Let $\xi$ 
be a leading exponent of $f$ along $(Q,v),$ then it follows from Lemma
\refer{l: q of Rv f} that $\xi$ is also a
leading exponent of $f^v$ along $(Q,e)$. Moreover, let $D\in\DX$, then
$(Df)^v=D^vf^v$ where $D^v:=\Ad(v)D\in\DXv$. Hence 
\begin{equation}
\naam{e: poly of Dvfv}
a^{\xi+\rho_Q} q_\xi(Q,e \asmid (Df)^v, \log a, m) =  
[\mu_Q^v(D^v) \gf](m a),
\end{equation}
for $m \in \MQgsp,\;a\in \AQq$,
where $\gf(ma) = a^{\xi+\rho_Q} q_\xi(Q,e\asmid f^v, \log a, m)$.
It follows from Lemma \refer{l: q of Rv f} that 
$\gf(ma) = a^{\xi+\rho_Q} q_\xi(Q,v\asmid f, \log a, m),$
and
$q_\xi(Q,e \asmid (Df)^v)=q_\xi(Q,v \asmid Df)$.
Now (\refer{e: poly of Df}) follows from (\refer{e: poly of Dvfv}) and 
(\refer{e: muQv}).
\qed
\begin{lemma}
\naam{l: action of DX on splitting}
Let $P \in \minparabs$ and assume that
$f \in \ExppolXptau.$ Let $S \subset \faqdc$ be a finite 
set as in Lemma \refer{l: splitting lemma}, and let $D \in \DX.$ 
Then $\Exp(P,v\asmid Df) \subset S -\N\gD$ for every $v \in\NKaq$
and, 
with notation as in Lemma \refer{l: splitting lemma}, 
\begin{equation}
\naam{e: D of fs}
(Df)_s = D(f_s).
\end{equation}
\end{lemma}
\proof
It follows immediately from
Lemma \refer{l: radial component applied to expansion}(b)
that $\Exp(P,v\asmid Df) \subset S -\N\gD$ and that
$\Exp(P,v\asmid D(f_s)) \subset s -\N\gD$ for $s\in S$.
Now (\refer{e: D of fs}) follows from Lemma \refer{l: splitting lemma}.
\qed

\section{Spherical eigenfunctions} 
In this section we assume that $(\tau, \Vtau)$ is a 
finite dimensional continuous representation of $K.$ 
Let $I$ be a cofinite ideal of the algebra 
$\DGH.$ Then by $\sphXpI$ 
we denote the 
space of $f \in \sphXp$ satisfying the system 
of differential equations
$$
Df = 0, \qquad (D \in I).
$$
\begin{rem}
\naam{r: extended results}
Many results of \bib{B87} that are formulated
for $\DGH$-finite $\tau$-spherical functions 
on $\spGH$ are actually valid for the bigger class
of $\DGH$-finite 
functions in $\sphXp$ as well, since their proofs
only involve behavior of functions and operators on
$\spXp.$ If such extended results are used in the
text, we  may give a reference to the present remark.
\end{rem}

\begin{rem}
\naam{r: right translate by v}
Let $v \in N_K(\faq).$
We recall from the text preceding  Lemma \refer{l: q of Rv f}
that right translation by $v$ induces a topological linear
isomorphism $R_v$ from  $\sphXp$ onto the space $\sphXvp.$ 
It maps the subspace of $\DGH$-finite functions onto 
the subspace of $\DXv$-finite functions. Thus,
if  $f \in \sphXp$ is a $\DX$-finite function, then the theory
of \bib{B87} may be applied to the $\DXv$-finite function $R_vf;$ the results
are then easily reformulated in terms of the function $f.$ 
\end{rem}

\begin{lemma}
\naam{l: DX finite in exppol}
Let $I \subset \DX$ be a cofinite ideal. Then 
$\sphXpI \subset \ExppolXptau.$ 
In particular, the elements of $\sphXpI$ are real analytic
functions on $\spXp$.
Moreover,
there exists a finite set $X_I\subset \faqdc$ 
such that $\ExpL(P,v\asmid f)\subset X_I,$ for all 
$f \in \sphXpI,$ $P\in \minparabs$ 
and $v \in \NKaq.$
\end{lemma}

\proof
Let $Q \in \minparabs.$ Applying Theorem 
2.5 of \bib{B87}, see Remark \refer{r: extended results},
we obtain that $f|\Aqp(Q)$ is given by a neatly converging
$\gD(Q)$-exponential polynomial expansion
for each $f \in \sphXpI$. Moreover, by Theorem 
2.4 of \bib{B87}, there exists a finite set $X_{I,Q,e}\subset \faqdc$,
such that $\ExpL(f|\Aqp(Q))\subset X_{I,Q,e}.$ 
Let $w \in \cW.$ Applying the above argument  to  $R_wf,$ 
cf.\ Remark  \refer{r: right translate by v}, we see, more generally,
that $\TdownQw f$ is given by the same type of expansion with
leading exponents in a finite set $X_{I,Q,w}\subset\faqdc$
independent of $f$. 
This implies that $f \in \ExppolXptau,$ with
$\ExpL(P,v\asmid f)\subset X_I:=\cup_{Q,w} X_{I,Q,w}$,
for all $P\in \minparabs$ and $v \in \cW.$ 
Finally, if $v \in \NKaq$ is arbitrary,
there exists $w\in\cW$, $m\in\KM$ and
$h\in N_{K\cap H}(\faq)$ such that $v=mwh$, and then 
$\ExpL(P,v\asmid f)=\ExpL(P,w\asmid f)\subset X_I$.
\qed

 \begin{cor}
Let $P \in \minparabs$ and let $\cW \subset \NKaq$ be a 
complete set of representatives
of $W/\WKH.$
Let $I$ be a cofinite ideal in $\DX.$ 
Then there exists a finite set $S= S_I$ satisfying the properties of 
Lemma \refer{l: splitting lemma} for every $f \in \CepXpI.$ Moreover, 
if $S_I$ is any such set, then $f_s \in \CepXpI$ 
for every $f\in \CepXpI$ and all $s \in S_I.$ 
\end{cor}
\proof
This is an immediate consequence of Lemmas 
\refer{l: DX finite in exppol}
and 
\refer{l: action of DX on splitting}.
\qed

The set $X_I$ in Lemma \refer{l: DX finite in exppol} 
can be described more explicitly if the ideal $I$ has codimension $1.$ 
Let $\fb$ be a maximal abelian subspace of $\fq$ containing $\faq,$  
let $\Sigma(\fb)$ be the restricted root system of $\fb$ in $\fg_\iC,$ 
and let $W(\fb)$ be the associated reflection group. 

Let $\gg$ be the  Harish-Chandra isomorphism from $\DX$  
onto the algebra $I(\fb)$ of $W(\fb)$-invariants in $S(\fb),$ 
see \bib{B91}, Sect.\ 2.
To an element  $\nu \in \fbdc$ we associate the character
$D \mapsto \gg(D\col \nu)$ of  $\DX$ and denote its kernel by $I_\nu.$ 
Then $I_\nu$ is 
an ideal of codimension one in $\DX;$ in fact, 
any codimension one ideal is of this form.

Let $W_0(\fb)$ 
be the normalizer of $\faq$ in $W(\fb).$ Then restriction
to $\faq$ induces an epimorphism from $W_0(\fb)$ onto $W,$ cf.\ \bib{B91}, 
Lemma 4.6. 
We put $\fb_\ik:= \fb \cap \fk.$ Then $\fb = \fbk \oplus \faq.$ 
Moreover, this decomposition is invariant under $W_0(\fb).$ 

\begin{lemma}
\naam{l: restriction on leading exponents}
There exists a finite subset $\cL = \cL_\tau$ of $\fbkdc$ with the following 
property.
Let $\nu \in \fbdc$ and $f \in \Ci(\spXp\col \tau\col I_\nu).$ 
Let $P \in \minparabs, v \in \NKaq$ and assume that 
$\xi \in \ExpL(P,v\asmid f).$ 
Then 
$$%\begin{equation}
%\naam{e: restriction leading exponent}
\nu \in W(\fb)(\cL + \xi + \rho_P).
$$%\end{equation}
\end{lemma}

The proof is based on the following result, which will be proved first.

\begin{lemma}
\naam{l: eigenfunctions on M1}
There exists a finite subset $\cL = \cL_\tau$ of $\fbkdc$ with the following 
property. Let $\nu \in \fbdc$ and $\gf\in\Ci(M_{1}/\HMone\col\tau)$,
and assume that
$$ 
\mu_P(D) \gf = \gg(D\col \nu) \gf
$$
for all $D\in\DGH$, where $\mu_P: \DGH\to \D(M_{1}/\HMone)$ is as defined
above Proposition \refer{p: radial deco with MQ}, with $P\in\minparabs$. 
Then $\gf|_{\Aq}$ is a linear combination of exponential polynomials of 
the form $a\mapsto p(\log a) a^{w\nu}$, where $p\in P(\faq)$ and where
$w\in W(\fb)$ satisfies $w\nu|_{\fbk}\in\cL$.
\end{lemma}

\proof
The algebra $\D(M/\HM)$ acts semisimply on  
$\Ci(M/\HM \col \tau),$ see \bib{B91}, Lemma 4.8;
let $\cL$ be the (finite) set of $\gL \in \fbkdc$ such that
the associated character of $\D(M/\HM)$ occurs.
We may assume that $\gf$ is a joint eigenfunction for $\D(M/\HM)$,
with eigenvalue character given by $\gL\in\cL$.
It follows that
$$ (D\gf)|_{\Aq}= \gg_{\iM_1}(D\col\gL+\dotvar)(\gf|_{\Aq})$$
for $D\in\DMoneH\simeq\D(M/\HM)\otimes S(\faq)$.
Here $\gg_{\iM_1}$ denotes the Harish-Chandra isomorphism 
from $\DMoneH$ into $S(\fb),$ defined as in \bib{BSmc}, above eq.~(2.11),
and $\gg_{\iM_1}(D\col\gL+\dotvar)\in S(\faq)$ is considered as a
differential operator on $\Aq$. Combining this identity with the 
assumption on $\gf$, the identity
$\gamma_{\iM_1}\after\mu_P=\gg$, and the surjectivity of
$\gamma:\DGH \to S(\fb)^{W(\fb)},$
it follows that
$$u(\gL+\dotvar)(\gf|_{\Aq})=u(\nu)\gf|_{\Aq}$$
for all $u\in S(\fb)^{W(\fb)}.$ Let $\tilde\gf\in C^\infty(\fb)$ be
defined by $\tilde\gf(X+Y)=e^{\gL(X)}\gf(\exp Y)$ for
$X\in\fbk$, $Y\in\faq$, then $u\tilde\gf=u(\nu)\tilde\gf$.
This implies that $\tilde\gf$ is a linear combination of exponential 
polynomials of the form $p\, e^{w\nu}$, where $p\in P(\fb)$ and 
$w\in W(\fb)$, see \bib{Helgason}, Thm.\ III.3.13. However, from the 
definition of $\tilde\gf$ it is readily seen that $w$ only contributes
if $w\nu|_{\fbk}=\gL$.
\qed

{\bf Proof of Lemma \refer{l: restriction on leading exponents}:}
We define the 
$\tauM$-spherical function 
$\gf: M_{1}/M_1\cap vHv^{-1} \simeq M/M\cap vHv^{-1} 
\times \Aq \to \Vtau$ by 
$$ 
\gf(ma) = a^{\rho_P + \xi} q_\xi(P,v \asmid f)(\log a, m).
$$
Then it follows from the equation $Df = \gg(D\col \nu) f$
and Lemma \refer{l: radial component applied to expansion} (c) applied
to $D - \gg(D\col \nu)$ in place of $D,$  
 that  
$$ 
\mu_{P,v}(D) \gf = \gg(D\col \nu) \gf.
$$ 
Since $\gf$ is $\tau$-spherical and non-zero, its restriction to $\Aq$ 
does not vanish. 

Let first $v=e$, and let $\cL$ be as in 
Lemma \refer{l: eigenfunctions on M1}.
It then follows immediately from that lemma
that there exists $w\in W(\fb)$ such that $w\nu|_{\fbk}\in \cL$
and $w\nu|_{\faq}=\xi+\rho_P$.

For general $v \in \NKaq$ we also obtain the result from 
Lemma \refer{l: eigenfunctions on M1}, by applying it to the function
$\gf^v:=\rho_{\tau,v^{-1}}\gf$. Indeed, it follows from the definition of 
$\mu_{P,v}$ that $\gf^v$
satisfies the assumption of the lemma. Hence there exists $w\in W(\fb)$
such that $w\nu|_{\fbk}\in\cL$ and $w\nu|_{\faq}=v^{-1}(\xi+\rho_P)$.
Let $v'\in W_0(\fb)$ be such that $v'Y=vY$ for all $Y\in\faq$,
then $\nu\in (v'w)^{-1}(v'\cL+\xi+\rho_P).$ 
\qed

We will also need a result on leading coefficients along non-minimal parabolic subgroups.

\begin{lemma}
\naam{l: D finiteness of leading coefficient}
Let $f \in C^\ep(\spXp\col \tau)$ be a $\DX$-finite function. Let $Q\in \allparabs,$ $v\in \NKaq$ 
and assume that $\xi \in \ExpL(Q,v\asmid f).$ 
Then the function $\gf: \spXoneQvp \to \Vtau$ defined
by 
$$
\gf(ma) = a^{\xi+\rho_Q} 
q_\xi(Q,v \asmid f, \log a, m)\qquad (m \in \spXQvp,\;a \in \AQq),
$$ 
is $\D(\spXoneQv)$-finite. 
\end{lemma}

\proof
Let $I$ be the annihilator of $f$ in the algebra $\DX.$ 
Then it follows from 
Lemma \refer{l: radial component applied to expansion} (c) 
that $\mu_{Q,v}(D) \gf = 0$ for all $D \in I.$ 
The algebra $\D(\spXoneQv)$ is a finite 
module over the image of the homorphism $\mu_{Q,v};$ see 
\bib{B91}, p.\ 342, and apply conjugation by $v$.
Hence $\mu_{Q,v}(I)$ generates a cofinite ideal in 
$\D(\spXoneQv).$ This establishes the result.
\qed

We end this section with a result that limits the asymptotic exponents occurring 
in discrete series representations to a countable set. 
Later 
we will apply this result to exclude the possibility 
of a `continuum of discrete series'
(see the proof of Lemma \refer{l: second step vanishing thm}).

 To
formulate the result we need 
to define asymptotic exponents for a $K$-finite rather than a $\tau$-spherical function.
We denote by $\dK$ the collection of equivalence classes of
irreducible continuous representations of $K.$ 
If $\types\subset \widehat K$ is a finite subset, then by $\Ci(\spXp)_\types$ we denote 
the space of smooth $K$-finite functions in $\Ci(\spXp)$ all of whose $K$-types
belong to $\types.$ By $\Vtypes:=C(K)_\types$ we denote the space of left $K$-finite continuous functions 
on $K$ all of whose left $K$-types belong to $\types.$ Moreover, by $\tau_\types$ we
denote the restriction of the right regular representation to $\Vtypes.$ 
If $f \in \Ci(\spXp)_\types,$ then the function $\spher_\types(f): \spX \to \Vtypes,$ 
defined by $\spher_\types(f)(x)(k) = f(kx)$ for $x \in \spXp, k\in K$ 
belongs to $\Ci(\spXp \col \tau_\types).$ 
The map $\spher: =\spher_\types$ is a topological linear isomorphism from 
$\Ci(\spXp)_\types $ onto $\Ci(\spXp \col \tau_\types),$ intertwining the $\DX$-actions
on these spaces. Moreover, $\spher$ maps the closed subspace $\Ci(\spX)_\types$ of globally
defined smooth functions onto the similar subspace $\Ci(\spX \col \tau_\types).$ 
We denote by $\Cep(\spXp)_\types$ the preimage of $\Cep(\spXp\col \tau_\types)$ under $\spher.$ 
It follows from Lemma \refer{l: DX finite in exppol}
that $\DX$-finite funcitons in $\Ci(\spXp)_\types$ 
belong to $\Cep(\spXp)_\types.$ 
Let $f \in \Cep(\spXp)_\types;$ then for $P \in \allparabs$ and $v \in \NKaq$ we 
define the set of exponents of $f$ along $(P,v)$ by 
$$
\Exp(P,v\asmid f):= \Exp(P, v\asmid \spher(f)).
$$
Note that this collection is the union for $k \in K$ and $m\in \spXPvp$  of the collections 
of exponents
occurring in the $\DP$-exponential polynomial expansions of $a \mapsto f(kamv).$ 

Let $\cC(\spX)$ denote the space of Schwartz functions on $\spX$, see 
\bib{BSmc}, Section 6, and let $\cA_2(\spX)_K$ denote the space of $K$-finite 
and $\DGH$-finite functions $f\in\cC(\spX).$ These functions are real
analytic and belong to $L^2(\spX)$, cf.\ \bib{B87}, Thm.\ 7.3.

\begin{lemma}
\naam{l: limitation on exponents of a Schwartz function}
Assume that the center of $G$ is compact. Then 
$$\{\xi\in\Exp(P,v\asmid f) \mid
{P \in \minparabs,v \in \NKaq,f \in \cA_2(\spX)_K}\}$$
is a countable subset of $\faqdc$.
\end{lemma}

\proof
Let $\discdual$ denote the set of equivalence classes of discrete series 
representations of the symmetric space $\spX.$  This set is 
countable, since $L^2(\spX)$ is a separable Hilbert space.
Given $\omega \in \discdual$ we denote by $L^2(\spX)_\omega$ the collection of 
functions $f \in L^2(\spX)$ whose closed $G$-span in $L^2(\spX)$ is equivalent 
to a finite direct sum of copies of $\omega.$ 
Let $\dK$ denote the countable set of equivalence classes of irreducible 
representations of $K.$ Given $\omega \in \discdual$ and $\gd \in \dK,$ we 
denote by $L^2(\spX)_{\omega, \gd}$ the collection of $K$-finite elements 
of type $\gd$ in $L^2(\spX)_\omega.$ It follows from \bib{B87}, Thm.\ 7.3,  
that  $L^2(\spX)_{\omega, \gd}$ is a subspace of $\cA_2(\spX)_K,$ and from 
\bib{B87b}, Lemma 3.9, that this subspace is finite dimensional.
On the other hand, let $f \in \cA_2(\spX)_K$, and let $V\subset L^2(\spX)$
denote the closed $G$-span of $f$. It follows from 
\bib{B87b}, Lemma 3.9, that $V$ is admissible. Since $V$ is finitely generated,
it must then be a finite direct sum of irreducible
representations. This implies that 
$f$ belongs to a finite direct sum of  spaces $L^2(\spX)_{\omega, \gd}.$ 
{}From the above
we conclude that $\cA_2(\spX)_K$ equals the following countable
algebraic direct sum:
\begin{equation}
\naam{e: cF as a direct sum} 
\cA_2(\spX)_K = \oplus_{\omega \in \discdual,\,\gd \in \dK }\;\;\; 
L^2(\spX)_{\omega, \gd}.
\end{equation}
Let $\omega \in \discdual$ and $\gd \in \dK.$ 
Then it follows from Lemma \refer{l: DX finite in exppol}
and the finite 
dimensionality of $L^2(\spX)_{\omega, \gd}$ that there exists a countable 
subset $\cE_{\omega, \gd} \subset \faqdc$ such that 
$$
\Exp(P, v\asmid f) \subset \cE_{\omega, \gd}
$$ 
for all $f \in L^2(\spX)_{\omega, \gd}, P \in \minparabs, v \in \NKaq.$ 
Combining this observation with (\refer{e: cF as a direct sum}), we obtain 
the desired result.
\qed

\section{Separation of exponents}
\naam{s: separation}
Let $Q\in\allparabs$.
In the next section we shall consider functions
$f_\gl\in\Cep(\spXp\col \tau)$, with parameter $\gl\in\faQqdc$, whose
exponents along $P\in\minparabs$ lie in sets of the form
$W\gl+S-\N\Delta(P)$, where $S\subset\faqdc$ is a finite set.
In general, given $\xi\in W\gl+S-\N\Delta(P)$, the elements
$s\in W/W_Q$ and $\eta\in S-\N\Delta(P)$, such that $\xi=s\gl+\eta$,
are not unique. In the present section we define a condition on $\gl$
that allows this unique determination for all $\xi$. In particular,
the condition is valid for generic $\gl\in\faQqdc$. 
We consider also the case where $P$ is non-minimal.

Let $P,Q \in \allparabs.$ We define the equivalence
relation 
$\simPQ$ on $W$ by 
\begin{equation}
\naam{e: equivalence relation PQ}
s \simPQ t \;\; \iff\;\; \forall \gl \in \faQqd:\;\; s\gl|_{\faPq} =  
t\gl|_{\faPq}.
\end{equation}
The associated quotient is denoted by $\WPQ.$ 
We note that the classes in $\WPQ$ are left
 $W_P$- and right $W_Q$-invariant.
Thus, $\WPQ$ may also be viewed as a quotient of $W_P\backslash W/W_Q.$ 

If $s,t \in W$ then one readily sees that $s \simPQ t \iff s^{-1} \simQP t^{-1}.$ 
Hence the anti-automorphism $s \mapsto s^{-1}$ of $W$ factors to
a bijection from $\WPQ$ onto $\WQP,$ which we denote by 
$\gs \mapsto \gs^{-1}.$

If $s \in W$ and $\gl \in \faQqdc,$ then 
the restriction $s\gl|_{\faPq}$ depends on $s$ through its class 
$[s]$ in $\WPQ.$ We therefore agree to write
$$
[s] \gl |_{\faPq} := s\gl |_{\faPq}.
$$

\begin{defi}
\naam{d: a zero sets newer}
For $S \subset \faqdc$ a finite subset, we define
$\faQqdzero(P,S)$ to be the subset of 
$\faQqdc$ consisting of elements $\gl$ such that,
for all $s_1,s_2 \in W,$ 
$$
(s_1\gl -  s_2\gl)|_{\faPq} \in  [S + (-S)]|_{\faPq} + \Z \DrP
 \;\;\;\implies \;\;\;s_1 \simPQ s_2.
$$
\end{defi}

\begin{lemma}
\naam{l: exponents disjoint}
Let $S\subset \faqdc$ be finite. Then, for $\gl \in \faQqdc,$
$$
W\gl|_{\faPq} + (S - \N \DP )|_{\faPq} = \bigcup_{\gs \in \WPQ } \; 
\left(\gs\gl|_{\faPq} + (S - \N \DP )|_{\faPq}\right).
$$
Moreover, the union is disjoint if and only if $\gl \in \faQqdzero(P,S).$
\end{lemma}

\proof
Straightforward.
\qed

\begin{lemma}
\naam{l: azero is full}
Let $Q,P \in \allparabs,$ and let $S$ be a finite subset of 
$\faqdc.$ Then $\faQqdzero(P,S)$ equals the complement of 
the union of a locally finite collection of proper 
affine subspaces in $\faQqdc.$ 
\end{lemma}

\proof
Let
$p: \faqdc \to \faPqdc$ denote  the map induced by restriction
to $\faPq.$ Let $\Pi$ be the complement of the diagonal in the set $\WPQ \times \WPQ.$ 
Then for every $\gs =(\gs_1, \gs_2) \in \Pi$ and every $\eta \in \faPqdc$ we write 
$\cA_{\gs, \eta} =\{\gl \in \faQqdc\mid p(\gs_1\gl - \gs_2\gl) = \eta\}.$ 
Note that $\cA_{\gs, 0}$ is a proper affine subspace of $\faQqdc.$ 
If $\gl \in \cA_{\gs, \eta}$ then $\cA_{\gs, \eta}$ equals
$ \gl + \cA_{\gs, 0};$ hence the set $\cA_{\gs, \eta}$ is either empty or a proper 
affine subspace. 

Let $\cA$ be the collection of subsets of the form 
$\cA_{\gs, \xi},$ for $\gs\in \Pi$ 
and $\xi  \in p(S + (-S)) + \Z\DrP.$ 
Then $\faQqdzero(P,S)$ equals the 
complement of $\cup\cA$ in $\faQqdc.$ Thus, it remains to show that 
the collection $\cA$ is locally finite.

Let $\cC$ be a compact subset of $\faQqdc$ and let $X$ be the collection of 
$\xi \in p(S + (-S)) + \Z\DrP$ such that $\cC \cap \cA_{\gs, \xi} \not= \emptyset$ 
for some $\gs\in \Pi.$ Then it suffices to show that $X$ is finite.

Let $\cC'\subset \faPqdc$ be the image of $\Pi \times \cC$ under the map 
$(\gs, \gl)\mapsto p(\gs_1 \gl - \gs_2 \gl).$ Then $X$ equals the intersection of
$\cC'$ with $p(S + (-S)) + \Z \DrP.$ The latter set is discrete since $S$ is finite, 
whereas the elements of $\DrP$ 
are linearly independent. It follows that $X$ is finite. 
\qed

\begin{rem}
In particular, it follows from the above lemma that $\faQqdczero(P,S)$ is 
a full open subset of $\faQqdc;$ see Appendix B for the notion of full.
\end{rem}

\begin{lemma}
\naam{l: WPQ as cosets}
Let $Q,P \in \allparabs.$ If either $\faQq$ or $\faPq$ has codimension
at most $1$ in $\faq,$ then the natural projection 
$W_P\bs W / W_Q \to \WPQ$ is a bijection.
\end{lemma}

\proof
It suffices to prove injectivity of the map.
Since $s \mapsto s^{-1}$ induces a bijection from $\WPQ$ onto $\WQP,$ 
it suffices to prove this when  $\faPq$ has codimension
at most $1.$ We assume the latter to hold.

For $s \in W,$ let $[s]$ denote its canonical image in $\WPQ.$ 
Assume that $s,t \in W$ and that $[s] = [t].$ Then for every $\gl \in \faQqd$ 
we have $s\gl = t\gl$ on $\faPq.$ If $\faPq = \faq,$ this implies
that $s =t$ on $\faQqd,$ hence $s W_Q = t W_Q,$ and since $W_P$ is trivial 
in this case, the proof is finished.
Thus, we may as well assume that $\faPq$ has codimension $1.$ Then there
exists a root $\ga \in \Sigma$ such that $\faPq = \ker \ga .$ 
Note that $W_P = \{1, s_\ga\}.$ 
For every $\gl \in \faQqd$ the Weyl group images $s\gl$ and $t\gl$ 
have equal length in $\faqd$ and equal image under the orthogonal 
projection to $\ga^\perp.$ Hence there exists a constant $\eta \in \{0,1\}$ 
such that
$\inp{s\gl}{\ga} = (-1)^\eta\inp{t\gl}{\ga}$ for all $\gl \in \faQqd.$
It follows that $s\gl = s_\ga^{\eta} t \gl$ for all $\gl \in \faQqd;$ 
hence $sW_Q= s_\ga^\eta  t W_Q,$ from which it follows in turn 
that $s$ and $t$ have the same image in $W_P\bs W / W_Q.$ 
\qed

In particular, if  $P$ is minimal, then 
the natural map 
$W/W_Q \to \WPQ$ is a bijection; we shall use it to identify 
the sets involved.

\section{Analytic families of spherical functions}
\naam{s: analytic families}
In this section we assume that $(\tau,\Vtau)$ is 
a finite dimensional continuous representation of $K.$ 
Let $Q \in \allparabs$ and let $Y$ be a finite subset of $\staQqdc,$ see 
(\refer{e: char staQq}).

In the following definition we introduce a space of analytic families of 
$\tau$-spherical functions that will play a crucial role in the rest of this 
paper.
\begin{defi}
\naam{d: anfamQY newer}
Let $Q,Y$ be as above and let $\Omega \subset \faQqdc$ be an open subset.
We define 
\begin{equation}
\naam{e: anfamQY new}
\anfamQY
\end{equation}
to be the space of $C^\infty$-functions 
$f: \Omega \times \spXp \to \Vtau$ 
satisfying the following conditions.
\begin{enumerate}
\itema
For every $\gl \in \Omega$ the function $f_\gl: x \mapsto f(\gl, x)$ 
belongs to $C^\infty(\spXp\col \tau).$
\itemb
There exists a constant $k \in \N,$ 
and, for every $P \in \minparabs$ and $v\in \NKaq,$  
a collection of functions
$q_{s, \xi}(P,v\asmid f)\in P_k(\faq)\otimes \cO(\Omega, \Ci(\Xov\col \tauM)),$ 
for $s \in W/W_Q$ and
$\xi \in - sW_Q Y + \N\DP,$ with the following property.
For all $\gl \in \Omega,$ $m \in \spXzerov$ and  $a\in \Aqp(P),$
\begin{equation}
\naam{e: expansion f from anfamQY new}
f_\gl(mav) = \sum_{s \in W/W_Q} a^{s\gl - \rho_P}
\sum_{\xi \in -sW_Q Y + \N \DP} a^{-\xi} q_{s, \xi}(P,v\asmid f, \log a)(\gl, m),
\qquad
\end{equation}
where the $\DP$-exponential polynomial series with coefficients in $\Vtau$ 
is neatly convergent on $\Aqp(P).$  
\itemc
For every $P\in \minparabs, v\in \NKaq$ and $s \in W/W_Q,$ the series 
$$ 
\sum_{\xi \in - s W_Q  Y + \N\DP}
a^{-\xi} q_{s, \xi}(P,v\asmid f, \log a)
$$
converges neatly on $\Aqp(P)$ as a $\DP$-exponential polynomial
series with coefficients in
 $\cO(\Omega, \Ci(\Xov\col\tauM)).$ 
\end{enumerate}
If $f\in C^\ep_{Q,Y}(\spXp\col \tau\col \Omega),$ we 
define the asymptotic degree of $f,$ denoted $\dega f,$ to be the smallest number
$k \in \N$ for which the above condition (b) is fulfilled.
\end{defi}

\begin{rem}
\naam{r: on defi anfamQY new}
We note that the space (\refer{e: anfamQY new}) depends on $Q$ through its
$\gs$-split component $\AQq.$ 
Moreover, from Lemma 
\refer{l: transformation of coeffs}
we see that in the above definition it suffices  
to require (b) and (c) 
for a fixed given $P \in \minparabs$ and for each $v$ 
in a  given set $\cW\subset \NKaq$ of representatives
for $W/\WKH.$ Alternatively, by the same lemma it  suffices to require (b) and (c)
for a fixed given $v\in \NKaq$ and arbitrary $P\in \minparabs.$
\end{rem}

\begin{lemma}
\naam{l: pointwise expansion of family} 
Let $f \in \anfamQY.$ Then $f_\gl \in \Cep(\spXp\col \tau)$
and 
\begin{equation}
\naam{e: exponents family}
\Exp(P,v\asmid f_\gl) \subset W(\gl + Y) - \rho_P - \N \DP
\end{equation}
for all $\gl \in \Omega,$ $P \in \minparabs,$ and
$v \in \NKaq.$ Moreover, let $\Omega':= \Omega \cap \faQqdczero(P, WY)$
(see Definition \refer{d: a zero sets newer}).
Then $\Omega'$ is open dense in $\Omega$ and
\begin{equation}
\naam{e: q-functions unique}
q_{s, \xi}(P,v\asmid f, X, \gl) = q_{s\gl -\rho_P - \xi}(P,v\asmid f_\gl, X) 
\end{equation}
for every $s \in W/W_Q$, $\xi \in - sW_Q Y + \N\DP,$ 
$X \in \faq$ and $\gl \in \Omega'.$ 
In particular, the functions
$q_{s, \xi}(P,v\asmid f)$ are uniquely determined.
\end{lemma}

\proof The first statement and (\refer{e: exponents family}) follow 
immediately from condition (b) in the above definition. The set
$\Omega'$ is open dense in $\Omega$ by Lemma \refer{l: azero is full},
and it follows from Lemmas \refer{l: exponents disjoint} and
\refer{l: WPQ as cosets} that if $\gl \in \Omega'$ then the sets
$s(\gl + W_Q Y) - \rho_P - \N\DP$, $s \in W/W_Q,$ are mutually disjoint.
Then (\refer{e: q-functions unique}) holds by uniqueness of asymptotics.
\qed

The following result shows that an element of $\anfamQY$ may 
be viewed as an analytic family of spherical functions.

\begin{lemma}
\naam{l: the family is analytic}
Let $f \in \anfamQY.$ Then $\gl \mapsto f_\gl$ is a holomorphic 
function on $\Omega$ with values in $C^\infty(\spXp\col\tau).$
\end{lemma}

\proof
Let $\cW\subset  \NKaq$ be a complete set of representatives for $W/\WKH.$ 
Note that for $v \in \cW$ the $\Vtau^{\KM \cap \vH}$-valued 
function $\TdownPv f_\gl$ 
on $\Aqp(P)$ is given by the series on the right-hand side of 
(\refer{e: expansion f from anfamQY new}) with $m=e.$ 
It follows from condition (c) of Definition \refer{d: anfamQY newer}
that $a \mapsto \TdownPv f_\gl (a)$ defines a smooth function on
$\Aqp(P)$ with values in $\cO(\Omega) \otimes \Vtau^{\KM \cap \vH}.$ 
According to
Appendix A, the function $\gl \mapsto \TdownPv f_\gl(\dotvar)$ is 
a holomorphic function on $\Omega$ 
with values in $\Ci(\Aqp(P), \Vtau^{\KM\cap \vH}).$ 
Hence $\gl \mapsto \TdownPcW (f_\gl)$ is a holomorphic function
on $\Omega$ with values in $\Ci(\Aqp(P), \oplus_{v \in \cW} \Vtau^{\KM \cap \vH}).$  
The conclusion of the lemma now follows by application of the isomorphism 
(\refer{e: the iso T down P cW}).
\qed

If $\Omega', \Omega$ are open subsets of $\faqdc$ with $\Omega' \subset \Omega,$ 
then restriction from $\Omega \times \spXp$ to $\Omega' \times \spXp$ 
obviously induces a linear map 
\begin{equation}
\naam{e: restriction map presheaf}
\rho^{\Omega}_{\Omega'}:  
C^\ep_{Q,Y}(\spXp\col\tau \col \Omega) \to C^\ep_{Q,Y}(\spXp\col\tau \col \Omega').
\end{equation}
Accordingly,  the assignment 
\begin{equation}
\naam{e: presheaf CepQY}
\Omega \mapsto C^\ep_{Q,Y}(\spXp\col\tau \col \Omega),
\end{equation}
defines a presheaf of complex linear spaces on $\faQqdc.$ 
Here we agree that (\refer{e: presheaf CepQY}) assigns the 
trivial space to the empty set.

The following lemma will be useful at a later stage.
\begin{lemma}
\naam{l: CepQY is sheaf}
Let $Q \in \allparabs$ and $Y \subset\staQqdc$ a finite subset.
\begin{enumerate}
\itema
If $\Omega', \Omega$ are open subsets of $\faQqdc$ with $\Omega' \neq \emptyset,$
$\Omega$ connected and $\Omega' \subset \Omega,$ then the restriction map
(\refer{e: restriction map presheaf}) is injective.
Moreover, $\dega( \rho^{\Omega}_{\Omega'} f) = \dega ( f)$ for all $f \in 
C^\ep_{Q,Y}(\spXp\col\tau \col \Omega).$ 
\itemb
The presheaf (\refer{e: presheaf CepQY}) is a sheaf.
\end{enumerate}
\end{lemma}

\proof
The injectivity of the restriction map follows by analytic continuation, 
in view of Lemma \refer{l: the family is analytic}. 
Let $f'=\rho^{\Omega}_{\Omega'} f$.
Let $P \in \minparabs,$ $v \in \NKaq,$ $s \in W/W_Q$ and $\xi \in -sW_Q \gl + \N\DP$
then it follows from (\refer{e: q-functions unique}) that
\begin{equation}
\naam{e: restriction of q} 
q_{s, \xi}(P, v\asmid f', \dotvar, \gl) = q_{s, \xi}(P, v\asmid 
f, \dotvar, \gl)
\end{equation}
for $\gl$ in a dense open subset of $\Omega'$, hence for all $\gl\in\Omega'$.
In particular this implies that the polynomial degree of the function
on the left-hand side of the equation is bounded by $\dega (f);$ hence
$\dega(f') \leq \dega (f).$ 
To prove the converse inequality, we note that the polynomial on the left-hand 
side of (\refer{e: restriction of q}) 
is of degree at most $k':= \dega (f') $ by the definition
of the latter number. Since $\Omega$ is connected, it follows
by analytic continuation that $\deg q_{s, \xi}(P, v\asmid f, \dotvar, \gl) 
\leq k'$ 
for all $\gl \in \Omega.$ Since this holds for all $P,v,\gs, \xi,$ 
it follows that $\dega (f ) \leq k'$ and we obtain (a).

Assertion (b) is equivalent with the assertion that the presheaf satisfies
the localization property (see \bib{mybook}, p.\ 9). This is
established in a straightforward manner, by using (a).
\qed

 We shall now  discuss the action of invariant differential
operators on families.
If $f$ is a family in $\CepQY(\spXp\col \tau\col \Omega),$ and $D \in \DX,$ then 
we define the family 
$Df: \Omega \times \spXp \to \Vtau$ 
by 
\begin{equation}
\naam{e: action D on family}
(Df)_\gl = D(f_\gl),\qquad (\gl \in \Omega).
\end{equation}
\begin{prop}
\naam{p: D on families}
Let $f \in \anfamQY.$ Then, for every $D \in \DX,$ 
the family $Df$ belongs to $\anfamQY;$ moreover, $\dega (Df) \leq \dega (f).$ 
\end{prop}

\proof
Let $D \in \DX.$ Then $g = Df$ 
is a smooth function $\Omega \times \spXp \to \Vtau;$ 
moreover, for $\gl \in \Omega$ the function $(Dg)_\gl = Df_\gl$ 
is $\tau$-spherical. Thus, $g$ satisfies 
condition (a) of Definition \refer{d: anfamQY newer} and it remains
to establish properties (b) and (c). In view of Remark \refer{r: on defi anfamQY new} 
it suffices to do this for $v =e$ and arbitrary $P \in \minparabs.$ 
Let $k := \dega f.$ 

It follows from condition (b) of Definition \refer{d: anfamQY newer}
that, for $\gl \in \Omega,$ the function $f_\gl$ belongs to 
$\Cep(\spXp\col \tau);$ 
moreover, its $(P,e)$-expansion is given by 
\begin{equation}
\naam{e: Pe expansion of f gl} 
f_\gl(ma) = \sum_{s \in W/W_Q} a^{s\gl - \rho_P}
\sum_{\xi \in -sW_Q Y + \N \gD(P)} a^{-\xi} q_{s,\xi}(P,e\asmid f, \log a)(\gl, m),
\end{equation}
for $a \in \Aqp(P)$  and $m\in M.$ 
Let $u:= \mu_P'(D) + \uplus$ be the element of $\cD_{1P}$ associated with $D$ as in 
Proposition  \refer{p: radial deco with MQ} with $P$ in place of $Q.$ 
In view of
Corollary \refer{c: cor on ringQ} its expansion $\ep(u),$ defined as in 
(\refer{e: the embedding ep}),
is the sum, as $i$ ranges over a finite index set $I,$ of series of the form 
$$
\ep(u)_i = \sum_{\nu \in \N\DP} a^{-\nu} \gf_{i,\nu}\otimes S_{i,\nu} \otimes u_{i,\nu} \otimes 
v_{i,\nu}.
$$
Here $\gf_{i,\nu} \in C^\infty(M_\gs),$ $S_{i,\nu} \in \End(\Vtau),$ 
$u_{i,\nu} \in U(\fm_{\gs})$ 
and $v_{i,\nu} \in U(\faq),$ 
and $\deg(u_{i,\nu}) + \deg(v_{i,\nu}) \leq d:=\deg(D)$ 
for 
all $i,\nu.$ 
By Lemma \refer{l: radial component applied to expansion},
the function 
$g_\gl$ belongs to $\Cep(\spXp\col \tau),$ for $\gl \in \Omega,$ and its 
$(P,e)$ expansion results from (\refer{e: Pe expansion of f gl}) 
by the formal application of the element $\ep(u).$ 
This 
gives, for $\gl \in \Omega,$ $m \in M$ and $a \in \Aqp(P),$
the neatly converging exponential polynomial expansion
$$
g_\gl(m a) = \sum_{s \in W/W_Q} a^{s\gl - \rho_P} \sum_{\eta \in -sW_Q Y + \N \DP} 
a^{-\eta} \tilde q_{s, \eta}(\log a)(\gl, m),
$$ 
where $\tilde q_{s, \eta}$ is given by the following finite sum
$$ 
\tilde q_{s, \eta}(X)(\gl, m) 
 := \sum_{i\in I} 
\sum_{{\nu \in \N \DP \atop {\xi \in -sW_Q Y + \N \DP}}\atop \nu + \xi = \eta}
\gf_{i, \nu}(m) S_{i, \nu}[\;
q_{s, \xi}(P,e\asmid f, X; T_{s\gl - \rho_P -\xi}(v_{i, \nu}), \gl, m ; u_{i,\nu})\;],
$$
for $\gl \in \Omega, $ $X \in \faq$ and $m \in M.$ 
Here
we have used Harish-Chandra's convention to indicate by a semicolon on the 
left or right-hand side of a Lie group variable the differentiation
on the corresponding side, with respect to that variable, 
by elements of the appropriate universal enveloping algebra. 
Moreover, 
given $\gg \in \faqdc$ we have denoted by $T_\gg$ the automorphism of $U(\faq)$ 
determined by $T_\gg(X) = X + \gg(X)$ for $X \in \faq.$ 

{}From the above formula  it readily follows  
that $\tilde q_{s, \eta}(X,\gl)$ is a smooth function of $(X,\gl)$ with 
values in $\Ci(M, \Vtau);$ moreover, it is polynomial in $X$ of degree at most
$k$ and holomorphic in $\gl\in \Omega.$ This establishes condition (b)
of Definition \refer{d: anfamQY newer} with $v=e,$ arbitrary $P\in \minparabs,$ 
and with 
$$
q_{s, \eta}(P,e\asmid g) = \tilde q_{s, \eta},
\qquad (s \in W/W_Q, \;\eta \in -sW_QY + \N \gD(P)).
$$  

For condition (c) we note that the series
\begin{equation} 
\naam{e: series for g}
\sum_{\eta \in -sW_Q Y + \N \DP} a^{-\eta} q_{s, \eta}(P,e\asmid g, \log a)
\end{equation}
arises from the series
\begin{equation} 
\naam{e: similar series for f}
\sum_{\xi \in -sW_Q Y + \N \DP} a^{-\xi} q_{s, \xi}(P,e\asmid f, \log a)
\end{equation}
by the formal application of $\ep(u)$ 
conjugated with multiplication by $a^{-s\gl + \rho_P}.$ {}From this we see that 
(\refer{e: series for g}) arises from (\refer{e: similar series for f})
by the formal application of the series
$$ 
\sum_{\nu \in \N \DP} a^{-\nu} \sum_{i \in I}  \gf_{i, \nu} \otimes S_{i, \nu} \otimes u_{i, \nu} 
\otimes v_{i, \nu}(\gl),
$$ 
with $v_{i, \nu}(\gl) = T_{s\gl - \rho_P} (v_{i, \nu}).$ 
We now observe that $\gl \mapsto T_{s\gl - \rho_P}|_{U_d(\faQq)}$ is a polynomial 
$\End(U_d(\faQq))$-valued function, of degree at most $d.$ 
Hence there exists a finite set $J$ and elements
$p_j \in P_d(\faQqd)$ and $T_j \in \End(U_d(\faQq)),$ for $j \in J,$ 
such that 
$$
T_{s\gl - \rho_P}|_{U_d(\faQq)} = \sum_{j \in J} p_j(\gl) T_j.
$$ 
Let $B_{i,\nu, j}$ be the continuous endomorphism of $\cO(\Omega, \Ci(M_{\gs}, \Vtau))$ defined
by 
$$ 
B_{i, \nu, j} \psi(\gl)(m) = p_j(\gl) \gf_{i, \nu}(m) S_{i,\nu}
[\psi(\gl)(m; u_{i,\nu})].
$$ 
Then the series (\refer{e: series for g}) arises from 
the formal application of 
the  series 
$$ 
\sum_{\nu \in \N \DP} a^{-\nu} \sum_{i\in I\atop j \in J} B_{i,\nu,j} \otimes T_j(v_{i,\nu}) 
$$ 
with coefficients in 
$\End(\cO(\Omega, C^\infty(M_{\gs}, \Vtau)))\otimes U(\faQq)$
to (\refer{e: similar series for f}), viewed as a 
series with coefficients in $\cO(\Omega, C^\infty(M_{\gs}, \Vtau)).$  
It follows from Lemmas \refer{l: formal application of Ufa} and 
\refer{l: formal application of hom valued series} that the resulting series 
is neatly convergent as a series 
on $\Aqp(P)$ with coefficients in $\cO(\Omega, C^\infty(M_{\gs}, \Vtau)).$
This establishes (c) with $v =e$ and arbitrary $P \in \minparabs.$ 
\qed

We will now describe the asymptotic behavior along walls for a family. 
If $P,Q\in\allparabs$
and $\gs \in \WPQ$ (see (\refer{e: equivalence relation PQ})),
then for every subset $Y \subset \faqdc$  we put 
\begin{equation}
\naam{e: defi gs cdot Y}
\gs \cdot Y :=
\{ s \eta|_{\faPq} \mid s \in W, [s] = \gs, \; \eta \in Y \}.
\end{equation}

\begin{thm}{\rm (Behavior along the walls).\ }
\naam{t: behavior along the walls for families new}
Let $Q\in \allparabs,$ $\Omega \subset \faQqdc$ a non-empty open subset and 
$Y \subset \staQqdc$  
a finite subset.  
Let $f \in \anfamQY$ 
and let $k= \dega(f).$ 

Let $P \in \allparabs$ and $v \in \NKaq.$ 
Then $\Exp(P,v\asmid f_\gl) \subset W(\gl + Y)|_{\faPq} -\rho_P- \N \DrP$ 
for every $\gl \in \Omega.$ Moreover, 
there exist unique 
functions
$$
q_{\gs, \xi}(P,v\asmid f) \in P_k(\faPq) \otimes \cO(\Omega, \Ci(\spXPvp\col \tauP)),
$$ 
for $\gs \in \WPQ$ and $\xi \in - \gs\cdot Y + \N \DrP,$ 
with the following property. For all
$\gl \in \Omega,$  $m \in \spXPvp$ and $a \in \APqp(R_{P,v}(m)^{-1}),$ 
\begin{equation}
\naam{e: series for fam along P}
f_\gl(mav) = \sum_{\gs \in \WPQ} a^{\gs\gl - \rho_P} 
\sum_{\xi \in - \gs\cdot Y + \N \DrP}
a^{-\xi}\, q_{\gs, \xi} (P, v\asmid f ,\log a)( \gl, m),
\end{equation} 
where the $\DrP$-exponential polynomial series with coefficients in 
$\Vtau$ is neatly convergent on $\APqp(R_{P,v}(m)^{-1}).$ 
In particular, if $\gl\in\Omega':=\Omega \cap \faQqdczero(P, WY)$
then
\begin{equation}
\naam{e: q of f in gl versus q of fgl new}
q_{\gs, \xi}(P, v\asmid f)(X,\gl) = 
q_{\gs \gl|_{\faPq }- \rho_P -\xi} (P, v\asmid f_\gl, X),
\end{equation}
for
$X\in \faPq.$

Finally, for each $\gs \in \WPQ$ and every
$R >1,$ 
the series
\begin{equation}
\naam{e: series with q gs xi}
\sum_{\xi \in - \gs \cdot Y + \N \DrP}
a^{-\xi} q_{\gs, \xi} (P, v \asmid f,\log a)
\end{equation} 
converges neatly on $\APqp(R^{-1})$ as a $\DrP$-exponential 
polynomial series with coefficients
in $\cO(\Omega, \Ci(\spXPvp[R]\col \tauP)).$ 
\end{thm}

\proof
Let $P \in \allparabs$ and let $v \in \NKaq.$
Fix a  minimal parabolic subgroup
${\Pmin} \in \minparabs,$ contained
in $P.$ 
Fix a set $\cW_{P,v}\subset \NKPaq$ of representatives
for $W_P/ W_P \cap v W_{K\cap H} v^{-1}.$ Then the natural map $\NKaq \to W$ 
induces an embedding  $\cW_{P,v}v \embeds W/\WKH.$ Therefore, we may fix 
a set $\cW \subset \NKaq$ of representatives for $W/\WKH$ containing 
$\cW_{P, v} v.$  
 
Fix $\gl \in \Omega$ for the moment. Then by 
Lemma \refer{l: pointwise expansion of family},
the function $f_\gl$ belongs to $\Cep(\spXp\col \tau),$ and 
$\Exp(P_1, w\asmid f_\gl) \subset W(\gl + Y) - \rho_{P_1} - \N \gD(P_1),$ 
for every $w \in \NKaq.$ 
According to
Theorem \refer{t: transitivity of asymptotics},
for every $u \in \cW_{P,v},$ 
the set $\Exp(P,v\asmid f_\gl)_{{\Pmin}, u}$ is 
contained in $\Exp({\Pmin}, uv\asmid f_\gl)|_{\faPq}.$ Hence, by 
of (\refer{e: Exp Q v f as union over index set cWQv}) with $P$ and $P_1$ 
in place of $Q$ and $P,$ respectively, we infer that 
\begin{eqnarray}
\Exp(P,v\asmid f_\gl) 
&\subset & [W(\gl + Y) - \rho_{{\Pmin}} - \N \gD({\Pmin})]|_{\faPq}
\nonumber\\ 
&=  & W(\gl + Y)|_{\faPq} - \rho_P - \N \DrP.
\naam{e: exponents f gl along P v in this set}
\end{eqnarray}

Notice that (\refer{e: q of f in gl versus q of fgl new}) is a consequence
of (\refer{e: series for fam along P}), by Lemma \refer{l: exponents disjoint}.
Therefore the functions $q_{\gs, \xi}(P,v\asmid f)$ are unique. We will
now establish their existence.

It follows form (\refer{e: exponents f gl along P v in this set})
that the elements of $\Exp(P,v\asmid f_\gl)$ are all of the 
form $\gs \gl|_{\faPq}   - \rho_P - \xi ,$ with $\gs \in \WPQ$ 
and $\xi \in -\gs\cdot Y + \N\DrP.$
Fix such elements $\gs$ and  $\xi.$ 
Then by transitivity of asymptotics, 
cf.\ Theorem \refer{t: transitivity of asymptotics},
we have,
for $u \in \cW_{P,v},$ $X \in \faPq,$ $m \in M$ and   $b \in \stAPqp(\starPmin),$  
that
\begin{equation}
\naam{e: first series for q gs f gl}
q_{\gs \gl|_{\faPq }- \rho_P -\xi} (P, v\asmid f_\gl, X, m bu)
=
\sum_{\zeta \in \Exp(P_1, uv \asmid f_\gl) \atop
\zeta|_{\faPq} = \gs \gl|_{\faPq} - \rho_P - \xi}
b^{\zeta} q_\zeta({\Pmin}, uv \asmid f_\gl, X + \log b, m),
\end{equation}
where the $\gD_P(P_1)$-exponential polynomial series 
in the variable $b$ converges neatly on $\stAPqp(\starPmin).$
It follows from condition (b) in Definition \refer{d: anfamQY newer} that,
for $\zeta \in \Exp(\Pmin, uv \asmid f_\gl),$ 
\def\etaormu{\mu}
\begin{equation}
\naam{e: formula for q xi in series}
q_\zeta({\Pmin}, uv \asmid f_\gl, X + \log b, m) = 
\sum_{{s \in W/W_Q \atop \etaormu \in -sW_QY + \N \gD(P_1)}\atop
      s \gl - \rho_{P_1} - \etaormu  = \zeta }
q_{s, \etaormu} (\Pmin, uv \asmid f, X + \log b)(\gl, m).
\end{equation}
Now assume that $\gl$ is contained in the full (cf.\ Lemma
\refer{l: azero is full}) subset $ \Omega'$ of $\Omega.$ Then, if
$s \in W$ and $\etaormu \in - sW_QY + \N\gD(\Pmin)$ satisfy 
$[s\gl -\rho_{\Pmin} - \etaormu]|_{\faPq} = \gs \gl|_{\faPq } - \rho_P -\xi$,
it follows that $[s] = \gs$ and $\etaormu|_{\faPq} = \xi,$ 
see Lemma \refer{l: exponents disjoint}. Hence, 
combining (\refer{e: first series for q gs f gl}) 
and (\refer{e: formula for q xi in series}) we infer that
for $\gl\in\Omega'$, $u \in \cW_{P,v},$ $X \in \faPq,$ $m \in M$ and 
$b \in \stAPqp(\starPmin),$ 
\begin{eqnarray}
\nonumber
\lefteqn{
q_{\gs \gl|_{\faPq }- \rho_P -\xi} (P, v\asmid f_\gl, X, m bu)=\quad}\\
\naam{e: series for q gs f gl}
&=& \label{e: guiding formula}
\sum_{s \in W/W_Q\atop [s]= \gs} b^{s\gl - \rho_{{\Pmin}}}
 \sum_{\etaormu \in  - sW_Q Y + \N \gD({\Pmin})\atop
\etaormu|_{\faPq} = \xi} 
b^{- \etaormu} q_{s, \etaormu}({\Pmin}, uv \asmid f, X + \log b, \gl)(m).
\end{eqnarray}
It will be seen below that each inner sum in
(\refer{e: guiding formula}) converges neatly, so that the separation
of terms by the outer sum is justified.
This formula will guide us towards 
the definition of the functions $q_{\gs, \xi}(P,v\asmid f).$ 

In the following we assume that $s \in W/W_Q$ and $[s]= \gs.$
For $w \in \cW$ 
we define the function 
$F_{s,w}: \Aqp({\Pmin}) \times \Omega \to \Vtau^{\KM\cap wHw^{-1}}$ by  
$$
F_{s,w}(a,\gl) = 
\sum_{\etaormu \in  - sW_Q Y + \N \gD({\Pmin})} 
a^{- \etaormu} q_{s, \etaormu}({\Pmin}, w \asmid f, \log a, \gl)(e),
$$
for $a \in \Aqp({\Pmin}),$ $\gl \in \Omega.$ 

The representation $\tilde \tau:= 1\otimes \tau$ of $K$ 
on the complete locally convex space $\cO(\Omega) \otimes \Vtau$ is smooth.
We shall apply the results of Section \refer{s: asymp walls}, with $\tilde \tau$ in place
of $\tau.$
The series defining $F_{s,w}$ is a $\gD({\Pmin})$-exponential polynomial series
with coefficients in $\cO(\Omega) \otimes \Vtau.$ By condition (c) 
of Definition \refer{d: anfamQY newer} it converges
neatly on $\Aqp({\Pmin});$ hence
$F_{s,w}$ may be viewed as an element of 
$C^\ep(\Aqp({\Pmin}), [\cO(\Omega) \otimes \Vtau]^{\KM \cap wHw^{-1}}).$ 
In view of the isomorphism (\refer{e: isomorphism of exppol}), 
there exists a unique 
function $F_s \in C^\ep(\spXp \col \tilde \tau)$ such that 
$\TdownPminw(F_s)(a) = F_s(aw) = F_{s,w}(a),$ 
for $w \in \cW$ and $a \in \Aqp({\Pmin}).$ 
{}From the  definition of $F_s$ it follows that
$\Exp({\Pmin},w\asmid F_s) \subset sW_Q Y - \N \gD({\Pmin}),$ for every $w \in \cW.$ 
Moreover, for every $w \in \cW$ and every $\etaormu \in - sW_QY + \N \gD({\Pmin}),$ 
\begin{equation}
\naam{e: q of Fs versus q of f}
q_{-\etaormu}({\Pmin}, w \asmid F_s,X,m)(\gl)  
= q_{s,\etaormu}({\Pmin}, w\asmid f,X)(\gl, m),
\end{equation}
for $X \in \faPq,$ $m \in \spX_{0,w}$ and $\gl \in \Omega.$ 
By transitivity of asymptotics, cf.\ Theorem  \refer{t: transitivity of asymptotics},
applied to $F_s,$ we have that
$\Exp(P, v \asmid F_s)_{{\Pmin}, u} \subset \gs\cdot Y - \N \DrP,$ 
for $u \in \cW_{P,v}.$ 
Moreover, by the same result it follows that,
for $\xi \in - \gs \cdot Y + \N \DrP,$ 
\begin{equation}
\naam{e: series for q minus xi FW}
q_{-\xi}(P, v\asmid F_s)(X, m bu) =
\sum_{\etaormu \in  - sW_Q Y + \N \gD({\Pmin})\atop
\etaormu|_{\faPq} = \xi} 
b^{- \etaormu} q_{-\etaormu}({\Pmin}, uv \asmid F_s, X + \log b,  m ),
\end{equation}
where the series on the right-hand side converges neatly as
a $\gD_P({\Pmin})$-exponential polynomial series in the variable
$b \in \stAPqp(\starPmin),$ with coefficients in
$\Ci(\spX_{0, uv}\col \tilde\tau_\iM).$ 
In particular, the asserted convergence of (\refer{e: guiding formula})
follows.

Substituting (\refer{e: q of Fs versus q of f}) in the right-hand side of 
(\refer{e: series for q gs f gl}) and using (\refer{e: series for q minus xi FW})
we find, for $\gl \in \Omega',$ that
\begin{equation}
\naam{e: q of f gl as s sum}
q_{\gs \gl|_{\faPq }- \rho_P -\xi} (P, v\asmid f_\gl, X, m bu)= 
\sum_{s \in W/W_Q\atop [s] =\gs} 
b^{s\gl- \rho_{{\Pmin}}} q_{-\xi}(P, v\asmid F_s)(X, m bu)(\gl).
\end{equation}
We
are now ready to define the functions $q_{\gs, \xi}(P,v\asmid f).$ 

Let $1$ denote the trivial representation of $K$ in $\C,$ 
and $1_P$ its restriction to $K_P.$
If $s \in W/W_Q,$ we define the function
$\gf_s \in \cO(\faQqdc,\Ci(\spXPvp\col 1_P))$ 
by
\begin{equation}
\naam{e: defi gf s}
\gf_s(\gl, kbu) = b^{s\gl-\rho_{{\Pmin}}},
\end{equation}
for 
$\gl \in \faQqdc,$  $u \in \cW_{P,v},$ $k \in K_P$ and $b \in \stAPqp(\starPmin).$ 
Moreover,
for $\gs \in \WPQ$ and $\xi \in - \gs\cdot Y + \N \DrP$  we 
define the function 
$q_{\gs, \xi}(P, v\asmid f): \faPq \times \Omega \to \Ci(\spXPvp\col\tauP)$ 
by 
\begin{equation}
\naam{e: q of f as sum of q of Fs new}
q_{\gs,\xi}(P,v\asmid f,X,\gl)(m)  = 
\sum_{s \in W/W_Q\atop [s] = \gs} \gf_s(\gl,  m) \,q_{-\xi}(P,v\asmid F_s, X,m)(\gl),
\end{equation}
for $X \in \faPq,$ $\gl \in \Omega$ and $m \in \spXPvp.$ 

If
$1 < R \leq \infty,$ then the locally convex space 
$\Ci(\spXPvp[R], \cO(\Omega) \otimes \Vtau)$ 
is naturally isomorphic with $\cO(\Omega, \Ci(\spXPvp[R], \Vtau)),$ see 
Appendix \refappa{}.
The isomorphism induces in turn a natural isomorphism of locally convex spaces
\begin{equation}
\naam{e: spaces cO Omega Ci}
\Ci(\spXPvp[R] \col \tilde\tauP)) \simeq \cO(\Omega, \Ci(\spXPvp[R] \col \tauP)).
\end{equation}
In particular, for $R=\infty,$ we obtain that
$\Ci(\spXPvp \col \tilde\tauP)$ is naturally isomorphic with 
$\cO(\Omega, \Ci(\spXPvp \col \tauP)).$ 
Thus, from (\refer{e: q of f as sum of q of Fs new})  
we deduce that $q_{\gs, \xi}(P, v\asmid f)$ 
is an element 
of $P_k(\faPq) \otimes \cO(\Omega, \Ci(\spXPvp\col \tauP)).$ 

Combining (\refer{e: q of f gl as s sum}), (\refer{e: defi gf s}) and  
(\refer{e: q of f as sum of q of Fs new}) we infer that 
(\refer{e: q of f in gl versus q of fgl new})
holds for $X\in \faPq,$ $\gl \in \Omega'.$
On the other hand, 
if $\gl \in \Omega',$ then it follows from 
(\refer{e: expansion f along Q v}) 
with $P$ and $f_\gl$ in place of $Q$ and $f,$
that, for $R > 1,$ $m \in \spXPvp[R]$ and $a \in \APqp(R^{-1}),$
\begin{equation}
\naam{e: series for f gl in proof}
f_\gl(mav) = \sum_{\gs \in \WPQ}a^{\gs \gl -\rho_P} 
\sum_{\xi \in -\gs \cdot Y + \N \DrP}
a^{-\xi} q_{\gs\gl|_{\faPq} - \rho_P - \xi}(P, v\asmid f_\gl, \log a)(m),
\end{equation}
where the series converges neatly on $\APqp(R^{-1}),$ as 
a $\DrP$-exponential polynomial series with coefficients in $V_\tau$ 
(use (\refer{e: exponents f gl along P v in this set}) 
and Lemma \refer{l: exponents disjoint}).
Substituting (\refer{e: q of f in gl versus q of fgl new}) in (\refer{e: series for f gl in proof}) we obtain 
the identity (\refer{e: series for fam along P})  for $\gl \in \Omega',$ 
$m \in \spXPvp[R]$ and $a \in \APqp(R^{-1}),$ with the convergence as asserted.

Thus, it remains to show that the identity (\refer{e: series for fam along P})
extends to all $\gl \in \Omega$ and that the final assertion of the theorem holds.
We will first establish the final assertion.

It follows from Theorem 
\refer{t: transitivity of asymptotics}  that the series
\begin{equation}
\naam{e: series with gf s q}
\sum_{\xi \in -sW_Q Y|_{\faPq} + \N \DrP} a^{-\xi} q_{-\xi}(P,v\asmid F_s,\log a)
\end{equation}
converges neatly on $\APqp(R^{-1})$ as a $\DrP$-exponential polynomial series  
with coefficients in the space (\refer{e: spaces cO Omega Ci}).
The series (\refer{e: series with q gs xi})
arises as the sum over $s \in W/W_Q$ with $[s] = \gs$ of the series 
in (\refer{e: series with gf s q}) multiplied by $\gf_s.$ Since multiplication
by $\gf_s$ induces a continuous linear endomorphism of the space 
(\refer{e: spaces cO Omega Ci}), this establishes the final assertion of the theorem.

{}From the final assertion it follows that, for every $R > 1,$ 
 the series on the right-hand side of (\refer{e: series for fam along P})
defines a holomorphic function of $\gl \in \Omega,$ for 
every $m \in \spXPvp[R]$ and $a \in \APqp(R^{-1}).$ For such $m,a$ 
the function $\gl \mapsto f_\gl(mav)$ is holomorphic in $\gl \in \Omega$ 
by Lemma \refer{l: the family is analytic}; hence the identity 
(\refer{e: series for fam along P}) extends to all $\gl \in \Omega,$ 
by density of $\Omega'$ in $\Omega.$ 
\qed

\begin{thm}{\rm (Transitivity of asymptotics).\ }
\naam{t: transitivity of asymptotics for families new}
Let $Q$, $\Omega$, $Y$, $f$, $P$ and $v$ be as in Theorem
\refer{t: behavior along the walls for families new}.
Let ${\Pmin} \in \minparabs$ be contained 
in $P.$ 
Let $\gs \in \WPQ$ and $\xi \in - \gs \cdot Y +\N \DrP.$ 
Then for every $X \in \faPq,$ all $u \in \NKPaq,$ 
$b\in \stAPqp(\starPmin),$ 
$m \in M$ 
and $\gl \in \Omega,$ 
\begin{equation}
\naam{e: series in transitivity for families}
q_{\gs,\xi}(P,v \asmid f, X)(\gl, mbu)
=
\sum_{s \in W/W_Q\atop [s] = \gs} b^{s\gl - \rho_{{\Pmin}}}
\sum_{\mu \in - sW_Q Y + \N \DPmin \atop \mu|_{\faPq} = \xi}
b^{-\mu}\; q_{s, \mu} ({\Pmin}, uv\asmid f, X + \log b)(\gl, m).
\end{equation}
Moreover, for every $s \in W/W_Q$ with $[s] = \gs$ and every
$X \in \faPq,$ the series
\begin{equation}
\naam{e: second series in transitivity for families}
\sum_{\mu \in - sW_Q Y + \N \DPmin \atop \mu|_{\faPq} = \xi}
b^{-\mu}\; q_{s, \mu} ({\Pmin}, uv\asmid f, X + \log b)
\end{equation}
converges neatly on $\stAPqp(\starPmin)$
as a $\DPmin$-exponential polynomial series in the variable $b$ 
with coefficients in $\cO(\Omega, \Ci(\spX_{0,uv}\col \tauM)).$
\end{thm}

\proof
Fix $u \in \NKPaq.$ Moreover,
we fix a set $\cW_{P,v}$ 
as in the beginning of the proof of 
Theorem \refer{t: behavior along the walls for families new}
such that it contains the element $u.$ 
We will also use the remaining notation of the proof of the mentioned 
theorem.

Using (\refer{e: q of Fs versus q of f}) we see that, via the natural 
isomorphism of $\cO(\Omega, \Ci(\spX_{0, uv}\col \tauM))$ with 
$\Ci(\spX_{0,uv} \col \tilde \tauM),$ the series 
(\refer{e: second series in transitivity for families})
may be identified with the  series
with coefficients in $\Ci(\spX_{0,uv} \col \tilde \tauM)$ 
that arises from the series on the right-hand side of 
(\refer{e: series for q minus xi FW}) 
by omitting the evaluation at $m.$  
The neat convergence of the latter series was noted already.
Moreover, the identity (\refer{e: series in transitivity for families})
follows by insertion of (\refer{e: series for q minus xi FW}) 
in the definition (\refer{e: q of f as sum of q of Fs new}) 
of $q_{\sigma,\xi}$.
\qed

The following result is an important consequence of 
`holomorphy of asymptotics.'

\begin{lemma}
\naam{l: holo of asymp}
Let $Q\in \allparabs,$  $Y\subset \staQqdc$ a finite subset and
$\Omega \subset \faQqdc$ a non-empty open subset. 
Let $f\in \CepQY(\spXp\col\tau \col \Omega)$
and let $P\in \allparabs,$ $v\in \NKaq,$ and 
$\gs\in \WPQ.$ 

Let $\xi \in -\gs \cdot Y + \N \DrP$ and assume that there exists
a  $\gl_0\in \faQqdczero(P, WY) \cap \Omega$ such that 
\begin{equation}
\naam{e: expression for gl zero is an exponent}
\gs \gl_0|_{\faPq} - \rho_P - \xi \in \Exp(P,v\asmid f_{\gl_0}).
\end{equation}
Then there exists a full open subset $\Omega_0$ of $\Omega$ such that 
$$
\gs \gl|_{\faPq} - \rho_P - \xi \in \Exp(P,v\asmid f_{\gl}),\qquad (\forall \gl \in \Omega_0).
$$
\end{lemma}

\proof
{}From (\refer{e: expression for gl zero is an exponent}) combined with 
(\refer{e: q of f in gl versus q of fgl new}) it follows that the 
$P_k(\faPq) \otimes \Ci(\spXPvp\col \tau_P)$-valued holomorphic
function $q: \gl \mapsto q_{\gs,\xi}(P,v\asmid f, \dotvar, \gl)$ 
does not vanish at $\gl= \gl_0.$
Hence there exists
a full open subset $\Omega_1 \subset \Omega$ such that $q(\gl)\neq 0$ for all $\gl \in \Omega.$
Let $\Omega_0:= \Omega_1 \cap \faQqdczero(P, WY),$ then the conclusion follows
by application of (\refer{e: q of f in gl versus q of fgl new}).
\qed

We end this section with a result describing the behavior of the functions
$q_{\gs,\xi}$  under the action of $\NKaq.$ 
Let $Q,P \in \allparabs$
and $u \in \NKaq,$ and put $P' =uPu^{-1}.$
The left multiplication by $u$ naturally induces
a  map $\WPQ \to W/\sim_{P' \mid Q},$ which we denote
by $\gs \mapsto u \gs.$ Moreover, the endomorphism  $\Ad(u^{-1})^*$ 
of $\faqdc$ restricts to a linear map $\faPqdc \to \fa_{P'\iq\iC}^*,$ 
which we denote by $\eta \mapsto u \eta.$ With these notations,
if $Y \subset \staQqdc$ is a finite subset and $\gs \in \WPQ,$ then
$$ 
u(\gs\cdot Y) = (u\gs)\cdot Y;
$$ 
see also (\refer{e: defi gs cdot Y}). For $v \in \NKaq,$ let the map 
$\rho_{\tau,u}: \Ci(\spXPvp\col \tau_P) \to 
\Ci(\spX_{P', uv, +}\col \tau_{P'})$ 
be defined by (\refer{e: defi rho tau u}) with $P$ in place of $Q.$ 

If $\Omega \subset \faQqdc$ is an open subset, let 
$\Ad(u^{-1})^* \otimes 1 \otimes \rho_{\tau,u}$ denote the naturally  induced map from
$P(\faPq) \otimes \cO(\Omega, \Ci(\spXPvp\col \tau_P))$ 
to 
$P(\fa_{P'\iq})  
\otimes \cO(\Omega, \Ci(\spX_{P', uv, +}\col \tau_{P'})).
$

\begin{lemma}
\naam{l: transformation holo coeffs}
Let $Q \in \allparabs,$ $Y \subset \staQqdc$ a finite subset
and $\Omega \subset \faQqdc$ a non-empty open subset.
Let $f\in C^\ep_{Q,Y}(\spXp\col \tau \col \Omega).$ 
If $P \in \allparabs$ and $u,v\in \NKaq,$ then for all $\gs \in \WPQ$ 
and $\xi \in \gs\cdot Y,$
$$ 
q_{u\gs, u\xi}(uPu^{-1}, uv \asmid f) = [\Ad(u^{-1})^* \otimes 1 \otimes \rho_{\tau,u}]
q_{\gs, \xi}(P, v\asmid f).
$$ 
\end{lemma} 

\proof
{}From combining (\refer{e: q of f in gl versus q of fgl new}) and
Lemma \refer{l: transformation of coeffs}
it follows that there exists a full open subset $\Omega_0$ of $\Omega$ such that,
for 
$\gl \in \Omega_0,$ 
$$ 
q_{u\gs, u\xi}(uPu^{-1}, uv \asmid f,\dotvar,\gl) 
= [\Ad(u^{-1})^* \otimes \rho_{\tau,u}]
q_{\gs, \xi}(P, v\asmid f,\dotvar,\gl).
$$
The result now follows by holomorphy of the above expressions in $\gl$ 
and density of $\Omega_0.$
\qed

\section{Asymptotic globality}
\naam{s: asymptotic globality}
In this section we introduce the notion of asymptotic globality 
of a spherical function on $\spXp$ and of an analytic family of such functions. 
We discuss properties needed in the statement and proof of the vanishing theorem 
in the next section.

\begin{defi}
\naam{d: globality property of exponents new}
Let  $P \in \allparabs$ and $v\in \NKaq.$ 
A function $f\in \Exppol(\spXp\col \tau)$ is said to be 
asymptotically global along $(P,v)$ at an element
$\xi \in \faPqdc$ if, for every $X \in \faPq,$  the $\Vtau$-valued smooth function
$q_\xi(P,v\asmid f,X)$ 
has a $C^\infty$-extension from $\spXPvp$ to $\spXPv.$ 
\end{defi}

\begin{rem}
\naam{r: after defi 7.1}
Since $q_\xi(P,v\asmid f,X)$ is polynomial in $X,$ with values in 
$\Ci(\spXPvp\col \tau_P),$ 
the requirement on $q_\xi$ implies that $q_\xi(P,v\asmid f)$ is a polynomial
$\Ci(\spXPv\col \tau_P)$-valued function on $\faPq.$ 

Note that for $P$ minimal the condition of asymptotic globality along 
$(P,v)$ is automatically fulfilled, since $\spXPvp = \spXPv.$ 

Finally, if $f \in \Cep(\spXp\col \tau),$ then $f$ is asymptotically 
global along $(G,e)$ at every exponent if and only if $f$ extends
smoothly to $\spX$ (use Remark \refer{r: 1.5 bis}).
\end{rem}

The property of asymptotic globality is preserved under the action
of $\DX$ in the following fashion.
If $P \in \allparabs,$ then by $\preceq_{\DrP}$ we denote the partial ordering
on $\faPqdc,$ defined as in (\refer{e: partial ordering preceq gD}), 
with $\faPq$ and $\DrP$ in place of $\fa$ and $\gD,$ 
respectively.

\begin{prop}
\naam{p: stability of globality new}
Let $f \in \Exppol(\spXp\col \tau)$ and $D \in \DGH.$ 
Let $P \in \allparabs,$ $v\in \NKaq$ and $\xi_0 \in \faPqdc.$ 
If
$f$ is asymptotically global along $(P,v)$ at 
every exponent $\xi \in \Exp(P,v\asmid f)$ with $\xi_0 \preceq_{\DrP} \xi,$ 
then $Df$ is asymptotically global along $(P,v)$ at $\xi_0.$ 
\end{prop}

\proof
Let   
$u := \mu_P'(D) + u_+$ be the element of $\cD_{1P}$ associated with $D$ as
in Proposition \refer{p: radial deco with MQ},
with $P$ in place of $Q.$  The key idea in the present proof 
is that $u$ has a $\DrP$-exponential polynomial expansion with coefficients 
that are globally defined smooth functions on $\MPgs,$ by
Cor. \refer{c: cor on ringQ}.
More precisely, the expansion $\ep(u)$ is a finite sum, as $i$ 
ranges over a finite index set $I,$ 
of terms of the form
$$ 
\ep(u)_i = \sum_{\nu \in \N\DrP} a^{-\nu} \gf_{i,\nu}\otimes S_{i,\nu} \otimes u_{i,\nu} \otimes 
v_{i,\nu}.
$$
Here $\gf_{i,\nu} \in C^\infty(\MPgs),$ $S_{i,\nu} \in \End(\Vtau),$ $u_{i,\nu} \in U(\fm_{P\gs})$ 
and $v_{i,\nu} \in U(\faPq),$ and $\deg(u_{i,\nu}) + \deg(v_{i,\nu}) \leq \deg(D)$ for 
all $i,\nu.$ 
By Lemma \refer{l: radial component applied to expansion}, $Df$ belongs to 
$\Cep(\spXp\col \tau)$
and its $(P,e)$-expansion results from the $(P,e)$-expansion of $f$  
by the formal application of the element $\ep(u).$  
Hence the asymptotic coefficient of $\xi_0$ 
is given by the finite sum
$$ 
q_{\xi_0}(P,e \asmid Df)(X , m) = \sum_{{\xi \in \Exp(P,e\asmid f) 
\atop \nu \in\N \DrP}\atop 
\xi - \nu = \xi_0} 
\sum_{i\in I} 
\gf_{i,\nu}(m) 
S_{i,\nu}[\, q_\xi(P,e\asmid f)(X; T_\xi(v_{i,\nu}), m ; u_{i,\nu})\,].
$$ 
Let now $f$ satisfy the hypothesis of the proposition.
The $\xi's$ occurring in the above sum 
belong to $\xi_0 + \N\DrP,$ hence satisfy $\xi_0 \preceq_{\DrP} \xi.$ 
By hypothesis, the associated coefficients $q_\xi(P,e\asmid f)$ all extend smoothly 
to $\faPq \times \MPgs,$ see Remark \refer{r: after defi 7.1}.
 Therefore, so does $q_{\xi_0}(P,e \asmid Df).$
This establishes the result for arbitrary $P \in \allparabs$ and the special choice $v=e.$ 
The result with general $v \in \NKaq$ now follows by application of Lemma
\refer{l: transformation of coeffs} (cf.\ Lemma 
\refer{l: transformation of globality} (a)).
\qed

We shall also introduce a notion of asymptotic globality for families 
from the space 
$\anfamQY$ introduced in the previous section, with $\Omega \subset \faQqdc$ 
an open subset.

\begin{defi}
\naam{d: s globality new}
Let $Q \in \allparabs,$  $Y $ a finite subset of 
$\staQqdc$ and $\Omega \subset \faQqdc$ a non-empty open subset. 
Let  $P \in \allparabs,$ $v \in \NKaq$ and $\gs \in \WPQ.$ 

We will say that a family $f \in \CepQY(\spXp\col \tau\col \Omega)$ 
is $\gs$-global along $(P,v),$ if there exists a dense open 
subset $\Omega_0$ of  $\Omega,$ 
such that, for every $\gl \in \Omega_0,$ the function $f_\gl$ 
is asymptotically global along $(P,v)$ at each exponent
$\xi \in \gs\gl|_{\faPq} + \gs\cdot Y -\rho_P -  \N\DrP.$
\end{defi}

\begin{rem}
\naam{r: globality independent of Y}
If $Y_1$ and $Y_2$ are finite subsets of $\staQqdc$ with $Y_1\subset Y_2,$ 
then obviously
$$ 
C^\ep_{Q, Y_1} (\spXp \col \tau \col \Omega) 
\subset 
C^\ep_{Q, Y_2} (\spXp \col \tau \col \Omega). 
$$
If $f$ belongs to the first of these spaces, then the condition of $\gs$-globality
along $(P,v)$ relative to $Y_1$ is equivalent to the similar condition relative 
to $Y_2.$ This is readily seen by using Lemmas
\refer{l: exponents disjoint} and \refer{l: azero is full}.
{}From this we see that the notion of $\gs$-globality along $(P,v)$ extends to the space
$$%\begin{equation}
%\naam{e: defi C ep Q without Y}
C^\ep_{Q} (\spXp \col \tau \col \Omega) :=
\bigcup_{Y\subset \staQqdc \; {\rm finite}}
C^\ep_{Q, Y} (\spXp \col \tau \col \Omega) 
$$%\end{equation}
\end{rem}

The
property of asymptotic globality for families is also stable under the action 
of $\DX.$ 
\begin{cor}
\naam{c: stability of globality}
Let $Q \in \allparabs,$  $Y $ a finite subset of 
$\staQqdc$ and $\Omega \subset \faQqdc$ a non-empty open subset.
Let  $P \in \allparabs,$ $v \in \NKaq$ and $\gs \in \WPQ.$ 

Let $f \in \CepQY(\spXp\col \tau\col \Omega)$ 
be $\gs$-global along $(P,v).$ Then for every $D \in \DX$ the family 
$Df \in \CepQY(\spXp\col \tau\col \Omega)$ is $\gs$-global
along $(P,v)$ as well.
\end{cor}

\proof
It
follows from Proposition \refer{p: D on families} that
$Df$ belongs to $\CepQY(\spXp\col \tau\col \Omega).$
According to  Theorem \refer{t: behavior along the walls for families new}, both sets
$\Exp(P, v\asmid f_\gl)$ and $\Exp(P, v \asmid Df_\gl)$ are contained
in the set $E_\gl:= W(\gl + Y)|_{\faPq} - \rho_P - \N \DrP,$ 
for every $\gl \in \Omega.$ 

Let $\Omega_0$ be as in Definition \refer{d: s globality new}. Then the set 
$\Omega_0': = \Omega_0 \cap \faQqdczero(P, WY)$ is open dense in $\Omega$ 
by Lemma \refer{l: azero is full}.
Let $\gl \in \Omega_0'$ and let 
$\xi_0 \in  \gs \gl|_{\faPq} + \gs \cdot Y  - \rho_P - \N \DrP.$
If $\xi \in \Exp(P,v\asmid f_\gl)$ satisfies $\xi_0 \preceq \xi,$ 
then $\xi \in \gs \gl |_{\faPq} + \gs \cdot Y - \rho_P - \N \DrP$ 
by Lemma \refer{l: exponents disjoint}. 
By hypothesis, $f_\gl$ is asymptotically global along $(P,v)$ at
the exponent $\xi.$ It now follows by application of 
Proposition \refer{p: stability of globality new} that $Df_\gl$ 
is asymptotically global along $(P,v)$ at $\xi_0.$ 
\qed

The following lemma describes the behavior of asymptotic
globality under the action of $\NKaq.$ 

\begin{lemma}
\naam{l: transformation of globality}
Let $P\in \allparabs$ and $u,v \in \NKaq$. Put $P'=u P u^{-1}$
and $v'=uv$.
\begin{enumerate}
\itema 
Let $f \in C^{\ep}(\spXp \col \tau )$ and $\xi\in\faPqdc$.
If $f$ is asymptotically global along $(P,v)$ at $\xi$, then
$f$ is asymptotically global along $(P',v')$ at $u\xi$.
\itemb 
Let $Q \in \allparabs,$ $\Omega \subset \faQqdc$ a non-empty open subset,
$Y\subset \staQqdc$ a finite subset,
$f \in C^{\ep}_{Q,Y} (\spXp \col \tau \col \Omega)$ and $\gs \in \WPQ.$
If $f$ is $\gs$-global along $(P,v),$ then
$f$ is $u\gs$-global along $(P',v')$. 
\end{enumerate}
\end{lemma}
\proof
{}From (\refer{e: defi rho tau u}) with $P$ in place of $Q$ 
it is readily seen that $\rho_{\tau, u}$ maps
$\Ci(\spXPv\col \tau_P)$ to $\Ci(\spX_{P',v'} \col \tau_{P'}).$ 
Then (a) and (b) follow immediately from Lemmas 
\refer{l: transformation of coeffs}
and \refer{l: transformation holo coeffs}, respectively. 
\qed

We end this section with the following result, which shows that the
globality condition is fulfilled for a certain natural class
of $\tau$-spherical functions. {}From
the text preceding Lemma \refer{l: restriction on leading exponents} we recall that 
$\fb$ is a maximal abelian subspace of $\fq$ containing
$\faq$ and that if $\mu \in \fbdc,$ then by $I_\mu$ we denote the kernel
of the character $\gg(\dotvar \col \mu)$ of $\DGH.$ Thus $I_\mu$ is an
ideal in $\DX$ of codimension one (over $\C$).

\begin{prop} 
\naam{p: global eigen implies as global}
Let $\mu\in\fbdc$ and let $f\in \cE(\spX\col\tau\col I_\mu)$. Then 
$f|_{\spXp}\in C^\ep(\spXp\col\tau)$. Moreover, this function
is asymptotically global along all pairs $(P,v)\in\allparabs\times\NKaq$
and at all exponents $\xi\in\faPqdc$.
\end{prop}

\proof
The first statement follows immediately from Lemma 
\refer{l: DX finite in exppol}. By Lemma 
\refer{l: transformation of globality} (a) it suffices to consider
$v=e$ and arbitrary $P\in\allparabs$. Let $\psi\in V_\tau$ be fixed.
Then it suffices to prove that the scalar valued function
$m\mapsto q_\xi(X,m):= \hinp{q_\xi(P,e|f,X,m)}{\psi}$
on $\spX_{P,+}$ has a $C^\infty$ extension to $\spX_P$, for each
$\xi\in\faPqdc$, $X\in\faPq$. It follows from 
Theorem \refer{t: expansion along the walls} that 
\begin{equation}
\naam{e: expansion of f(ma)psi}
\hinp{f(ma)}{\psi}=\sum_{\xi\in Y-\N\Delta_r(P)}a^\xi q_\xi(\log a,m).
\end{equation}
On the other hand, it follows from \bib{B91}, Lemma 12.3, that
\bib{B91}, Thm.\ 12.8 can be applied to the $K$-finite function
$F\colon x\mapsto \hinp{f(x)}{\psi}$. By uniqueness of asymptotics
(see Lemma \refer{l: uniqueness of asymp} and its proof) the expansion 
(\refer{e: expansion of f(ma)psi}) coincides with that of 
\bib{B91}, Thm.\ 12.8. We conclude that, in the notation of loc.\ cit.,
$q_\xi(X,m)=p_{\mu|_{\faq},\xi}(P|F,m,X)$ for all $X\in\faPq$, 
$m\in\spX_{P,+}$. The function $x\mapsto p_{\mu|_{\faq},\xi}(P|F,x,X)$
is smooth on $G$.
{}From this the smooth extension of $q_\xi(X,m)$
follows immediately.\qed

\section{A vanishing theorem}
\naam{s: vanishing thm}
In this section we  formulate and prove the vanishing theorem.
We assume that $Q$ is a $\gs$-parabolic subgroup containing $\Aq.$ 

As before, let $\fb$ be a maximal abelian subspace of $\fq$ containing
$\faq.$ By  $\staQq$ and $\stbQ$ 
we denote the orthocomplements of $\faQq$ in $\faq$ and $\fb,$ respectively.
Let $\fbk: = \fb \cap \fk;$ then   
$$%\begin{equation}
%\naam{e: defi star fb Q}
\stbQ = \fbk \oplus \staQq.
$$%\end{equation} 
We  write $\DQmaps$ for the collection of functions
$\gdmap: \stbQdc \to \N$ with finite support $\supp \gdmap.$ 
For $\gd \in \DQmaps$ we put 
$$
|\gd|:=  \sum_{\nu \in \supp \gd} \gd(\nu).
$$
For $\gdmap \in \DQmaps$ and $\gl \in \faQqdc$ we define
the ideal $\Igdgl$ in $\DGH$  as the following product of ideals 
\begin{equation}
\naam{e: ideal I gl gdmap}
\Igdgl := \prod_{\nu \in \supp{\gdmap}} (I_{\nu+ \gl})^{\gdmap(\nu)}.
\end{equation}
If $\gdmap = 0,$ this ideal is understood to be the full ring $\DGH.$ 
Being a product of cofinite ideals in the Noetherian ring $\DGH,$ the 
ideal $\Igdgl$ is  cofinite.

\begin{defi}
\naam{d: defi cE Q Y revised new}
Let    $\Omega \subset \faQqdc$ be a non-empty open subset and  
$\gd\in\DQmaps.$ For every finite subset  $Y \subset \staQqdc$ 
we define 
\begin{equation}
\naam{e: anfamQY gd} 
\cEQY(\spXp\col \tau \col \Omega \col \gdmap) 
\end{equation}
to be the space of families 
$f \in \anfamQY$ (cf.\ Def.\ \refer{d: anfamQY newer}) such that
for every $\gl \in \Omega$ the function $f_\gl: x \mapsto f(\gl, x)$ 
is annihilated by the cofinite ideal (\refer{e: ideal I gl gdmap}).
Moreover, we define 
$$
\FQOmegad:= 
\bigcup_{Y\subset \staQqdc \;\;{\rm finite}}\;\; 
\cE_{Q,Y}(\spXp\col \tau \col \Omega\col\gdmap).
$$
\end{defi}
Note that the space (\refer{e: anfamQY gd}) depends on $Q$ through its
$\gs$-split component $\AQq.$ 
If $\nu \in \stbQdc,$ we denote by $\gd_\nu$ the characteristic 
function of the set $\{\nu\}.$ Then $\gd_\nu \in \DQmaps.$ Moreover,
if $\gd \in \DQmaps$ and $\nu \in  \supp \gd,$ then $\gd - \gd_\nu \in \DQmaps$ 
and $|\gd - \gd_\nu| = |\gd| -1.$

\begin{lemma}
\naam{l: stability family under D new}
Let  $f \in \cEQY(\spXp\col \tau \col \Omega \col \gdmap).$ 
\begin{enumerate}
\itema
If $D \in \DX$ then $Df \in \cEQY(\spXp\col \tau \col \Omega \col \gdmap).$ 
\itemb
If $D \in \DX$ and $\nu \in \supp \gd$ then the function
 $g: \Omega \times \spXp \to \Vtau$ 
defined by 
\begin{equation}
\naam{e: defi g as D minus gamma f}
g(\gl,x): = [D - \gg(D\col \nu + \gl)]f_\gl(x),\qquad (\gl \in \Omega, \; x \in \spXp),
\end{equation}
belongs to $\cEQY(\spXp\col \tau \col \Omega \col \gdmap- \gd_\nu).$ 
\end{enumerate}
\end{lemma}
\proof
Let $D \in \DX.$ 
By Proposition \refer{p: D on families}, the family $Df$ belongs to
$C^{\ep}_{Q,Y}(\spXp\col \tau\col \Omega).$ Moreover, if $\gl\in\Omega$ and 
$D ' \in \Igdgl,$ then $D' (Df)_\gl = D'D f_\gl = D D' f_\gl = 0$ 
and we see that assertion (a) holds.

The function $\gl \mapsto \gg(D\col \nu + \gl)$ is polynomial on $\faQqdc,$ 
hence holomorphic
on $\Omega$ and it follows that $G: (\gl, x) \mapsto \gg(D\col \nu +\gl) f(\gl, x)$ 
belongs to $C^\ep_{Q,Y}(\spXp\col \tau\col \Omega).$ Hence 
$g = Df - G$ belongs to the latter space as well. 
Furthermore, if $D' \in I_{\gd - \gd_\nu, \gl},$ then 
$D'':= D' (D - \gg(D\col \nu + \gl)) \in \Igdgl,$ and we see that
$D'g_\gl = D'' f_\gl = 0.$ Hence (b) holds.
\qed

\begin{rem}
It follows from Lemma \refer{l: stability family under D new} (a)
that (\refer{e: action D on family}) defines a representation of $\DX$ 
in
$\cEQ(\spXp\col \tau \col\Omega\col \gdmap),$ leaving the subspaces 
$\cEQY(\spXp\col \tau \col \Omega\col\gdmap)$
invariant.
\end{rem}

\begin{lemma}
\naam{l: locally killed by I gl gd}
Let $Q \in \allparabs,$ $\gd\in \DQmaps$ 
and $\Omega$ a connected non-empty open subset of $\faQqdc.$ 
Assume that $f \in C^\ep_Q(\spXp\col \tau\col \Omega).$
If $f_\gl$ is annihilated by $\Igdgl$ for $\gl$ in a non-empty open 
subset $\Omega'$ of $\Omega,$ then 
$f \in \cEQ(\spXp\col \tau \col \Omega\col \gd).$
\end{lemma}
\proof 
Fix a finite subset $Y\subset \staQqdc$ such that 
$f \in C^{\ep}_{Q,Y}(\spXp\col \tau\col \Omega).$ 
We proceed by induction on $|\gd|.$ 

First, assume that  $|\gd| = 0.$ Then $I_{\gd, \gl} = \DX$ for all $\gl$ 
and hence $f|_{\Omega'\times\spXp}=0.$
Since $\Omega$ is connected, this implies that $f = 0,$ 
see Lemma \refer{l: the family is analytic}.

Next assume that $|\gd| = k \geq 1$ and assume the result 
has already been established for all $\gd \in \DQmaps$ with 
$|\gd| < k.$ 
Fix $\nu \in \supp \gd$ and put $\gd' = \gd - \gd_\nu,$ then $|\gd'| < k.$ 
Let $D \in \DX$ and define $g$ as in (\refer{e: defi g as D minus gamma f}). 
Then $g \in C^\ep_{Q,Y}(\spXp \col \tau \col \Omega)$, as seen
in the proof of Lemma \refer{l: stability family under D new}.
On the other hand, it follows from (b) of that lemma
that $g|_{\Omega'\times\spXp}\in \cEQY(\spXp\col \tau \col \Omega' \col \gd')$.
Hence $g \in \cEQ(\spXp\col \tau \col \Omega\col \gd')$ 
by the induction hypothesis. Fix $\gl \in \Omega.$ Then it follows,
for $D' \in I_{\gd', \gl},$ that $D'(D - \gg(D\col \nu + \gl)) f_\gl = D'g_\gl = 0.$ 
Since $D$ was arbitrary, we conclude that $f_\gl$ is annihilated
by the ideal $I_{\gd', \gl} I_{\nu+\gl} = \Igdgl.$ 
\qed

We define the following subset of $\allparabs,$ consisting of the  
parabolic subgroups whose $\gs$-split rank is of codimension one,
$$ 
\parone:= \{P \in \allparabs \mid \dim \faq / \dim \faPq  = 1\}.
$$ 
\begin{defi}
\naam{d: family for the vanishing thm newer}
Let $Q \in \allparabs,$ $\Omega\subset \faQqdc$ 
a non-empty open
subset and $\gd \in \DQmaps.$ 
By $\FQOmegadP$ we denote the space of functions 
$f \in \FQOmegad$ 
satisfying the following condition.

\hbox{\hspace{-12pt}
\vbox{\vspace{-2mm}
\begin{enumerate}
\item[]
For every $s \in W$ and every $P \in \parone$ with $ s (\faQq)  \not\subset \faPq,$ 
the family $f$ is $[s]$-global along $(P,v),$ for all $v \in \NKaq;$ 
here $[s]$ denotes the image of $s$ in $W/{\sim_{P|Q}}= W_P \bs W / W_Q.$
\end{enumerate}
\vspace{-2mm}}}
\noindent
If $Y \subset \staQqdc$ is a finite subset, we define 
$$
\FQYOmegadP: = \FQYOmegad \cap \FQOmegadP.
$$ 
\end{defi}
\begin{rem}
\naam{r: W Pga Q new}
Note that $\FQOmegadP$ depends on $Q$ through its $\gs$-split component $\faQq.$ 

The equality $W/{\sim_{P|Q}} = W_P \bs W/ W_Q$  follows from 
Lemma \refer{l: WPQ as cosets}. 
Note that the condition $s(\faQq) \not \subset \faPq$ on $s$
factors to a condition on its class in $W_P\bs W/ W_Q$. 
\end{rem}

The following result reduces the globality condition of Definition
\refer{d: family for the vanishing thm newer} to a condition involving 
a smaller set of $(s, P) \in W \times \parone.$ Its formulation requires
some more notation. 

Let
$\gD$ be a fixed basis for the root system $\Sigma,$ let $\Sigma^+$ be
the associated system of positive roots and $\faqp$ the associated open positive
chamber. Let $\stPo$ be the unique element of $\minparabs$ with $\gD(\stPo) = \gD.$
A $\gs$-parabolic subgroup $Q$ is said to be standard if it contains
$\stPo;$ of course then $Q \in \allparabs.$
Given  such a $Q,$ we write $\Delta_Q$ for the subset of $\gD$ consisting 
of the roots vanishing 
on $\faQq$ and $\Delta(Q)$ for its complement. 

If $\ga$ is any root in $\gD,$ we write
$
\fn_\ga
$ 
for the sum of the root spaces $\fg_\gb$ where $\gb$ ranges over the set
$\Sigma^+\setminus \N \ga.$  
Moreover, we put $N_\ga: = \exp(\fn_\ga)$ and write $M_{1\ga}$ for 
the centralizer in $G$ of the root hyperplane $\ker \ga.$ Then $P_\ga = M_{1\ga} N_\ga$ 
is the standard parabolic subgroup with $\Delta_{P_\ga}=\{\ga\}$.
We write $P_\ga = M_\ga A_\ga N_\ga$ and $P_\ga = M_{\gs \ga} A_{\ga \iq} N_\ga$ 
for the Langlands and $\gs$-Langlands decompositions of $P_\ga,$ respectively.
Accordingly, $\fagaq = \ker \ga$ and $\stagaq = (\ker \ga)^\perp.$
Finally, we write $W_\ga = W_{P_\ga}$ for the centralizer of $\ker \ga$ in $W.$

\begin{lemma}
\naam{l: minimal condition for glob new}
Let $Q\in \allparabs$ be a standard parabolic subgroup,
$\Omega \subset \faQqdc$ a non-empty open subset,
$\gd \in \DQmaps$ and $f \in \cE_Q(\spXp\col \tau \col \Omega\col \gd).$ 
Then $f$ belongs to $\cE_Q(\spXp\col \tau \col \Omega\col \gd)_\glob$
if and only if the following condition is fulfilled.

\hbox{\hspace{-12pt}
\vbox{\vspace{-2mm}
\begin{enumerate}
\item[]
For every $s \in W$ and every $\ga \in \gD$ with $s^{-1}\ga|_{\faQq} \neq 0,$ 
the family $f$ is $[s]$-global along $(\stPga,v),$ for all $v \in \NKaq;$ 
here $[s]$ denotes the image of $s$ in $W/{\sim_{\stPga|Q}}= W_\ga \bs W / W_Q.$
\end{enumerate}
\vspace{-2mm}}}
\end{lemma}

\proof
We must show that the condition of Definition \refer{d: family for the vanishing thm newer}
is fulfilled if and only if 
the above condition holds.
For this we first observe that for $\ga \in \gD$ and $s \in W,$ 
$$
s^{-1}\ga|_{\faQq} \neq 0 \iff s(\faQq) \not\subset \fa_{\ga\iq}.
$$ 
The `only if part' is now immediate. 
For the `if part', assume that the above condition is fulfilled.  
Let $(s, P) \in W \times \parone$ be such that $s(\faQq) \not\subset \faPq.$ 
There exist $\ga \in \gD$ and $t \in W$ such that $tP t^{-1} = P_\ga.$
It follows that  $ts(\faQq) \not\subset t\faPq = \ker \ga,$ hence $(ts)^{-1} \ga = \ga \after (ts) $ is not 
identically zero on $\faQq.$ 
From the
hypothesis it now follows
that $f$ is $[ts]$-global along $(t P t^{-1}, v),$ for all $v \in \NKaq.$ 
By Lemma \refer{l: transformation of globality} 
it follows that $f$ is $[s]$-global along $(P, w),$ for all 
$w \in \NKaq.$  
\qed

\begin{lemma}
\naam{l: stability family with globality under D}
Let $Q \in \allparabs,$  $\Omega \subset \faQqdc$ a non-empty open subset 
and $\gd\in \DQmaps.$ 
Then the space $\FQOmegadP$ is $\DX$-invariant. Moreover,  $\FQYOmegadP$ is a $\DX$-submodule,
for every finite subset $Y \subset \staQqdc.$ 
\end{lemma}

\proof 
This follows from combining the $\DX$-invariance of the space
$\FQYOmegad$  with Proposition \refer{p: stability of globality new}.
\qed

\begin{defi}\naam{d: Q-distinguished}
Let $Q\in \allparabs.$
An open subset $\Omega$ of $\faQqdc$ will be called 
$Q$-distinguished if it is connected and if 
for every $\ga \in \Sigma(Q)$ the function
$\gl \mapsto \inp{\Re \gl}{\ga}$ is not bounded from above
on $\Omega.$
\end{defi}

In particular, a connected open dense subset of $\faQqdc$ is 
$Q$-distinguished.
In the following theorem we assume that $\QW \subset \NKaq$ is a complete 
set of representatives for $W_Q\backslash W / \WKH.$ 

\begin{thm}
{\rm (Vanishing theorem).\ }
\naam{t: vanishing theorem new}
Let $Q \in \allparabs$
and $\gdmap \in \DQmaps.$ 
Let $\Omega \subset \faQqdc$ 
be a $Q$-distinguished open subset and let 
$f \in \cE_Q(\spXp\col\tau\col\Omega\col \gd)_\glob.$ 
Assume that there exists a non-empty open subset 
$\Omega' \subset \Omega$ such that, for each $v \in \QW, $ 
\begin{equation}
\naam{e: hypothesis vanishing thm new}
\gl - \rho_Q \notin \Exp(Q,v \asmid f_\gl), \qquad (\gl \in \Omega').
\end{equation}
Then $f = 0.$ 
\end{thm}

The proof of this theorem will be given after the following lemmas on which it is 
based. We may and shall assume that $Q$ is standard. Thus, $Q$ contains
the minimal standard $\gs$-parabolic subgroup $P_0$ which will 
be denoted by $P$ in the rest of this section.

\begin{lemma}
\naam{l: first step vanishing thm}
Let $\Omega \subset \faQqdc $ be a  non-empty connected open subset,
 $\gd \in \DQmaps$ and assume that $|\gd| =1.$ Let $Y \subset \staQqdc$ be a finite subset and
 let $f \in \cE_{Q,Y}(\spXp\col\tau\col\Omega\col \gd).$
 Moreover, let $v \in \NKaq$ and 
 assume that
 there
exist $t \in W_Q,\, \eta \in Y,$  $\mu \in \N\Delta$ and $u \in \NKQaq$ 
such that
\begin{equation}
\naam{e: gl s eta mu in Exp new}
\gl  + t\eta - \rho - \mu \in \Exp (P, uv \asmid f_\gl)
\end{equation}
for  $\gl$ in some non-empty open subset of $ \Omega.$ 
Then there exists a full open subset $\Omega_0 \subset \Omega$ such that 
$$
\gl - \rho_Q \in \Exp(Q,v \asmid f_\gl), \qquad (\gl \in \Omega_0).
$$
\end{lemma}

\proof
Let $\nu\in\stbQdc$ be the unique element such that $\supp\gd=\{\nu\}$.
Fix $t, \eta, \mu$ and $u$ with the mentioned property.
Replacing $\mu$ by a $\preceq_{\gD}$-smaller
element if necessary 
we may in addition assume that $\mu$ is $\preceq_{\gD}$-minimal subject to the condition
that (\refer{e: gl s eta mu in Exp new}) holds for 
 $\gl$ in some non-empty open subset of $\Omega.$ 
By holomorphy of asymptotics, see
Lemma \refer{l: holo of asymp},
it follows that (\refer{e: gl s eta mu in Exp new}) holds
for $\gl$ in a full open subset $\Omega'$ of $\Omega.$ 
Moreover, using the minimality of $\mu$ and applying 
Lemma \refer{l: exponents disjoint} we see that 
for every $\gl$ in the full open subset $\Omega_0:= \Omega' \cap \faQqdzero(P, WY)$ of $\Omega,$
$$
\gl + t\eta - \rho - \mu \in \ExpL(P,uv\asmid f_\gl).
$$
Since $f_\gl$ is annihilated by $\Igdgl = I_{\nu + \gl},$ this implies, 
in view of Lemma \refer{l: restriction on leading exponents},
that there exists a finite subset $\cL \subset \fbkdc$ such 
that
$$
\nu + \gl \in 
W(\fb) ( \cL + \gl + t \eta - \mu),\qquad (\gl \in \Omega_0).
$$ 
For  $\gL_0 \in \cL, \, w \in W(\fb)$ 
we define 
$\Omega_0(\gL_0, w)$ to be the set 
of $\gl \in \Omega_0$ satisfying 
\begin{equation}
\naam{e: defining relation for Omega 1 gL zero w}
\nu + \gl = w(\gL_0 + \gl + t \eta -\mu).
\end{equation}
The union of these sets, as $\gL_0 \in \cL,$ $w \in W(\fb),$ equals $\Omega_0.$ By finiteness
of the union, we may select $\gL_0$ and $w$ such that $\Omega_0(\gL_0,w)$ 
has a non-empty interior in $\Omega_0.$ Since $\Omega_0(\gL_0,w)$ 
is also the intersection of $\Omega_0$ with an affine linear subspace of $\fbdc,$ 
it must  be all of $\Omega_0.$ 
Hence for all $\gl_1,\gl_2 \in \Omega_0$ 
we have $w(\gl_1 - \gl_2) = \gl_1 - \gl_2.$ Since $\Omega_0$ is 
a non-empty open subset of $\faQqdc$ this implies that
$w$ belongs to $W_Q(\fb),$ the centralizer of $\faQq$ in $W(\fb).$
{}From (\refer{e: defining relation for Omega 1 gL zero w}) we now deduce that
$-w\mu = \nu -w\gL_0 -w t \eta.$ The expression on the right-hand side of this
equality 
has zero restriction to $\faQq.$ Therefore, so has $w\mu,$ 
and we conclude that also $\mu|_{\faQq} = 0.$ 
Combining
this fact with (\refer{e: gl s eta mu in Exp new})
and transitivity of asymptotics, see Theorem \refer{t: transitivity of asymptotics},
we
conclude that 
$$
\gl - \rho_Q =
[\gl + t\eta - \rho - \mu]|_{\faQq} \in \Exp(Q,v\asmid f_\gl),
$$ 
for
all $\gl \in \Omega_0.$ 
\qed

For the formulation of the next lemma, we need the following definition.

\begin{defi} 
\naam{d: special W set}
Let $\Omega\subset\faQqdc$ and $s_0\in W$ be given.
The subset $W(\Omega,s_0)$ of $W$ is defined as follows. Let $s'\in W$.
Then $s'\in W(\Omega,s_0)$ if and only if 
there exists a chain $s_1, \ldots, s_k=s'$ of elements in $W$, 
with $s_j s_{j-1}^{-1}=s_{\ga_j}$ 
a simple reflection, such that the following condition 
(\refer{e: condition of unboundedness}) holds
for each of the pairs $(s,\alpha)=(s_{j-1},\ga_j)\in W\times\gD$, 
$j=1,\ldots,k$.
\begin{equation}
\naam{e: condition of unboundedness}
\text{If $s^{-1}\ga|_{\faQq}\neq0$ then $\gl\mapsto\Re\inp{s\gl}{\ga}$ is
not bounded from below on $\Omega$.}
\end{equation}
\end{defi}

Notice that if $\Omega$ is dense in $\faQqdc$, then $W(\Omega,s_0)=W$
for all $s_0\in W$. Indeed, (\refer{e: condition of unboundedness}) is 
then fulfilled by all
elements $\ga\in\gD$. Hence, in order to verify the conditions
of Definition \refer{d: special W set} for $s'\in W$ arbitrary,
we may choose as $s_{\ga_1},\ldots,s_{\ga_k}$ 
the elements in a reduced expression 
$s's_0^{-1}=s_{\ga_k}\cdots s_{\ga_1}$.

\begin{lemma}
\naam{l: second step vanishing thm}
Let $\Omega\subset \faQqdc$ be a non-empty connected open subset, 
$Y \subset \staQqdc$ a finite subset, and $\gd \in \DQmaps.$
Let $f \in \cE_{Q,Y}(\spXp\col\tau\col\Omega\col \gd)_\glob$ 
and $s \in W.$ Assume that there exist
$t \in W_Q,\, \eta \in Y$, $\mu \in \N\Delta$ and
$w \in \NKaq$
such that 
\begin{equation}
\naam{e: exponents in Exp P w f new}
s\gl + st \eta - \rho - \mu \in \Exp(P,  w\asmid f_\gl),
\end{equation}
for all $\gl$ in some non-empty open subset of $\Omega.$
Then for every $s_1\in W(\Omega,s)$ 
there exist
$t_1 \in W_Q,\, \eta_1 \in Y,$  $\mu_1 \in \N\Delta$ and $w_1 \in \NKaq,$
such that 
\begin{equation}
\naam{e: one exponents in Exp P w f new}
s_1\gl + s_1t_1 \eta_1 - \rho - \mu_1 \in \Exp(P, w_1\asmid f_\gl),
\end{equation}
for all $\gl$ in a full open subset of $\Omega.$ 
In particular, if $\Omega$ is dense in $\faQqdc$, then the above conclusion
holds for every $s_1\in W$. 
\end{lemma}

\proof
In the proof we will frequently use the following consequence of
Lemma \refer{l: holo of asymp}, based on holomorphy of asymptotics.
If $s_1\in W,$ 
$t_1 \in W_Q,\, \eta_1 \in Y,$  $\mu_1 \in \N\Delta$ and $w_1 \in \NKaq,$
then (\refer{e: one exponents in Exp P w f new}) holds for $\gl$ in a full open
subset of $\Omega$ as soon as it holds for a fixed $\gl$ in the full open subset 
$\Omega \cap \faQqdczero(P, WY)$  of $\Omega.$ We now turn to the proof.

If $s_1 = s,$ or more generally, if $s_1\in sW_Q$,
then the conclusion readily follows by the previous remark.
By Definition \refer{d: special W set} we now  see that it suffices to
prove the lemma for $s_1 = s_\ga s,$ with $\ga\in \Delta$ such that 
(\refer{e: condition of unboundedness}) holds.
There are two cases to consider, namely that $s^{-1} \ga |_{ \faQqd}$ 
equals zero or not.
In the first case, $s_1 = s s_{s^{-1} \ga}\in sW_Q$ and the conclusion
is valid.
We may thus assume that we are in the second case, i.e., $s_1=s_\ga s$ with
\begin{equation}
\naam{e: s inv ga in min gsQ new}
s^{-1} \ga |_{ \faQqd}\neq 0.
\end{equation}
We will complete the proof by showing that the following assumption leads to a contradiction.

{\bf Assumption:} for all $t_1 \in W_Q,$ $\eta_1 \in Y,$ $\mu_1\in \N \Delta$ and
$w_1 \in \NKaq$ there exists
no non-empty open subset $\Omega'$ of $\Omega$ such that 
(\refer{e: one exponents in Exp P w f new}) holds for $\gl \in \Omega'.$

Let $\expset$ be the set of elements $(st\eta - \mu)|_{\fagaq}$ with
$t \in W_Q, \eta \in Y, \mu \in \N \Delta$ such that 
(\refer{e: exponents in Exp P w f new}) holds 
for  $\gl$ in a non-empty open subset of $\Omega$,
for some $w\in \NKaq$. Then $\expset$ is a non-empty subset 
of $\fagaqdc$ contained in a set of the form $X - \N \Dr({\stPga}),$ with
$X \subset \fagaqdc$ finite. Hence we may select $t \in W_Q, \eta \in Y$ and
$\mu \in \N \gD$ such that $(st\eta -\mu)|_{\fagaq}$ is 
$\preceq_{\Dr({\stPga})}$-maximal in $\expset.$ 
According to the first paragraph of the proof,
there exists $w\in\NKaq$ such that
(\refer{e: exponents in Exp P w f new}) is
valid for  $\gl$ in a full open subset $\Omega_0$ of $\Omega.$ 
For $\gl \in \Omega_0$ we put 
$$
\xi(\gl) = [s\gl + st \eta - \rho - \mu]|_{\fagaq}.
$$
Then by  transitivity of asymptotics, see 
Theorem \refer{t: transitivity of asymptotics},
it follows that 
$$
\xi(\gl) \in \Exp({\stPga}, w\asmid f_\gl)
$$ 
for $\gl \in \Omega_0.$ In the following we shall investigate the  
coefficient of the expansion of $f_\gl$ along
$({\stPga}, w),$ for $\gl \in \Omega_0,$ given by 
$$
\gf_\gl (m) := q_{\xi(\gl)}({\stPga},w\asmid f_\gl, \dotvar, m).
$$
Here $\gf_\gl$ is a non-trivial $\tauQga$-spherical function 
on $\spXgawp$ with values in $P_k(\fagaq),$ for $k=\dega f$, see
Thm.\ \refer{t: expansion along the walls} (b). 

It follows from (\refer{e: s inv ga in min gsQ new}) and 
the asymptotic globality assumption on $f,$ see 
Lemma \refer{l: minimal condition for glob new},
that actually $\gf_\gl$ extends
to a smooth function on $\spXgaw,$ for every $\gl$ in an dense open subset
$\Omega_0'$ of $\Omega_0.$ 
This observation will play a crucial role at a later stage of this proof. 

Let 
$$
\Omega_1  := \Omega_0' \cap \faQqdczero(P, WY) \cap \faQqdczero(\stPga, WY). 
$$
The second and third set in this intersection 
are full open subset of $\faQqdc,$ 
see Lemma \refer{l: azero is full}.
Hence $\Omega_1$ is a dense open subset of $\Omega.$
We claim that for $\gl \in \Omega_1$ the following holds.
If $s' \in W,\, t' \in W_Q, \,\eta' \in Y,$ $\mu' \in \N\gD$  and $w'\in \NKaq$
are such that
\begin{equation}
\naam{e: hypothesis claim on exponents new}
\left\{
\begin{array}{l}
s'\gl + s't'\eta' -\rho -\mu' \in \Exp(P, w'\asmid f_\gl) 
\text{and}\\
\xi(\gl) \preceq_{\Dr({\stPga})} (s'\gl + s't'\eta' -\rho -\mu')|_\fagaq, 
\end{array}\right.
\end{equation}
then 
\begin{equation}
\naam{e: conclusion claim on exponents}
s' \in  sW_Q  \text{and} (s't'\eta' - \mu')|_{\fagaq} = (st \eta - \mu)|_{\fagaq}.
\end{equation}
To prove the claim, let $s',t',\eta',\mu',w'$ satisfy (\refer{e: hypothesis claim on exponents new}).
Then there exists a $\nu \in \N\gD({\stPga})$ 
such that $s'\gl + s't'\eta' -\rho -\mu'- \nu$ and $s\gl + st \eta -\rho - \mu$ 
have the same restriction $\xi(\gl)$ to $\fagaq.$ By the definition of $\Omega_1$ 
this implies that $s'$ and $s$ define the same class in $W/\sim_{\stPga|Q},$ 
see Lemma \refer{l: exponents disjoint}. 
The latter set equals $W_\ga\backslash W/W_Q,$ by Lemma \refer{l: WPQ as cosets},
hence $s'$ belongs to $s_\ga s W_Q=s_1W_Q$ or to $ s W_Q.$ In the first case
it follows that $s' \gl= s_1 \gl,$ hence 
$s_1 \gl + s_1 t''\eta' -\rho -\mu' \in \Exp(P,  w'\asmid f_\gl)$
for some $t''\in W_Q$. This assertion 
then 
holds for $\gl$ in a full open subset of $\Omega_1,$ contradicting the above
assumption.

It follows that we are in the second case $s' \in s W_Q,$
hence $s' = s t''$ for some $t'' \in W_Q.$ 
The element 
$(s't'\eta' - \mu')|_{\fagaq} = (st''t'\eta' - \mu')|_{\fagaq}$ 
therefore belongs to $\expset;$ from (\refer{e: hypothesis claim on exponents new})
it follows that it dominates
the maximal element $(st \eta -\mu)|_{\fagaq},$ hence is equal to that
element. This implies (\refer{e: conclusion claim on exponents}), hence establishes
the claim.

It follows from the above claim that, for $\gl \in \Omega_1,$ 
the exponent $\xi(\gl)$ is 
actually a leading exponent of $f_\gl$ along $({\stPga}, w).$ 
To see this, let $\gl \in \Omega_1$ and let $\xi \in \Exp({\stPga}, w \asmid f_\gl)$ 
be an exponent with  $\xi(\gl) \preceq_{\Dr({\stPga})} \xi.$
Then, in view of Theorem \refer{t: transitivity of asymptotics},
there exist $s' \in W,\, t' \in W_Q, \,\eta' \in Y$ and $\mu' \in \N\gD$ 
such that the element $s'\gl + s't'\eta' -\rho -\mu'$ restricts to $\xi$ on $\fagaq$ 
and belongs to $\Exp(P, w'\asmid f_\gl)$ 
for some $w' \in \cW.$ It now follows from 
the claim established above that $\xi = \xi(\gl).$ 

Thus, we see that $\xi(\gl)$ is a leading exponent indeed.
Consequently, by 
Lemma \refer{l: D finiteness of leading coefficient} 
the function $\gf_\gl$ is $\D(\spXgaw)$-finite, for every $\gl \in \Omega_1.$ 
We proceed by investigating the exponents of its expansion.
 
Select a complete set $\cW_{\ga, w}$  of representatives for 
$W_\ga / (W_\ga \cap W_{K \cap wHw^{-1}})$ in $\NKaq.$ 
We put $\starP = P \cap M_\ga.$ Then by transitivity 
of asymptotics, cf.\ Theorem \refer{t: transitivity of asymptotics}, 
we see that for  
the set of $(\starP, u)$-exponents of $\gf_\gl,$ as $u\in \cW_{\ga,w},$ 
the following inclusion
holds:
$$
\Exp(\starP, u\asmid \gf_\gl) \subset \{ \xi|_{\stagaq}\setmid \xi \in \Exp(P, uw\asmid f_\gl)
\;\;\; \xi | \fagaq= \xi(\gl) | \fagaq \}.
$$
Hence, for  $\gl\in \Omega_1,$ 
every exponent in $\Exp(\starP,u \asmid \gf_\gl)$ 
is of the form 
$(s' \gl + s't'\eta' - \rho - \mu')|_{\stagaq} $ with
$s' \in W, \;t' \in W_Q, \;\eta' \in Y$ and $\mu' \in \N\Delta$ 
satisfying 
$$
\left\{
\begin{array}{l}
s'\gl + s't'\eta' -\rho -\mu' \in \Exp(P, u w\asmid f_\gl),\\{}  
[s'\gl + s't'\eta' - \rho - \mu']|_{\fagaq} 
= \xi(\gl)|_{\fagaq}. 
\end{array}
\right.
$$
It follows from the claim established above that 
(\refer{e: conclusion claim on exponents}) holds.  

We have thus shown that for every $\gl \in \Omega_1$  the exponents in 
$\Exp(\starP,u \asmid \gf_\gl)$ 
are  of the form 
$ (s \gl + st'\eta' - \rho - \mu' )|_{\stagaq}$ with
$t' \in W_Q,\; \eta' \in Y,\;\mu' \in \N\Delta$ satisfying
$$
[st'\eta' - \mu'] |_{\fagaq} = [st \eta - \mu ]|_{\fagaq}.
$$
{}From this it follows that
the restriction $\mu'|_{\fagaq}$ of the $\mu'$ occurring 
runs through a finite subset of 
$\N\Dr(\stPga)= \N[\gD\setminus\{\ga\}]|_\fagaq,$ independent of $\gl$. 
Hence there exists a finite subset $S' \subset \N\gD$ 
such that $\mu'$ runs through $S' - \N \ga.$ We thus see
that there exists  
a finite subset $S \subset \stagaqdc$ 
such that, for every $\gl \in \Omega_1,$ 
\begin{equation}
\naam{e: inclusions for Exp starP u gf}
\cup_{u \in \cW_{\ga,w}}\;
\Exp(\starP,u\asmid \gf_\gl) \subset  s\gl |_{\stagaq} + S - \N\ga.
\end{equation}
{}From (\refer{e: condition of unboundedness}) and
(\refer{e: s inv ga in min gsQ new}) it now follows that we may select
a non-empty open 
subset $\Omega_2$ of the dense  open subset $\Omega_1$ of $\Omega$ 
such that, for every 
$\gl \in \Omega_2,$ each $u \in \cW_{\ga,w}$ 
and all $\xi \in \Exp(\starP, u\asmid \gf_\gl),$ 
$$
\inp{\Re \xi + {}^* \rho }{\ga} < 0.
$$
Since $\gf_\gl$ is $\D(\spXgaw)$-finite this implies that
$\gf_\gl$ is square integrable on $\spXgaw,$  
see \bib{B87}, Thm. 6.4 with $p =2;$ hence $\gf_\gl$ a Schwartz function
for $\gl \in \Omega_2,$
see \bib{B87}, Thm.\ 7.3.

On 
the other hand, from (\refer{e: s inv ga in min gsQ new}) it follows that 
the linear map $\gl \mapsto s\gl |_{\stagaq}$ is surjective
from $\faQqdc$ onto $\stagaqdc.$ 
Therefore, the set
$\{s\gl|_{\stagaq}\mid \gl \in \Omega_2\}$ has a non-empty interior
in $\stagaqdc.$ Combining this observation 
with (\refer{e: inclusions for Exp starP u gf}) we infer that there exists a non-empty 
open subset $\Omega_3 \subset \Omega_2,$ such that the sets 
$\cup_{u \in \cW_{\ga,w}}\;
\Exp(\starP,u\asmid \gf_\gl),$ for $\gl \in \Omega_3,$ are mutually disjoint.
Now these sets are non-empty, since $\gf_\gl \neq 0,$ for $\gl\in \Omega_3.$
Therefore,  the union of these sets, as $\gl \in \Omega_3,$ is 
uncountable. This contradicts Lemma 
\refer{l: limitation on exponents of a Schwartz function},
applied to the space $\spXgaw.$ 
\qed

\begin{lemma}
\naam{l: e in W set}
Assume that $\Omega\subset\faQqdc$ is $Q$-distinguished.
Then $e\in W(\Omega,s_0)$ for all $s_0\in W$.
\end{lemma}

\proof Let $k=l(s_0)$ denote the length of $s_0$, and let 
$s_0=s_{\ga_1}\cdots s_{\ga_k}$ be a reduced expression for $s_0$.
Put $s_j=s_{\ga_j}\cdots s_{\ga_1}s_0=s_{\ga_{j+1}}\cdots s_{\ga_k}$ 
for $j=1,\ldots,k$, then $s_k=e$. We claim that 
(\refer{e: condition of unboundedness}) holds for each pair 
$(s,\ga)=(s_{j-1},\ga_j)$. Since $l(s_j)=l(s_{j-1})-1$, the root
$s_{j-1}^{-1}\ga_j$ must be negative. Hence the restriction of
this root to $\faQq$ is zero or belongs to $-\gS(Q)$. Now 
(\refer{e: condition of unboundedness})
follows immediately from Definition \refer{d: Q-distinguished}.\qed

{\bf Proof of Theorem \refer{t: vanishing theorem new}:\ }
We prove the result  by induction on $|\gdmap|.$ 
If $\gd = 0,$ then for $\gl \in \faQqdc$ the ideal 
$\Igdgl$ equals $\DGH;$ 
hence $\cE_Q(\spXp\col \tau \col \Omega \col \gd)_\glob = 0$ and the result follows.

Let now $|\gd| =1,$ let  
$f \in \cE_Q(\spXp\col \tau\col \Omega\col \gd)_\glob$ and let 
(\refer{e: hypothesis vanishing thm new}) be fulfilled for all $v \in \QW.$  
Assume that  $f \neq 0.$ We will show that this assumption leads to a contradiction. 
There exists a finite subset $Y\subset \staQqdc$ 
such that $f \in \cEQY(\spXp\col \tau\col \Omega\col \gd)_\glob$ and
a $\gl_0 \in \Omega \cap \faQqdczero(P, WY)$ 
such that $f_{\gl_0} \neq 0.$ Let $\cW$ be a complete set of representatives
of $W/\WKH$ in $\NKaq$ containing $\QW.$ Then 
$\Exp(P, w\asmid f_{\gl_0}) \neq \emptyset$ for some $w \in \cW.$ 
In view of (\refer{e: exponents family}) 
it  follows that there exist 
$s \in W,$ $t\in W_Q,$ $\eta \in Y$ and $\mu \in \N \gD,$
such that 
\begin{equation}
\naam{e: exponent for gl is gl zero}
s \gl+ st \eta - \rho -\mu \in \Exp(P, w\asmid f_\gl),
\end{equation}
for $\gl = \gl_0.$ From Lemma \refer{l: holo of asymp} 
it follows that (\refer{e: exponent for gl is gl zero}) is valid for $\gl$ 
in a full
open subset of $\Omega.$  
By Lemmas \refer{l: second step vanishing thm} and \refer{l:  e in W set} 
this implies that there exist 
$t_1 \in W_Q, $ $\eta_1 \in Y,$ $\mu_1 \in \N\gD$ and $w_1 \in \NKaq,$ such that 
$\gl + t_1 \eta_1 - \rho -\mu_1 \in \Exp(P, w_1 \asmid f_\gl)$ 
for $\gl$ in a full open subset of $\Omega.$ 
Let $v\in\QW$ be the representative of $W_Qw_1W_{K\cap H}$.
By Lemma \refer{l: first step vanishing thm} 
it follows that 
$\gl- \rho_Q \in \Exp(Q, v\asmid f_\gl)$ 
for $\gl$ in a full open subset $\Omega_0$ of $\Omega.$ 
Since $\Omega_0 \cap \Omega'$ is non-empty, we obtain 
a contradiction with (\refer{e: hypothesis vanishing thm new}). 

Now suppose that $|\gdmap| = k > 1,$ and assume that the result 
has already been established for $\gdmap \in \DQmaps$ with $|\gdmap| < k.$ 
Fix $\nu \in \supp (\gdmap)$ and put 
 $\gdmap'  = \gdmap - \gdmap_\nu.$ Then $\gdmap' \in \DQmaps;$ 
moreover, $|\gd_\nu| =1$ and $|\gd'| = k-1.$ 
Fix any $D \in \DGH$ and define the family $g$ by 
(\refer{e: defi g as D minus gamma f}).
Then $g \in \cE_{Q}(\spXp\col \tau\col \Omega\col \gd')$ 
by Lemma \refer{l: stability family under D new}.
Moreover, it readily follows from 
Lemma \refer{l: stability family with globality under D} 
that the family $g$ belongs to 
$\cE_{Q}(\spXp\col \tau\col\Omega\col  \gd')_\glob.$ 

For $\gl \in \Omega$ and $v \in \NKaq$ we have 
\begin{equation}
\naam{e: inclusion exponents F gl}
\Exp(Q,v\asmid g_\gl) \subset \Exp(Q, v \asmid f_\gl) - \N \Sigma_r(Q),
\end{equation}
in view of Lemma \refer{l: radial component applied to expansion} (b).
Moreover, by hypothesis we have the following inclusion, for every 
$\gl \in \Omega',$ 
\begin{equation}
\naam{e: inclusion exponents f gl}
\Exp(Q,v\asmid f_\gl) \subset [W(\gl + Y)|_{\faQq}  - \rho_Q - \N \Sigma_r(Q)]
\setminus \{\gl - \rho_Q\}.
\end{equation}
Combining (\refer{e: inclusion exponents F gl}) and 
(\refer{e: inclusion exponents f gl}) we infer that $\Exp(Q,w\asmid g_\gl)$ 
does not contain $\gl - \rho_Q$ for
$\gl \in \Omega'$ and every $w \in \NKaq.$  
Consequently, the family $g$ satisfies the hypotheses 
of Theorem \refer{t: vanishing theorem new}. 
Since $|\gdmap'| = k-1,$ it follows from the induction hypothesis 
that $g =0.$ Since $D$ was arbitrary, we see that $f_\gl$ is 
annihilated by $I_{\gdmap_\nu, \gl},$ for every $\gl \in \Omega.$ 
Hence $f$ belongs to $\cE_{Q}(\spXp\col \tau\col\Omega\col \gdmap_\nu)_\glob.$ 
Since $|\gdmap_\nu| =1 < k,$ it now follows from the induction hypothesis that 
$f=0.$
\qed

The following result is also based on Lemma 
\refer{l: second step vanishing thm}.

\begin{cor}
\naam{c: variant of vanishing theorem}
Let $\Omega\subset \faQqdc$ be a connected dense open subset, 
$Y \subset \staQqdc$ a finite subset, and $\gd \in \DQmaps.$
Let $f \in \cE_{Q,Y}(\spXp\col\tau\col\Omega\col \gd)_\glob$
and $s_1\in W$. If 
$$
(s_1\gl + WY - \rho - \N\gD)\cap \Exp(P, w\asmid f_\gl)=\emptyset,
$$
for all $\gl$ in a non-empty open subset of $\Omega$ 
and for all $w\in N_K(\faq)$, then $f=0$.
\end{cor}

\proof Assume that $f\neq0$. Then there exists an element
$\gl\in \Omega \cap \faQqdczero(P, WY)$ such that $f_\gl\neq0$,
and then 
\begin{equation}
\naam{e: not an inclusion}
s\gl + st \eta - \rho - \mu \in \Exp(P,  w\asmid f_\gl)
\end{equation}
for some $s\in W$, $t \in W_Q,\, \eta \in Y$, $\mu \in \N\Delta$ and
$w\in N_K(\faq)$. 
As remarked in the beginning of the proof of 
Lemma \refer{l: second step vanishing thm},
(\refer{e: not an inclusion})
then holds for all $\gl$ in a full open subset of $\Omega$. 
Hence Lemma \refer{l: second step vanishing thm} applies;
its final statement contradicts the present assumption for $s_1$.\qed

Finally in this section we will show that for a family in 
$\cE_{Q}(\spXp\col\tau\col\Omega\col \gd)$ that allows a smooth
extension to $\spX$, the hypothesis of asymptotic globality can be
left out in the vanishing theorem. Let
$$\cE_{Q}(\spX\col\tau\col\Omega\col \gd)=
\{f\in\cE_{Q}(\spXp\col\tau\col\Omega\col \gd)\mid
f_\gl\in C^\infty(\spX\col\tau), \gl\in\Omega\}.$$

\begin{cor}
Let $Q \in \allparabs$ and $\gdmap \in \DQmaps.$ 
Let $\Omega \subset \faQqdc$ 
be a $Q$-distinguished open subset and let 
$f \in \cE_Q(\spX\col\tau\col\Omega\col \gd).$ 
Assume that there exists a non-empty open subset 
$\Omega' \subset \Omega$ such that, for each $v \in \QW, $ 
$$
\gl - \rho_Q \notin \Exp(Q,v \asmid f_\gl), \qquad (\gl \in \Omega').
$$
Then $f = 0.$ 
\end{cor}

\proof As in the proof of Theorem \refer{t: vanishing theorem new}
we proceed by induction on $|\gd|$. If $|\gd|=0$ the result is trivial.
If $|\gd|=1$ it follows from Proposition 
\refer{p: global eigen implies as global} that
$\cE_Q(\spX\col\tau\col\Omega\col \gd)\subset
\cE_Q(\spXp\col\tau\col\Omega\col \gd)_\glob$, and then the result
follows directly from Theorem \refer{t: vanishing theorem new}.

Now suppose that $|\gd|=k>1$, and assume that the result has already been 
established for all $\gd\in D_Q$ with $|\gd|<k$. Let $\gd'$ and $g$
be as in the proof of Theorem \refer{t: vanishing theorem new}.
Then it is easily seen that $g\in \cE_Q(\spX\col\tau\col\Omega\col \gd')$.

For the rest of the proof we can now proceed exactly as in the proof of
Theorem \refer{t: vanishing theorem new}.\qed

\section{Laurent functionals}
\naam{s: Laurent functionals}
In this section we define Laurent functionals and describe their
actions on suitable spaces of meromorphic functions.  

Throughout this section, $V$ will be a finite dimensional 
real linear space, equipped with a (positive definite) inner product $\inp{\dotvar}{\dotvar}.$ 
Its complexification $V_\biC$ is equipped with the complex bilinear 
extension of this inner product.

Let $X$ be a (possibly empty) finite set of non-zero elements of $V.$ 
At this stage we allow proportionality between elements of $X.$ 
By an $X$-hyperplane in $V_\biC,$ we mean an affine hyperplane of the form
$H = a + \ga_{\iC}^\perp,$ with $a \in V_\biC,$ $\ga \in X.$ The hyperplane is
called real if $a$ can be chosen from $V,$ or, equivalently, if it is 
the complexification of a real hyperplane from $V.$ 
A locally finite collection of $X$-hyperplanes in $V_\biC$ is called
an $X$-configuration in $V_\biC.$ It is called real if all its
hyperplanes are real.

If $a \in V_\biC,$ we denote the collection of $X$-hyperplanes in $V_\biC$ 
through $a$ by $\Hyp(a, X) = \Hyp(V_\biC, a, X).$ 
If $E$ is a complete locally convex space, then by 
$\Mer(a, X,E) = \Mer(V_\biC, a, X,E)$ 
we denote the ring of germs of $E$-valued meromorphic functions at $a$ whose singular 
locus at $a$ is contained in $\Hyp(a,X).$ Here and in the following we will suppress 
the space $E$ in the notation if $E = \C.$ Thus, $\Mer(a,X) = \Mer(a,X,\C).$ 

Let $\N^X$ denote the set of maps $X \to \N.$ If $d \in \N^X,$ 
we define the polynomial function $\pi_{a,d} = \pi_{a,X,d}: V_\biC \to \C$ by 
\begin{equation}
\naam{e: defi pi a X d} 
\pi_{a,d}(z) =\prod_{\xi \in X} \inp{\xi}{z -a}^{d(\xi)}, \qquad (z \in V_\biC).
\end{equation}
If $X=\emptyset$ then $\N^X$ has one element which we agree to denote by $0$.
We also agree that $\pi_{a,0}=1$.
Let $\cO_a(E)= \cO_a(V_\biC, E)$ denote the ring of germs of $E$-valued holomorphic
functions at $a.$ Then 
$$
\Mer(a,X,E) = \cup_{d \in \N^X} \;\; \pi_{a,d}^{-1} \cO_a(E).
$$ 
In the following we shall identify $S(V)$ with the algebra of constant 
coefficient holomorphic differential operators on $V_\biC$ in the usual way;
in particular an element
$v \in V$ corresponds to the operator 
$\gf \mapsto v\gf(z) = \left.\frac{d}{d\tau}\right|_{\tau=0} \gf(z + \tau v).$
\begin{defi}
\naam{d: Laurent functional at a point}
{\rm (Laurent functional at a point)\ }
An $X$-Laurent functional at $a$ is a linear functional $\Lau: \Mer(a,X) \to \C$ 
such that for every $d \in \N^X$ there exists an element  $u_d \in S(V)$ such that 
\begin{equation}
\naam{e: Lau and u d}
\Lau \gf = u_d(\pi_{a,d} \gf)(a),
\end{equation} 
for all $\gf \in \pi_{a,d}^{-1} \cO_a.$ 
The space of all Laurent functionals at $a$ is denoted by $\Mer(a, X)^*_\laur =
\Mer(V_\biC, a, X)^*_\laur.$ 
\end{defi}

\begin{rem}
Obviously, the string  $(u_d)_{d \in \N^X}$ of elements from  $S(V)$ is 
uniquely determined by the requirement (\refer{e: Lau and u d}). We shall denote 
it by $u_\cL.$ 

If $E$ is a complete locally convex space, then $X$-Laurent functionals 
at $a$ may naturally be viewed as linear maps from $\Mer(a,X,E)$ to $E$.
Indeed, let $\Lau \in \Mer(a, X)^*_\laur$ and let 
$u_\Lau = (u_d)_{d \in \N^X}$ be the associated string of elements from $S(V).$ 
If  $\gf \in \pi_{a,d}^{-1}\cO_a(E)$ then  $\Lau \gf$ 
is given by formula (\refer{e: Lau and u d}).
\end{rem}

Let $T_a: z \mapsto z + a$ denote translation by $a$ in $V_\biC.$ 
Then $T_a$ maps $\Hyp(0,X)$ bijectively onto $\Hyp(a,X).$ 
Pull-back under $T_a$ induces an isomorphism of rings $T_a^*: \gf \mapsto \gf\after T_a$ 
from $\cO_a$ onto $\cO_0.$ Therefore, pull-back under $T_a$ also induces
an isomorphism of rings $T_a^*: \Mer(a, X) \to \Mer(0,X).$ 
By transposition we obtain an isomorphism of linear spaces 
$T_{a*}: \Mer(0, X)^* \to \Mer(a,X)^*.$ It is readily seen that $T_a^*(\pi_{a,d}) 
=\pi_{0,d}$ for every $d \in \N^X.$ {}From the definition of Laurent functionals
it now follows that $T_{a*}$ maps $\Mer(0,X)^*_\laur$ isomorphically 
onto $\Mer(a,X)^*_\laur.$ Moreover, 
$$%\begin{equation}
%\naam{e: translation of u Lau}
u_{T_{a*}\Lau} = u_{\Lau}
$$%\end{equation} 
for all $\Lau \in \Mer(0,X)^*.$ 

Let $X'$ be another finite collection of non-zero elements of $V.$ We say 
that $X$ and $X'$ are proportional if $\Hyp(0,X) = \Hyp(0,X').$ 

\begin{lemma}
\naam{l: Laurent functionals and proportional X}
Let $X,X'$ be proportional finite subsets of $V\setminus\{0\}$ 
and let $a \in V_\biC.$ Then $\Mer(a, X) = \Mer(a, X')$ and
$\Mer(a,X)^*_\laur = \Mer(a,X')^*_\laur.$ 
\end{lemma}

\proof It is obvious that $\Mer(a, X) = \Mer(a, X')$. 
Let $\Lau \in \Mer(a,X)^*=\Mer(a, X')^*$, and assume that
$\Lau \in \Mer(a,X')^*_\laur.$ Let $(u_{d'})_{d'\in\N^{X'}}$
be the associated string. Let $d\in\N^X$. 
Then, by proportionality, there exists $d'\in\N^{X'}$ and 
$c\in\R\setminus\{0\}$ such that $\pi_{a,X,d}=c\pi_{a,X',d'}.$
Let $u_d=c^{-1}u_{d'}$, then (\refer{e: Lau and u d}) follows immediately.
This shows that $\Lau \in \Mer(a,X)^*_\laur$ 
and  establishes the inclusion $\Mer(a,X')^*_\laur \subset \Mer(a,X)^*_\laur.$ 
The converse inclusion is proved similarly.
\qed

Following 
the method 
of \bib{BSres}, Sect.\ 1.3,
we shall now give 
a description of the space of strings 
$u_\Lau,$ as $\Lau \in \Mer(a,X)^*_\laur$.

Put $\varpi_d := \pi_{0,d}$ and equip the space $\N^X$ with the partial ordering 
$\preceq$ defined by $d' \preceq d$ if and only if $d'(\xi) \leq d(\xi)$ for
every $\xi \in X.$ 
If $d'\preceq d$ then we define $d -d'$ componentwise as suggested by the notation.
In \bib{BSres}, Sect.\ 1.3, we defined the linear space
$S_{\leftarrow}(V,X)$ as follows.
Let $d,d' \in \N^X$ with $d' \preceq d.$ If $u \in S(V),$ then by the Leibniz
rule there exists a unique $u' \in S(V)$ such that
$$%\begin{equation}
%\naam{e: condition strings}
u(\varpi_{d - d'}\gf)(0) = u'(\gf)(0), \qquad (\gf \in \cO_0).
$$%\end{equation}
We denote the element $u'$ by $j_{d',d}(u).$ The map $j_{d',d}: S(V) \to S(V)$ thus defined 
is linear. Note
that it only depends on $d - d';$ note also that, for $d,d',d'' \in \N^X$ with 
$d'' \preceq d' \preceq d,$
$$%\begin{equation}
%\naam{e: composition j maps}
j_{d'',d'}\after j_{d', d} = j_{d'', d}.
$$%\end{equation}
We now define $\SprojVX$ as the linear space of strings $(u_d)_{d \in \N^X}$ in 
$S(V)$ such that $j_{d',d}( u_d) = u_{d'}$ for all $d,d \in \N^X$  with $d' \preceq d.$ 
Thus, this space is the projective limit:
$$ 
\SprojVX  = \lim_{\leftarrow} (S(V), j_{\cdot}).
$$
The natural map $\SprojVX \to S(V)$ that maps a string to its $d$-component
is denoted by $j_d.$   

\begin{lemma}
\naam{l: iso Laurent functionals with Sproj}
The map $\Lau \mapsto u_\Lau$ is a linear isomorphism from $\Mer(a,X)^*_\laur$ onto 
$S_\leftarrow(V,X).$ 
\end{lemma}

\proof 
See \bib{BSfi}, Appendix B, Lemma B.2.
\qed

\begin{lemma}
\naam{l: extension to Laurent functional}
Let $a \in V_\biC,$ $d \in \N^X$ and $u \in S(V).$ 
Then there exists a Laurent functional 
$\Lau \in \Mer(a, X)^*_\laur$ such that $(u_\Lau)_d = u.$ 
\end{lemma}
\proof
See \bib{BSres}, Lemma 1.7.
\qed

\begin{rem}
\naam{r: evaluation is a Laurent functional}
In particular, it follows that for each $a\in V_\C$ there 
exists a Laurent functional $\cL\in\cM(a,X)^*_\laur$ such that
$\cL\gf=\gf(a)$ for all $\gf\in\cO_a$. Note however, that this
functional is not unique, unless $X=\emptyset$.
\end{rem}

\begin{lemma} 
\naam{l: annihilator of annihilator}
Let $\cM(a,X)^{*\cO}_\laur$ denote the annihilator
of $\cO_a$ in $\cM(a,X)^*_\laur$. Then all functions $\gf$ in
$\cM(a,X)$, that are annihilated by $\cM(a,X)^{*\cO}_\laur$,
belong to $\cO_a$.
\end{lemma}

\proof We may assume that $a=0$. Let $\gf\in\cM(0,X)$ 
and assume that $\gf\not\in\cO_0$. Then there exist elements
$d, d'\in\N^X$ and $\xi\in X$ such that $\pi_{0,d'}=\xi\pi_{0,d}$
and $\pi_{0,d'}\gf\in\cO_0$ but $\pi_{0,d}\gf\not\in\cO_0$. 
Here we have written $\xi$ also for the function $z\mapsto
\inp{\xi}{z}$ on $\Vc$.
Since $\pi_{0,d'}\gf$ is not divisible by $\xi$, its
restriction to $\xi^\perp=\xi^{-1}(0)$ does not vanish.
Hence there exists $u\in S(\xi^\perp)$ such that 
$u(\pi_{0,d'}\gf)(0)\neq 0$. 
By Lemma \refer{l: extension to Laurent functional} there exists an element 
$\cL\in\cM(a,X)^*_\laur$ such that the $d'$ term of $u_\Lau$ is $u$. Then
$\cL\gf=u(\pi_{0,d'}\gf)(0)\neq0$.
However, for each $\psi\in\cO_0$
we have $\cL\psi=u(\pi_{0,d'}\psi)(0)=[\xi u(\pi_{0,d}\psi)](0)=0$.
Hence $\cL\in\cM(a,X)^{*\cO}_\laur$.\qed

We extend the notion of a Laurent functional as follows.
The disjoint union of the spaces $\Mer(a,X)^*_\laur$ as $a \in V_\biC$ is denoted 
by $\Mer(*,X)^*_\laur = \Mer(V_\biC, *, X)^*_\laur.$ 
By a section of $\Mer(*, X)^*_\laur$ 
we mean a map $\Lau: V_\biC \to \Mer(*,X)^*_\laur$ with $\Lau_a \in \Mer(a,X)^*_\laur$ 
for all $a \in V_\biC.$  
The closure of the set $\{a \in V_\biC\mid \Lau_a \neq 0\}$ is called the support 
of $\Lau$ and denoted by $\supp(\Lau).$ 

\begin{defi}
{\rm (Laurent functional)}
An $X$-Laurent functional on $V_\biC$ is a finitely supported section of 
$\Mer(*, X)^*_\laur.$
The set of $X$-Laurent functionals is denoted by $\Mer(V_\biC, X)^*_\laur$
and equipped with the obvious structure of a linear space.

Is $S$ is a subset of $V_\biC,$ we define the space 
$\Mer(S, X)^*_\laur = \Mer(V_\biC, S, X)^*_\laur$ 
by 
$$
\Mer(S, X)^*_\laur = \{ \Lau \in \Mer(V_\biC, X)^*_\laur \mid \supp \Lau \subset S \}
$$
and call this the space of $X$-Laurent functionals on $V_\biC$ supported in $S.$ 
\end{defi}

\begin{rem}
\naam{r: Laurent functionals at a as subspace}
Note
that, for $a \in V_\biC,$ the map $\Mer(\{a\}, X)^*_\laur \to \Mer(a, X)^*_\laur,$ 
defined by $\Lau \mapsto \Lau_a,$ is a linear isomorphism. Accordingly we shall
view $\Mer(a, X)^*_\laur$ as a linear subspace of $\Mer(V_\biC, X)^*_\laur.$ 
In this way $\Mer(S, X)_\laur^*$ becomes identified with 
the algebraic direct sum of the linear spaces $\Mer(a, X)_\laur^*,$ as $a \in S,$ 
for $S$ any subset of $V_\biC.$ Accordingly,
if $\Lau \in \Mer(V_\biC, X)^*_\laur,$ then $\Lau_a \in 
 \Mer(a, X)^*_\laur \subset  \Mer(V_\biC, X)^*_\laur$ for $a \in V_\biC,$ and
$$
\Lau = \sum_{a \in \supp \Lau} \Lau_a.
$$ 
\end{rem}

\begin{lemma}
\naam{l: global Laurent functionals and proportional X}
Let $X$ and $X'$ be proportional finite subsets of $V\setminus \{0\}.$
Then $$\Mer(V_\biC, X)^*_\laur = \Mer(V_\biC, X')^*_\laur.$$ 
\end{lemma}
\proof 
This is an immediate consequence of 
Lemma \refer{l: Laurent functionals and proportional X} and the above definition.
\qed

We proceed by discussing the action of a Laurent functional on meromorphic 
functions.
Let $E$ be a complete locally convex space and $\Omega\subset V_\biC$ an open 
subset.
If  $a \in \Omega,$ 
then by $\Mer(\Omega, a, X, E)$ we denote the space of meromorphic functions 
$\gf: \Omega \to E$ whose germ $\gf_a$ at $a$ belongs to $\Mer(a,X,E).$ 
If $S\subset \Omega,$ we define 
$$ 
\Mer(\Omega, S, X, E) := \cap_{a \in S} \Mer(\Omega, a, X, E).
$$ 
Finally, we write $\Mer(\Omega, X, E )$ for $\Mer(\Omega, \Omega, X, E).$ 
In particular, $\Mer(V_\biC, X, E)$ denotes the space of functions
$\gf \in \Mer(V_\biC, E)$ with singular locus $\sing(\gf)$ contained in an $X$-configuration.

There is a natural pairing $\Mer(S,X)^*_\laur \times \Mer(\Omega, S, X, E) \to E,$ given
by
\begin{equation}
\naam{e: pairing Lau and functions}
\Lau \gf = \sum_{a \in \supp \Lau} \Lau_a \gf_a.
\end{equation} 
\hide{The pairing naturally induces a linear map $\Mer(S, X)^*_\laur \to \Mer(\Omega, S, X, E)^*.$ 
Although we do not need the following result in the rest of the paper, we give
a proof for the sake of completeness.}
\begin{lemma}
Let $S \subset \Vc$ be arbitrary, and let  $\Omega$ be an open
subset of $\Vc$ containing $S.$ Then the pairing given by 
(\refer{e: pairing Lau and functions}) for $E=\C$
induces a linear embedding
$$
\Mer(S,X)^*_\laur \embeds \Mer(\Omega, S, X)^*.
$$
\end{lemma}

\proof
Let $\Lau \in \Mer(S,X)^*_\laur$ and assume that $\Lau = 0$ on
$\Mer(\Omega, S, X).$ We may assume that $S=\supp \Lau$.
For every $a \in S$ we write $u^a = (u^a_d)_{d \in \natX}$ for the string
determined by~$\Lau_a.$ 

Select $b \in S.$ Then it suffices to prove
that $\Lau_b = 0.$ 
Fix $d \in \natX$ and $\phi \in \cO_b.$ Then it suffices to show that
$u^b_d(\phi)(b) = 0.$  

For every $a \in S\setminus\{b\}$ we may select 
$d(a) \in \natX$ such that $\pi_{a,d(a)}\pi_{b,d}^{-1}$ is holomorphic at
$a.$ Moreover, we put $d(b) =d.$ 
For $a \in S$ there exists a unique  $v_a \in S(V)$ such that
for all $f \in \cO_a$ we have 
$$
v_a(f)(a) = u^a_{d(a)}(\pi_{a,d(a)}\pi_{b,d}^{-1}f)(a). 
$$ 
We note that $v_b = u_{d}^b.$ 
We may now apply the
lemma below,
with $E_a = \C v_a,$ for $a \in S,$ and, finally with $\xi_a = 0$ if $a \neq b$ and
with $\xi_b$ defined by $\xi_b(v_b)  = v_b(\phi)(b).$
Hence there exists a polynomial function  $\psi$ on $\Vc$ such that 
$v_a(\psi)(a) = 0$ for all $a \in S\setminus \{b\},$ and such that
$v_b(\psi)(b) = v_b(\phi)(b).$ 

Define $\gf =\pi_{b,d}^{-1}\psi .$ Then $\gf \in \Mer(\Omega, S, X).$ 
Hence $\cL\gf = 0.$ On the other hand, 
\begin{eqnarray*}
\cL\gf = \sum_{a \in S} \cL_a\gf_a &=&
\sum_{a \in S} \cL_a(\pi_{a,d(a)}^{-1} \pi_{a,d(a)}\pi_{b,d}^{-1} \psi)
\\
&=&
\sum_{a \in S} u^a_{d(a)}(\pi_{a,d(a)}\pi_{b,d}^{-1} \psi)(a)
=
\sum_{a \in S} v_a( \psi)(a)
=
v_b(\psi)(b) = u^b_d(\phi)(b).
\end{eqnarray*}
It follows that $u^b_{d}(\phi)(b) =0.$  
\qed

\begin{lemma}
Let $S \subset \Vc$ be a finite set. Suppose that for every $a \in S$ a  
finite dimensional complex linear subspace $E_a \subset S(V)$ together with 
a complex linear functional $\xi_a \in E_a^*$ is given. 
Then there exists a polynomial function $\psi$ on $\Vc$  such 
that $u\psi(a) = \xi_a(u)$ for every $a \in S$ and all $u \in E_a.$ 
\end{lemma}

\proof
This result is well known.
\qed
We proceed by discussing the push-forward of a Laurent functional by an injective linear mapping.
Let $V_0$ be a real linear space and $\iota: V_0 \to V$ an injective linear map.
We assume that no element of $X$ is orthogonal to $\iota(V_0)$.
We equip $V_0$ with the pull-back of the inner product of $V$ under $\iota$ 
and denote the corresponding transpose of $\iota$ by $p.$ 
Then $X_0 := p(X)$ consists of non-zero elements.
We denote the complex linear extensions of $\iota$ and $p$ by the same symbols. 
Then, if $H \subset V_\biC$ is an $X$-hyperplane, its preimage
$\iota^{-1}(H)$ is an $X_0$-hyperplane of $V_{0\biC}.$ 

Let $a_0 \in V_{0\biC}$ and put $a =\iota(a_0).$ 
Then pull-back by $\iota$ induces a 
natural algebra homomorphism $\iota^*: \cO_a(V_\biC) \to \cO_{a_0}(V_{0\biC}).$ 
On the other hand, pull-back by $p$ induces 
a natural algebra homomorphism $p^*: \cO_{a_0}(V_{0\biC}) \to \cO_a(V_\biC).$ 
{}From $p\after \iota= I_{V_0}$ it follows that $\iota^* \after p^*  = I$ on 
$\cO_{a_0}(V_{0\biC}),$ hence $\iota^*$ is surjective. 

If $d: X \to \N$ is a map, then we write $p_*(d)$ for the map 
$X_0 \to \N$ defined by 
$$
p_*(d) (\xi_0) = \sum_{\xi \in X, p(\xi) = \xi_0}  d(\xi).
$$ 
One readily verifies that for every $d: X \to \N$ we have
\begin{equation}
\naam{e: pull back of pi}
\iota^* (\pi_{a,X,d}) = \pi_{a_0, X_0, p_*(d)}.
\end{equation}
Let $E$ be a complete locally convex space.
Then it follows that pull-back by $\iota$ induces a linear
map 
\begin{equation}
\naam{e: pull back by iota}
\iota^*: \Mer(V_\biC,a, X, E) \to \Mer(V_{0\biC},a_0, X_0, E).
\end{equation}
\begin{lemma}
\naam{l: injectivity pull back by iota}
The linear map $\iota^*$ in (\refer{e: pull back by iota}) is surjective.
\end{lemma}

\proof
Let $d_0 : X_0 \to \N$ be a map. Then one readily checks that there 
exists a map $d: X \to \N$ such that $d_0 = p_*(d).$ {}From this it
follows that
$$
\pi_{a_0,X_0, d_0}^{-1} \cO_{a_0}(V_{0\biC},E) 
=
\iota^*  (\pi_{a,X, d}^{-1})\; \iota^*p^*( \cO_{a_0}(V_{0\biC},E) )
\subset 
\iota^*( \pi_{a,X,d}^{-1} \cO_a(V_\biC,E) ).
$$ 
where the first equality follows from (\refer{e: pull back of pi}).
\qed

The pull-back map $\iota^*$ in (\refer{e: pull back by iota}) with $E = \C$ 
has a transpose $\iota_*: \Mer(V_{0\biC},a_0, X_0)^* \to  \Mer(V_\biC,a, X)^* $ 
which is injective by Lemma \refer{l: injectivity pull back by iota}.

\begin{lemma}
\naam{l: iota under star}
The map $\iota_*$ maps  $\Mer(V_{0\biC},a_0, X_0)^*_\laur$ injectively
into $\Mer(V_\biC,a, X)^*_\laur.$
\end{lemma}

\proof 
Let
$\Lau \in \Mer(V_{0\biC},a_0, X_0)^*_\laur.$ Then it suffices to 
show that $\iota_* \Lau$ belongs to the space 
$\Mer(V_\biC,a, X)^*_\laur.$

We first note that 
$\iota : V_0 \to V$ has a unique extension  to 
an algebra homomorphism $\iota_*: S(V_0) \to S(V).$ One readily verifies
that $u [\iota^*(\gf)] = \iota^*(\iota_*(u)\gf)$ for every 
$\gf \in \cO_a(V_\biC)$ 
and all $u \in S(V_0).$
Let $d$ be a map $X \to \N.$ Then there exists a $u_d \in S(V_0)$ 
such that
$\Lau = \ev_{a_0}\after u_d \after \pi_{a_0,X_0, p_*(d)}$ on 
$\pi_{a_0,X_0, p_*(d)}^{-1}\cO_{a_0}(V_{0\biC});$
here
$\ev_{a_0}$ denotes evaluation
at the point $a_0.$
Put $v_d = \iota_*(u_d).$ Then, for $\gf \in \cO_a(V_\biC),$ 
$$ 
\iota_*(\Lau) [\pi_{a,X,d}^{-1} \gf]
=
\Lau[\iota^*(\pi_{a,X,d})^{-1} \iota^* \gf)]
=
\Lau [\pi_{a_0,X_0,p_*d}^{-1} \iota^* \gf] 
=
\iota^*(v_d \gf)(a_0) = v_d \gf(a).
$$
Hence $\iota_*(\Lau) = \ev_a \after v_d \after\pi_{a,X,d}$ on
$\pi_{a,X,d}^{-1} \cO_a(V_\biC)$ and we see that 
$\iota_*(\Lau)\in\Mer(V_\biC, a, X)^*_\laur.$ 
\qed

There exists a unique linear 
map $\iota_*: \Mer(V_{0\biC}, X_0)^*_\laur \to \Mer(V_\biC, X)^*_\laur$ 
that restricts to the map $\iota_*$ of 
Lemma \refer{l: iota under star}
for every $a_0 \in V_{0\biC},$ see Remark 
\refer{r: Laurent functionals at a as subspace}.
Clearly,
$\supp(\iota_* \Lau) = \iota(\supp(\Lau)),$ for every $\Lau \in \Mer(V_{0\biC}, X_0)^*_\laur.$

On the other hand, if $E$ is a complete locally convex space, 
$\Omega \subset V_\biC$ open subset
and $S \subset \iota^{-1}(\Omega)$ a subset, then pull-back by $\iota$ 
induces a natural map 
$\iota^*: \Mer(\Omega, \iota(S), X, E) \to \Mer(\iota^{-1}(\Omega), S, X_0, E).$ Moreover,
if $\Lau \in \Mer(V_{0\biC}, S, X_0)^*_\laur$ and $\gf \in \Mer(\Omega, 
\iota(S), X, E),$ then
\begin{equation}
\naam{e: iota push forward and pull back}
\iota_*(\Lau)\gf = \Lau[\iota^*\gf].
\end{equation}

We end this section with a discussion of the multiplication 
by a meromorphic function and the application of a differential operator
to a Laurent functional.

First, assume that $a \in \Vc$ and that $\psi\in \Mer(a, X).$ 
Then multiplication by $\psi$ induces a linear endomorphism of 
$\Mer(a, X),$ which we denote
by $m_\psi.$ The transpose of this linear endomorphism is denoted by 
$m_\psi^*: \Mer(a, X)^* \to \Mer(a, X)^*.$ It readily follows from the definition
of $X$-Laurent functionals at $a$ that $m_\psi^*$ leaves the space $\Mer(a,X)^*_\laur$ 
of those functionals invariant. 

Let now $S \subset \Vc$ be a finite subset, let  $\Omega \subset \Vc$ be an open subset
containing $S$ and let $\psi \in \Mer(\Omega, S, X).$ If 
$\Lau \in \Mer(\Vc, S, X)^*_\laur,$ 
we define the Laurent functional $m_\psi^*(\Lau)\in\Mer(\Vc,S,X)^*_\laur$ 
by 
$$ 
m_\psi^*(\Lau) = \sum_{a \in S} m_{\psi_a}^*(\Lau_a).
$$ 
On the other hand, multiplication
by $\psi$ induces a linear endomorphism of $\Mer(\Omega, S, X),$ and it 
is immediate from the definitions that 
\begin{equation}
\naam{e: previous commutative diagram}
m_\psi^*(\Lau)(\gf)=\Lau(\psi\gf)
\end{equation}
for $\gf\in \Mer(\Omega,S,X)$.

\begin{lemma} 
\naam{l: diff of Laur}
Let $v\in S(V)$, then $v\gf\in\cM(a,X)$ for all $\gf\in\cM(a,X)$, 
and the transpose $\partial^*_v$ of the endomorphism
$\partial_v\colon \gf\mapsto v\gf$ of $\cM(a,X)$
leaves $\cM(a,X)^*_\laur$ invariant.
\end{lemma}

\proof We may assume $v\in V$. Let $d\in \N^X$ and define $d'\in \N^X$
by $d'(\xi)=d(\xi)+1$ for all $\xi\in X$. Then 
$\pi_{a,d}$ divides $v(\pi_{a,d'})$, and hence 
$$\pi_{a,d'}v\gf=v(\pi_{a,d'}\gf)-v(\pi_{a,d'})\gf\in\cO_a$$
for all $\gf\in \pi_{a,d}^{-1}\cO_a$.
Thus $\partial_v\gf=v\gf\in\pi_{a,d'}^{-1}\cO_a$ for 
$\gf\in \pi_{a,d}^{-1}\cO_a$.

Let now $\cL\in\cM(a,X)^*_\laur$, and let $u=u_\cL\in\SprojVX$.
Then for $d,d'$ and $\gf$ as above
$$\partial_v^*\cL(\gf)=\cL(v\gf)=u_{d'}(\pi_{a,d'}v\gf)(a)
=u_{d'}v(\pi_{a,d'}\gf)(a)-u_{d'}(v(\pi_{a,d'})\gf)(a).$$
Each term on the right hand side of this equation has the form
$u'(p\gf)(a)$ with $u'\in S(V)$ and $p$ a polynomial which is divisible by 
$\pi_{a,d}$. Hence, by the Leibniz rule, $\partial_v^*\cL(f)$ has the 
required form $u''(\pi_{a,d}\gf)(a)$, where $u''\in S(V)$.\qed

For $\cL\in\cM(V_\C,X)^*_\laur$ and $v\in S(V)$ we now define 
$\partial_v^*\cL\in\cM(V_\C,X)^*_\laur$ by 
$$\partial_v^*\cL=\sum_{a\in\supp\cL} \partial_v^*\cL_a.$$ It
is immediately seen that $\partial_v^*\cL(\gf)=\cL(\partial_v\gf)$
for each $\gf\in\cM(\Omega,\supp\cL,X)$, where $\Omega$ is an 
arbitrary open neighborhood of $\supp\cL$.

\section{Laurent operators}
\naam{s: Laurent operators}
In this section we discuss Laurent operators, originally introduced in \bib{BSres}, Section 5,
in the slightly different context of meromorphic functions with values in a complete locally
convex space, whose singular locus is contained in an $X$-configuration.

Let $V$ and $X$ be as in the previous section, 
let $\Hyp$ be an $X$-configuration
and let  $E$ be a complete locally convex space.

We define $\Mer(V_\biC, \Hyp, E)$ to be the space of meromorphic functions
$\gf: V_\biC \to E$ whose singular locus is contained in $\cup \Hyp.$ 
If $\Hyp$ is real, we put $\Hyp_V = \{H \cap V \mid H \in \Hyp\}.$ 
Then $\Mer(V_\biC, \Hyp) = \Mer(V_\biC, \Hyp, \C)$ 
equals the space $\Mer(V, \Hyp_V)$ introduced in \bib{BSres}.

It is convenient to select a minimal subset $\subX$ of $X$ that is proportional to
$X.$ Then for every $X$-hyperplane $H \subset V_\biC$ there exists a unique
$\ga_H \in \subX$ and a unique first order polynomial $l_H$ of the form
$z \mapsto \inp{\ga_H}{z} - c,$ with $c \in \C,$ such that $H = l_H^{-1}(0).$ 
Note that a different choice of $\subX$ causes only a change of $l_H$ by a non-zero factor.

Let $\N^\Hyp$ denote the collection of maps $\Hyp \to \N.$ 

\begin{rem}
\naam{r: convention about d}
If $d \in \N^\Hyp,$ then for convenience we agree to write $d(H) = 0$ for any $X$-hyperplane 
$H$ not contained in $\Hyp.$ 
\end{rem}

If $\omega \subset V_\biC$ is a bounded subset and $d \in \N^\Hyp$ 
we define the polynomial function $\pi_{\omega, d}: V_\biC \to \C$ 
by
\begin{equation}
\naam{e: defi pi omega d}
\pi_{\omega, d}  =  \prod_{H \in \Hyp\atop  H  \cap \omega \not= \emptyset} l_H^{d(H)}
\end{equation}
Note that a change of $\subX$ only causes this polynomial to 
be multiplied by a positive
factor. Let $\Mer(V_\biC, \Hyp, d, E)$ 
be the collection of meromorphic functions $\gf\in \Mer(V_\biC, E)$ such that
$\pi_{\omega, d} \gf\in \cO(\omega, E)$ 
for every bounded open subset $\omega \subset V_\biC.$ 
We equip the space $\Mer(V_\biC, \Hyp, d, E)$ with the weakest locally convex topology  
such that for every bounded open subset $\omega \subset V_\biC$ the map 
map $\gf \mapsto \pi_{\omega, d}\gf$ is continuous into $\cO(\omega, E).$ 
This topology is complete; moreover, it is Fr\'echet if $E$ is Fr\'echet. 

We now note that  
\begin{equation}
\naam{e: mer hyp E as union}
\Mer(V_\biC, \Hyp , E) = \cup_{d \in \N^\Hyp} \Mer(V_\biC, \Hyp , d, E).
\end{equation} 
We equip $\N^\Hyp$ with the partial ordering $\preceq$ defined by $d' \preceq d$ 
if and only if $d'(H) \leq d(H)$ for all $H \in \Hyp.$ If $d,d'$ are elements of $\N^\Hyp$
with  $d' \preceq d$ then $\Mer(V_\biC, \Hyp , d', E) \subset \Mer(V_\biC, \Hyp , d, E)$
and the inclusion map $i_{d', d}$ is continuous. Thus, the inclusion maps form a directed 
family and from (\refer{e: mer hyp E as union}) we see that the space $\Mer(V_\biC, \Hyp, E)$ may 
be viewed as the direct limit of the spaces $\Mer(V_\biC, \Hyp , d, E).$ Accordingly 
we equip $\Mer(V_\biC, \Hyp ,E)$ with the direct limit locally convex topology.

By an $X$-subspace of $V_\biC$ we mean any non-empty intersection of 
$X$-hyperplanes; we agree that $\Vc$ itself is also an $X$-subspace.
We denote the set of such affine subspaces by $\cA = \cA(V_\biC, X).$ 
For $L \in \cA$ there exists a unique real linear subspace $V_L \subset V$ such that 
$L = a + V_{L\biC}$ for some $a \in V_\biC.$ The intersection 
$V_{L\biC}^\perp \cap L$ consists of a single point, called the central point of $L;$ 
it is denoted by $c(L).$ The space $L$ is said to be real if $c(L) \in V;$ 
this means precisely that $L$ is the complexification of an affine subspace of 
$V.$ Translation  by $c(L)$ 
induces an affine isomorphism from 
$V_{L\biC}$ onto $L.$ Via this isomorphism we equip $L$ with the structure of a complex
linear space together with a real form  that is equipped with an inner product.

If $L \in \cA,$ the collection of $X$-hyperplanes containing $L$ is finite; we denote
this collection by $\Hyp(L, X).$ Moreover, we put $X(L):= X \cap V_L^\perp$ and
$\subX(L):= \subX \cap V_L^\perp.$ {}From the definition of $\subX$ it 
follows that the map $H \mapsto \ga_H$ is a bijection from $\Hyp(L, X)$ onto
$\subXbL.$ Accordingly we shall identify the sets $\N^{\HypbLX}$ and $\N^{\subXbL}.$ 
If $\Hyp$ is any $X$-configuration and $d \in \N^{\Hyp},$ 
we define the polynomial function
$q_{L,d}$ by 
$$%\begin{equation}
%\naam{e: defi q L d}
q_{L,d}: = \prod_{H \in \Hyp(L,X)} l_H^{d(H)},
$$%\end{equation}
see also Remark \refer{r: convention about d}.
Let $X_r$ be the orthogonal projection of $X\setminus X(L)$ onto $V_L;$ then
$X_r$ is a finite set of non-zero elements. Its image in $L$ under translation
by $c(L)$ is denoted by $X_L.$ 
If $\Hyp$ is  an $X$-configuration in $V_\biC,$  then
the collection 
$$ 
\Hyp_L := \{H\cap L\mid H \in \Hyp,\; \emptyset \subsetneqq H\cap L \subsetneqq L\}
$$ 
is an $X_L$-configuration in $L;$ here $L$ is viewed as a complex linear space
in the way described above.
 
We now assume that $L \in \cA$ and
that $\Hyp$ is an $X$-configuration in $V_\biC.$ In accordance with
\bib{BSres}, Sect. 1.3,
a linear map $R: \Mer(V_\biC, \Hyp) \to \Mer(L, \Hyp_L)$ is called a 
Laurent operator 
if for every $d \in \Hyp^\N$ there exists an element $u_d \in S(V_L^\perp)$ 
such that
\begin{equation}
\naam{e: defi Laurent operator by u d}
R\gf = u_d(q_{L,d} \gf)|_L \text{for all} \gf \in \Mer(V_\biC, \Hyp, d).
\end{equation}
The space of such Laurent operators is denoted by $\Laur(V_\biC, L, \Hyp).$ 

Assume now in addition that $\Hyp$ contains $\Hyp(L,X).$ 
Then as in loc.\ cit.\ it is seen that, for $R \in \Laur(V_\biC, L, \Hyp)$ and 
$d \in \N^\Hyp$,
the element $u_d \in S(V_L^\perp)$ such that 
(\refer{e: defi Laurent operator by u d}) holds, is uniquely determined.
Moreover, it only depends on the restriction of $d$ to $\Hyp(L,X),$  and 
the associated string $u_R:= (u_d\mid d \in \natHypbLX)$ belongs to 
$\Sproj(V_L^\perp, \subXbL).$ 
As in to \bib{BSres}, Lemma 1.5,
the map $R \mapsto u_R$ defines a linear isomorphism 
\begin{equation}
\naam{e: iso Laurent operators with Sproj}
\Laur(V_\biC, L, \Hyp) \simeq \Sproj(V_L^\perp, \subXbL).
\end{equation}
If $E$ is a complete locally convex space, 
and $R \in \Laur(V_\biC, L, \Hyp)$ a Laurent operator,
we may define a linear operator $R_E$ from $\Mer(V_\biC, \Hyp, E)$ 
to $\Mer(L, \Hyp_L, E)$ by the formula (\refer{e: defi Laurent operator by u d}), for 
$\gf \in  \Mer(V_\biC, \Hyp, d, E)$ and with 
$u_d$ equal to the $d$-component of $u_R.$ We shall often denote $R_E$ by $R$ as well.

\begin{rem} Here we note that the algebraic tensor product 
$\Mer(V_\biC, \Hyp) \otimes E$ naturally embeds onto a subspace of 
$\Mer(V_\biC, \Hyp, E)$ which is dense. Thus, $R_E$ is the unique continuous
linear extension of $R \otimes I_E.$ However, we shall not need this.
\end{rem}

\begin{lemma}
\naam{l: continuity Laurent operator}
Let $L \in \cA$ and let $\Hyp$ be an $X$-configuration in 
$V_\biC$ containing $\Hyp(L,X).$ 
Let $R \in \Laur(V_\biC, L, \Hyp).$ 
Then for every $d \in \N^\Hyp$ there exists a $d' \in \N^{\Hyp_L}$ 
with the following property. For every complete locally convex space $E$ 
the operator $R_E$ maps $\Mer(V_\biC, \Hyp, d, E)$ continuously into
the space $\Mer(L, \Hyp_L, d', E).$ 
\end{lemma}

\proof 
This is proved in a similar fashion as \bib{BSres}, Lemma 1.10.
\qed

We shall now relate Laurent operators to the Laurent functionals
introduced in the previous section.
Let $\subXr$ be a minimal subset of $X_r$ 
subject to the condition that it be proportional to $X_r.$ 
Let $\subXL$ be its image in $L$ under translation by $c(L).$ 
Thus, with respect to the linear 
structure of $L,$ the set $\subXL$ is an analogue for the pair 
$(L, X_L)$
of the set $\subX$ for the pair
$(V,X).$

\begin{lemma}
\naam{l: iso Laurent functionals and operators}
Let $L\in \aff$ and let $\Hyp$ be an $X$-configuration in $V_\biC$ 
containing
$\Hyp(L,X).$
Let $E$ be a complete locally convex space.
\begin{enumerate}
\itema
If $\gf \in \Mer(V_\biC, \Hyp, E),$ then for 
$w \in L\setminus \cup \Hyp_L$ the function $z \mapsto \gf(w + z)$ 
is meromorphic on $\VLperpc,$ with a germ at $0$ that belongs to 
$\Mer(\VLperpc,0, X(L),E).$ 
\itemb
If  $\Lau \in \Mer(\VLperpc, 0, X(L))^*_\laur $ is an $X(L)$-Laurent functional in 
$\VLperpc,$ supported at the origin, 
then for $\gf \in \Mer(V_\biC, \Hyp,E)$ the function 
\begin{equation}
\naam{e: defi Lau star}
\Lau_*\gf :  w \mapsto \Lau(\gf(w +  \dotvar ))
\end{equation}
belongs to the space $\Mer(L, \Hyp_L,E).$ 
The operator $\Lau_*: \Mer(V_\biC, \Hyp) \to \Mer(L, \Hyp_L),$ 
defined by (\refer{e: defi Lau star}) 
for $E =\C,$ 
is a Laurent operator. 
\itemc
The map $\Lau \mapsto \Lau_*,$ defined by (\refer{e: defi Lau star})
for $E = \C,$ is an isomorphism from the space  
$\Mer(\VLperpc, 0, X(L))^*_\laur$ onto the space $\Laur(V_\biC, L,\Hyp).$ 
This isomorphism corresponds with the identity 
on $\Sproj(\VLperp,\subXbL),$ via the isomorphisms
of Lemma \refer{l: iso Laurent functionals with Sproj} and 
eq.\ (\refer{e: iso Laurent operators with Sproj}).
\end{enumerate}
\end{lemma}

\proof
See \bib{BSfi}, Appendix B, Lemma B.3.
\qed

\begin{rem}
In the formulation of (c) we have used that the spaces  
$\Mer(\VLperpc, 0, X(L))^*_\laur$ and $\Mer(\VLperpc, 0, \subX(L))^*_\laur$ 
are equal, see Lemma \refer{l: Laurent functionals and proportional X}.
\end{rem}

We now assume that $\Hyp$ is an $X$-configuration, and that $L \in \cA.$ 
If $a \in \VLperpc,$ then by $\Hyp_L(a)$ we denote the collection
of hyperplanes $H'$ in $L$ for which there exists a $H \in \Hyp$ 
such that $H' = L \cap [(-a) + H].$ Thus, $\Hyp_L(a) = (T_{-a}\Hyp)_L$ and we see that 
$\Hyp_L(a)$ is an $\XL$-configuration.
If $S \subset \VLperpc$ is a finite subset, 
then
\begin{equation}
\naam{e: defi Hyp L S}
\Hyp_L(S) = \cup_{a \in S} \Hyp_L(a)
\end{equation}
is an $\XL$-configuration in $L$ as well. 
The corresponding set of regular points in $L$ 
equals
$$
L\setminus \cup \Hyp_L(S) = 
\{w \in L\mid \forall a \in S \,\forall H \in \Hyp: \;\; 
a + w \in H \implies a + L \subset H \}.
$$ 
\begin{cor}
\naam{c: continuity of Laustar}
Let $L \in \cA$ and let $\Hyp$ be an $X$-configuration. 
Let $S \subset \VLperpc$ be a finite subset and let $E$ be a complete locally
convex space.
\begin{enumerate}
\itema
For every $\gf \in \Mer(V_\biC, \Hyp,E)$ and each $w \in L\setminus\cup \Hyp_L(S),$ 
there exists an open  neighborhood $\Omega$ of $S$ in $\VLperpc$ such that the
 function
$\gf(w + \dotvar): z \mapsto \gf(w + z)$ belongs to $\Mer(\Omega,  X(L),E).$ 
\itemb
Let $\Lau \in \Mer(\VLperpc, X(L))^*_\laur$ be a Laurent functional supported at $S.$ 
For every $\gf \in \Mer(V_\biC, \Hyp,E)$ the function 
$\Lau_*\gf: L\setminus\cup \Hyp_L(S) \to E$ 
defined by 
\begin{equation}
\naam{e: defi Laustar}
\Lau_*\gf(w) :=  \Lau(\gf(w + \dotvar)) 
\end{equation} 
belongs to $\Mer(L, \Hyp_L(S),E).$ 
Finally, $\Lau_*$ is a continuous linear map from $\Mer(V_\biC, \Hyp,E)$ to 
$\Mer(L, \Hyp_L(S), E).$ In fact, for every $d \in \N^\Hyp$ there exists
a $d' \in \N^{\Hyp_L(S)},$ independent of $E,$ such that
$\Laustar$ maps $\Mer(V_\biC, \Hyp,d,E)$ continuously into $\Mer(L, \Hyp_L(S), d', E).$ 
\end{enumerate}
\end{cor}

\proof
It
suffices to prove the result for $S$ consisting of a single point $a.$ 
Applying a translation by $-a$ if necessary, we may as well assume that $a = 0.$ 
Then $\Hyp_L(S) = \Hyp_L(0) = \Hyp_L.$ 
Let $\Hyp'$ be the union of $\Hyp$ with $\Hyp(L, X).$
Then $\Mer(V_\biC, \Hyp,E )\subset \Mer(V_\biC, \Hyp', E)$
and $(\Hyp')_L=\Hyp_L=\Hyp_L(S)$, hence assertions (a) and (b) of 
Lemma \refer{l: iso Laurent functionals and operators} 
with $\Hyp'$ in place of $\Hyp$ imply assertion (a) and (b),
except for the final statement about the continuity.

For the final statement of (b), we note that by 
Lemma  \refer{l: iso Laurent functionals and operators}(b),
$\Lau_*$ is a Laurent operator 
$ \Mer(V_\biC, \Hyp') \to \Mer(L, \Hyp_L(S)).$
Let $d: \Hyp \to \N$ be a map. We extend $d$ to  $\Hyp'$ 
by triviality on $\Hyp'\setminus \Hyp.$ Then according to Lemma
\refer{l: continuity Laurent operator} there exists a map $d': \Hyp_L(S) \to \N$ such that
for any complete locally convex space $E$ the map 
$$
\Lau_*: \Mer(V_\biC, \Hyp', d, E) \to \Mer(L, \Hyp_L(S), d', E) 
$$ 
is continuous linear. Since $d$ is zero on $\Hyp'\setminus \Hyp,$
the first of these spaces equals $\Mer(V_\biC, \Hyp, d, E)$ and the asserted
continuity follows.
\qed

\begin{lemma} 
\naam{l: eval of Laustar is a Lau}
Let $L$, $\cH$, $S$ and $\cL$ be as in Cor.\ \refer{c: continuity of Laustar},
and fix $w\in L\setminus\cup\cH_L(S)$. There exists a Laurent
functional (in general not unique) $\cL'\in\cM(V_\C,X)^*_\laur$,
supported in $w+S$, such that $\cL'\gf=\cL(\gf(w+\,\cdot\,))$
for all $\gf\in\cM(V_\C,\cH)$.\end{lemma}

\proof As in the proof of Cor.\ \refer{c: continuity of Laustar} 
we may assume that $S=\{0\}$. Let $\tilde\cH=\cH\cup\cH(w,X)$.
Then $\cL_*\colon\gf\mapsto\cL(\gf(w+\,\cdot\,))$ is a Laurent operator
in $\Laur(V_\C,L,\tilde\cH)$, according to Lemma 
\refer{l: iso Laurent functionals and operators} (b).
On the other hand, it follows from
Lemma \refer{l: extension to Laurent functional} (see Remark 
\refer{r: evaluation is a Laurent functional}) that there exists a 
(in general not unique) $X_L$-Laurent 
functional $\cL''$ on $L$ such that $\psi(w)=\cL''(\psi_w)$ 
for each $\psi\in\cO_w(L)$.
The functional $\psi\mapsto\cL''(\psi_w)$ is defined 
for $\psi\in\Mer(L,\tilde\cH_L)$, and it may be viewed as a 
Laurent operator in $\Laur(L,\{w\},\tilde\cH_L)$, which we
denote by the same symbol $\cL''$ (see \bib{BSfi} Appendix, Remark B.4). 
It now follows from \bib{BSres}, Lemma 1.8 that the composed map 
$\cL''\after\cL_*$ belongs to $\Laur(V_\C,\{w\},\tilde\cH)$ 
and hence by \bib{BSfi} Appendix, Remark B.4
it is given by an $X$-Laurent functional $\cL'$, supported at $w$.
In particular, for $\gf\in\cM(V_\C,\cH)$ we have from Lemma 
\refer{l: iso Laurent functionals and operators} (b) that 
$w\mapsto\cL(\gf(w+\,\cdot\,))$ is holomorphic in a neighborhood of $w$, 
hence its evaluation at $w$ is obtained from the application of $\cL''$ to it. 
Thus $\cL(\gf(w+\,\cdot\,))=\cL'\gf$
for $\gf\in\cM(V_\C,\cH)$.\qed

Recall from Section \refer{s: Laurent functionals} that
$\Mer(V_\biC, X, E)$ is the union of the spaces
$\Mer(V_\biC, \Hyp, E)$ with $\Hyp$ an $X$-configuration. 

\begin{lemma}
\naam{l: Lau under star new}
Let
$L \in \cA$ and let $\Lau \in \Mer(\VLperpc, X(L))^*_\laur$ be a Laurent
functional. Then for any complete locally convex space $E$ 
there exists a unique linear operator 
$$ 
\Lau_*: \Mer(V_\biC, X, E) \to \Mer(L, X_L, E)
$$ 
that coincides on the subspace $\Mer(V_\biC, \Hyp,E)$ 
with the operator $\Lau_*$ defined in 
Corollary \refer{c: continuity of Laustar}, for every $X$-configuration $\Hyp$ in $V_\biC.$ 
\end{lemma}

\proof
Let  $\Hyp_1$ and $\Hyp_2$ be two  $X$-configurations. Let $S = \supp(\Lau)$ 
and let, for $j=1,2,$ 
the continuous linear operator 
$\Lau^j_*: \Mer(V_\biC, \Hyp_j,E) \to \Mer(L, \Hyp_{jL}(S),E)$ be
defined as in Corollary \refer{c: continuity of Laustar} 
with $\Hyp_j$ in place of $\Hyp.$ 
Then it suffices to show that $\Lau_*^1$ and $\Lau_*^2$ coincide on the 
intersection of $\Mer(V_\biC, \Hyp_1,E)$ and $\Mer(V_\biC, \Hyp_2,E).$ 
That intersection equals $\Mer(V_\biC, \Hyp_1 \cap \Hyp_2,E ).$
Let $\gf$ be a function in the latter space,
then from the defining formula (\refer{e: defi Laustar}) 
it follows that $\Lau_{*}^1 \gf = \Lau_*^2\gf$ 
on the intersection of the sets $L\setminus \cup \Hyp_{jL}(S),$ for $j=1,2.$ 
This implies 
that  $\Lau_{*}^1 \gf$ and $\Lau_*^2\gf$ coincide as elements of 
$\Mer(L).$ 
\qed

We end this section with another useful consequence.
\begin{lemma}
\naam{l: diagonal action of Laurent functional}
Let $\Lau \in \Mer(\VLperpc, X(L))^*_\laur.$ Let the finite subset $\widetilde X$ 
of $V\times V\setminus \{(0,0)\}$ be defined by 
$\widetilde X = (X \times \{0\}) \cup 
(\{0\} \times X).$ 
If $\Phi \in \Mer(V_\biC \times V_\biC, \widetilde X),$ then 
$$
\Psi: (w_1, w_2) \mapsto \Lau(\Phi(\dotvar + w_1, \dotvar + w_2))
$$ 
defines a function in $\Mer(L\times L, \widetilde X_L),$ 
where
 $\widetilde X_L = (X_L \times\{c(L)\}) \cup (\{c(L)\} \times X_L).$ 
In particular, the  pull-back of $\Psi$ under the diagonal embedding 
$j: L \to L \times L$ belongs to the space $\Mer(L, X_L).$ 
\end{lemma}

\proof
Equip $\VLperp \times \VLperp$ with one half times the direct sum
inner product.
Then the diagonal embedding $\iota: z \mapsto (z,z)$ is an isometry of 
$\VLperp$ into $\VLperp \times \VLperp.$ Its adjoint is the map 
$p: (z_1, z_2) \mapsto \frac12 (z_1 + z_2)$ from $\VLperp \times \VLperp$ onto $\VLperp.$ 
The intersection $\widetilde X(L) := \widetilde X \cap (\VLperp\times \VLperp)$ 
equals $(X(L) \times \{0\} )\cup (\{0\} \times X(L)).$ Its image under $p$ is given by 
$\widetilde X(L)_0 = \frac12 X(L).$ Thus, according to Lemma 
\refer{l: global Laurent functionals and proportional X},
the space of 
$\widetilde X(L)_0$-Laurent functionals on $\VLperpc$ is equal to the space of  
$X(L)$-Laurent functionals on $\VLperpc.$ Hence, 
according to Lemma \refer{l: iota under star} and the remark following 
its proof, we have an associated 
push-forward map $\iota_*$ from $\Mer(\VLperpc, X(L))^*_\laur$ to 
$\Mer(\VLperpc \times \VLperpc, \widetilde X(L))^*_\laur.$ 

For generic $w_1, w_2 \in L$ we define the meromorphic function $\Phi^{(w_1, w_2)}$
on $\VLperpc \times \VLperpc$ by $\Phi^{(w_1, w_2)}(z_1, z_2) = \Phi(w_1 + z_1 , w_2 + z_2).$ 
The definition of $\Psi$ may now be rewritten as 
$\Psi(w_1, w_2) = \Lau [\iota^*(\Phi^{(w_1, w_2)})].$
By  (\refer{e: iota push forward and pull back}) 
it follows that 
$\Psi(w_1, w_2) = \iota_*(\Lau)(\Phi^{(w_1, w_2)}),$ or, equivalently, 
in the notation
of Lemma \refer{l: Lau under star new},
$$ 
\Psi = [\iota_*(\Lau)]_*\Phi.
$$ 
We now observe that $\widetilde X_L = (\widetilde X)_{L \times L}.$ 
Hence it follows by application of 
Lemma \refer{l: Lau under star new}.
that $\Psi\in\Mer(L\times L,\tilde X_L)$.
There exists an $\tilde X_L$-configuration $\tilde\Hyp$ in
$L\times L$ such that $\Psi\in\Mer(L\times L,\tilde\Hyp)$.
Any hyperplane $\tilde H\in\tilde\Hyp$ is of the form
$\tilde H=H\times L$ or $\tilde H=L\times H$, with $H$ an $X_L$-hyperplane
in $L$. In both cases $j^{-1}(\tilde H)=H$. It now follows that
$j^{-1}(\tilde\Hyp)$ is an $X_L$-configuration in $L$,
and that $j^*\Psi\in\cM(L,X_L)$.
\qed

\section{Analytic families of a special type}
\naam{s: spec fam new}
In
this section we introduce a space $\cEhypQgd$ of analytic families 
of $\DX$-finite $\tau$-spherical functions whose singular locus is
a $\gS$-configuration. The definition of this space is motivated 
by the fact that 
it contains the families obtained from applying
Laurent functionals to Eisenstein integrals related
to a minimal $\gs$-parabolic subgroup, as we shall see in
the following sections, and by the fact that
the vanishing theorem is applicable, see 
Theorem \refer{t: special vanishing theorem}.

In this section we fix a choice $\Sigma^+$ of positive roots
for $\Sigma$ and denote by $P_0$ the associated minimal 
standard $\gs$-parabolic subgroup.

\begin{defi}
\naam{d: Cephyp}
Let $Q \in \allparabs$ and let $Y \subset \staQqdc$ be a finite subset.
We define 
\begin{equation}
\naam{e: CephypQY}
\CephypQY
\end{equation}
to be the space of functions
$f:\faQqdc \times \spXp \to \Vtau,$ meromorphic 
in the first variable, for which there exist a constant $k \in \N,$ 
a $\gS_r(Q)$-hyperplane configuration $\Hyp$ in $\faQqdc$ and a function
$d: \cH \to \N$ such that
the following conditions are fulfilled.
\begin{enumerate}
\itema
The function $\gl \mapsto f_\gl$ belongs to 
$\Mer(\faQqdc, \Hyp, d, \Ci(\spXp\col \tau)).$
\itemb
For every $P \in \minparabs$ and $v \in \NKaq$ there exist functions 
$
q_{s,\xi}(P,v\asmid f)$ in $ P_k(\faq) \otimes 
\Mer(\faQqdc, \Hyp, d, \Ci(\spXzerov\col \tauM)),$
for $s \in W/W_Q$ and $\xi \in -sW_Q Y + \N \Delta(P),$ with the 
following property.
For all $\gl \in \faQqdc\setminus \cup \Hyp,$  $m \in \spXzerov$ and $a \in \Aqp(P),$
\begin{equation}
\naam{e: expansion Cephyp} 
f_\gl(mav) = \sum_{s \in W/W_Q} a^{s\gl - \rho_P} \sum_{\xi  \in -sW_Q Y + \N \gD(P)}
a^{-\xi}\, q_{s,\xi}(P,v\asmid f, \log a)( \gl, m),
\end{equation} 
where the $\DP$-exponential polynomial series of each inner sum
converges neatly on $\Aqp(P).$
\itemc
For every $P \in \minparabs,$ $v \in \NKaq$ and $s\in W/W_Q,$ the series
$$
\sum_{\xi \in -sW_Q Y + \N\DP} a^{-\xi} q_{s, \xi}(P,v\asmid f,\log a)
$$ 
converges neatly on $\Aqp(P),$ as an exponential polynomial series
with coefficients in the space
$\Mer(\faQqdc, \Hyp, d, \Ci(\spXzerov \col \tauM)).$ 
\end{enumerate}
Finally, we define
\begin{equation}
\naam{e: Cephyp}
\Cephyp: = C^{\ep,\hyp}_{P_0, \{0\}}(\spXp\col \tau).
\end{equation} 
\end{defi}

\begin{rem}
\naam{r: relation between hypfam and fam}
Note the analogy between the above definition and Definition
\refer{d: anfamQY newer}. In fact, let  $\Omega = \faQqdc\setminus \cup \Hyp,$
then it follows immediately from the definitions
that the restriction of $f$ to $\Omega \times \spXp$  
belongs to $C^\ep_{Q, Y}(\spXp\col \tau \col \Omega).$
Moreover, it follows from Lemma \refer{l: pointwise expansion of family}
that the functions $q_{s,\mu}(P,v\asmid f)$ introduced 
above are unique, and that the notation used here is consistent 
with the notation in Definition
\refer{d: anfamQY newer}.
The precise relation between the definitions is given in 
Lemma \refer{l: relation specfam and fam second} below.
\end{rem}

\begin{rem}
\naam{r: second on Cephyp}
In analogy with Remark \refer{r: on defi anfamQY new} we note that 
the space (\refer{e: CephypQY}) depends on $Q$ through its $\gs$-split component
$\AQq.$ 
Moreover, it suffices in the above definition to require conditions 
(b) and (c) for
a fixed $P \in \minparabs$ and all $v$ in a given set $\cW \subset \NKaq$ 
of representatives 
for $W/\WKH.$ Alternatively, it suffices to require those conditions 
for a fixed given $v\in \NKaq$ and each $P \in \minparabs.$ 

Finally, we note that $\fa_{P_0 \iq} = \faq,$ hence ${}^*\fa_{P_0} = \{0\}.$ 
Thus, if $Q = P_0,$ we only need to consider the finite set 
$Y = \{0\}.$ This explains the limitation in (\refer{e: Cephyp}).
\end{rem}

It follows from Remark \refer{r: relation between hypfam and fam}
that the following definition of the notion of 
asymptotic degree is in accordance 
with the definition of the similar notion in Definition \refer{d: anfamQY newer}.

\begin{defi}
\naam{d: asymptotic degree Cephyp}
Let $f \in \CephypQY.$ 
We define the asymptotic degree of $f,$ denoted  $\dega(f),$ to be 
the smallest integer $k$ for which there 
exist $\Hyp, d$ such that the conditions of Definition \refer{d: Cephyp}
are fulfilled.
Moreover, we denote by $\Hyp_f$ the smallest $\gS_r(Q)$-configuration
in $\faQqdc$ 
such that the conditions of 
Definition \refer{d: Cephyp} are fulfilled with $k = \dega(f)$ and for 
some $d: \Hyp_f \to \N.$ These choices being fixed,
we denote by $d_f$ the $\preceq$-minimal 
map $\Hyp_f \to \N$ for which the conditions of the definition are fulfilled.
Finally, we put $\rega (f) := \faQqdc\setminus \cup\Hyp_f.$ 
\end{defi}

If $Q \in \allparabs,$ we denote
by $\Sigma_{r0}(Q)$ the set of indivisible roots in $\gS_r(Q),$ i.e., the roots
$\ga \in \gS_r(Q)$ with $]0, 1]\ga \cap \gS_r(Q) = \{\ga \}.$ 
Moreover, we put $\Sigma^+_0 = \gS_{r0}(P_0).$ 
Let $\Hyp$ be a $\gS_r(Q)$-configuration in $\faQqdc$  and $d: \Hyp \to \N$ a map.
If $\omega \subset \faQqdc$ 
is a bounded subset, we define  $\pi_{\omega,d}$ 
as in (\refer{e: defi pi omega d}) with $V = \faQqd,$ $X = \gS_r(Q)$ and 
$X^0 = \gS_{r0}(Q).$ 

\begin{lemma}
\naam{l: relation specfam and fam second}
Let $Q \in \allparabs,$ $Y \subset \staQqdc$ a finite subset,
$\Hyp$ a $\gS_r(Q)$-configuration
in $\faQqdc$ and $d \in \N^\Hyp.$ 
Assume that $f \in \Mer(\faQqdc, \Ci(\spXp\col \tau)).$ 
Then the following two conditions
are equivalent.
\begin{enumerate}
\itema
The function $f$ belongs to $\CephypQY$ 
and satisfies $\Hyp_f \subset \Hyp$ and $d_f \preceq d.$ 
\itemb
For every non-empty bounded open subset $\omega \subset \faQqdc,$ 
the function $ f_{\pi_{\omega, d}}: (\gl, x) \mapsto \pi_{\omega, d}(\gl) f(\gl, x),$ 
$\omega \times \spXp \to \Vtau$ 
belongs to $C_{Q,Y}^\ep(\spXp\col \tau \col \omega).$
\end{enumerate}
Moreover, if one of the above equivalent conditions is fulfilled,
then for every non-empty bounded open subset $\omega \subset \faqdc$ and  
all $P \in \minparabs,$ $v \in \NKaq,$ $s \in W/W_Q$ and $\xi \in -sW_Q Y + \N\DP,$ 
\begin{equation}
\naam{e: relation q of fpi and q}
q_{s, \xi}(P,v\asmid f_{\pi_{\omega,d}}) = \pi_{\omega, d}\, q_{s, \xi} (P,v\asmid f),
\end{equation}
where on the right-hand side we have identified $\pi_{\omega, d}$ with the 
function $1 \otimes \pi_{\omega, d} \otimes 1$ in $P(\faq) \otimes \cO(\omega) \otimes 
\Ci(\spXzerov \col \tau).$ 
\end{lemma}

\proof
Assume that (a) holds and that $\omega \subset \faQqdc$ is a non-empty bounded 
open subset. Put  $\pi = \pi_{\omega,d}$ and $f_\pi=f_{\pi_{\omega,d}}.$
It follows from Definition \refer{d: Cephyp} (a) that 
$f_\pi: \omega \times \spXp \to \Vtau$ 
is smooth and that $f_{\pi\gl}$ is $\tau$-spherical for every $\gl \in \omega.$ 
Thus, it remains to verify conditions (b) and (c) of Definition 
\refer{d: anfamQY newer} for $f_\pi.$ 
Let $P \in \minparabs$ and $v \in \NKaq.$ For $s \in W/W_Q$ and 
$\xi \in -sW_QY + \N \DP$ we define
$$ 
q_{s,\xi}'(P, v \asmid f_\pi, X, \gl, m) := \pi(\gl) q_{s, \xi} (P, v\asmid f, X, \gl, m).
$$ 
Then conditions (b) and (c) of Definition \refer{d: anfamQY newer},
with $k = \dega f$ 
and with $q_{s,\xi}'$ in place of
$q_{s, \xi},$
follow from the similar conditions of Definition \refer{d: Cephyp}. 
Thus, it follows that 
$f_\pi \in C^\ep_{Q,Y}(\spXp \col \tau \col \omega)$ 
and that (\refer{e: relation q of fpi and q}) holds for
all $P \in \minparabs,$ $v \in \NKaq,$ $s \in W$ and $\xi \in -sW_QY +\N \DP.$ 

Now assume that (b) holds, then it suffices to show that (a) holds.
Let $\omega$ be a bounded non-empty 
open subset of $\faQqdc.$  Then it follows
from Definition \refer{d: anfamQY newer}  that the function
$f_\pi=f_{\pi_{\omega,d}}: \omega \times \spXp \to \Vtau$ is smooth; moreover, 
from condition (a) of the mentioned definition it follows that $f_{\pi,\gl}$
is $\tau$-spherical for every $\gl \in \omega.$ 
Hence the map $\gl \mapsto f_\pi$ belongs to 
$\cO(\omega, \Ci(\spXp\col \tau)).$ 
Since $\omega$ was arbitrary, this implies  that $\gl \mapsto
f_\gl$ belongs to $\Mer(\faQqdc, \Hyp, d, \Ci(\spXp\col \tau)).$ 
Hence $f$ satisfies condition (a) of Definition 
\refer{d: Cephyp}. 
Let now $P \in \minparabs$ and $v \in \NKaq.$ Then 
it remains to establish conditions (b) and (c) of 
that definition.

If $\omega$ is a non-empty bounded open subset 
of $\faQqdc,$ then obviously the restriction to $\omega\setminus \cup \Hyp $ 
of the function $f_\pi$ belongs
to $C^\ep_{Q,Y}(\spXp\col \tau\col \omega\setminus \cup\Hyp).$ Moreover,
since $\pi_{\omega, d}$ is nowhere zero on $\omega \setminus \cup \Hyp,$ it 
follows from division by $\pi_{\omega, d}$ 
that the restriction $f|_{(\omega\setminus \cup \Hyp) \times \spXp}$ belongs
to  $C^\ep_{Q,Y}(\spXp\col \tau\col \omega\setminus \cup\Hyp).$ 
Hence, in view of  
Lemma \refer{l: CepQY is sheaf}, the function 
$f$ belongs to 
$C^\ep_{Q,Y}(\spXp\col \tau\col \Omega),$ 
where $\Omega:= \faQqdc \setminus \cup \Hyp.$  Let $k = \dega f.$ 

It follows from the division by $\pi_{\omega,d}$, that 
for every $s \in W$ and $\xi \in -sW_QY +\N \DP,$ 
$$%\begin{equation}
%\naam{e: q of f omega}
\pi_{\omega,d}(\gl) q_{s, \xi}(P,v\asmid f, \dotvar, \gl)
=
q_{s, \xi}(P, v\asmid f_\pi, \dotvar, \gl),
\qquad (\gl \in \omega\setminus \cup \Hyp).
$$%\end{equation}
In particular, the function $(X,\gl) \mapsto 
\pi_{\omega,d}(\gl) q_{s, \xi}(P,v\asmid f, X, \gl)$ belongs 
to the space $P_k(\faq) \otimes \cO(\omega, \Ci(\spXzerov \col \tauM)).$ Since 
$\omega $ is arbitrary, this implies that $f$ satisfies 
condition (b) of Definition \refer{d: Cephyp}. 

{}From condition  (c) of Definition \refer{d: anfamQY newer} with 
$f_\pi$ 
and $\omega$ 
in place of $f$ and $\Omega,$ respectively, it  follows that, for $s \in W,$ 
the series 
$$
\sum_{\xi\in -sW_QY +\N\DP} a^{-\xi} 
\pi_{\omega, d}(\gl)\,q_{s, \xi}(P, v\asmid f ,\log a , \gl)
$$
converges neatly on $\Aqp(P)$ as a $\DP$-exponential polynomial series  
with coefficients in $\cO(\omega, \Ci(\spXzerov \col \tau)).$ 
Since $\omega$ was arbitrary, it follows from the definition of the
topology on $\Mer(\faQqdc, \Hyp, d, \Ci(\spXzerov \col \tauM))$ 
(see Section \refer{s: Laurent operators})
that $f$ satisfies condition (c) of Definition \refer{d: Cephyp}.
\qed

\begin{lemma}
Let $f\in \CephypQY$ and $D \in \DX.$ Then $Df \in \CephypQY.$ 
Moreover, $\Hyp_{Df} \subset \Hyp_f,$ $d_{Df}\preceq d_f$ 
and $\dega Df \leq \dega f.$ 
\end{lemma}

\proof
This follows from a straightforward combination of Lemma
\refer{l: relation specfam and fam second}
with Proposition \refer{p: D on families}.
\qed

If $f\in \CephypQY,$ then by Remark \refer{r: relation between hypfam and fam}
the function
$f$ belongs to $C^\ep_{Q,Y}(\spXp\col \tau\col \Omega),$ 
with $\Omega = \rega f.$ 
Let $k = \dega f.$ 
For $P \in \allparabs,$ $v \in \NKaq,$ $\gs \in \WPQ$ and 
$\xi \in -\gs \cdot Y + \N\DrP,$ let 
$q_{\gs, \xi}(P, v \asmid f) \in P_k(\faPq) \otimes \cO(\Omega, \Ci(\spXPvp\col \tau_P))$ 
be the function defined in Theorem \refer{t: behavior along the walls for families new}.

\begin{lemma}
\naam{l: asymptotics along walls for CephypQY}
Let $Q\in \allparabs$ and $Y \subset \staQqdc$ a finite subset.
Assume that  
$f \in \CephypQY$ and put $k = \dega f.$   
Let $P \in \allparabs$ and $v \in \NKaq.$ 
Then, for every 
$\gl \in \rega f,$ 
the set $\Exp(P, v\asmid f_\gl)$ is contained in $W(\gl+ Y)|_{\faPq} - \rho_P - \N\DrP.$ 
Moreover, let $\gs \in \WPQ.$ Then
\begin{enumerate}
\itema
for every $\xi \in -\gs\cdot Y + \N \DrP,$ 
$$ 
q_{\gs, \xi}(P,v \asmid f) 
\in P_k(\faPq) \otimes \Mer(\faQqdc, \Hyp_f, d_f, \Ci(\spXPvp\col\tauP));
$$
\itemb
for every $R > 1,$ the series 
$$%\begin{equation}
%\naam{e: series along wall of special family new}
\sum_{\xi \in -\gs \cdot Y + \N \DrP} a^{-\xi} q_{\gs, \xi}(P,v\asmid f,\log a)
$$%\end{equation} 
converges neatly on $\APqp(R^{-1})$ 
as a $\DrP$-exponential polynomial series with coefficients
in $\Mer(\faqdc, \Hyp_f, d_f, \Ci(\spXPvp[R]\col \tauP)).$
\end{enumerate}
\end{lemma}

\proof
Let $\Omega=\rega f$.
Then $f\in C^\ep_{Q, Y}(\spXp\col \tau \col \Omega)$.
It follows from Theorem \refer{t: behavior along the walls for families new}
that the assertion about the $(P,v)$-exponents of $f_\gl$ holds.
That (a) and (b) hold can be seen as in the last part of the proof of
Lemma \refer{l: relation specfam and fam second}, with
the reference to Definition \refer{d: anfamQY newer} replaced
by reference to Theorem \refer{t: behavior along the walls for families new}.
\qed

The following definition is the analogue for $\CephypQY$ of 
Definitions \refer{d: defi cE Q Y revised new}
and \refer{d: family for the vanishing thm newer}. 

\begin{defi}
\naam{d: defi cEhyp Q Y gd}
Let $Q \in \allparabs$ and $\gd \in \DQmaps.$ Then for $Y \subset \staQqdc$ a finite 
subset we define
$$
\cEhypQYgd
$$ 
to be the space of functions $f \in \CephypQY$ 
(see Definition \refer{d: Cephyp})
such that, for all 
$\gl \in \rega(f),$ 
the function $f_\gl: x \mapsto f(\gl, x)$ is annihilated by the cofinite ideal 
$\Igdgl.$ 
Moreover, we define
$$
\cEhypQgd: = \bigcup_{Y \subset \staQqdc \;{\rm finite}} \cEhypQYgd.
$$ 

The spaces
$$
\cEhypQYgd_\glob,\quad \cEhypQgd_\glob
$$
are defined to be the spaces of functions $f$ in $\cEhypQYgd$, resp.\ 
$\cEhypQgd$, for which the condition in Definition 
\refer{d: family for the vanishing thm newer}
is satisfied by the restriction to $\Omega=\rega f$. 

Finally, we define $$\cE^\hyp_0(\spXp\col\tau\col\gd):=
\cE^\hyp_{P_0}(\spXp\col\tau\col\gd),\quad 
\cE^\hyp_0(\spXp\col\tau\col\gd)_\glob:=
\cE^\hyp_{P_0}(\spXp\col\tau\col\gd)_\glob$$
for $\gd\in D_{P_0}$.
\end{defi}

\begin{rem}
\naam{r: only local annihilation needed}
Combining Lemmas \refer{l: relation specfam and fam second} and 
\refer{l: locally killed by I gl gd} we 
see that, in the above definition of $\cEhypQYgd$, it suffices to  require
that $\Igdgl$ annihilates $f_\gl$ for $\gl$ in a non-empty open 
subset of $\rega(f).$ 
\end{rem}

We now come to a special case of the vanishing theorem that will
be particularly useful in the following.
Let $\QW \subset \NKaq$ be a complete set of representatives for $W_Q\bs W/\WKH.$

\begin{thm}
{\rm (A special case of the vanishing theorem)\ }
\naam{t: special vanishing theorem}
Let $Q \in \allparabs$  
and let $\delta\in D_Q$. 
Let $f \in \cEhypQgd_\glob$ and let $\Omega'$ be a non-empty open subset of  
$\rega f.$ If
$$ 
\gl - \rho_Q \notin \Exp(Q,u\asmid f_\gl)
$$
for each $u \in \QW$ and all $\gl \in \Omega',$ 
then $f = 0.$
\end{thm}

\proof
Put $\Omega = \rega(f).$ 
It follows immediately from the definitions that 
the restriction $f_\Omega$ of $f$ to $\Omega$ is a family in 
$\cE_Q(\spXp\col \tau \col \Omega\col \gd)_\glob.$ 
Moreover, being the complement of a locally finite collection 
of hyperplanes, $\Omega$ is  $Q$-distinguished in $\faQqdc.$
It follows that $f_\Omega$ satisfies all hypothesis of 
Theorem \refer{t: vanishing theorem new}; hence $f_\Omega = 0$ 
and hence $f=0$.
\qed

\section{Action of Laurent functionals on analytic families}
\naam{s: act Laufunc}
Let $Q\in \allparabs$ be fixed. 
We shall discuss the application of a Laurent functional 
$\Lau \in \Mer(\staQqdc, \Sigma_Q)^*_\laur,$  
to families $f\in\Cephyp.$ More precisely, we want to set up natural
conditions on $f$ under which the family obtained from applying $\Lau$
to $f$ belongs to $\cEhypQgd_\glob$, so that Theorem 
\refer{t: special vanishing theorem} is applicable.

Given a $\gS$-configuration $\Hyp$ 
in $\faqdc$ and a finite subset $S \subset \staQqdc$ we define 
the $\Sigma_r(Q)$-configuration 
$\Hyp_Q(S) = \Hyp_{\faQqdc}(S)$ as in (\refer{e: defi Hyp L S}), 
with $V = \faqdc,$ $X = \Sigma,$ and $L = \faQqdc.$  Thus, for 
$\nu \in \faQqdc$ we have
$$ 
\nu\notin\cup \Hyp_Q(S) \iff [\; \forall \gl \in S \,\forall H \in \Hyp :\;\;\;  
\gl + \nu \in H \implies \gl + \faQqdc \subset H\;].
$$
We recall from  Lemma \refer{l: Lau under star new} that a Laurent functional
$\Lau \in \Mer(\staQqdc, \Sigma_Q)_\laur^*$ induces a 
linear operator 
\def\lcspace{U}
\begin{equation}
\naam{e: Laustar in root context}
\Lau_*: \Mer(\faqdc, \Sigma, \lcspace) \to \Mer(\faQqdc, \Sigma_r(Q), \lcspace),
\end{equation}
for any complete locally convex space $\lcspace.$ 

\begin{lemma}
\naam{l: Laustar for roots and d prime}
Let $\Lau \in \Mer(\staQqdc, \Sigma_Q)_\laur^*$ and put $Y = \supp \Lau.$ 
Let $\Hyp$ be a $\gS$-configuration in $\faqdc,$ and let $\Hyp'=\Hyp_Q(Y)$.
Then for every 
map $d: \Hyp \to  \N$ there exists a map $d': \Hyp' \to \N$ such that, for every complete
locally convex space $\lcspace,$ the linear map (\refer{e: Laustar in root context})
restricts to a continuous linear operator
$$ 
\Lau_*: \Mer(\faqdc, \Hyp,d, \lcspace) \to \Mer(\faQqdc, \Hyp', d', \lcspace),
$$ 
\end{lemma}
\proof 
This follows immediately from Corollary \refer{c: continuity of Laustar}.
\qed

For the formulation of the next result it will be convenient 
to introduce a particular linear map. Let  $\Lau \in \Mer(\staQqdc, \Sigma_Q)_\laur^*$ 
and let $\gl_0 \in Y:= \supp \Lau.$
Let $\Lau_{\gl_0}\in \Mer(\staQqdc, \Sigma_Q)_\laur^*$
 be the Laurent functional supported at $\gl_0,$ 
defined as in Remark \refer{r: Laurent functionals at a as subspace}, 
and let $\lcspace$ be a complete locally convex space.
If $P \in \allparabs$ and 
$s \in W_P\bs W,$ then we define the linear operator $\Lau_{\gl_0 *}^\Ps$ from 
$\Mer(\faqdc, \Sigma, \lcspace) $ 
into $C(\faPq, \Mer(\faQqdc, \Sigma_r(Q), \lcspace) )$  
by the formula
\begin{equation}
\naam{e: defi Laustar s} 
\Lau_{\gl_0*}^\Ps\gf (X, \nu) = e^{-s(\gl_0 + \nu)(X)} 
\Lau_{\gl_0*}[ e^{s(\dotvar)(X)} \gf(\dotvar) ](\nu),
\end{equation}
for $\gf \in \Mer(\faqdc, \Hyp, \lcspace),$  
$X \in \faPq$ and  
$\nu \in \faQqdc\setminus\cup\Hyp_Q(Y).$ 

If $f \in \Cephyp,$ then $f,$ viewed as the function 
$\gl \mapsto f_\gl,$ belongs to the complete
locally convex space $\Mer(\faqdc, \Hyp_f, d_f, \Ci(\spXp\col \tau)).$ 
Accordingly, 
\begin{equation}
\naam{e: Laustar f in merhyp}
\Lau_* f \in \Mer(\faQqdc, \Hyp', d', \Ci(\spXp\col \tau)),
\end{equation} 
where $\Hyp'=\Hyp_{fQ}(Y)$ and $d': \Hyp'\to \N$ is associated with $\Lau, \Hyp_f$ and $d_f$ as in Lemma
\refer{l: Laustar for roots and d prime}.  We note that by definition
\begin{equation}
\naam{e: defi Laustar of family}
\Laustar f(\nu, x) =  \Lau [ f(\dotvar + \nu, x) ],\qquad 
(\nu \in \faQqdc\setminus \cup \Hyp', \;  x \in \spXp).
\end{equation}

\begin{prop}
\naam{p: Lau to hyp family}
Let $Q\in \allparabs$ and let  $\Lau \in \Mer(\staQq, \Sigma_Q)^*_\laur$ 
be a Laurent functional with support contained in the finite subset
$Y \subset \staQqdc.$ 
Assume that $f \in \Cephyp,$ and let $k=\dega f$.
\begin{enumerate}
\itema
The  function $\Laustar f,$ defined as in 
(\refer{e: defi Laustar of family}), belongs
to the space $\CephypQY.$ Moreover, $\Hyp_{\Laustar f}\subset\Hyp'=
\Hyp_{fQ}(Y)$
and $\dega \Laustar f \leq k+ k',$ 
with $k' \in \N$ a constant only depending on $\Lau, \Hyp_f$ and $d_f.$ 
\itemb
Let $P \in \allparabs, v \in \NKaq.$ Then, for $\gs \in \WPQ$ and
$\xi \in -\gs \cdot Y + \N \DrP,$ 
\begin{equation}
\naam{e: sum for q of Laustar f}
q_{\gs, \xi}(P, v\asmid \Lau_*f, X, \nu) 
=
\sum_{\gl \in Y}
\sum_{s \in W_P\bs W, \,  [s] = \gs \atop  s\gl|_{\faPq} + \xi \in \N\DrP} 
\Lau_{\gl*}^\Ps \left[
q_{s, s\gl|_{\faPq} + \xi} (P,v\asmid f)(X,\dotvar)\right] (\nu, X), 
\end{equation}
for all $X \in \faPq$ and  $\nu \in \faQqdc\setminus \cup \Hyp'.$ 
In particular,
\begin{eqnarray*}
&&\Exp(P,v\asmid (\Laustar f)_\nu)\\
&&\quad\subset
\{s(\nu+\gl)|_{\faPq}-\rho_P-\mu\mid s\in W, \gl\in Y, \mu\in\N\Delta_r(P),
q_{s,\mu}(P,v\asmid f)\neq0\}.
\end{eqnarray*} 
\end{enumerate}
\end{prop}
\begin{rem}
\naam{r: sum with empty index set}
Note that the index set of the inner sum in (\refer{e: sum for q of Laustar f})
may be empty. We agree that such a sum  should be interpreted
as zero.
\end{rem}
The following lemma prepares for the proof of the proposition.
\begin{lemma}
\naam{l: continuity Laustar s new}
Let $\Lau \in \Mer(\staQqdc, \Sigma_Q)_\laur^*$ be a Laurent functional with
support contained in the finite set $Y \subset \staQqdc.$ 
Let $\cH$ be a $\Sigma$-configuration in $\faqdc$ and $d: \cH \to \N$ a map.
Let $\Hyp'=\Hyp_Q(Y)$ and
$d':  \Hyp' \to \N$ be as in Lemma \refer{l: Laustar for roots and d prime}.
There exists a 
natural number $k' \in \N$ with the following property.

For every $\gl_0 \in Y,$  every $P \in \allparabs,$ 
each  $s \in W_P\bs W$ and 
any complete locally convex space
$\lcspace,$ the operator $\Lau_{\gl_0*}^\Ps$ restricts to a continuous linear map
$$
\Lau_{\gl_0*}^\Ps :\;\; \Mer(\faqdc, \Hyp, d, \lcspace) \to P_{k'}(\faPq)\otimes 
 \Mer (\faQqdc, \Hyp', d', \lcspace).
$$
\end{lemma}

\proof
For a fixed $X \in \faPq,$ multiplication by the holomorphic 
function
$e^{s(\dotvar)(X)}: \faqdc \to \C$ yields a continuous linear endomorphism 
of the space $\Mer(\faqdc, \Hyp, d, \lcspace);$ similarly, multiplication
by the holomorphic function $e^{-s(\gl_0 + \dotvar)(X)}: \faQqdc \to \C$ 
yields a continuous linear endomorphism of $\Mer(\faQqdc, \Hyp', d', \lcspace).$ 
It now follows from (\refer{e: defi Laustar s}) that for a fixed
$X \in \faPq,$ the function $\Lau_{\gl_0*}^\Ps\gf (X)$ belongs
to the space $\Mer(\faQqdc, \Hyp', d', \lcspace)$ and depends continuously
on $\gf.$ 
Thus, it remains to establish the 
polynomial dependence on $X.$ 
 
For any $\Sigma$-hyperplane $H \subset \faqdc$ we denote by 
$\ga_H$ the root from $\Sigma_0^+$ such that $H$ is a translate 
of $\ga_{H\iC}^\perp.$  Let $\Sigma_{Q,0}^+:= \Sigma_Q \cap \Sigma_0^+$ 
and let $d_0: \Sigma_{Q,0}^+ \to \N$ be defined 
by $d_0(\ga) = d(\ga^\perp + \gl_0);$ thus 
$d_0(\ga) = 0$ if $\ga^\perp + \gl_0 \notin \Hyp.$ 
We define $\pi_0 = \pi_{\gl_0, d_0}$ as in (\refer{e: defi pi a X d}) with
$\staQqd, \gl_0, \Sigma_{Q,0}^+$ and $ d_0$ in place of $V, a, X$ 
and $d,$ respectively.
If $\gf \in \Mer(\faqdc, \Hyp, d, \lcspace),$ then for 
$\nu \in \faQqdc\setminus \cup\Hyp',$ the germ 
of the function $\gf^\nu: \gl \mapsto \gf(\gl + \nu)$ at $\gl_0$ 
belongs to $\pi_0^{-1} \cO_{\gl_0}(\staQqdc, \lcspace)$. 
Hence there exists
a constant coefficient differential operator $u_0 \in S(\staQqd),$ independent 
of $\lcspace,$ such that 
\begin{equation}
\naam{e: action of Laustar as dif op}
\Lau_{\gl_0 *} \gf (\nu)  =  u_0[ \pi_0(\dotvar)\gf(\dotvar + \nu) ](\gl_0),     
\qquad (\nu \in \faQqd\setminus \cup \Hyp'),
\end{equation}
for any $\gf \in \Mer(\faqdc, \Hyp, d, \lcspace).$  
Inserting (\refer{e: action of Laustar as dif op}) in (\refer{e: defi Laustar s}) 
we find that 
\begin{eqnarray*}
\Lau_{\gl_0*}^\Ps \gf (X,\nu) 
&=& e^{-s(\gl_0 + \nu)(X)} u_0[ e^{s(\dotvar +\nu )(X)}
 \pi_0(\dotvar) \gf(\dotvar + \nu) ](\gl_0)\\
&=& e^{-s(\gl_0)(X)} u_0[ e^{s(\dotvar)(X)}
 \pi_0(\dotvar) \gf(\dotvar + \nu) ](\gl_0).
\end{eqnarray*}
By application of the Leibniz rule it finally follows that 
this expression is polynomial  in the variable $X$ of degree 
at most $k' := {\rm order}(u_0).$ 
\qed
\medno
{\bf Proof of Proposition \refer{p: Lau to hyp family}:\ }
By
linearity we may assume that 
$\supp \Lau$ consists of a single point $\gl_0 \in \staQqdc.$ 
Let $\Hyp=\Hyp_f$ and $d = d_f$, and let
$d':\Hyp'\to\N$ and $k'\in\N$ be associated as in Lemmas 
\refer{l: Laustar for roots and d prime} and
\refer{l: continuity Laustar s new}. 
We will establish parts (a), (b) and (c) of Definition \refer{d: Cephyp}
for $\Laustar f$ with $k$, $\Hyp$ and $d$ replaced by $k+k'$, 
$\Hyp'$ and $d'$.    
Note that part (a) was observed already in (\refer{e: Laustar f in merhyp}).
Put $\Omega:= \faQqdc\setminus \cup \Hyp'.$ Then, in particular, 
the function $\Laustar f: \Omega \times  \spXp \to \Vtau$ is smooth.

We will establish parts (b) and (c) of Definition \refer{d: Cephyp} 
by obtaining an exponential polynomial 
expansion for $(\Laustar f)_\nu,$ for $\nu \in \Omega$,
along $P\in\minparabs$. However, having the proof of
(\refer{e: sum for q of Laustar f}) in mind, we assume only 
$P \in \allparabs$ at present.
Let $v \in \NKaq.$ 
Then $f \in C^\ep_{P_0, \{0\}}(\spXp\col \tau\col \faqdc\setminus \cup \Hyp)$ 
by Remark \refer{r: relation between hypfam and fam}. Hence 
by Lemma \refer{l: asymptotics along walls for CephypQY} and 
(\refer{e: series for fam along P}) we obtain,
for $\gl \in \faqdc\setminus\cup\cH,$ 
\begin{equation}
\naam{e: f as sum fs}
f(\gl, mav) = \sum_{s \in W_P \bs W}
f_{s}(\gl, a, m),\qquad (m \in \spXPvp,\; a \in \APqp(R_{P,v}(m)^{-1})),
\end{equation}
where the functions $f_{s}$ on the right-hand side are defined by
\begin{equation}
\naam{e: series for fs}
f_{s}(\gl,a,m ) =  a^{s \gl -\rho_P } \sum_{\mu \in \N\DrP}
a^{-\mu}  q_{s,\mu} (P,v\asmid f)( \log a , \gl , m).
\end{equation}
Here the functions
$q_{s,\mu}(P,v\asmid f)$ belong to the space
$P_k(\faPq) \otimes \Mer(\faqdc, \Hyp, d, 
\Ci(\spXPvp\col \tauP)).$  
By Lemma \refer{l: asymptotics along walls for CephypQY} (b),
for every $R>1$ the series 
in (\refer{e: series for fs}) converges neatly
on $\APqp(R^{-1})$
as a series with coefficients in 
$\Mer(\faqdc, \Hyp, d, C^\infty(\spXPvp[R]\col \tauP)).$ 
By (\refer{e: defi Laustar s}) we have, for 
$\nu \in \Omega,$ $m \in \spXPvp[R]$ and $a \in \APqp(R^{-1})$
$$\Laustar(f_s)(\nu,a,m)= a^{s(\gl_0+\nu)-\rho_P}
\Lau_{\gl_0*}^\Ps [\sum_{\mu \in \N\DrP}
a^{-\mu}  q_{s,\mu} (P,v\asmid f) ( \log a , \dotvar , m)](\log a,\nu).$$
It follows from 
Lemma \refer{l: continuity Laustar s new} that
$\Lau_{\gl_0*}^\Ps$ may 
be applied term by term to the series. Moreover, 
the resulting series is neatly convergent on
$\APqp(R^{-1})$ as a $\DrP$-exponential polynomial series with coefficients
in $\Mer(\faQqdc, \Hyppr, \dpr, \Ci(\spXPvp[R]\col\tauP)).$

The application of $\Laustar$ 
thus leads to the following identity, 
\begin{equation}
\naam{e: series for Laustarfs}
\Laustar(f_s)(\nu, a, m) = a^{s(\gl_0 + \nu) - \rho_P}
\sum_{\mu \in \N\DrP} a^{-\mu} q_{s,\mu}^\Lau (P,v\asmid  f)(\log a, \nu, m),
\end{equation}
where the function  
$q_{s,\mu}^\Lau(P,v\asmid f): \faPq \times \Omega \to \Ci(\spXPvp\col \tauP)$ is given 
by  
\begin{equation}
\naam{e: defi q Lau s mu new}
q_{s,\mu}^\Lau(P,v\asmid f)(\log a, \nu)
=
\Lau_{\gl_0*}^\Ps [q_{s,\mu}(P,v\asmid f,\log a, \dotvar)](\log a,\nu).
\end{equation}
Using Lemma \refer{l: continuity Laustar s new} we deduce that
$$
q_{s,\mu}^\Lau(P,v\asmid f)
\in
P_{k+ k'}(\faPq) \otimes 
\Mer(\faQqdc, \Hyppr, \dpr, \Ci(\spXPvp\col \tauP)).
$$

Combining (\refer{e: series for Laustarfs}) with 
(\refer{e: f as sum fs}) we obtain
an  exponential polynomial expansion along $(P,v)$ for the 
$\tau$-spherical function $(\Laustar f)_\nu$ as  
\begin{equation}
\naam{e: series for Laustar f nu}
(\Laustar f)_\nu (mav) = \sum_{s \in W_P\bs W}
a^{s(\gl_0 + \nu) -\rho_P}
\sum_{\mu \in \N\DrP}
a^{-\mu} q_{s,\mu}^\Lau(P,v\asmid f)(\log a, \nu,m).
\end{equation}
If $s \in W_P\bs W$ and $\nu \in \faQqdc,$ then $s\nu|_{\faPq} = [s]\nu|_{\faPq},$ 
where $[s]$ denotes the class of $s$ in $\WPQ.$ It follows that the series 
in (\refer{e: series for Laustar f nu}) may be rewritten as
$$
\sum_{\gs \in \WPQ} a^{\gs \nu - \rho_P}
\sum_{s \in W_P\bs W, [s] = \gs\atop\mu \in \N\DrP}
a^{s \gl_0 -\mu} q_{s,\mu}^\Lau(P,v\asmid f)(\log a, \nu,m).
$$ 
The exponents $s\gl_0 - \mu$ as $s \in W_P\bs W,$ $[s]  =\gs$ and $\mu \in \N \DrP,$
are all of the form $-\xi,$ with $\xi \in -\gs\cdot\{\gl_0\} + \N\DrP.$ 
Thus, we see that, for $\nu \in \Omega,$ $m\in \spXPvp[R]$ and 
$a \in \APqp(R^{-1}),$   
\begin{equation}
\naam{e: series for Laustar f along P v}
(\Laustar f)_\nu (mav) = \sum_{\gs \in \WPQ}
a^{\gs\nu -\rho_P}
\sum_{\xi \in -\gs\cdot\{\gl_0\} + \N\DrP} a^{-\xi}\;\widetilde 
q_{\gs, \xi}(\log a, \nu, m)
\end{equation}
with 
\begin{eqnarray}
\naam{e: q gs xi of Laustar f}
 \widetilde q_{\gs, \xi}
&= &
\sum_{s \in W_P\bs W,\, [s] = \gs \atop s\gl_0|_{\faPq} + \xi  \in \N\DrP}
 q_{s,s\gl_0|_{\faPq} + \xi }^\Lau(P,v\asmid f)\\
&\in&
P_{k+ k'}(\faPq) \otimes 
\Mer(\faQqdc, \Hyppr, \dpr, \Ci(\spXPvp\col \tauP)).
\nonumber
\end{eqnarray}
{}From what we said earlier about the convergence 
of the series in (\refer{e: series for Laustarfs}), it follows that, for every $R>1,$ the
inner series on the right-hand side of (\refer{e: series for Laustar f along P v})
converges neatly on $\APqp(R^{-1})$ as a 
$\DrP$-exponential polynomial series with
coefficients in the space $\Mer(\faQqdc, \Hyp', d', \Ci(\spXPvp[R] \col\tauP)).$ 

If $P$ is minimal, then $\spXPvp[R] = \spXzerov$
and we see that $\Laustar f$ satisfies conditions (b) and (c) 
of Definition \refer{d: Cephyp} with 
$q_{\gs,\xi}(P,v\asmid \Laustar f) = \tilde q_{\gs, \xi}$ 
for $\gs\in\WPQ\,=W/W_Q$.
This establishes part (a) of the proposition.

For general $P$ we now see that the functions $\tilde q_{\gs, \xi}$ 
introduced above coincide with functions $q_{\gs, \xi}(P,v\asmid \Laustar f)$ 
introduced in Theorem \refer{t: behavior along the walls for families new}. 
Finally, combining (\refer{e: q gs xi of Laustar f}) and 
(\refer{e: defi q Lau s mu new}) we see that we have 
established part (b) of the proposition as well.
\qed

\begin{lemma}
\naam{l: Lau to family of eigenfunctions new}
Let $\gd\in D_{P_0}$ and $f \in \cE_0^\hyp(\spXp\col \tau \col \gd).$ 
Let $Q\in \allparabs$ and $\Lau \in \Mer(\staQqdc, \Sigma_Q)^*_\laur,$ and
put $Y = \supp \Lau.$   
There exists a $\gd' \in \DQmaps$ such that 
$$
\Lau_*f \in \cE_{Q,Y}^\hyp(\spXp\col \tau\col \gd').
$$  
\end{lemma}
\proof
It follows from Proposition \refer{p: Lau to hyp family} 
that $\Laustar f \in C^{\ep, \hyp}_{Q,Y}(\spXp\col\tau).$ 
Moreover, $\rega{\Laustar f} \supset \Omega = \faQqdc\setminus \Hyp_{fQ}(Y).$ 
Then in view of  Definition \refer{d: defi cEhyp Q Y gd} and Remark
\refer{r: only local annihilation needed}
it suffices to establish the existence of  
a $\gdmap' \in \DQmaps$ such that, for every $\nu \in \Omega,$ the function
$(\Lau_*f)_\nu$ is annihilated by the cofinite ideal
$I_{\gd',\nu}.$ 

By linearity we may assume that $\supp \Lau$ consists of a single 
point $\gl_0 \in \staQqdc.$ Then $\Lau = \Lau_{\gl_0}.$ 

Let $\pi_0,u_0$ be as in the proof of Lemma \refer{l: continuity Laustar s new}.
Then from (\refer{e: action of Laustar as dif op}) we see that 
$$
(\Laustar f)_\nu( x) = u_0 [\pi_0(\dotvar) f(\dotvar + \nu, x)](\gl_0) ,
$$ 
for $x \in \spXp,\, \nu \in \Omega.$ 
Moreover, since $(\gl, x) \mapsto \pi_0(\gl)f_{\gl + \nu}(x)$ 
is smooth in a neighborhood of $\{\gl_0\} \times \spXp,$ 
it follows that, for $D \in \DX,$ $\nu \in \Omega$ and $x \in \spXp,$ 
\begin{equation}
\naam{e: D on Laustar f}
D(\Laustar f)_\nu (x)  =
u_0[ \pi_0 (\dotvar) D (f_{\dotvar + \nu})(x) ](\gl_0).
\end{equation}
Put $l = {\rm order}(u_0)$ and define $\gd'\in D_Q$ by
$\supp\gd'=\{\gl_0\}+\supp\gd$ and $\gd'(\gl_0+\gL)=\gd(\gL)+l$ for
$\gL\in\supp\gd$. It suffices to prove the following.
Let elements $D^\gL_i\in\DX$ be given for $i=1,\dots,\gd(\gL)+l$,
for each $\gL\in\supp\gd$, and define the differential operator
\begin{equation}
\naam{e: defi Dnu}
D_\nu:=\prod_{\gL\in\supp\gd}\prod_{i=1}^{\gd(\gL)+l}
(D^\gL_i-\gg(D^\gL_i,\gl_0+\gL+\nu))\in\DX
\end{equation}
for $\nu\in\faQqdc$.
Then $D_\nu$ annihilates $(\Laustar f)_\nu$ for each $\nu\in\Omega$.

It follows from (\refer{e: defi Laustar of family}) and
(\refer{e: D on Laustar f}) that
\begin{equation}
\naam{e: Dnu on Laustar f}
D_\nu(\Laustar f)_\nu(x)= u_0[\pi_0(\dotvar) D_\nu f_{\dotvar+\nu}(x)](\gl_0),
\end{equation}
where the dots indicate a variable in $\staQqdc$.
We write each factor in $D_\nu$ as
\begin{eqnarray*}
&&D^\gL_i-\gg(D^\gL_i,\gl_0+\gL+\nu)\\
&&\qquad=[D^\gL_i-\gg(D^\gL_i,\dotvar+\gL+\nu)]+
[\gg(D^\gL_i,\dotvar+\gL+\nu)-\gg(D^\gL_i,\gl_0+\gL+\nu)],
\end{eqnarray*}
also with variables in $\staQqdc$ indicated by dots.
Inserting this into (\refer{e: defi Dnu}) and
(\refer{e: Dnu on Laustar f}) we obtain an
expression for $D_\nu(\Laustar f)_\nu(x)$ as a sum of terms
each of the form
\begin{equation}
\naam{e: summands of Dnu Laustar f}
u_0[\pi_0(\dotvar) 
\prod_{\gL\in\supp\gd}p^\gL(\dotvar)D^\gL(\dotvar)
f_{\dotvar+\nu}(x)](\gl_0),
\end{equation}
where
$$D^\gL(\gl)=\prod_{i\in S_\gL}[D^\gL_i-\gg(D^\gL_i,\gl+\gL+\nu)]$$
and
$$p^\gL(\gl)=\prod_{i\in S_\gL^c}
[\gg(D^\gL_i,\gl+\gL+\nu)-\gg(D^\gL_i,\gl_0+\gL+\nu)]$$
with $S_\gL$ a subset of $\{1,\ldots,\gd(\gL)+l\}$
and $S^c_\gL$ its complement in this set.
On the one hand, if $S_\gL$ has fewer than 
$\gd(\gL)$ elements for some $\gL$, there are at least
$l+1$ factors in the corresponding product $p^\gL$. Since each of these 
factors vanish at $\gl_0$, it follows from the Leibniz rule that then 
(\refer{e: summands of Dnu Laustar f}) vanishes. On the other hand,
if for each $\gL$ the set $S_\gL$ has at least $\gd(\gL)$ elements, 
then the differential operator $\prod_{\gL}D^\gL(\gl)$ 
annihilates $f_{\gl+\nu}$, again causing
(\refer{e: summands of Dnu Laustar f}) to vanish. It follows
that $D_\nu(\Laustar f)_\nu(x)=0$.
\qed

In the following definition we introduce a  notion of asymptotic
globality that is somewhat stronger than the one 
in Definition \refer{d: s globality new}. It is motivated by the fact that
it carries over by the application of Laurent functionals, as we shall see
in Proposition \refer{p: transference of s globality}

\begin{defi}
\naam{d: holomorphic s globality new}
Let $Q \in \allparabs,$ and let $Y\subset \staQqdc$ be finite.
Let  $P \in \allparabs,$ $v \in\NKaq$ and $\gs \in \WPQ.$
\begin{enumerate}
\itema Let $\Omega \subset \faQqdc$ be an open subset.
A family
$f \in C^\ep_{Q,Y}(\spXp\col \tau \col\Omega)$ 
is called holomorphically $\gs$-global along
$(P,v)$ if there exists a full open subset $\Omega^*$ of $\faQqdc$ such 
that, for every  $\xi \in  -\gs\cdot Y + \N\DrP,$
the function $\gl \mapsto q_{\gs,\xi}(P,v\asmid f,\dotvar)(\gl)$ is a 
holomorphic $P_k(\faPq) \otimes \Ci(\spXPv \col \tauP)$-valued function on 
$\Omega^*\cap\Omega,$ for some $k \in \N.$  
\itemb A family $f \in C^{\ep,\hyp}_{Q,Y}(\spXp\col \tau)$ is called 
holomorphically $\gs$-global along $(P,v)$ if its restriction 
to $\Omega=\rega f$ is holomorphically $\gs$-global along $(P,v)$, according to (a).
\end{enumerate}
\end{defi}

It is easily seen that the property of holomorphic globality
according to (a) of the above definition implies the globality
in Definition \refer{d: s globality new}.
We have the following analogue of Lemma \refer{l: transformation of globality}, 
describing 
how the property of holomorphic globality transforms under
the action of $\NKaq.$ 

\begin{lemma}
\naam{l: transformation of holomorphic globality}
Let $Q$, $Y$, $P$, $v$ and $\gs$ be as above,
and let $f \in C^{\ep,\hyp}_{Q,Y} (\spXp \col \tau).$ 

If $f$ is holomorphically $\gs$-global along $(P,v),$ then
$f$ is holomorphically $u\gs$-global along $(uPu^{-1}, uv),$ for every $u \in \NKaq.$ 
\end{lemma}

\proof
The proof is completely analogous to the proof
of Lemma \refer{l: transformation of globality},
involving an application of Lemma
\refer{l: transformation holo coeffs}.
\qed

\begin{prop}
\naam{p: mero H globality newer}
Let $Q \in \allparabs,$ 
$Y\subset \staQqdc$ a finite subset and let
$P \in \allparabs,$ $v \in\NKaq$ and $\gs \in \WPQ.$
Let $f \in \CephypQY$ and put $\Hyp=\Hyp_f$, $d=d_f$ and $k=\dega f$.

The family $f$ is holomorphically $\gs$-global along 
$(P,v)$ if and only if, for every element $\xi \in - \gs\cdot Y + \N\DrP,$
the function $\gl \mapsto q_{\gs, \xi}(P,v\asmid f,\dotvar)(\gl)$ 
belongs to the space $\Mer(\faQqdc, \Hyp, d, P_k(\faPq) \otimes
\Ci(\spXPv\col \tauP)).$ 
\end{prop}

\proof
The `if'-statement is obvious. Assume that $f$ is holomorphically 
$\gs$-global along $(P,v)$, and let $\xi \in -\gs \cdot Y + \N\DrP.$ 
According to Lemma \refer{l: asymptotics along walls for CephypQY}, 
the function 
\begin{equation}
\naam{e: globality of q}
\gl \mapsto q_{\gs, \xi}(P,v\asmid f, \dotvar , \gl)
\end{equation}
belongs to the space
\begin{equation}
\naam{e: space to which q belongs}
\Mer(\faQqdc, \Hyp,d, P_k(\faPq)\otimes \Ci(\spXPvp\col \tauP)).
\end{equation}

Let $\Omega=\rega(f)$ and let
$\Omega^*$ be a full open subset of $\faQqdc$ satisfying the properties of 
Definition \refer{d: holomorphic s globality new} (a) for the restriction
of $f$ to $\Omega$.
Then the function (\refer{e: globality of q}) not only belongs to the space
(\refer{e: space to which q belongs}), but also to the space
$\cO(\Omega^*\cap\Omega, P_l(\faPq) \otimes \Ci(\spXPv\col \tauP)),$ for some $l \in \N.$ 
In particular we see that this is true with $l = k.$ 

Let now $X \in\faPq$ be fixed.
Then it suffices to show that the function
(\refer{e: globality of q}), with $X$ substituted for the dot, belongs to
the space
$\Mer(\faQqdc, \Hyp, d,  \Ci(\spXPv\col \tauP)).$ To prove the latter,
we fix an arbitrary bounded non-empty open set $\omega \subset \faQqdc$ and 
put  $\pi := \pi_{\omega, d},$ see above Lemma
\refer{l: relation specfam and fam second}. 
Then the function
$F: \omega \times \spXPvp \to \Vtau,$ defined by
$$
F(\gl, m) = \pi(\gl)\, q_{\gs, \xi}(P,v\asmid f, X, \gl)(m)
$$ 
is $C^\infty$ and holomorphic in its first variable.
Moreover, let $\omega_0$ be the full open subset 
$\omega \cap \Omega^*\cap\Omega$ 
of $\omega.$ Then by what we said above,
the restricted function $ F|_{\omega_0 \times \spXPvp}$ 
admits a smooth extension to the manifold $\omega_0 \times \spXPv.$ 
It now follows from Corollary \refer{c: aux smooth extension} that 
$F$ has a unique smooth extension
to $\omega \times \spXPv;$ 
this extension is holomorphic in the first variable.
It follows that the function 
$\gl \mapsto \pi(\gl) q_{\gs, \xi}(P,v\asmid f, X, \gl)$ 
belongs to $\cO(\omega, \Ci(\spXPv\col \tauP)).$ 
Since $\omega$ was arbitrary, this completes the proof.
\qed

\begin{prop}
\naam{p: transference of s globality}
Let $f \in \Cephyp,$ let $Q \in \allparabs$ 
and let  $\cL$ be a Laurent functional in $\cM(\staQqdc, \Sigma_Q)^*_\laur.$ 
Put $Y = \supp \Lau.$
Let $P \in \allparabs,$ $ v \in \NKaq$ and  $\gs \in \WPQ.$

If $f$ is holomorphically $s$-global along $(P,v)$ for every $s \in W_P\bs W$ 
with $[s] = \gs,$ 
then  $\cL_* f \in \CephypQY$ is holomorphically 
$\gs$-global along $(P,v).$ 
\end{prop}

\proof
It follows from Proposition  \refer{p: Lau to hyp family} (a) 
that $\cL_* f \in \CephypQY.$ 
Assume that $f$ satisfies the globality assumptions. Then 
it remains to establish the assertion on $\gs$-globality for $\Laustar f$.

Let $k = \dega f.$ 
Let $\Hyp = \Hyp_f,$ $ d = d_f$ and $\Hyp' = \Hyp_Q(Y).$ 
Moreover, let $d': \Hyp' \to \N$ be associated with these data 
as in Lemma \refer{l: Laustar for roots and d prime} and let $k' \in \N$ 
be associated as in Proposition \refer{p: Lau to hyp family} (a).
According to the latter proposition, the set 
$\Omega' = \faQqdc\setminus\cup \Hyp'$ 
is contained in $\rega(\Laustar f)$.

Let
$\xi \in - \gs\cdot Y + \N \DrP.$ Moreover, let $s \in W_P\bs W$ 
be such that $[s] = \gs$ 
and let $\gl_0 \in Y$ be such that $\eta:=s\gl_0|_{\faPq} + \xi $ 
belongs to $\N \DrP.$ 
Then by Proposition \refer{p: mero H globality newer}, the function 
$$
\gl \mapsto q_{s,\eta}(P,v\asmid f, \dotvar,\gl)
$$ 
belongs to $\Mer(\faQqdc, \Hyp, d, P_k(\faPq) \otimes \Ci(\spXPv \col \tauP)).$ 
Using Lemma \refer{l: continuity Laustar s new} with $C^\infty(\spXPv\col \tauP)$ 
in place of $\lcspace,$ we see that, for $X \in \faPq,$ the function
$$
\gf_X := \Lau^\Ps_{\gl_0 *} [q_{s, \eta}(P,v\asmid f, X, \dotvar)]
$$ 
belongs to $\Mer(\faQqdc, \Hyppr, d', P_{k'}(\faPq) \otimes \Ci(\spXPv\col \tauP)).$ 
Moreover, it depends on $X \in \faPq$ as a polynomial function of degree at most $k.$ 
It follows that the function $(\nu, X) \mapsto \gf_X(\nu)(X)$ 
belongs to the space
\begin{equation}
\naam{e: Mer prime pol k kpr}
\Mer(\faQqdc, \Hyppr, d', P_{k + k'}(\faPq) \otimes \Ci(\spXPv\col \tauP)).
\end{equation}
Each term in the finite sum (\refer{e: sum for q of Laustar f}) is of this form. Hence the function
$$
(\nu, X) \mapsto q_{\gs, \xi}(P, v\asmid \Laustar f , X, \nu)
$$
belongs to the space (\refer{e: Mer prime pol k kpr}) as well.
This holds for all $\xi \in - \gs\cdot Y + \N\DrP.$ Therefore the 
restriction of $\Laustar f$ to $\rega(\Laustar f)$ satisfies Definition 
\refer{d: holomorphic s globality new} (a)
with $\Omega^*=\Omega'$.
\qed

The following definition is an analogue of the final part of
Definition \refer{d: defi cEhyp Q Y gd}, 
replacing the globality condition by a condition of
holomorphic globality.
\begin{defi}
\naam{d: family for special case of the vanishing thm newer}
Let $Q\in \allparabs$ and let $\gd \in \DQmaps.$ We define 
$$
\cEhypQgdhglob
$$ 
to be the space of functions  $f \in \cEhypQgd$ satisfying the following 
condition.

\hbox{\hspace{-12pt}
\vbox{\vspace{-3mm}
\begin{enumerate}
\item[]
For each $s \in W$ and every $P \in \parone$ with $s(\faQq) \not\subset \faPq,$ 
the family $f$ is holomorphically $[s]$-global along $(P,v),$ for all 
$v \in \NKaq;$ 
here $[s]$ denotes the image of $s$ in $W/{\sim_{P|Q}}= W_P \bs W / W_Q.$
\end{enumerate}}}
\vspace{-3mm}

\noindent
If $Y \subset \staQqdc$ is a finite subset, we define 
$$
\cEhypQYgdhglob =  \cEhypQYgd \cap \cEhypQgdhglob.
$$ 
\end{defi}

It is easily seen that $\cEhypQgdhglob\subset\cEhypQgd_\glob$.
As in Lemma \refer{l: minimal condition for glob new} the above
condition allows a reduction to a smaller set of $(s,P).$ 

\begin{lemma}
\naam{l: minimal condition for hglob new}
Let $Q\in \allparabs$ be standard, let 
$\gd \in \DQmaps$ and $f \in \cE^\hyp_Q(\spXp\col \tau\col \gd).$ 
Then $f$ belongs to $\cE^\hyp_Q(\spXp\col \tau \col \gd)_\hglob$ if and
only if the following condition is fulfilled.

\hbox{\hspace{-12pt}
\vbox{\vspace{-2mm}
\begin{enumerate}
\item[]
For each $s \in W$ and every $\ga \in \gD$ with $s^{-1} \ga|_{\faQq}  \neq 0,$ 
the family $f$ is holomorphically $[s]$-global along $(\stPga,v),$ 
for all $v \in \NKaq.$ 
\end{enumerate}
\vspace{-2mm}}}
\end{lemma}

\proof
The proof is similar to the proof of 
Lemma \refer{l: minimal condition for glob new},
involving Lemma \refer{l: transformation of holomorphic globality}
instead of Lemma \refer{l: transformation of globality}.
\qed

We now come to the main result of this section, which provides 
a source of functions to which the vanishing theorem
(Theorem \refer{t: special vanishing theorem}) can be applied.

\begin{thm}
\naam{t: source of functions by Lau new}
Let $\gd\in D_{P_0}$ and 
$f \in \cEzerohyp(\spXp\col \tau\col\gd)_\hglob,$ let $Q \in \allparabs$ 
be a standard
$\gs$-parabolic subgroup
and let $\Lau \in \Mer(\staQqdc, \Sigma_Q)^*_\laur.$ 
Put $Y = \supp \Lau.$
Then there exists a $\gd' \in \DQmaps$ such that 
$$
\Laustar f \in \cE^\hyp_{Q,Y}(\spXp\col \tau \col \gd')_\hglob.
$$ 
\end{thm}

\proof
{}From Lemma \refer{l: Lau to family of eigenfunctions new}
it follows that $\cL_*f$ is 
a family in $\cE^\hyp_{Q,Y}(\spXp \col \tau\col \gdmap')$ for some 
$\gdmap' \in \DQmaps.$
Let $s \in W$ and $\ga \in \gD$ be such that 
$s^{-1}\ga |_{\faQq} \neq 0.$  
Then every $t\in W_\ga s W_Q$ also satisfies the condition
$t^{-1}\ga |_{\faQq} \neq 0;$ hence $t(\faQq) \not\subset \fa_{\ga\iq}.$ 
Thus, from the hypothesis it follows that
$f$ is holomorphically $W_\ga t$-global along $(P_\ga,v)$
for every $t$ in the double coset $W_\ga s W_Q.$ 
According to Lemma \refer{l: WPQ as cosets}, see also Remark \refer{r: W Pga Q new},
the latter set equals the class $[s]$ of $s$ for 
the equivalence relation $\sim_{{\stPga}|Q}$ in $W.$ 
It now follows from Proposition
\refer{p: transference of s globality} 
that $\cL_*f$ is holomorphically $[s]$-global 
along $({\stPga},v).$
We conclude that $\cL_* f$ satisfies the conditions 
of Lemma \refer{l: minimal condition for hglob new},
hence  belongs
to $\cE^\hyp_{Q,Y}(\spXp \col \tau \col\gdmap')_\hglob.$ 
\qed

\section{Partial Eisenstein integrals}
\naam{s: partial Eisenstein integrals}
In this section we will define partial Eisenstein integrals and
show that they belong to the families of eigenfunctions
introduced in the previous section.

We start by recalling some properties of Eisenstein integrals.
Let $P\in \minparabs$ be a minimal $\gs$-parabolic subgroup.
Let $(\tau, \Vtau)$ be a finite dimensional unitary representation of $K.$ 
Let $\cW \subset \NKaq$ be a fixed set of representatives
for $W/\WKH.$ Following
\bib{BSmc}, eq.\ (5.1), we define the
complex linear space $\oC= \oCtau$ as the following formal direct sum
of finite dimensional linear spaces
\begin{equation}
\naam{e: defi oCtau}
\oC:= \oplus_{w \in \cW}\;\; \Ci(\spX_{0,w}\col \tauM).
\end{equation}
Every summand in the above sum, as $w \in \cW,$ is a finite dimensional subspace
of the Hilbert space $L^2(\spX_{0,w}, \Vtau);$ here the $L^2$-inner product 
is defined relative to the normalized $M$-invariant measure of the compact
space $\spX_{0,w} = M/ M \cap wHw^{-1}$ and the Hilbert structure of $\Vtau.$ 
Thus, every summand is a finite dimensional Hilbert space of its own right.
The formal direct sum $\oC$ is equipped with the direct sum inner product,
turning (\refer{e: defi oCtau}) into an orthogonal direct sum.

For $\psi \in \oC,$ $\gl \in \faqdc$ and $x \in \spX,$ 
the Eisenstein integral $E(\psi\col\gl \col x) = E(P \col \psi \col \gl \col x)$ 
and its normalized version 
$\nE(\psi\col \gl \col x)= \nE(P\col \psi\col \gl \col x)$ are defined as in 
\bib{BSmc}, \S{} 5. The Eisenstein integrals are $\tau$-spherical functions
of $x,$ depend meromorphically on $\gl$ and linearly on $\psi.$  We 
view $\nE(\gl \col x) := \nE(\dotvar \col \gl \col x)$ 
(and similarly its unnormalized version)
as an element 
of $\Hom(\oC, \Vtau)\simeq \Vtau \otimes \oC^*.$ Thus, for generic $\gl \in \faqdc,$ 
$\nE(\gl)$ is a $\tau \otimes 1$-spherical function on $\spX.$  The connection between 
the unnormalized and the normalized Eisenstein integral is now given by the identity
\begin{equation}
\naam{e: relation unnormalized and normalized Eis}
\nE(\gl\col x ) = E(\gl\col x ) \after C(1 \col \gl)^{-1},\qquad (x \in\spX),
\end{equation}
for generic $\gl \in \faqdc.$ 
Here $C(1 \col \gl) := C_{P|P}(1 \col \gl)$ is a meromorphic $\End(\oC)$-valued
function of $\gl \in \faqdc;$ see \bib{BSmc}, p.\ 283.
  
The Eisenstein integral is $\DGH$-finite. In fact, we recall from \bib{BSmc}, eq.\ (5.11),
that there exists a homomorphism $\mu$ from $\DGH$ to the algebra of 
$\End(\oC)$-valued polynomial functions  on $\faqdc$ such that
$$%\begin{equation}
%\naam{e: D on nE}
D \nE(\gl) = [I \otimes \mu(D\col \gl)^*]  \nE(\gl), \qquad (D \in \DGH).
$$%\end{equation}
It now follows from Lemma \refer{l: DX finite in exppol} that,
for generic $\gl \in \faqdc,$ the Eisenstein integral $\nE(\gl)$ 
belongs to $C^\ep(\spXp \col \tau \otimes 1).$ It therefore has expansions
of the form (\refer{e: expansion f on PqPw with m}). 
These expansions have been determined
explicitly in  \bib{BSexp}. We 
recall some of the results of that 
paper.

In \bib{BSexp}, eq.\ (15), we define a function
$\Phi_{P}(\gl \col \dotvar)$ on $\Aqp(P)$ by an exponential polynomial 
series with coefficients in $\End(\VtauMKH)$ of the form
\begin{equation}
\naam{e: expansion Phi gl}
\Phi_{P}(\gl \col a) = 
a^{\gl - \rho_P} \sum_{\nu \in \Delta(P)} a^{-\nu} \Gamma_{{P},\nu}(\gl),\qquad 
(a \in \Aqp(P)).
\end{equation}
Note that here ${P}$ replaces the $Q$ of \bib{BSexp}, Sect. 5; also, in
\bib{BSexp} we suppressed the $Q$ in the notation. The coefficients
in the expansion (\refer{e: expansion Phi gl}) 
are defined by recursive relations (see \bib{BSexp}, eq.\ (18) 
and Prop.\ 5.2); it follows from these that the coefficients
depend meromorphically on $\gl,$ and that the expansion (\refer{e: expansion Phi gl})
converges to a smooth function on $\Aqp(P),$ 
depending meromorphically on $\gl.$ In fact, we have the following stronger result.

Let $\PiSRaq$ be the collection of polynomial
functions $\faqdc \to \C$ that can be written as finite 
products of linear factors of the form $\gl \mapsto \inp{\gl}{\ga} - c,$ with 
$\ga  \in \Sigma$ and $c \in \R.$  
For $R \in \R,$ we define the set 
$$
\faqd({P},R):= \{ \gl \in \faqdc \mid \Re\inp{\gl}{\ga} < R \;\; \forall\ga \in \Sigma(P)\}.
$$ 
\begin{lemma}
\naam{l: neat convergence of series Phi}
Let $R \in \R.$ Then there exists a polynomial function $p \in \PiSRaq$ such that
the functions $ p \Gamma_{{P},\nu},$ for $\nu \in \N\Delta(P),$ 
are all regular on $\faq({P},R).$ 
Moreover, if $p$ is a polynomial function with the above property, then 
the series 
\begin{equation}
\naam{e: series with Gamma}
\sum_{\nu \in \N\Delta(P)} a^{-\nu} p(\dotvar) \Gamma_{{P},\nu}(\dotvar)
\end{equation} 
converges neatly on $\Aqp(P)$ as an exponential series with coefficients
in $\cO(\faqd({P},R)) \otimes \End(\VtauMKH).$ In particular, 
the function $(a,\gl) \mapsto p(\gl)\Phi_{P}(\gl \col a)$ is smooth
on $\Aqp(P) \times \faqd({P},R),$ and in addition holomorphic in its 
second variable.
\end{lemma}

\proof
Let $p_R$ be the polynomial function 
described in \bib{BSexp}, Thm.\ 9.1.
As in the proof of that theorem, it follows from the 
estimates in \bib{BSexp}, Thm.\ 7.4, that the power series
$$ 
\Psi(\gl\col z) = \sum_{\nu \in \N \gD(P)} z^{-\nu} p_R(\gl)\Gamma_{{P},\nu}(\gl) 
$$ 
converges absolutely locally uniformly in the variables $z \in D^{\Delta(P)}$ and
$\gl \in \faqd({P},R).$  Here we have used the notation of Sect.\ 1 of the present paper.
Since $p_R(\gl) \Phi_{P}(\gl\col a) = a^{\gl - \rho_P} \Psi(\gl\col \uz(a)),$ 
for $a \in \Aqp(P),$  
this implies all assertions of the lemma with $p_R$ in place of $p.$ 

This is not immediately good enough, since $p_R$ is a finite product of 
linear factors of the form $\gl \mapsto \inp{\gl}{\nu} -c,$ with 
$\nu \in \N\gD(P)$ and $c \in \R,$ see \bib{BSexp}, the equation preceding Lemma 7.3.
To overcome this, we invoke 
\bib{BSexp}, Prop.\  9.4. It follows from that result 
and its proof that there exists a $p \in \PiSRaq$ such that 
$p\Gamma_{P, \nu}$ is regular on $\faqd({P},R),$ 
for every $\nu \in \N {\Delta(P)}.$
Let $p$ be any polynomial with this property, and let $\bpr p $ be the least common 
multiple of $p$ and $p_R.$ Then all assertions of the lemma hold with $\bpr p$ in
place of $p.$ Let $q$ be the quotient of $\bpr p$ by $p.$ 
Denote the image of the linear endomorphism $m_q: \gf \mapsto q \gf$ of  $\cO(\faqd({P},R))$ 
by $\cF,$ and equip this space with the locally convex topology
induced from $ \cO(\faqd({P},R)).$  It follows from an easy application
of the Cauchy integral formula that $m_q$ is a topological linear isomorphism
from $\cO(\faqd({P}, R))$ onto $\cF;$ see also \bib{BSmc}, Lemma 20.7. 
 As said above, all assertions of the lemma hold with $\bpr p$ in place of $p;$ 
on the other hand, by the hypothesis the series 
(\refer{e: series with Gamma}) with $\bpr p$ in place of $p$ has coefficients in $\cF.$ 
Applying the continuous linear map $m_q^{-1}$ to that series, 
we infer that all assertions of the lemma are true
with the polynomial $q^{-1} \bpr p = p.$ 
\qed
 
Following
\bib{BSexp}, Sect. 11, we define  the function 
$\Phi_{{P},w}: \faqdc\times \Aqp(P) \to \End(\VtauMKwH),$ for  $w \in \cW,$ by 
\begin{equation}
\naam{e: defi Phi P w} 
\Phi_{{P},w}(\gl\col a) = \tau(w) \after \Phi_{w^{-1}{P} w}(w^{-1} \gl 
\col w^{-1} aw) \after \tau(w)^{-1}.
\end{equation}
Following \bib{BSmc}, p.\ 283, we define normalized $C$-functions 
$\nC(s \col \gl) = \nC_{{P}|{P}}(s\col \gl),$ for $s \in W,$  by 
\begin{equation}
\naam{e: defi normalized c function} 
\nC(s\col \gl)= C(s\col \gl) \after C(1 \col \gl)^{-1};
\end{equation}
these are $\End(\oC)$-valued meromorphic functions of $\gl \in \faqdc.$ 
{}From (\refer{e: relation unnormalized and normalized Eis}) and \bib{BSexp}, eq.\ (54),
we now obtain the following description of  the normalized Eisenstein integral 
in terms of the functions $\Phi_{{P},w}.$
Let $\psi \in \oC$ and $w \in \cW.$ Then, for $a \in \Aqp(P),$ 
\begin{equation}
\naam{e: deco nE in Phi P w}
\nE(\gl \col aw)\psi = \sum_{s \in W} \Phi_{{P},w}(s \gl \col a)
[\nC(s\col \gl)\psi]_w(e),
\end{equation}
as a meromorphic identity in $\gl \in \faqd.$ 

{}From 
(\refer{e: defi Phi P w}) and (\refer{e: expansion Phi gl}) it follows that, for $w \in \cW,$ the function
$\Phi_{{P},w}$ 
is given by the series
\begin{equation}
\naam{e: series for Phi P w}
\Phi_{{P},w}(\gl\col a) = a^{\gl - \rho_P} \sum_{\nu \in \N{\Delta(P)}} a^{-\nu} \Gamma_{{P},w,\nu}(\gl),
\end{equation} 
with coefficients
\begin{equation}
\naam{e: defi of Gamma P w}
\Gamma_{{P},w,\nu}(\gl)  = 
\tau(w) \after \Gamma_{w^{-1}{P} w, w^{-1}\nu}(w^{-1} \gl)\after \tau(w)^{-1}.
\end{equation}
We now have the following result on the convergence of the series
(\refer{e: series for Phi P w}).

\begin{cor}
\naam{c: neat convergence of Phi P w}
Let $w \in \cW.$ Then there exists a  locally finite real $\Sigma$-hyperplane configuration
$\Hyp = \Hyp_w$ in $\faqdc$ and a map $d = d_w: \Hyp \to \N,$ such that 
the functions $\Gamma_{{P},w,\nu}$ belong to $\Mer(\faqdc, \Hyp, d, \End(\VtauMKwH)),$ 
for every $\nu \in \N\gD(P).$ Moreover, the series
\begin{equation}
\naam{e: series with Gamma P w}
\sum_{\nu \in \N \gD(P)} a^{-\nu} \Gamma_{{P},w,\nu}
\end{equation}
converges neatly on $\Aqp(P)$ as an exponential polynomial series
with coefficients in the space $\Mer(\faqdc,\Hyp, d, \End(\VtauMKwH)).$ In particular,
the function $\gl \mapsto \Phi_{{P},w}(\gl\col  \dotvar)$ belongs to the space
$\Mer(\faqdc,\Hyp, d, \Ci(\Aqp(P)) \otimes \End(\VtauMKwH)).$
\end{cor}

\proof
For $w =1$ the assertion of the corollary follows immediately 
from Lemma \refer{l: neat convergence of series Phi}.
For arbitrary $w \in \cW$ it then follows by application of 
(\refer{e: defi of Gamma P w}).
\qed

For $s \in W$ we define the so called partial Eisenstein integral 
$E_{+,s}( \gl)= E_{+,s}({P}\col\gl)$ as the $\tau \otimes 1$-spherical 
function  $\spXp \to \Vtau \otimes \oC^*$ determined by 
\begin{equation}
\naam{e: defi partial Eis}
E_{+,s}(\gl\col a w) \psi = \Phi_{{P},w}(s\gl \col a) [\nC(s\col \gl)\psi]_w(e),
\end{equation}
for $\psi \in \oC,\; w \in \cW,\;a \in \Aqp(P)$ and generic $\gl \in \faqdc$ 
(use the isomorphism (\refer{e: the iso T down P cW})).
It follows from Corollary \refer{c: neat convergence of Phi P w}
that $E_{+,s}$ is a meromorphic $C^\infty(\spXp\col\tau\otimes 1)$-valued
function on $\faqdc$.
By sphericality it follows from
(\refer{e: deco nE in Phi P w}) and (\refer{e: defi partial Eis}) that 
\begin{equation}
\naam{e: splitting of Eis}
\nE( \gl) = \sum_{s \in W} E_{+,s}( \gl)\quad\text{on}\quad \spXp.
\end{equation}
It follows from the definitions and the isomorphism (\refer{e: isomorphism of exppol})
 that, for generic $\gl \in \faqdc,$ the function 
$E_{+,s}( \gl)\psi$ belongs to $C^\ep(\spXp \col \tau \otimes 1)$ 
for each $\psi\in \oC$. Moreover, 
\begin{equation}
\naam{e: exponents Eps}
\Exp(P,v\asmid E_{+,s}(\gl)\psi)\subset s\gl - \rho_P -\N \gD(P),
\end{equation} 
for every $v \in \cW$ and hence also for every $v \in \NKaq.$ 
Thus, we see that (\refer{e: splitting of Eis}) is 
the splitting of Lemma \refer{l: splitting lemma} applied to  
the Eisenstein integral.
We abbreviate $E_+(\gl) = E_{+,1}( \gl).$ Then from (\refer{e: defi partial Eis}) and (\refer{e: defi normalized c function}) we see that
$$%\begin{equation}
%\naam{e: formula Ep in Phi}
E_+(\gl)(aw)\psi = \Phi_{{P},w}(\gl\col a)\psi_w(e),
$$%\end{equation}
for $\psi \in \oC,$ $w \in \cW,$ $a \in \Aqp(P)$ and generic $\gl \in \faqdc.$
Moreover, the following holds as a meromorphic identity in $\gl \in \faqdc$
\begin{equation}
\naam{e: Esp in terms of Ep and C}
E_{+,s} (\gl\col x) = E_+(s\gl\col x) \,\nC(s\col \gl). 
\end{equation}

In the next lemma we will need the following notation.
If $\gL \in \fbkdc,$ we denote by $\oC[\gL]$ the subspace of $\oC$ 
consisting of elements $\psi$ 
satisfying $\mu(D\col \gl)\psi = \gg(D\col \gL + \gl)\psi$ for all $D \in \DGH,$ 
$\gl \in \faqdc.$ 
We recall from \bib{BSmc}, eq.\ (5.14), that $\oC$ is a finite direct sum
$$
\oC = \oplus_\gL\; \oC[\gL],
$$  
where $\gL$ ranges over a finite 
subset $L_\tau$ of $\fbkdc.$ 
For each $\gL\in\fbkdc,$ we denote by $\specfamgL$
the space $\cE^\hyp_{0}(\spXp \col \tau \col\gdmap)$
(see Definition \refer{d: defi cEhyp Q Y gd})
where $\gd\in D_P$ is the characteristic function of $\{\gL\}$.

\begin{lemma}
\naam{l: f sub s belongs to specfam}
Let $P \in \minparabs,$ $t \in W$ and $\psi \in \oC[\gL],$ 
where $\gL\in \fbkdc$. Define the family 
$f=f_{\{t\}}:\faqdc \times \spXp \to \Vtau,$ 
by
$$
f(\gl, x) = \Ept(P\col\gl\col x)\psi.
$$ 
Then $f\in\cE^\hyp_{0}(\spXp \col \tau \col\gL)$ 
and $\dega f = 0$.
\end{lemma}

\proof According to Definition \refer{d: defi cEhyp Q Y gd} and 
Remark \refer{r: only local annihilation needed},
in order to prove that $f\in\cE^\hyp_{0}(\spXp \col \tau \col\gL)$
we must establish that $f\in C^{\ep,\hyp}_0(\spXp\col\tau)$ and that 
$f_\gl$ is annihilated by $I_{\gL+\gl}$ for $\gl$ in a non-empty
open subset of $\rega f$.

We first assume that $t=1.$ 
Then $f(\gl,x) = E_{+}(\gl\col x)\psi.$
It follows immediately from \bib{BSexp} Cor.\ 9.3 and the hypothesis on 
$\psi$ that $f_\gl$ is annihilated by $I_{\gL+\gl}$ for generic $\gl\in\faqdc$.
We will now show that 
$f \in C^{\ep,\hyp}_0(\spXp\col \tau).$ 
Let
$\Hyp$ be the union of the hyperplane configurations
 $\Hyp_w,\, w \in \cW,$ of 
Corollary \refer{c: neat convergence of Phi P w},
and let $d: \Hyp \to \N$ be defined by 
$ 
d = \max_{w \in \cW} \,d_w 
$
(see Remark  \refer{r: convention about d}).
Then for every complete locally convex space $\lcspace,$ the spaces
$\Mer(\faqdc, \Hyp_w, d_w, \lcspace)$ are included in the space $\Mer(\faqdc, \Hyp, d, \lcspace),$ 
with continuous inclusion maps.
Hence for each $w\in\cW$
the series  (\refer{e: series with Gamma P w}) converges as
a $\gD(P)$-exponential polynomial series on $\Aqp(P),$ with coefficients in 
the space
$\Mer(\faqdc, \Hyp, d, \End(\Vtau^{\KM \cap wHw^{-1}})).$ 
Moreover, 
the function $\gl \mapsto \Phi_{{P},w}(\gl \col \dotvar)$ 
is contained in $\Mer(\faqdc,\Hyp, d, \Ci(\Aqp(P)) \otimes \End(\VtauMKwH))$.

On the other hand, from (\refer{e: defi partial Eis}) 
and (\refer{e: defi normalized c function}) with $s =1,$ 
it follows that 
\begin{equation}
\naam{e: TdownPw f as Phi}
\TdownPw (f_\gl)(a) = f(\gl, aw) = \Phi_{P,w}(\gl \col a) \psi_w(e),
\end{equation}
for all $w \in \cW,$ $a \in \Aqp(P)$ and $\gl \in \faqdc\setminus \cup \Hyp.$ 
Hence the function $\gl \mapsto \TdownPw(f_\gl)$ belongs
to the space $\Mer(\faqdc, \Hyp, d, \Ci(\Aqp(P),\VtauMKwH)).$ 
In view of the isomorphism (\refer{e: the iso T down P cW}), it now follows
that 
the function $\gl \mapsto f_\gl$ belongs
to $\Mer(\faqdc, \Hyp,d, \Ci(\spXp\col \tau)).$ 
This establishes condition (a) of 
Definition \refer{d: Cephyp}, with $Q = P_0$ and $Y = \{0\}.$

The evaluation map $\psi \mapsto \psi(e)$ 
is a linear isomorphism 
from $\Ci(\spXzerow\col \tauM)$ onto $\VtauMKwH.$ Thus, for $w \in \cW$ 
and $\nu \in \N \DP$ we may define a
function $\tilde q_{1,\nu}(P,w\asmid f) : 
\faqdc \to \Ci(\spXzerow\col \tauM)$ 
by 
\begin{equation}
\naam{e: tilde q one as Gamma}
\tilde q_{1,\nu}(P,w\asmid f , \gl, e) = \Gamma_{P,w,\nu}(\gl)\psi_w(e),
\end{equation}
for $\gl \in \faqdc.$ 
Then $\tilde q_{1,\nu}(P,w \asmid f)
 \in \Mer(\faqd, \Hyp, d, \Ci(\spXzerow\col \tauM)).$ 
Moreover, from what we said earlier about the convergence 
of the series (\refer{e: series with Gamma P w}), it follows that,
for $w \in \cW,$  the series
$$
\sum_{\nu \in \N\DP} a^{-\nu} \tilde q_{1, \nu}(P, w\asmid f)
$$ 
converges neatly as a $\DP$-exponential polynomial series on $\Aqp(P),$ 
with coefficients in $\Mer(\faqdc, \Hyp, d, \Ci(\spXzerow\col \tauM)).$ 

{}From (\refer{e: TdownPw f as Phi}), (\refer{e: series for Phi P w}) 
and (\refer{e: tilde q one as Gamma}) 
it follows by sphericality that, for
 $w \in \cW,$ $\gl \in \faqdc \setminus \cup \Hyp,$ $m \in \spXzerow$
and $a \in \Aqp(P),$ 
$$ 
f_\gl(maw) = a^{\gl -\rho_P} \sum_{\nu \in \N \gD(P)} 
a^{-\nu} \tilde q_{1, \nu} (P,w\asmid f)(\gl, m),
$$ 
This establishes assertions 
(b) and (c) of Definition \refer{d: Cephyp} with a fixed $P,$  arbitrary 
$v\in \cW,$ and, for $\nu \in \N\gD(P),$ $X\in\faq$,
$$
q_{s, \nu}(P,v \asmid f,X) = \left\{ 
\begin{array}{lcl}
\tilde q_{1, \nu}(P, v \asmid f)& \text{for} &s=1;\\
0& \text{for} &s \in W\setminus \{1\}.
\end{array}
\right.
$$ 
In view of Remark \refer{r: second on Cephyp} we have shown that 
$f \in C^{\ep,\hyp}_0(\spXp\col \tau).$ Moreover,
$\dega f=0$.
This completes the proof for $t =1.$

Let now $t \in W$ be arbitrary and 
let $\tilde t \in W_0(\fb)$ be such that $\tilde t |_{\faq} = t;$ see
the text preceding 
Lemma \refer{l: restriction on leading exponents}.
{}From (\refer{e: Esp in terms of Ep and C}) we see that
\begin{equation}
\naam{e: expression f t in Ep and nC}
f(\gl, x) = E_+(t\gl\col x) \nC(t\col \gl) \psi.
\end{equation} 
It follows from \bib{BSmc}, Lemma 20.6, that there exists
a $\Sigma$-configuration $\Hyp'$ in $\faqdc$ and a map $d': \Hyp' \to \N,$ 
such that
\begin{equation}
\naam{e: C-function in MerSigma}
\nC(t\col \dotvar)\in\Mer(\faqdc, \Hyp', d', \End(\oC)). 
\end{equation}
{}From \bib{BSmc}, eq.\ (5.13),
it follows that $\nC(t \col \gl)$ 
maps $\oC[\gL]$ into $\oC[\tilde t \gL].$ 
Fix a basis $\psi_1,\ldots, \psi_s$ for $\oC[\tilde t \gL].$ 
Then there exist unique functions $c_j \in \Mer(\faqdc, \Hyp', d')$ 
such that 
\begin{equation}
\naam{e: nC t gl in components}
\nC(t\col \gl) \psi = \sum_{j=1}^r c_j(\gl) \psi_j.
\end{equation} 
For $1\leq j \leq r$ we define the family $g_j: \faqdc \times \spXp \to \Vtau$ by 
\begin{equation}
\naam{e: defi g j} 
g_j(\gl, x) = E_+(\gl \col x) \psi_j.
\end{equation}
Then by the first part of the proof, each $g_j$ belongs to 
$\cEstarzero(\spXp\col\tau\col\tilde t \gL).$ 
Moreover, for every $1 \leq j \leq r,$ 
the family $g_j$  satisfies the conditions
of Definition \refer{d: Cephyp} with $Q = P_0$ and $Y = \{0\},$ with
$\Hyp$ and $d$ as in the first part of the proof,
and with $k =0.$ 

For $1 \leq j \leq r$ we define the family $f_j: \faqdc \times \spXp \to \Vtau$ by 
$f_j(\gl, x) = g_j(t\gl, x).$ Then we readily see
that $f_j$ satisfies the conditions of Definition \refer{d: Cephyp} with 
$t^{-1}\Hyp$ and $d\after t$ in place of $\Hyp$ and $d,$ respectively,
and with $k =0.$ Hence $f_j \in C^{\ep,\hyp}_0(\spXp\col \tau).$ 
Since $I_{\tilde t \gL +  t\gl} = I_{ \gL+ \gl}$ we see that 
$f_j \in \cEstarzero(\spXp\col \tau\col \gL)$. Moreover,
$\dega f_j=0$.

Combining (\refer{e: expression f t in Ep and nC})
and (\refer{e: nC t gl in components}) with  
(\refer{e: defi g j}) and the definition of $f_j,$ 
we find that 
$$ 
f(\gl, x) 
 = \sum_{j=1}^r c_j(\gl) f_j(\gl, x).
$$
Let $\Hyp'' = t^{-1}\Hyp \cup \Hyp'$  
and define $d'': \Hyp'' \to \N$ 
by $d''(H) = d(tH) + d'(H)$ (see Remark \refer{r: convention about d}).
Then by linearity it readily follows that $f$ satisfies all
conditions of Definition \refer{d: Cephyp}, with $k=0$ and with 
$\Hyp''$ and  $d''$ in place of $\Hyp$ and $d,$ 
respectively. Hence $f \in C^{\ep, \hyp}_0(\spXp\col \tau)$ 
and $\dega f = 0.$ 
Moreover, for generic $\gl,$ 
$f_\gl$ is annihilated by $I_{\gL +\gl},$
and hence $f\in\cEstarzero(\spXp\col \tau\col \gL)$. \qed

\begin{cor}
\naam{c: Lau to partial Eis}
Let assumptions be as in Lemma \refer{l: f sub s belongs to specfam}
and let $Q$ be a $\gs$-parabolic subgroup.
Let $\Lau \in \Mer(\staQqdc, \Sigma_Q)^*_\laur$.
Then $\Laustar f \in \cE^\hyp_{Q,Y}(\spXp\col \tau \col \gd)$
for $Y = \supp \Lau$ and $\gd$ a suitable element in $\DQmaps$. 
Moreover, $$\Exp(P,v\asmid (\Laustar f)_\nu)\subset t(\nu+Y)-\rho_P-\N\gD(P)$$
for $v\in \NKaq$ and $\nu\in\rega \Laustar f$.
\end{cor}
\proof This follows immediately from Lemmas
\refer{l: f sub s belongs to specfam} and 
\refer{l: Lau to family of eigenfunctions new},
and from (\refer{e: exponents Eps})
combined with the final statement in Proposition 
\refer{p: Lau to hyp family} (b).
\qed

\begin{lemma}
\naam{l: globality for the full Eisenstein}
Let $\psi \in \oC[\gL]$ where $\gL\in \fbkdc$.
Then the family $f:\faqdc \times \spXp \to \Vtau,$ defined by
$$
f(\gl, x) = \nE(\stPo\col \gl \col x) \psi 
$$ 
belongs to $\specfamgL.$ Moreover, $\dega f=0$ and
for all $P \in \allparabs, v \in \NKaq$ 
and every 
$s \in W_P\bs W,$ the family $f$ is holomorphically $s$-global along $(P,v).$ 
\end{lemma}
\proof
The function $f$ equals the sum, for $t \in W,$ of the functions $f_{\{t\}}$ defined 
in Lemma \refer{l: f sub s belongs to specfam}, with $\stPo$ in place of $P.$ 
Hence $f \in \specfamgL$ and $\dega f = 0.$
Moreover, for each $\gl\in\rega f$, the function
$f_\gl$ is asymptotically global along all pairs $(P,v)$ 
by Proposition \refer{p: global eigen implies as global}.
Thus, it remains to prove the assertion on holomorphic globality.
In view of Lemma \refer{l: transformation of holomorphic globality},
 it suffices to do this for arbitrary 
$P \in \allparabs$ and the special value $v = e.$ 

In the rest of this proof we shall use  notation of the paper \bib{BSft}.
According to \bib{BSft}, Lemma 14, there exists a locally finite collection $\Hyp$ of
$\Sigma$-hyperplanes such that $\gl \mapsto f_\gl$ is holomorphic 
on $\Omega_0: =  \faqdc\setminus \cup\Hyp,$ 
with values in $\Ci(\spX\col \tau).$ 
According to the same mentioned lemma it follows that
$f \in \cE_*(G/H, \Vtau, \Omega_0).$
According to \bib{BSft}, p.\ 562, Cor.\ 1,  for generic $\gl \in \faqdc$ 
the function $f_\gl$ has an asymptotic expansion of the form
\begin{equation}
\naam{e: expansion Eis in terms of p}
f_\gl(x\exp tX) \sim \sum_{s \in W_P\bs W\atop \nu \in \N\DrP} 
p_{P,\nu}(f_\gl\col s\col \gl)(x)\, e^{(s \gl -\rho_P - \nu)(tX)} \;\; (t \to \infty)
\end{equation}
for $X \in \faPq$ at every $X_0\in \fa_{P\iq}^+.$ 
Proposition 10 of \bib{BSft} is valid with $\cE_*(G/H, \Vtau, \Omega_0)$ 
in place of $\cE_*(G/H, \gL, \Omega_0),$ by the remarks in the beginning of 
\bib{BSft}, Sect.\ 12.
In particular, there exists a full open subset 
$\ppfaqdc$ of $\faqdc$ such that, for all $s \in W_P\bs W$ and 
$\nu \in \N\DrP,$ the coefficient
$p_{P,\nu}(f_\gl\col s\col \gl)$ is holomorphic as a $C^\infty(G,\Vtau)$-valued
function of $\gl$ on the full open set $\Omega_0 \cap \ppfaqdc.$ 

On the other hand, since $f \in \cE^\hyp_0(\spXp\col \tau\col\gL),$ and
$\dega f = 0,$ 
the expansion (\refer{e: expansion Cephyp})
holds, with $k = 0$ and $Y = \{0\},$ 
for all $\gl \in \Omega:= \rega f.$ 
Thus, if $\gl \in  \Omega \cap \fa_{\iq\iC}^{*0}(P, \{0\})$ is generic,
then it follows from comparing the expansions (\refer{e: expansion Eis in terms of p})
and (\refer{e: expansion Cephyp}), and 
using Lemma \refer{l: exponents disjoint} and uniqueness 
of asymptotics (see the proof of 
Lemma \refer{l: uniqueness of asymp}), that 
\begin{equation}
\naam{e: equality of q and p of Eis}
q_{s,\nu}(P,e \asmid f, X, \gl)(m) = p_{P,\nu}(f_\gl\col s\col \gl)(m),
\end{equation}
for all $s \in W_P\bs W,$ $\nu \in \N \DrP,$ 
$X \in \faPq$ and $m \in M_{P,+};$ here we have
written $M_{P,+}$ for the preimage of $\spX_{P,e,+}$ in $M_P.$

By analytic continuation the equality 
(\refer{e: equality of q and p of Eis}) 
holds for all $\gl$ in the full, hence connected, open 
subset  $\Omega':= \Omega \cap \Omega_0 \cap \ppfaqdc$ of $\faqdc.$ 
In particular it follows
that $\gl \mapsto q_{s,\nu}(P,e \asmid f, \gl)$ is holomorphic 
on $\Omega'$ as a function with values in 
$P_0(\faPq) \otimes \Ci(\spX_{P,e}\col \tauP),$ for all
$s \in W_P\bs W$ and $\nu \in \N \DrP.$ This establishes the 
assertion on holomorphic globality, 
see Definition \refer{d: holomorphic s globality new}.
\qed

\begin{lemma}
\naam{l: s globality for partial Eisenstein}
Let $\gL\in \fbkdc$, $\psi \in \oC[\gL],$  $S \subset W$ and define 
$ f_S:\faqdc \times \spXp \to \Vtau$ 
 by 
$$
f_S(\gl,x) := \sum_{s \in S} E_{+,s}(\stPo\col \gl \col x)\psi.
$$ 
Then the family $f_S$ belongs to $\specfamgL.$ Moreover, let $t \in W$ and
$\ga \in \Delta$, and assume that either $W_\ga t\subset S$
or $W_\ga t\cap S=\emptyset$, where $W_\ga = \{1, s_\ga\}.$ Then  
the family $f_S$ is holomorphically $W_\ga t$-global along 
$({\stPga}, v),$ for every $v \in \NKaq.$ 
\end{lemma}
\proof
The first assertion is an immediate consequence of Lemma
\refer{l: f sub s belongs to specfam} with $\stPo$ in place of $P.$ 

Let $v \in \NKaq.$
It follows from (\refer{e: exponents Eps}) and Theorem
\refer{t: transitivity of asymptotics} that
$$\Exp(P_\ga,v\asmid f_{\{s\}\gl})\subset s\gl|_{\fagaq}-\rho_\ga-\N\Dr(\stPga)$$
for each $s\in W$.
For  $\gl$ in the full open subset $\fa_{{\stPga} q \iC}^{*0}(\stPo, \{0\})$ 
of $\fagaqdc$ the sets $s\gl|_{\fagaq} - \rho_\ga - \N\Dr({\stPga})$ 
are mutually disjoint for different $[s]=W_\ga s$ from $W_\ga \bs W,$
see Lemmas \refer{l: exponents disjoint} and \refer{l: WPQ as cosets}.
Hence  
\begin{equation}
\naam{e: q sub t is zero}
q_{[t], \xi}({\stPga},v \mid f_{\{s\}}) = 0,
\end{equation}
for all $s\in W\setminus W_\ga t$ and all $\xi \in \Dr({\stPga}).$

First assume that $W_\ga t \cap S = \emptyset.$ 
Then 
it follows from (\refer{e: q sub t is zero}) that 
$q_{[t],\xi}({\stPga}, v \mid f_S) = 0$ 
for all $\xi \in \Dr({\stPga}).$ Hence $f_S$ is holomorphically $[t]$-global 
along $({\stPga}, v)$.

Next assume that $W_\ga t \subset S.$ Let $S^c=W\setminus S$.
Then $f_S=f_W-f_{S^c}$, and it follows from Lemma 
\refer{l: globality for the full Eisenstein}
and what was just proved, that $f_S$ is holomorphically $[t]$-global along 
$({\stPga}, v).$ 
\qed

If $Q \in \allparabs$ is standard, then we define the subset $W^Q$ of $W$ by 
\begin{equation}
\naam{e: defi W^Q}
W^Q = \{ s \in W \mid s(\gD_Q) \subset \Sigma^+ \}
\end{equation} 
It is well known, see e.g. \bib{Carter}, Thm.\ 2.5.8, that the multiplication 
map $W^Q \times W_Q \to W$ is bijective. Moreover, if $s\in W^Q$ and 
$t \in W_Q,$ then $l(st) = l(s) + l(t);$ 
here $l: W \to \N$ denotes the length function relative to $\gD.$ In particular
this means that
$W^Q$ consists of the minimal length representatives in $W$ of the cosets in 
$W/W_Q.$ 

\begin{lemma}
\naam{l: W Q and left Wga invariance}
Let $s \in W,$ $\ga \in \gD$ and assume that 
$s^{-1}\ga |_{\faQq} \neq 0.$ Let $t \in W_Q$. Then 
$s \in W^Qt$ if and only if $s_\ga s \in W^Qt$.
\end{lemma}
\proof
The hypothesis $s^{-1}\ga |_{\faQq} \neq 0$ is also satisfied by
the elements $s_1 = s t^{-1}$ and $s_2=s_\ga s t^{-1}$.
Hence we need only prove the implication
$s \in W^Q\Rightarrow s_\ga s \in W^Q$.

Assume that $s \in W^Q.$ Then $s(\gD_Q) \subset \Sigma^+.$ 
{}From the hypothesis it follows that  $s^{-1}\ga \notin \gD_Q,$ 
hence $\ga \notin s(\gD_Q).$ Since $\ga$ is simple, it follows
that $s_\ga (s(\gD_Q)) \subset \Sigma^+.$ Hence $s_\ga s \in W^Q$.
\qed

\begin{cor}
\naam{c: Lau on partial WQ Eis}
Let $\psi \in \oC[\gL]$ where $\gL\in\fbkdc$ and let $Q \in \allparabs$ be 
a standard parabolic subgroup. Fix $t \in W_Q,$ and 
let the family $f: \faqdc \times \spXp \to \Vtau$  be defined
by 
$$
f(\gl, x) = \sum_{s \in W^Q} E_{+,st}(\stPo\col\gl\col x)(\psi).
$$
Then $f \in \specfamgL_\hglob.$ 
If $\Lau \in \Mer(\staQqdc, \Sigma_Q)^*_\laur,$ then the 
family  
$\cL_* f$ belongs to the space $\cEhypQYgdhglob,$ 
where $Y = \supp \Lau,$ and where $\gdmap$ is a suitable element of $\DQmaps.$ 
\end{cor}

\proof
Let $S = W^Qt.$
Then,  $f = f_S,$ where we have used the notation 
of Lemma \refer{l: s globality for partial Eisenstein}.  
It follows from the mentioned lemma that $f \in \specfamgL$.
Moreover, let $s \in W$ and $\ga \in \gD$ be such that 
$s^{-1} \ga |_{\faQq} \neq 0.$ 
Then it follows from Lemma
\refer{l: W Q and left Wga invariance} that either
$W_\ga s \subset S$ or $W_\ga s \cap S$ is empty.
Hence it follows from  Lemma \refer{l: s globality for partial Eisenstein} 
that  $f$ is holomorphically $W_\ga s$-global along $(P_\ga, v),$ for every 
$v \in \NKaq.$ Thus $f \in \specfamgL_\hglob$ by Lemma 
\refer{l: minimal condition for hglob new}. 
The remaining assertion now follows from 
Theorem \refer{t: source of functions by Lau new}.
\qed

\section{Asymptotics of partial Eisenstein integrals}
\naam{s: asymptotics of partial Eisenstein integrals}
Let $P\in \minparabs$ and let 
$Q$ be a $\gs$-parabolic subgroup containing $P.$ 
For the application of the asymptotic vanishing theorem, 
Theorem \refer{t: vanishing theorem new},
in the next section we need
to determine the coefficient of the leading exponent in the $(Q,v)$-expansion
of the Eisenstein integral $\nE(P\col \gl),$ for every $v \in \NKaq.$ 
To formulate a result in this direction,
we need some additional notation.

Let $v \in \NKaq$ and select
a complete set of representatives $\cW_{Q,v}$ 
in $\NKQaq$ for
$W_Q/ W_{K_Q \cap vHv^{-1}}.$ 
We define $\oC(Q,v) = \oC(Q,v,\tau)$ to be the analogue of the space
$\oC$ for the data $\spXoneQv, \tauQ.$ Thus 
\begin{equation}
\naam{e: dir sum oC Q v}
\oC(Q,v) = \oplus_{u \in \cW_{Q,v}} \Ci(M/M \cap uv H (uv)^{-1} \col \tau)
\end{equation}
with an orthogonal direct sum. Note that $\oC(Q,v)$ is also the analogue of
$\oC$ for the data $\spXQv, \tauQ.$ 

One readily checks that the map $\cW_{Q,v} \to W/\WKH$  given by 
$u \mapsto \Ad(uv)|\faq$ is injective. Hence we may extend $\cW_{Q,v}v $ 
to a complete set $\cW \subset \NKaq$ of representatives for 
$W/\WKH.$ If $w \in \cW,$ then  $w \in \cW_{Q,v}v \iff wv^{-1} \in K_Q.$ 
With such choices made we have a natural isometric embedding 
$
\iQv : \oC(Q,v) \embed \oC,
$
defined by
\begin{equation}
\naam{e: defi i Q v}
(\iQv\psi)_w  = \left\{
\begin{array}{ll}
\psi_{wv^{-1}}&\text{if} w \in \cW_{Q,v} v,\\
0 &\text{otherwise.}
\end{array}
\right.
\end{equation}
The adjoint of the embedding $\iQv$ is denoted by 
$
\pr_{Q,v} : \oC \to \oC(Q,v).
$
It is given by the following formula, for $\psi \in \oC,$ 
\begin{equation}
\naam{e: formula for pr Q v}
(\pr_{Q,v}\psi)_u = \psi_{uv}, \qquad (u \in \cW_{Q,v}).
\end{equation}

The normalized Eisenstein integral  associated with the data 
$\spXoneQv, \tauQ$ and $\starP:= P \cap \MoneQ$ is denoted
by 
$\nE(\spXoneQv \col \starP \col \nu),$ for $\nu \in\faqdc.$ 
Similarly, the partial Eisenstein integrals associated with these data
are denoted by
$\Eps(\spXoneQv \col \starP\col \nu),$ for $s \in W_Q$ and 
$\nu \in \faqdc.$ Note that all of these are $(\tauQ \otimes 1)$-spherical 
smooth functions on $\spXoneQvp$ with values in 
$\Hom(\oCQv, \Vtau)\simeq \Vtau \otimes \oCQv^*.$

\begin{prop}
\naam{p: q Q v of nE}
Let $P \in \minparabs,$ $Q \in \allparabs$ and assume that $Q\supset P.$ 
Let $v \in \NKaq$, and choose $\cW_{Q,v}$, $\cW$
as above such that $\cW_{Q,v} \subset \cW v^{-1}$. 
Let $\psi \in \oC$ and let the family $f:\faqdc \times \spX \to \Vtau$
be defined by 
$$ 
f(\gl, x) = \nE(P \col \gl \col x) \psi.
$$ 
Then, for $\gl \in \faqdc$ generic, and for all  $X \in \faQq$ and $m \in \spXQvp,$ 
\begin{equation}
\naam{e: q gl for nE}
q_{\gl|_{\faQq} - \rho_Q}(Q,v \asmid f_\gl,X,m) = 
\nE(\spXoneQv \col \starP \col \gl \col m)\,\prQv \psi.
\end{equation}
\end{prop}

\proof
We first assume that 
$v = e.$ Then $\spXoneQv = \spX_{1Q,e} = \MoneQ/ \MoneQ \cap H.$ 
Moreover, the set $\cW_Q:= \cW_{Q,e}$ is contained in $\cW.$ 
{}From \bib{BSft}, p.\ 563, Thm.\ 4, it follows that
$$
q_{\gl|_{\faQq} - \rho_Q}(Q,e\asmid f_\gl , X, m) = 
\nE(\spX_{1Q,e} \col \starP \col \gl \col m)\, \pr_Q \psi,
$$ 
for generic $\gl \in \faqdc$ and all $X \in \faQq$ and $m \in \spX_{Q,e,+}.$ 
Here $\pr_Q$ is the natural projection map from $\oC$ onto 
$\oC_Q(\tau) = \oplus_{v \in \cW_Q} \Ci(M/M\cap vHv^{-1}\col \tauM),$ 
see \bib{BSft}, pp.\ 544 and 547. Thus, $\pr_Q$ equals the map $\pr_{Q,e}$ defined above
and it follows that (\refer{e: q gl for nE}) holds with $v=e.$ 
To establish the result for arbitrary $v\in \NKaq,$ we first need a lemma.
\medbreak
{}From  Remark \refer{r: extreme cases subspaces} 
we recall that $\spX_v = \spX_{1G,v} = G/vHv^{-1}.$
The set $\cW_{G,v}:= \cW v^{-1}$ is a complete set of representatives
for $W/W_{K \cap vHv^{-1}}.$ Accordingly, the analogue $\oC(G,v) = \oC(G,v,\tau)$ 
of the space $\oC$ is given by 
(\refer{e: dir sum oC Q v}) with $G$ in place of $Q.$ 
The associated map $\rmi_{G,v}: \oC(G,v) \to  \oC$ is now 
a bijective isometry; moreover, its adjoint 
$\pr_{G,v}$ is its two-sided inverse.

We recall from the end of Section \refer{s: asymp walls} that 
right translation by $v$ induces a topological linear
isomorphism $R_v$ from $\Ci(\spX\col \tau)$ onto 
$\Ci(\spX_{v} \col \tau).$ In the following lemma 
we will relate the right translate of $\nE(P\col \gl)$ to the normalized 
Eisenstein integral associated with $\spX_{v},$ $\cW v^{-1}$ and $P.$

\begin{lemma}
\naam{l: Rv of nE}
Let $\psi \in \oC.$ Then, for generic $\gl \in \faqdc,$  
\begin{equation}
\naam{e: Rv of nE}
R_v (\nE(\spX\col P\col \gl)\psi) = \nE(\spX_{v} \col P \col \gl) 
[\pr_{G,v}\psi].
\end{equation}
The formula remains valid if  the normalized Eisenstein integrals 
are replaced by their unnormalized versions. 
\end{lemma}

\proof
We first prove the formula for the unnormalized Eisenstein integrals.
Let $\gl \in \faqdc$ be such that 
$\inp{\Re \gl + \rho_P }{\ga} < 0$ for all $\ga \in \Sigma(P).$ 
Define the function $\tilde\psi(\gl): G \to V_\tau$ as in \bib{BSft}, Eq.\ (19).
Then $E(P\col \gl \col x) \psi = \int_K \tau(k)\, \tilde\psi(\gl \col k^{-1} x) \, dk.$ 
Hence $E(P\col \gl \col x v)\psi = \int_K \tau(k) \tilde\psi_{G,v}(\gl\col k^{-1} x) \, dk,$ 
where $\tilde\psi_{G,v}(\gl\col x) = \tilde\psi(\gl\col xv).$ One now readily checks 
that $\tilde \psi_{G,v}(\gl)$ is the analogue of $\tilde\psi(\gl),$ associated
with the data $\spX_{v}, \cW v^{-1}$ and with the element 
$\psi_{G,v} := \pr_{G,v}\psi $ of  $\oC(G,v).$ {}From
this we obtain the equality (\refer{e: Rv of nE}) for
the present $\gl.$ For general $\gl,$ the result follows by
meromorphic continuation. 

Let $Q \in \minparabs.$ Then it follows, by application of Lemma \refer{l: q of Rv f}
and
the definition of the $\bf c$-functions (cf.\ \bib{BSft}, \S{} 4), that, for every $s \in W,$ each $u \in \cW v^{-1}$ 
and generic $\gl \in \faqdc,$ we have  
$[C_{Q|P}(\spX\col s \col \gl) \psi]_{uv} = 
[C_{Q|P}(\spX_{v}\col s\col \gl) \pr_{G,v} \psi]_u.$ In other words,
$$ 
\pr_{G,v} \after  C_{Q|P}(\spX\col s \col \gl) = C_{Q|P}(\spX_{v}\col s\col \gl)\after 
 \pr_{G,v}.
$$ 
The proof is completed by combining this equation, after substitution of  $P$ and $1$ 
for $Q$ and $s,$ respectively, 
with the unnormalized version of 
(\refer{e: Rv of nE}) and the definitions of the normalized Eisenstein integrals 
(cf.\ \bib{BSft}, eq.\ (49)).
\qed
\medbreak\noindent
{\bf Completion of the proof of Prop.\ \refer{p: q Q v of nE}.\ }
Let $v \in \NKaq$ be arbitrary. Then from Lemmas \refer{l: q of Rv f},
\refer{l: Rv of nE} and equation
(\refer{e: q gl for nE}) with $\spX_{v},$ $e$ and $\pr_{G,v}\psi$ 
in place of $\spX,$  $v$ and $\psi,$ respectively,
it follows that, for $X \in \faQq$ and $m \in \spXQvp,$  
\begin{eqnarray*}
q_{\gl|_{\faQq} - \rho_Q}(Q,v \asmid f_\gl , X, m) &=&
q_{\gl|_{\faQq} - \rho_Q}(Q, e \asmid R_v(f_\gl) , X, m) \\
&=& \nE(\tilde \spX_{1Q,e} \col \starP \col \gl \col m) \, \tilde \pr_{Q,e} \pr_{G,v}\psi. 
\end{eqnarray*}
In the last expression the two tildes over objects indicate that the analogues of 
the objects for the symmetric space $\spX_{v}$ are taken.
We now observe that $\tilde \spX_{1Q,e}$ equals  the space 
$\MoneQ / \MoneQ \cap e vHv^{-1} e = \spX_{1Q,v}.$ 
Hence, to establish (\refer{e: q gl for nE}), it suffices to show that $\tilde \pr_{Q,e} \pr_{G,v}\psi 
=\pr_{Q,v} \psi.$ 
For this we note that $\tilde \pr_{Q,e}$ is the projection
from $\oC(G,v)$ onto the sum of the components
parametrized by the elements $u$ of $M_{1Q} \cap \cW v^{-1} = \cW_{Q,v}.$ 
Moreover, for $u \in \cW_{Q,v},$  
$$
[\tilde \pr_{Q,e} \pr_{G,v}\psi]_u = [\pr_{G,v}\psi]_{u} =
\psi_{uv} 
=
[\pr_{Q,v}\psi]_u.
$$
\qed
The result just proved generalizes to partial Eisenstein integrals.
\begin{prop}
\naam{p: q Q v of sum partial eis}
Let $P \in \minparabs.$  
Let $\psi \in \oC,$ let 
$S \subset W$ and let the family $f = f_S $ 
be defined by 
$$ 
f(\gl, x) = \sum_{s \in S} \Eps(P\col \gl \col x) \psi,
$$ 
see Lemma \refer{l: s globality for partial Eisenstein}.
Assume that $Q \in \allparabs$ contains $P$ and that $v \in \NKaq.$ 
Then, for generic $\gl \in \faqdc,$ and all  $X \in \faQq$ and $m \in \spXQvp,$
\begin{equation}
\naam{e: q gl for sum partial Eis}
q_{\gl|_{\faQq} - \rho_Q}(Q,v \asmid f_\gl,X,m) = 
\sum_{s \in S \cap W_Q} 
\Eps(\spXoneQv \col \starP \col \gl\col m)\,\prQv \psi.
\end{equation}
In particular, if $S\cap W_Q=\emptyset$ then 
$\gl|_{\faQq} - \rho_Q\notin \Exp(Q,v\asmid f_\gl)$.
\end{prop}

\proof
For $S = W$ this result is precisely  
Prop.\ \refer{p: q Q v of nE}. 
We shall use transitivity of asymptotics to derive the result for arbitrary $S$ 
from it.

It suffices to prove the 
above identity for $m = bu \in \spXQvp,$ 
with  $u \in \NKQaq$ and $b \in \stAQqp(\starP)$ arbitrary.

According to Lemma \refer{l: f sub s belongs to specfam} 
and Remark \refer{r: relation between hypfam and fam},
the function $f_S$ belongs to 
$C^{\ep}_{0,\{0\}}(\spXp\col \tau\col \Omega),$ 
for the full open subset $\Omega: = \rega{f_S}$ of $\faqdc.$ 

Hence, according to Theorem \refer{t: transitivity of asymptotics for families new}
with 
$\stPo, Q$  and $P$ 
in place of $Q, P$ and $\Pmin,$ respectively, 
for $\gl \in \faqdc$ generic the following holds, with 
$[1]$ the class of $1 \in W$ in $W/\sim_{Q|\stPo} = W_Q\bs W,$ 
\begin{eqnarray*}
q_{\gl|_{\faQq} - \rho_Q}(Q,v \asmid f_{S\gl},X,bu)
&=& q_{[1],0}(Q,v\asmid f_S, X)(\gl, bu) \\
& = &
\sum_{s \in W_Q}
\sum_{\mu \in \N\Delta_Q(P)} b^{s \gl - \rho_P - \mu} 
q_{s,\mu}(P, uv\asmid f_S, X + \log b)(\gl, e).
\end{eqnarray*}
Now, for all $s,t\in W$, $\mu\in\N\gD$ and $v\in\NKaq$ it follows from
(\refer{e: exponents Eps}) and Lemma \refer{l: exponents disjoint}
that $q_{s,\mu}(P,v\asmid f_{\{t\}})=0$ if $s\neq t$. Hence
$$ 
q_{s,\mu}(P, v\asmid f_S ) =
\left\{
\begin{array}{ll}
q_{s,\mu}(P, v\asmid f_W)& \text{if} s \in S, \\
0 & \text{otherwise.}
\end{array}
\right.
$$ 
Thus,
we obtain that 
\begin{equation}
\naam{e: formula for q one zero of fS}
q_{\gl|_{\faQq} - \rho_Q}(Q,v \asmid f_{S\gl},X,bu)
=
\sum_{s \in S \cap W_Q}
\sum_{\mu \in \N\Delta_Q(P)} b^{s \gl - \rho - \mu} 
q_{s,\mu}(P, uv\asmid f_W, X + \log b)(\gl, e).
\end{equation} 
This equation is valid for any subset $S$ of $W;$ in particular,
it 
is valid for $S = W.$ Using (\refer{e: q gl for nE}) we now obtain that,
for any $u \in \NKQaq$ and all $b \in \stAQqp(\starP),$
\begin{equation}
\naam{e: Eisenstein integral for spXQv}
\nE(\spXoneQv\col \starP \col \gl\col bu)\,\prQv \psi
=
\sum_{s \in W_Q}
\sum_{\mu \in \N\Delta_Q(P)} b^{s \gl - \rho - \mu} 
q_{s,\mu}(P, uv\asmid f_W, X + \log b)(\gl, e).
\end{equation}
This is the $\DQP$-exponential polynomial expansion of the Eisenstein integral
along $(\starP, u).$ In view of (\refer{e: splitting of Eis})
and the remark following (\refer{e: exponents Eps}), with
$\spX_{1Q,v}$ in place of $\spX$,
we infer from (\refer{e: Eisenstein integral for spXQv}) that,
for each $s \in W_Q,$ and every $u\in \NKQaq$ and $b \in \stAQqp(\starP),$  
\begin{equation}
\naam{e: expansion for partial Eisenstein integral along Q}
\Eps(\spXoneQv \col \starP \col \gl\col bu)\prQv \psi
=
\sum_{\mu \in \N\Delta_Q(P)} b^{s \gl - \rho - \mu} 
q_{s,\mu}(P, uv\asmid f_W, X + \log b)(\gl, e).
\end{equation} 
Finally, (\refer{e: q gl for sum partial Eis}) with $m = bu$ follows 
by combining (\refer{e: formula for q one zero of fS}) 
and (\refer{e: expansion for partial Eisenstein integral along Q}).
\qed

We end this section with a 
generalization of Proposition \refer{p: q Q v of sum partial eis}, 
involving the application of a Laurent functional.

\begin{prop}
\naam{p: q of Laustar of part Eis}
Let assumptions be as in Prop.\ \refer{p: q Q v of sum partial eis}
and let 
$\cL \in \Mer(\staQqdc, \Sigma_Q)^*_\laur.$ 
Then the family $\Lau_*f$ defined by $\Lau_*f(\nu,x) = \Lau[f(\dotvar + \nu , x)],$
for generic $\nu \in \faQqdc$ and $x \in \spXp,$ belongs to
 $\cEhypQYgd,$ 
with $Y = \supp \Lau$ and for a suitable $\gd \in \rmD_Q.$ 

Moreover, for generic $\nu \in \faQqdc$ 
and all $X \in \faQq$ and $m \in \spXQvp,$
\begin{equation}
\naam{e: formula for q of Laustar f}
q_{\nu- \rho_Q}(Q,v\asmid (\cL_* f)_\nu , X, m) = 
\cL[\sum_{s \in S \cap W_Q} 
\Eps(\spXoneQv \col \starP \col  \dotvar + \nu  \col m)\,\prQv \psi].
\end{equation}
In particular, if $S\cap W_Q=\emptyset$ then 
$\nu - \rho_Q\notin \Exp(Q,v\asmid (\cL_* f)_\nu)$.
\end{prop}

\proof 
The first assertion follows from Cor.\ \refer{c: Lau to partial Eis}. 
For the 
second assertion, we note that $\Lau_*f \in \CepQY(\spXp\col \tau \col \Omega),$ 
where $\Omega$ is the full open subset $\faQqdc\setminus\cup \Hyp_{\Laustar f}$ 
of  $\faQqdc,$ see Remark \refer{r: relation between hypfam and fam}.
The set $\Omega_*: = \Omega \cap \faQqdczero(P,\{0\})$ is a full open subset
of $\faQqdc.$ Moreover, from (\refer{e: q of f in gl versus q of fgl new})
it follows that, 
for $\nu \in \Omega_*,$ 
\begin{equation}
\naam{e: q nu min rho as q one zero}
q_{\nu - \rho_Q}(Q,v\asmid (\Lau_* f)_\nu , X) = q_{[1],0}(Q,v\asmid \Lau_* f, X)(\nu),
\quad (X \in \faQq);
\end{equation}
here $[1]$ denotes the image of the identity element of $W$ in $\WQQ.$ 
The expression on the right-hand side of the above equation is given 
by (\refer{e: sum for q of Laustar f}), with $P=Q, \gs = [1] \in \WQQ$ and $\xi = 0.$ 
Note that an element $s \in W$ satisfies $[s] = [1]$ if and only if $s \in W_Q.$ 
It follows from this that $[1]\cdot Y = \{0\}.$ 
Hence from (\refer{e: sum for q of Laustar f}) and (\refer{e: defi Laustar s})
we conclude, with $\bar 1$ denoting
the image of $1 \in W$ in $W_Q\bs W,$  
\begin{eqnarray}
q_{[1],0}(Q,v\asmid \Lau_* f, X)(\nu) 
&=& 
\sum_{\gl \in Y} \Lau_{\gl *}^{Q, \bar 1} [q_{\bar 1, 0}(Q,v\asmid f)(X,\dotvar)](\nu,X)
\nonumber\\
\naam{e: q one zero as Laustar of q one zero}
&=&
\sum_{\gl \in Y} 
e^{-(\gl + \nu)(X)}\Lau_{\gl*}
[ e^{(\dotvar)(X)}q_{\bar 1, 0}(Q,v\asmid f)(X, \dotvar)](\nu).
\end{eqnarray}
for $X \in \faQq$ and generic $\nu \in \faQqdc.$ From
$(\gl + \nu)(X) = \nu(X)$ we deduce that the last expression 
in (\refer{e: q one zero as Laustar of q one zero}) equals
$\sum_{\gl \in Y} \Lau_{\gl*}
              [q_{\bar 1, 0}(Q,v\asmid f)(X, \dotvar)](\nu).$ 
Hence from (\refer{e: q nu min rho as q one zero}) and 
(\refer{e: q one zero as Laustar of q one zero}) we obtain
\begin{equation}
\naam{e: q as Laustar q}
q_{\nu - \rho_Q}(Q,v\asmid (\Lau_* f)_\nu , X)
=
\Lau_{*}[q_{\bar 1, 0}(Q,v\asmid f)(X, \dotvar)](\nu).
\end{equation}
It follows from (\refer{e: q gl for sum partial Eis})
and (\refer{e: q of f in gl versus q of fgl new})
that, for $X \in \faQq,$ $m \in \spXQvp,$ 
\begin{equation}
\naam{e: q bar one zero as Eis} 
q_{\bar 1, 0}(Q,v\asmid f)(X, \gl, m) =
\sum_{s \in S \cap W_Q} 
\Eps(\spXoneQv \col \starP \col \gl\col m)\,\prQv \psi,
\end{equation} 
as a meromorphic identity in $\gl \in \faqdc.$
The equality (\refer{e: formula for q of Laustar f})
now follows by combining
 (\refer{e: q as Laustar q})
with 
(\refer{e: q bar one zero as Eis}).
\qed

\section{Induction of relations}
\naam{s: induction of relations}
After the preparations of the previous sections
we are now able to apply the vanishing theorem, 
Theorem \refer{t: special vanishing theorem}, to families
obtained from applying Laurent functionals to partial Eisenstein integrals.
This will lead to what we call induction of relations.

We retain the notation of the previous section.
Moreover, we assume that $Q \in \allparabs$ is a standard parabolic subgroup.
Thus $\starstPo: = M_Q \cap \stPo$ is the standard
minimal $\gs$-parabolic subgroup of $M_Q,$ relative to the positive
system $\Sigma_Q^+:= \Sigma_Q \cap \Sigma.$ 

We assume that $\QW$ is a complete
set of representatives in $\NKaq$ for the double coset space $W_Q\bsl W/\WKH.$  
We also assume that for each $v \in \QW$ 
a set $\cW_{Q,v}$ as above (\refer{e: dir sum oC Q v}) is chosen.
Then one readily verifies that
\begin{equation}
\naam{e: cW as disjoint union over QW}
\cW = \cup_{v \in \QW} \;\;\cW_{Q,v} v \qquad \text{(disjoint union).}
\end{equation} 
is a complete set of representatives for $W/\WKH$ in $\NKaq.$ 
Combining this with (\refer{e: defi i Q v}) and (\refer{e: formula for pr Q v})
we find that  
$$
\sum_{v \in \QW} \iQv\after \pr_{Q,v} = I_{\oC}.
$$
Combining
(\refer{e: cW as disjoint union over QW}) with 
(\refer{e: defi i Q v}) and (\refer{e: formula for pr Q v}), 
it also follows,  for $u,v \in \QW,$ that
\begin{equation}
\naam{e: prQu after iQv}
\pr_{Q,u} \after \rmi_{Q,v} = \left\{
\begin{array}{ll}
I_{\oC(Q,v)} & \text{if} u=v;\\
0 &\text{otherwise.}
\end{array}
\right.
\end{equation}

\begin{thm}
\naam{t: induction of relations, new}
{\rm (Induction of relations)\ }
Let 
$\cL_t\in \cM(\staQqdc, \Sigma_Q)^*_\laur \otimes \oC$ be given
for each $t\in W_Q$.
If, for each $v \in \QW$, 
\begin{equation}
\naam{e: hypo ind rels with pr, new}
\sum_{t\in W_Q}\cL_t[E_{+,t}(\XQv\col\starstPo\col\dotvar\col m) 
\after \pr_{Q,v} ] 
=
0,\qquad (m \in \XQvp)
\end{equation}
then for each $s\in W^Q$ the following holds 
as a meromorphic identity in $\nu \in \faQqdc:$ 
\begin{equation}
\naam{e: conclusion ind rels with pr, new}
\sum_{t\in W_Q}\cL_t[\;
E_{+,st}(\spX\col\stPo\col  \dotvar + \nu \col x)\;]
= 0, \qquad (x \in \spXp).
\end{equation}

Conversely, if the identity (\refer{e: conclusion ind rels with pr, new}) 
holds for some $s\in W^Q$ and all
$\nu$ in a non-empty open subset of $\faQqdc$, then 
(\refer{e: hypo ind rels with pr, new}) holds 
for each $v \in \QW$.
\end{thm}

\proof
Define for each $w\in W$ the family
$g_w\colon \faQqdc \times \spXp \to \Vtau \otimes \oC^*$  by 
$$
g_w(\nu,x) = \cL_t[\;E_{+,st}(\spX\col P_0\col  \dotvar+\nu \col x)\;]
$$ 
for generic $\nu\in \faQqdc$ and every $x \in \spXp;$ 
the elements $s\in W^Q, t\in W_Q$ are determined by
the unique product decomposition $w=st$ (see below (\refer{e: defi W^Q})).

It follows from Cor.\ \refer{c: Lau to partial Eis}, that there
exist $\gd_w \in \rmD_Q$ such that 
$g_w \in  \cE_{Q,Y_w}^\hyp(\spXp \col \tau\col \gd_w)$;
here $Y_w = \supp \Lau_t,$ where $t\in W_Q$ is determined as above.
If we put $Y = \cup Y_w$ and $\gd = \max(\gd_w),$ then
$g_w$ belongs to $\cEQY^\hyp(\spXp\col \tau\col \gd)$ 
for all $w\in W$. Moreover, for generic $\nu\in\faQqdc$,
\begin{equation}
\naam{e: 3}
\Exp(P_0,v\asmid (g_w)_\nu)\subset w(\nu+Y)-\rho-\N\Delta.
\end{equation}
In view of Proposition \refer{p: q of Laustar of part Eis} it also follows 
for $X\in\faQq$, $m \in \XQvp$ and generic $\nu\in \faQqdc$ that
\begin{equation}
\naam{e: 4}
q_{\nu - \rho_Q}(Q,v\asmid (g_t)_\nu, X, m) = 
\cL_t[E_{+,t}(\XQv\col\starstPo\col\dotvar+\nu\col m) 
\after \pr_{Q,v} ] \qquad (t\in W_Q), 
\end{equation}
and 
\begin{equation}
\naam{e: 4'}
q_{\nu - \rho_Q}(Q,v\asmid (g_w)_\nu, X, m) = 0
\qquad (w\notin W_Q).
\end{equation}

According to Cor.\ \refer{c: Lau on partial WQ Eis} 
the family $\sum_{s\in W^Q} g_{st}$
belongs to the space
$\cE_{Q,Y}^\hyp(\spXp\col \tau \col \gd)_\glob$ for each $t\in W_Q$.
Hence so does the family 
$g=\sum_{w\in W} g_w=\sum_{t\in W_Q,s\in W^Q} g_{st}$.
Moreover, by (\refer{e: 4}) and (\refer{e: 4'})
$$
q_{\nu - \rho_Q}(Q,v\asmid (g)_\nu, X, m) = 
\sum_{t\in W_Q}\cL_t[E_{+,t}(\XQv\col\starstPo\col\dotvar+\nu\col m) 
\after \pr_{Q,v} ] \qquad (m \in \XQvp).
$$
{}From Theorem \refer{t: special vanishing theorem}
we now see that (\refer{e: hypo ind rels with pr, new})
holds for each $v \in \QW$ if and only if $g=0$.

On the other hand, let $g^s=\sum_{t\in W_Q} g_{st}$ for $s\in W^Q$.
It follows from (\refer{e: 3}) that 
$$\Exp(P_0,v\asmid (g^s)_\nu)
\subset s\nu+WY-\rho-\N\Delta.$$
Since the latter sets are mutually disjoint as $s$ runs over
$W^Q$, for $\nu$ in a full open subset
(see Lemma \refer{l: exponents disjoint}), we conclude that
for such $\nu$,
$$(s\nu+WY-\rho-\N\Delta)\cap\Exp(P_0,v\asmid g_\nu) 
=\Exp(P_0,v\asmid (g^s)_\nu).$$
Hence $g=0$ implies that $g^s=0$
for each $s\in W^Q$. Conversely
it follows from Corollary \refer{c: variant of vanishing theorem} that 
$g=0$ if $g^s=0$
for some $s\in W^Q$.
The theorem follows immediately.
\qed

\begin{cor}
\naam{c: ind rels with i, new}
Let  $v \in \QW$ and let 
$\cL_t\in\cM(\staQqdc, \Sigma_Q)^*_\laur \otimes \oC(Q,v)$  be given
for each $t\in W_Q$. If
\begin{equation}
\naam{e: hypo ind rels with i, new}
\sum_{t\in W_Q}
\cL_t[E_{+,t}(\XQv\col \starstPo\col \dotvar\col m)  ] 
=0,\qquad (m \in \XQvp)
\end{equation}
then for each $s\in W^Q$ the following holds 
as a meromorphic identity in $\nu \in \faQqdc:$
\begin{equation}
\naam{e: conclusion ind rels with i, new}
\sum_{t\in W_Q}\cL_t[\;
E_{+,st}(\spX\col\stPo\col  \dotvar + \nu \col x)\after\iQv\;]
= 0,\qquad (x\in \spXp).
\end{equation}

Conversely, if the identity (\refer{e: conclusion ind rels with i, new}) 
holds for some $s\in W^Q$ and all
$\nu$ in a non-empty open subset of $\faQqdc$, then 
(\refer{e: hypo ind rels with i, new}) holds.
\end{cor}

\proof
For $t\in W_Q$ we define the functional 
$\Lau_t^\circ\in \Mer(\staQqdc,\Sigma_Q)^*_\laur \otimes \oC$ 
by
$\Lau_t^\circ = [ I \otimes \iQv ] (\Lau_t).$ 
Then for $F \in  \Mer(\staQqdc,\Sigma_Q) \otimes \oC^*$ 
we have 
\begin{equation}
\naam{e: Lau t circ}
\Lau_t^\circ F = \Lau_t [ F(\dotvar)\,\iQv].
\end{equation}
Let $u \in \QW.$ Then from
(\refer{e: prQu after iQv}) and (\refer{e: hypo ind rels with i, new}) 
we deduce
that (\refer{e: hypo ind rels with pr, new}) holds with $u$ and $\Lau_t^\circ$ 
in place of $v$ and $\Lau_t,$ respectively.
It follows that (\refer{e: conclusion ind rels with pr, new}) holds with $\Lau_t^\circ$ 
in place of $\Lau_t.$ In view of (\refer{e: Lau t circ}) 
this implies (\refer{e: conclusion ind rels with i, new}).
The converse statement is seen similarly.
\qed

Another useful formulation of the principle of induction of relations
is the following.

\begin{cor}
\naam{c: ind rels version 3, new}
Let $v \in \QW.$ Let $\Lau_t\in \Mer(\staQqdc, \Sigma_Q)^*_\laur$ 
and $\gf_t \in \Mer(\faqdc,\Sigma) \otimes \oC(Q,v)$ be given
for each $t\in W_Q$.
Assume that 
\begin{equation}
\naam{e: hypo ind rels with i and nu, new}
\sum_{t\in W_Q}\Lau_t[E_{+,t}(\spXQv\col \starstPo\col \dotvar \col m) 
\gf_t(\dotvar + \nu)]
=0,\qquad(m \in \spXQvp)
\end{equation}
for generic $\nu \in \faQqdc.$ Define 
$\psi_t = (I \otimes \iQv)\gf_t \in \Mer(\faqdc, \Sigma) \otimes \oC,
$ 
for $t\in W_Q.$
Then, for each $s\in W^Q$,
$$%\begin{equation}
%\naam{e: concl ind rels with i and nu, new}
\sum_{t\in W_Q}\Lau_t [E_{+,st}(\spX\col\stPo\col \dotvar + \nu \col x) 
\psi_t (\dotvar + \nu )]
=
0,\qquad(x\in \spXp)
$$%\end{equation}
as an identity
of $\Vtau$-valued meromorphic functions in the variable $\nu \in \faQqdc.$ 
\end{cor}

\proof
Let $\Hyp$ be a $\Sigma$-configuration such that $\sing(\gf_t)\subset\cup\Hyp,$ 
for each $t\in W_Q.$ Moreover, let 
$Y=\cup_{t\in W_Q}\supp\Lau_t \subset \staQqdc$. Fix  $t\in W_Q.$ Let
$\Hyp':= \Hyp_{\faQqdc}(Y)$ be the $\Sigma_r(Q)$-configuration
in $\faQqdc$ defined as in Corollary \refer{c: continuity of Laustar}. 
Let $\nu \in \faQqdc \setminus \cup \Hyp';$ then the function
$\gf_t^\nu: \gl \mapsto \gf_t(\gl + \nu)$ belongs to 
$\Mer(\staQqdc, Y, \Sigma_Q).$ 
It follows from (\refer{e: previous commutative diagram})
that the functional
$\Lau_t^\nu \in \Mer(\staQqdc, \Sigma_Q)^* \otimes \oC(Q,v)$ defined by 
$$ 
\Lau_t^\nu[F(\dotvar)] := \Lau_t[F(\dotvar)\gf_t(\dotvar + \nu)],
$$
for $F \in \Mer(\staQqdc, \Sigma_Q) \otimes \oC(Q,v)^*,$
is  a $\oC(Q,v)$-valued
$\Sigma_Q$-Laurent functional on $\staQqdc.$ 
The hypothesis (\refer{e: hypo ind rels with i and nu, new}) may be rewritten 
as (\refer{e: hypo ind rels with i, new}) with $\Lau_t^\nu$ in place of
$\Lau_t,$ for each $t\in W_Q.$ By application of Corollary 
\refer{c: ind rels with i, new} we therefore obtain, for
$\nu \in \faQqdc\setminus \cup \Hyp',$ that 
\begin{equation}
\naam{e: conclusion ind rels with i nu, new}
\sum_{t\in W_Q}\Lau_t [\;
E_{+,st}(\spX\col P_0\col \dotvar + \mu\col x)\, \psi_t(\dotvar + \nu) \;]
=
0
\end{equation}
as an identity of $\Vtau$-valued meromorphic functions in the variable 
$\mu \in \faQqdc.$ According to Lemma 
\refer{l: diagonal action of Laurent functional} the expression
in this equation defines a meromorphic $\Vtau$-valued 
function on $\faQqdc \times \faQqdc$ whose restriction to the diagonal
is a meromorphic function on $\faQqdc.$ Thus, if we substitute $\nu$ for $\mu$ 
in (\refer{e: conclusion ind rels with i nu, new}), we obtain an identity of 
$\Vtau$-valued meromorphic functions in the variable $\nu \in \faQqdc.$ 
\qed

\begin{cor}
\naam{c: induction of relations}
Let $\cL_1,\cL_2 \in \cM(\staQqdc, \Sigma_Q)^*_\laur \otimes \oC$. 
If, for each $v \in \QW$, 
\begin{equation}
\naam{e: hypo ind rels with pr}
\cL_1[ \nEp(\XQv\col \starstPo\col\dotvar\col m) \after \pr_{Q,v} ] 
=
\cL_2[
\nE(\XQv\col\starstPo\col \dotvar \col m) \after \pr_{Q,v} ],
\qquad(m \in \XQvp)
\end{equation}
then the following holds 
as a $\Vtau$-valued meromorphic identity in $\nu \in \faQqdc:$ 
\begin{equation}
\naam{e: conclusion ind rels with pr}
\cL_1[\;
\sum_{s \in W^Q} 
\Eps(\spX\col\stPo\col  \dotvar + \nu \col x)\;]
= \cL_2 [\;
\nE(\spX\col\stPo\col  \dotvar + \nu \col x) \;],\qquad(x\in\spXp).
\end{equation}
In particular, for regular values of $\nu,$ the expression 
on the left extends smoothly in the variable $x$ to all of $X.$ 

Conversely, if the identity (\refer{e: conclusion ind rels with pr})
holds for $\nu$ in a non-empty open subset of $\faQqdc$, then
(\refer{e: hypo ind rels with pr}) holds for each $v\in\QW$.
\end{cor}

\proof It follows from (\refer{e: splitting of Eis}) that 
$\nE(\XQv\col\starstPo\col \gl)=\sum_{t\in W_Q}
 E_{+,t}(\XQv\col \starstPo\col\gl).$
Define $\cL_t \in \cM(\staQqdc, \Sigma_Q)^*_\laur \otimes \oC$
for $t\in W_Q$ as follows. If $t=e$ then $\cL_t:=\cL_2-\cL_1$, 
and otherwise $\cL_t:=\cL_2$. Then the hypothesis
(\refer{e: hypo ind rels with pr, new}) in Theorem 
\refer{t: induction of relations, new} 
follows from (\refer{e: hypo ind rels with pr}). 
Hence the conclusion (\refer{e: conclusion ind rels with pr, new})
holds for each $s\in W^Q$. By summation over $s$ this implies that
\begin{equation}
\naam{e: equivalent of conclusion ind rels with pr}
\sum_{s\in W^Q}\sum_{t\in W_Q} 
\Lau_t[E_{+,st}(\spX\col\stPo\col  \dotvar + \nu \col x)]=0,
\qquad(x\in \spXp)
\end{equation}
which, by the definition of the operators $\Lau_t$ is equivalent
to (\refer{e: conclusion ind rels with pr}).

For the converse, let $g^s(\nu,x)$ denote the expression in
(\refer{e: conclusion ind rels with pr, new}), as in the proof of 
Theorem \refer{t: induction of relations, new}, with $\Lau_t$ 
specified as above. Then it was seen in the mentioned proof that
if the sum $g$ of the $g^s$ vanishes then so does each $g^s$
separately. Now (\refer{e: conclusion ind rels with pr})
implies (\refer{e: equivalent of conclusion ind rels with pr})
which exactly reads that $g=0$. Thus 
(\refer{e: conclusion ind rels with pr, new}) holds for each $s\in W^Q$,
so that the converse statement in Theorem 
\refer{t: induction of relations, new} can be applied.
\qed

The result just proved allows a straightforward 
corollary similar to Corollary \refer{c: ind rels with i, new},
in which the maps $\iQv$ are used instead of the maps $\pr_{Q,v}$.
We omit the details. The following result is derived from Corollary
\refer{c: ind rels version 3, new} in exactly the same way as the first
part of Corollary \refer{c: induction of relations} was 
derived from Theorem \refer{t: induction of relations, new}.

\begin{cor}
\naam{c: ind rels version 3}
Let $v \in \QW.$ Let $\Lau_1, \Lau_2\in \Mer(\staQqdc, \Sigma_Q)^*_\laur$ be 
$\Sigma_Q$-Laurent functionals on $\staQqdc,$ 
and let $\gf_1,\gf_2 \in \Mer(\faqdc,\Sigma) \otimes \oC(Q,v).$ 
Assume that 
$$%\begin{equation}
%\naam{e: hypo ind rels with i and nu}
\Lau_1(\Ep(\spXQv\col \starstPo\col \dotvar \col m) \gf_1(\dotvar + \nu))
=
\Lau_2(\nE(\spXQv\col \starstPo\col \dotvar\col m)\gf_2(\dotvar +\nu)),
$$%\end{equation}
for all $m \in \spXQvp$ and generic $\nu \in \faQqdc.$ Define 
$\psi_j = (I \otimes \iQv)\gf_j \in \Mer(\faqdc, \Sigma) \otimes \oC,
$ 
for $j=1,2.$
Then, for every $x\in \spXp,$
$$%\begin{equation}
%\naam{e: concl ind rels with i and nu}
\Lau_1 (\sum_{s \in W^Q} \Eps(\spX\col\stPo\col \dotvar + \nu \col x) \psi_1 (\dotvar + \nu ))
=
\Lau_2(\nE(\spX\col\stPo\col \dotvar + \nu \col x) \psi_2(\dotvar + \nu )),
$$%\end{equation}
as an identity
of $\Vtau$-valued meromorphic functions in the variable $\nu \in \faQqdc.$ 
\end{cor}

\begin{cor}  Let $v\in\QW$ and let $\psi_t\in
\cM(\staQqdc, \Sigma_Q)\otimes \oC(Q,v)$ be given
for each $t\in W_Q$. Let $\gl_0\in \staQqdc$.
Assume that for each $m\in\XQvp$, the meromorphic 
$V_\tau$-valued function on $\faQqdc$, given by
$$
\gl\mapsto\sum_{t\in W_Q}E_{+,t}(\XQv\col\starstPo\col\gl\col m) 
\psi_t(\gl), 
$$
is regular at $\gl_0$. Then for $s\in W^Q$, 
$x\in\spXp$ and generic $\nu \in \faQqdc$ 
the meromorphic function
\begin{equation}
\naam{e: function sum of Epst}
\gl\mapsto\sum_{t\in W_Q}
E_{+,st}(\spX\col\stPo\col  \gl+ \nu \col x)\iQv\psi_t(\gl)
\end{equation}
is also regular at $\gl_0$.
\end{cor}

\proof The function in (\refer{e: function sum of Epst}) has a germ
at $\gl_0$ in $\Mer(\staQqdc,\gl_0,\gS_Q)$. By 
Lemma \refer{l: annihilator of annihilator} it suffices to show
that it is annihilated by $\Mer(\staQqdc,\gl_0,\gS_Q)^{*\cO}_\laur$.
Let $\cL\in\Mer(\staQqdc,\gl_0,\gS_Q)^{*\cO}_\laur$ and
define $\cL_t\in\cM(\staQqdc, \Sigma_Q)^*_\laur \otimes \oC$
for $t\in W_Q$ by $\Lau_t=m^*_{\psi_t}\Lau$, see 
(\refer{e: previous commutative diagram}). The desired conclusion
now follows from Corollary \refer{c: ind rels with i, new}.\qed

We shall now give an equivalent formulation of the induction
of relations. We call it the lifting principle. For the group case
a similar principle was formulated by Casselman, 
see \bib{Arthur}, Thm.\ II.4.1, however with Eisenstein integrals
that carry a different normalization. 

\def\Asp{\cA_\laur}

\begin{defi}
\naam{d: defi Aspace}
The space $\Asp(\spXp\col\tau)$ is defined as
the space of functions
$$x\mapsto \Lau[E_+(\stPo\col\dotvar\col x)]\in V_\tau$$
where $\Lau\in\cM(\faqdc,\gS)^*_\laur\otimes\oC$.
It is a linear subspace of $C^\infty(\spXp\col\tau).$
\end{defi}

It follows from Corollary \refer{c: Lau to partial Eis} with $Q=G$ that
$\Asp(\spXp\col\tau)$ consists of $\DGH$-finite functions
in $C^\ep(\spXp\col\tau)$.  

\begin{rem} Let $\Lau\in\cM(\faqdc,\gS)^*_\laur\otimes\oC$.
Then $\Lau[\gf(\dotvar)E_+(\stPo\col\dotvar)]\in\Asp(\spXp\col\tau)$
for all $\gf\in\cM(\faqdc,\gS)$ (see (\refer{e: previous commutative diagram})).
In particular, it follows from (\refer{e: C-function in MerSigma}) that
$\nC(s\col\dotvar)\in \cM(\faqdc,\gS)\otimes\End(\oC).$
Hence it follows from the identity 
(\refer{e: Esp in terms of Ep and C})
that $\Lau[E_{+,s}(\stPo\col\dotvar)]\in\Asp(\spXp\col\tau)$
for each $s\in W$. Moreover, by similar reasoning it can be seen that
the space $\Asp(\spXp\col\tau)$ does not depend on the choice of 
$P_0\in\minparabs$.
\end{rem}

\begin{rem} Let $\gl_0\in\faqdc$ and 
$\gf\in\cM(\faqdc,\gS)\otimes\oC$, and assume that
$\gl\mapsto E_+(\stPo\col\gl)\gf(\gl)$ is regular at $\gl_0$.
Then the function 
$x\mapsto u[E_+(\stPo\col\gl\col x)\gf(\gl)]|_{\gl=\gl_0}$
belongs to $\Asp(\spXp\col\tau)$ for each $u\in S(\faqd)$
(see the previous remark and Lemma \refer{l: diff of Laur}).
Moreover, it follows easily from the definition of 
$\cM(\faqdc,\gS)^*_\laur$ that $\Asp(\spXp\col\tau)$ is spanned 
by functions of this form. 
\end{rem}

\begin{thm} 
\naam{t: lifting principle}
{\rm (Lifting principle)}
Let $Q\in\allparabs$ be a standard parabolic subgroup,
and let $s\in W^Q$ be fixed. 
\begin{enumerate}
\itema There exists for each $v\in\QW$ a unique linear map 
$$F_{+,s,v}\colon \Asp(\spX_{Q,v,+}\col\tau_Q)
\to\cM(\faQqdc,\gS_r(Q),C^\infty(\spXp\col\tau))$$ with the following
property. If 
$\gf\in\Asp(\spX_{Q,v,+}\col\tau_Q)$ is given by
\begin{equation}
\naam{e: 5}
\gf(m)=
\sum_{t\in W_Q}\Lau_t[E_{+,t}(\XQv\col\starstPo\col\dotvar\col m)]
\qquad(m\in \XQvp),
\end{equation}
for some $\Lau_t\in\cM(\staQqdc, \Sigma_Q)^*_\laur \otimes \oC(Q,v)$, 
$t\in W_Q$, then
\begin{equation}
\naam{e: 6}
F_{+,s,v}(\gf)(\nu,x)=\sum_{t\in W_Q}\Lau_t[\;
E_{+,st}(\spX\col\stPo\col  \dotvar + \nu \col x)\;\iQv\;]
\end{equation}
for $x \in \spXp$ and generic $\nu\in\faQqdc$. 
\itemb
The function $x\mapsto F_{+,s,v}(\gf,\nu,x)$ belongs to $\Asp(\spXp\col\tau)$
for generic $\nu$.
\itemc
The map $$F_{+,s}\colon\quad \oplus_{v\in\QW}\Asp(\spX_{Q,v,+}\col\tau_Q)
\to\cM(\faQqdc,\gS_r(Q),C^\infty(\spXp\col\tau)),$$ given by 
$F_{+,s}(\gf)=\sum_v F_{+,s,v}\gf_v$,
is injective.
\end{enumerate}
\end{thm}

\proof The uniqueness is clear from Definition \refer{d: defi Aspace}.
We use (\refer{e: 5}) and (\refer{e: 6}) as the definition of $F_{+,s,v}$; 
the fact that $F_{+,s,v}(\gf)$ is well defined for all 
$\gf\in\Asp(\spX_{Q,v,+}\col\tau_Q)$ is equivalent with the first
statement in Theorem \refer{t: induction of relations, new} 
(see also Corollary \refer{c: ind rels with i, new}). Once the definition makes sense,
it is easily seen that $F_{+,s,v}(\gf)$ depends linearly on $\gf$. 
That $F_{+,s,v}(\gf,\nu)\in\Asp(\spXp\col\tau)$ for generic $\nu$
is seen from Lemma \refer{l: eval of Laustar is a Lau}.
Finally, the injectivity of $F_{+,s}$ is equivalent with the
final statement of Theorem \refer{t: induction of relations, new}.\qed

\begin{rem} Note that 
with $\gf_v=\nE(\XQv\col\starstPo\col\gl)$ for each $v\in\QW$ we obtain
$$\sum_{t\in W_Q}E_{+,st}(\spX\col\stPo\col  \gl + \nu\col x )\;\iQv=
F_{+,s,v}(\gf_v,\nu,x),
$$
for $x \in \spXp$, and hence by summation over $v$ and $s$
$$
\nE(\spX\col\stPo\col  \gl + \nu \col x)=
\sum_{s\in W^Q} F_{+,s}(\gf,\nu,x).
$$
\end{rem}

\begin{rem} In \bib{BSfi}, Definition 10.7, we define the generalized
Eisenstein integral $E^\circ_F(\psi\col\nu)\in\Ci(\spX\col\tau)$ for
$\psi\in\cC_F$, $\nu\in\fa_{F\iq\iC}^*$ (with the notation of {\it loc.\ cit.}).
By comparison with Theorem \refer{t: lifting principle}
for $Q=P_F$ it is easily seen that
$E^\circ_F(\psi\col\nu\col x)=F_{+,1}(\psi,\nu,x)$
for $x\in\spXp$.
\end{rem}

\section{Appendix A: spaces of holomorphic functions}
\naam{s: certain function spaces}
If
$\Omega$ is a complex analytic manifold, then by $\cO(\Omega)$ we denote the space
of holomorphic and by $\Mer(\Omega)$ the space of meromorphic functions on $\Omega.$ 

If $V$ is a complete locally convex (Hausdorff) space, we say that a function $\gf: \Omega \to V$ is 
holomorphic if for every $a \in \Omega$ there exists a holomorphic coordinatisation $z = (z_1, \ldots, z_n)$ 
at $a$ such that in a neighborhood of $a$ the function $\gf$ is expressible as
a converging $V$-valued power series in the coordinates $z.$ The space of such holomorphic 
functions is denoted by $\cO(\Omega, V).$ We equip this space with a 
locally convex topology as follows. Let $\cP$ be a separating collection of 
continuous seminorms
for $V.$ For every $p \in \cP$ and every compact set $K \subset \Omega$ we define
the seminorm $\nu_{K,p}$ on $\cO(\Omega, V)$ by $\nu_{K,p}(\gf) = \sup_{K} p\after \gf.$ 
This collection of seminorms is separating
hence equips $\cO(\Omega, V)$ with a locally convex topology.
Note that this topology is independent of the choice of $\cP.$ Moreover, it is complete; it is 
Fr\'echet if $V$ is a Fr\'echet space. 

We recall that $\cO(\Omega, V)$ is a closed subspace of $\Ci(\Omega, V).$ 
Indeed, if $\bar\partial$ denotes the anti-linear part of exterior differentiation,
then $\cO(\Omega, V)$ is the kernel of $\bar\partial$ in $\Ci(\Omega, V).$

A densely defined function $\gf: \Omega \to V$ is called meromorphic 
if for every $a \in \Omega$ there exists an open neighborhood 
$U$ of $a,$  
and a function $\psi \in \cO(U)\setminus\{0\}$ such that $\psi \gf \in \cO(U,V).$ 
As usual, meromorphic functions are considered to be equal if 
they coincide on a dense open subset. The space of $V$-valued meromorphic 
functions on $\Omega$ is denoted by $\Mer(\Omega, V).$ 
If $\gf$ is an $V$-valued
meromorphic function on $\Omega$ we define $\reg(\gf)$ to be the largest open 
subset $U$ of $\Omega$ for which $\gf|_U$ coincides (densely) with an element of $\cO(U, V).$ 
The complement $\sing(\gf) = \Omega\setminus \reg(\gf)$ is called the singular locus of $\gf.$

\begin{lemma}
\naam{l: natural isomorphism of function spaces}
Let $X$ be a $C^\infty$ and $\Omega$ a complex analytic manifold.
Let $V$ be a complete locally convex space.

Let $\cF$ be the locally convex space of $C^\infty$-functions
$X \times \Omega \to V$ that are holomorphic in the second variable.
Given $f \in \cF$ and $x \in X,$ we define the function
${}_1f(x) : \Omega \to V$ by ${}_1f(x)(z) = f(x, z).$ Given $z \in \Omega$ 
we define the function ${}_2f(z): X \to V$ by ${}_2f(z)(x) = f(x, z).$ 
\begin{enumerate}
\itema
The map $f \mapsto {}_1f$ defines a natural isomorphism
of locally convex spaces from $\cF$ onto $\Ci(X, \cO(\Omega, V)).$ 
\itemb
The map $f \mapsto {}_2f$ defines a natural isomorphism 
of locally convex spaces from $\cF$ onto $\cO(\Omega, \Ci(X,V)).$  
\end{enumerate}
In particular, the above maps lead to a natural isomorphism
$$
\Ci(X, \cO(\Omega, V)) \simeq \cO(\Omega, \Ci(X, V)).
$$
\end{lemma}

\proof
The above isomorphisms are valid with $\cO$ replaced by $\Ci$ everywhere.
This is a well known fact, and basically a 
straightforward consequence of the definitions, though somewhat tedious to check.
The isomorphisms with $\cO$ are seen to be valid by showing that the appropriate 
kernels of the operator $\bar\partial$ correspond. 
Checking this involves a local application of the multivariable Cauchy integral formula.
\qed

\section{Appendix B: removable singularities}
\naam{s: appendix b}
We discuss a variation on the idea of removable singularities
for holomorphic functions that is particularly useful for application 
in the present paper.

A subset $T$ of a finite dimensional complex analytic manifold  $\Omega$ will be called thin
if for every $\gl\in \Omega$ there exists a connected open neighborhood $U$ and 
a non-zero holomorphic function $\gf \in \cO(U)$ 
such that $T \cap U \subset  \gf^{-1}(0)$, see
\bib{GuRo}, p.\ 19.
An open subset $U$ of $\Omega$ will be called full if its complement is thin.
It is clear that a full subset of $\Omega$ is dense in $\Omega$. 
Note that the union of finitely many thin subsets is thin again; accordingly,
the intersection of finitely many full open subsets of $\Omega$ is again a full open
subset. Obviously any union of full open subsets
is a full open subset.
Note also that if $\Omega$ is connected, then every full open subset of $\Omega$ 
is connected (\bib{GuRo} p.\ 20).

\begin{lemma}
Let $j: V \to W$ be an injective continuous linear 
map of complete locally convex Hausdorf spaces,
and let $F$ be a $W$-valued holomorphic function
on a complex analytic manifold $\Omega.$ 
Assume that there exists a full open subset $\Omega_0$ of $\Omega$ 
and a holomorphic function $G_0: \Omega_0 \to V$ such that 
such that $F = j \after G_0$ on $\Omega_0.$ 
Then there exists a unique holomorphic map $G : \Omega \to V$ 
such that $j\after G = F.$ 
\end{lemma}

\proof
Clearly the result is of a local nature 
in the $\Omega$-variable, so that we may assume that
$\Omega$ is a connected open subset of $\C^n,$ for some $n \in \N.$ 
Moreover, we may as well assume that $\Omega_0 = \Omega\setminus \gf^{-1}(0),$ 
with $\gf \in \cO(\Omega)$ a non-zero holomorphic function.

Fix $\gl_0 \in \Omega.$ Since $\gf$ is non-zero, the function
$z\mapsto \gf(\gl_0+z\mu)$, defined on a neighborhood of 0 in $\C$,
is non-zero for some $\mu\in\C^n\setminus \{0\}$. Being holomorphic,
this function then takes the value 0 in isolated points. Hence we may
choose $\mu$ such that $\gl_0+z\mu\in\Omega_0$ for
$0<|z|\le 1$. By compactness there exists an open neighborhood $N_0$ 
of $\gl_0$ in $\Omega$ such that $\gl + z \mu \in 
\Omega$  for all $\gl \in N_0$
and  $z \in \C$ with $|z|\le 1,$ and such that
$\gl + z \mu \in\Omega_0$  for  $|z| =1.$ 
By the Cauchy integral formula we  have:
\begin{equation}
\naam{e: Cauchy int}
F(\gl) = \frac{1}{2\pi i} \int_{\partial D} F(\gl + z \mu) \frac{dz}{z}.
\end{equation}
Here $\partial D$ denotes the boundary of the unit circle 
in $\C,$ equipped with the orientation induced by the complex structure,
i.e., the counter clockwise direction. 
 
Note that the $W$-valued (or $V$-valued) integration is well defined, since
$W$ (or $V$) is complete locally convex. In the integrand 
of (\refer{e: Cauchy int}) the function
$F(\gl + z \mu)$ may be replaced by 
$ j(G_0(\gl + z \mu)).$ Using that 
$j$ is continuous linear we then obtain that
\begin{equation}
\naam{e: F is j F zero}
F(\gl) = j ( G(\gl) ),
\end{equation}
where 
$$
G(\gl) := \frac{1}{2 \pi i} \int_{\partial D}  
G_0(\gl + z \mu) \; \frac{dz}{z}\qquad (\gl\in N_0).
$$
Clearly $G: N_0 \to V$ is a holomorphic function;
moreover, it is uniquely determined by equation
(\refer{e: F is j F zero}), since $j$ is injective. This implies
that the local definition of $G$ is independent of the particular
choice of $\mu.$ Moreover, it also follows
from (\refer{e: F is j F zero}) and the injectivity of
$j$ that all local definitions 
match and determine a global holomorphic function $G: \Omega \to V$ 
satisfying our requirement.
\qed

\begin{cor}
\naam{c: aux smooth extension}
Let $\Omega_0$ be a full open subset of a complex analytic manifold  $\Omega$ 
and let $X_0$ be a dense open subset of a $C^\infty$-manifold
$X.$ Moreover, let $F: \Omega \times X_0 \to \C$ be 
a $C^\infty$ function that is holomorphic in its first variable,
and assume that its restriction to $\Omega_0 \times X_0$ has 
a smooth extension to $\Omega_0 \times X.$  Then the 
function $F$ has a unique smooth extension to $\Omega \times X.$ 
Moreover, the extension
is holomorphic in its first variable. 
\end{cor}

\proof
As in the proof of the above lemma we may as well assume 
that $\Omega$ is an open subset of $\C^n,$ for some $n.$ 

Let $V = C^\infty(X)$ and $W= C^\infty(X_0)$ be equipped with 
the usual Fr\'echet topologies. Restriction to $X_0$ 
induces an injective continuous linear map $j: V \to W.$ 

By 
Lemma \refer{l: natural isomorphism of function spaces}(b)
we see that the function 
$\widetilde F: \Omega \to W ,$ 
defined by $\widetilde F(z) = F(z, \dotvar)$ 
is holomorphic. Let $G_0$ be the extension
of $(z, x) \mapsto  F(z, x)$ to
a smooth map $\Omega_0 \times X \to \C.$ Then by density and continuity
the function $G_0$ satisfies the Cauchy-Riemann equations
in its first variable. Hence it is holomorphic in its 
first variable, and it follows that the function
$\widetilde G_0: \Omega_0 \to V$ defined by $\widetilde G_0(z) = G_0(z, \dotvar)$ 
is holomorphic. {}From the definitions given 
we obtain that $ \widetilde F = j \after \widetilde G_0$ on $\Omega_0.$ 
By the above lemma
there exists a unique holomorphic function $\widetilde G: \Omega \to V$
such that $\widetilde F = j \after \widetilde G.$ 
The function $G: (z, x) \mapsto \widetilde G(z)(x)$ 
is the desired extension of $F.$ 
\qed

     %references
\def\adritem#1{\hbox{\small #1}}
\def\distance{\hbox{\hspace{3.5cm}}}
\def\apetail{@}

\def\adderik{\vbox{
\adritem{E.P. van den Ban}
\adritem{Mathematisch Instituut}
\adritem{Universiteit Utrecht}
\adritem{PO Box 80 010}
\adritem{3508 TA Utrecht}
\adritem{Netherlands}
\adritem{E-mail: ban{\apetail}math.uu.nl}
}
}

\def\addhenrik{\vbox{
\adritem{H. Schlichtkrull}
\adritem{Matematisk Institut}
\adritem{K\o benhavns Universitet}
\adritem{Universitetsparken 5}
\adritem{2100 K\o benhavn \O}
\adritem{Denmark}
\adritem{E-mail: schlicht@math.ku.dk}
}
}

\vfill
\hbox{\vbox{\adderik}\vbox{\distance}\vbox{\addhenrik}}

    %our addresses
\end{document}